Escuela Superior Técnica de Ingeniería
Departamento de Ingeniería de Sistemas y Automática

Doctoral Thesis

# Implementation of MPC in embedded systems using first order methods

**Pablo Krupa García**

---

Supervised by:
Daniel Limon Marruedo
Teodoro Alamo Cantarero

Seville, June 2021

# Contents

# Acknowledgements

First and foremost, I want to thank my supervisors, Professors Daniel Limón and Teodoro Álamo. The passion and energy of Daniel, along with his good-hearted disposition and ever-present willingness to discuss my problems and questions, hooked me in and kept me going. It was in his office, in 2015, that this journey began when I met him to ask about a graduate thesis he was offering. I found myself, a solid hour and a half later, having received a private lecture on MPC simply because I showed some interest. Looking back at that meeting, I have come to realize that it was one of the most life-changing moments of my life. I wouldn't be here without the support and the knowledge he has shared with me since then. Teodoro's attitude towards science and research, his strive for excellence, and his seemingly endless knowledge have inspired me and made me grow throughout these years. His constant search of a better algorithm and a better proof, and his obsession with using the proper punctuation marks after every single equation, have, for better or worse, been deeply drilled into my mind. I will always fondly remember the hours spent in his office discussing optimization algorithms and MPC; oftentimes with me struggling to keep up with his line of thought, but finally managing to do so through his sheer perseverance. I don't know if I would have the patience. I am deeply in debt with both of them; more, I think, than they are aware of.

I would also like to thank all the other members of the *GEPOC* research group; those currently there and those now working elsewhere. They have always been helpful and kind. One cannot ask for a better work environment. I would also like to acknowledge Nilay Saraf, who is a pleasure to work with. I know few people as willing as he to work the extra hour to get a better scientific result.

A huge thanks to my family, friends and loved ones. I will not list them all for fear of leaving someone out, but there are a few I feel deserve a special mention. To my parents and brother, for the numerous times they have asked about my progress during these years, even though they may not understand half of what I'm talking about. Likewise, I would like to thank my friends Manuel and Isa, who are always supportive (and as insistent as humanly possible). In particular, I would like to thank Manuel for his huge effort in trying to understand my work and all the interest he has shown over the years. Thanks to Laura, who has also always shown a great interest in what I do; going so far as to promise that she will read this dissertation. I'm holding you to your word. A deep and special thanks to Paula, for always being there and encouraging me, particularly during these last few months while I was writing this dissertation.

Finally, thanks to all the people I worked with in the *Mitsubishi Electric Research Laboratories*. In particular, a special thanks to Dr. Claus Danielson, who welcomed and helped me during my stay and from whom I learnt a lot.

Pablo Krupa García,
Seville, June 2021.

# Notation, conventions and definitions

We list here the notation, basic definitions, conventions and well-known mathematical background that we use throughout the paper. Additionally, section-specific definitions and notation may be presented throughout the manuscript when needed. The definitions and conventions shown here are standard in the literature, and therefore no proofs nor specific references are provided. We refer the reader to [1, 2, 3] as general references containing most of the definitions and conventions stated here.

### Spaces and general set notation

Let $\mathcal{C}$ and $\mathcal{D}$ be two sets. We use the standard set notations $\mathcal{C} \subset \mathcal{D}$ ($\mathcal{C}$ is a strict subset of $\mathcal{D}$), $\mathcal{D} \subseteq \mathcal{D}$ ($\mathcal{C}$ is a subset or equal to $\mathcal{D}$), $\mathcal{C} \cap \mathcal{D}$ (intersection of $\mathcal{C}$ and $\mathcal{D}$) and $\mathcal{C} \cup \mathcal{D}$ (union of $\mathcal{C}$ and $\mathcal{D}$). If $x$ is an element of set $\mathcal{C}$, we write $x \in \mathcal{C}$. The empty set is denoted by $\emptyset$.
The set of real numbers is denoted by $\mathbb{R}$. The set of extended real numbers, i.e., $\mathbb{R}$ extended with $-\infty$ and $+\infty$, is denoted by $\overline{\mathbb{R}}$. We denote by $\mathbb{R}_{\geq 0}$ and $\mathbb{R}_{>0}$ the set of non-negative and (strictly) positive real numbers, respectively. We use analogous notations for the non-positive ($\mathbb{R}_{\leq 0}$) and negative ($\mathbb{R}_{<0}$) real numbers as well as for their extended real number counterparts ($\overline{\mathbb{R}}_{>0}$, $\overline{\mathbb{R}}_{\geq 0}$, $\overline{\mathbb{R}}_{<0}$, $\overline{\mathbb{R}}_{\leq 0}$). We denote by $[a, b]$ the set of real numbers $x$ satisfying $a \leq x \leq b$. We use a rounded bracket (instead of a square one), to indicate strict inequality, e.g., $(a, b]$ denotes the set of real numbers $x$ satisfying $a < x \leq b$.
The set of integer numbers is denoted by $\mathbb{Z}$. For integers $i, j \in \mathbb{Z}$ satisfying $i \leq j$, we denote $\mathbb{Z}_i^j \doteq \{ x \in \mathbb{Z} : i \leq x \leq j \}$. We use the same notation used with real numbers to indicate the sets of positive ($\mathbb{Z}_{>0}$), non-negative ($\mathbb{Z}_{\geq 0}$), negative ($\mathbb{Z}_{<0}$) and non-positive ($\mathbb{Z}_{\leq 0}$) integer numbers.

### Vectors and matrices

We denote by $\mathbb{R}^n$ the set of $n$-dimensional real vectors. For $x \in \mathbb{R}^n$, we denote by $x_{(i)}$ its $i$th component. The inner product (or *dot product)* of two vectors $x, y \in \mathbb{R}^n$ is denoted by $\langle x, y \rangle = \sum_{i=1}^n x_{(i)} y_{(i)}$. All vectors are considered column vectors unless otherwise specified. For vectors $x_1 \in \mathbb{R}^{n_1}, x_2 \in \mathbb{R}^{n_2}, \dots, x_N \in \mathbb{R}^{n_N}$, with $N \in \mathbb{Z}_{>0}$, we denote by $(x_1, x_2, \dots, x_N)$ the column vector formed by the concatenation of column vectors $x_1$ to $x_N$. Vector $(x_1, x_2, \dots, x_N) \in \mathbb{R}^{n_1} \times \mathbb{R}^{n_2} \times \cdots \times \mathbb{R}^{n_N}$ can also be viewed as the Cartesian product of vectors $x_1$ to $x_N$. The transposed of a vector $x \in \mathbb{R}^n$ is denoted by $x^\top$. For vectors $x, y \in \mathbb{R}^n$, the notations $x \leq y$, $x > y$, etc. are to be taken componentwise.
We denote by $\mathbb{R}^{m \times n}$ the space of $m$ by $n$ dimensional real matrices. For $M \in \mathbb{R}^{m \times n}$, we denote by $M_{(i,j)}$ its $(i, j)$th component, $M^\top$ its transposed and $M^{-1}$ its inverse (if $M$ is non-singular). We say that $M$ is diagonal if $M_{(i,j)} = 0$ whenever $i \neq j$. The space of diagonal matrices in $\mathbb{R}^{n \times n}$ is denoted by $\mathbb{D}^n$. We denote by $\mathbb{S}^n_{++}, \mathbb{S}^n_{+} \subset \mathbb{R}^{n \times n}$ the spaces of positive definite and positive semi-definite matrices in $\mathbb{R}^{n \times n}$, respectively. We denote their diagonal counterparts by $\mathbb{D}^n_{++} \subset \mathbb{S}^n_{++}$ and $\mathbb{D}^n_{+} \subset \mathbb{S}^n_{+}$. Given a symmetric matrix $M$, its maximum and minimum eigenvalues are given by $\lambda_{\max}(M)$ and $\lambda_{\min}(M)$, respectively. For matrices and/or scalars $M_1, M_2, \ldots, M_N$ (not necessarily of the same dimension), with $N \in \mathbb{Z}_{>0}$, we denote by $\text{diag}(M_1, M_2, \ldots, M_N)$ the diagonal concatenation of $M_1$ to $M_N$.
We denote by $\mathbf{1}_{n \times m} \in \mathbb{R}^{n \times m}$ and $\mathbf{0}_{n \times m} \in \mathbb{R}^{n \times m}$ the matrices of all ones and all zeros, respectively. The vectors of ones and zeros in $\mathbb{R}^n$ are denoted by $\mathbf{1}_n$ and $\mathbf{0}_n$, respectively. The identity matrix of dimension $n$ is denoted by $\mathbf{I}_n$. We may drop the dimensions of the matrices/vectors if they are clear from the context, and instead simply write $\mathbf{1}$, $\mathbf{0}$ or $\mathbf{I}$.

### Norms

Given $a \in \mathbb{R}$, $|a|$ denotes its absolute value. Given $x \in \mathbb{R}^n$, we denote: the standard Euclidean norm by $\|x\|_2 \doteq \sqrt{\langle x, x \rangle}$; the $M$-weighted Euclidean norm by $\|x\|_M \doteq \sqrt{\langle x, Mx \rangle}$ (provided that $M \in \mathbb{S}^n_{++}$); the $\ell_1$-norm by $\|x\|_1 \doteq \sum_{i=1}^n |x_{(i)}|$; the $\ell_\infty$-norm by $\|x\|_\infty \doteq \max_{i \in \mathbb{Z}_1^n} |x_{(i)}|$; and $\|x\|$ to represent any vector norm. The spectral norm of a matrix $M \in \mathbb{R}^{n \times m}$ is denoted by $\|M\| \doteq \sqrt{\lambda_{\max}(M^\top M)}$.

**Definition N.1** (Dual norm). Given any norm $\|\cdot\| : \mathbb{R}^n \to \mathbb{R}_{\geq 0}$, its *dual norm* $\|\cdot\|_* : \mathbb{R}^n \to \mathbb{R}_{\geq 0}$ is defined as $\|x\|_* \doteq \sup_{z \in \mathbb{R}^n} \{ \langle x, z \rangle \,:\, \|z\| \leq 1 \}$.

If $\|\cdot\|$ is an $M$-weighted Euclidean norm $\|x\| = \sqrt{\langle x, Mx \rangle}$, where $M \in \mathbb{S}^n_{++}$, then its dual norm is given by $\|x\|_* = \sqrt{\langle x, M^{-1}x \rangle}$.

### Sequences

A sequence of elements $x_j$ indexed by $j \in \mathcal{J} \subseteq \mathbb{Z}$ is denoted by $\{x_j\}$, where the set $\mathcal{J}$ will be clear from the context. Additionally, we use $\{x_j\}_{\geq 0}$ if $\mathcal{J} = \mathbb{Z}_{\geq 0}$, and $\{x_j\}_{>0}$ if $\mathcal{J} = \mathbb{Z}_{>0}$. If the index set $\mathcal{J}$ is finite, then we denote by a bold $\mathbf{x}$ the Cartesian product of $\{x_j\}$. We often times append the index term $j$ as a superscript, instead of a subscript, i.e., we write $x^j$ to denote element $j$ of an ordered sequence. In particular, we do this to express the iterates of some algorithms. For signals (or sequences) that evolve over time we use the following specific notation. Let $x \in \mathbb{R}^n$ be a vector whose value changes over time, which we index by letter $t$. We denote by $x(t)$ the value of $x$ at time instant $t$.

### Set definitions

**Definition N.2** (Convex set)**.** A set $\mathcal{C} \subseteq \mathbb{R}^n$ is *convex* if $\alpha x + (1-\alpha) y \in \mathcal{C}$, $\forall x, y \in \mathcal{C}$, $\forall \alpha \in [0,\, 1]$.

**Definition N.3** (Affine set)**.** A set $\mathcal{C} \subseteq \mathbb{R}^n$ is said to be *affine* if it contains all the lines that pass through pairs of points $x, y \in \mathcal{C}$ with $x \neq y$.

**Definition N.4** (Closed, open, bounded and compact sets)**.** A non-empty set $\mathcal{C} \in \mathbb{R}^n$ is *closed* if the limit of every converging sequence $\{x_k\}$, with $x_k \in \mathcal{C}$, is also contained in $\mathcal{C}$. It is *open* if its complement $\{\, x \in \mathbb{R}^n \,:\, x \notin \mathcal{C} \,\}$ is closed. It is *bounded* is there exists $c \in \mathbb{R}_{>0}$ such that $\|x\| \leq c$, $\forall x \in \mathcal{C}$. It is *compact* if it is closed and bounded.

**Definition N.5** (Affine hull)**.** The *affine hull* of a set $\mathcal{C}$, which is denoted by $\operatorname{aff}(\mathcal{C})$, is the intersection of all the affine sets containing $\mathcal{C}$.

**Definition N.6** (Interior and relative interior of a set)**.** A point $x \in \mathbb{R}^n$ is an *interior point* of a set $\mathcal{C} \subseteq \mathbb{R}^n$ is there exists an open sphere $S \doteq \{\, z \in R^n \,:\, \|z - x\|_2 < \epsilon \,\}$ that is contained in $\mathcal{C}$. The set of all interior points of $\mathcal{C}$ is called the *interior*, and is denoted by $\operatorname{int}(\mathcal{C})$. A point $x \in \mathbb{R}^n$ is a *relative interior point* of a set $\mathcal{C} \subseteq \mathbb{R}^n$ is there exists an open sphere $S$ for some $\epsilon \in \mathbb{R}_{>0}$ such that $S \cap \operatorname{aff}(\mathcal{C}) \subset \mathcal{C}$. The set of all relative interior points of of $\mathcal{C}$ is called the *relative interior*, and is denoted by $\operatorname{ri}(\mathcal{C})$.

**Definition N.7** (Ellipsoid)**.** For a given $P \in \mathbb{S}^n_{++}$, $c \in \mathbb{R}^n$, $r \in \mathbb{R}_{>0}$ we denote by

$$\mathcal{E}(P, c, r) \doteq \{\, x \in \mathbb{R}^n \,:\, \|x - c\|_P^2 \leq r^2 \,\}$$

the ellipsoid centered at $c$ with radius $r$ and geometry determined by $P$.

**Definition N.8** (Admissible invariant set)**.** Consider a discrete-time autonomous system $z(t+1) = f(z(t))$, where $t$ is the current discrete-time instant and $f : \mathbb{R}^n \to \mathbb{R}^n$ describes the system dynamics, that is subject to $z(t) \in \mathcal{Z} \subseteq \mathbb{R}^n$, $\forall t$. An admissible invariant set of the system is a set $\Omega \subseteq \mathcal{Z}$ such that

$$z(t) \in \Omega \implies z(t+1) \in \Omega.$$

### Functions

We use the standard notation for the derivative and partial derivative of a differentiable real-valued function $f$. That is, $\dfrac{df}{dx}$ is the derivative of $f(x)$ with respect to $x$ (we may also use the shorthand $f'$) and $\dfrac{\partial f}{\partial x}$ is the partial derivative of $f(x, y, \dots)$ with respect to $x$. For time derivatives we may use the shorthands $\dot{f} \doteq \frac{df}{dt}$ and $\ddot{f} \doteq \frac{d^2 f}{dt^2}$. We use the notation $\left.\frac{df(x)}{dx}\right|_{\hat{x}}$ to indicate the value of the derivative of $f(x)$ evaluated at $\hat{x}$.
The gradient of a differentiable function $f : \mathbb{R}^n \to \mathbb{R}$ is denoted by $\nabla f$. We denote by $\nabla f(x)$ its value at $x \in \mathbb{R}^n$.

**Definition N.9** (Lipschitz continuity)**.** A function $f : \mathcal{C} \to \mathbb{R}^m$, with $\mathcal{C} \subseteq \mathbb{R}^n$ is *Lipschitz continuous* if there exists a constant $L \in \mathbb{R}_{>0}$ such that

$$\|f(x) - f(y)\|_2 \leq L\|x - y\|_2, \ \forall x, y \in \mathcal{C}.$$

This condition will also be satisfied for any $\hat{L} > L$. We call the smallest value for which it is satisfied the *Lipschitz constant*.

**Definition N.10** (Extended real-valued function)**.** A function $f$ is said to be a *real valued function* if it maps $f : \mathbb{R}^n \to \mathbb{R}$, and it is said to be an *extended real-valued function* it it maps $f : \mathbb{R}^n \to \overline{\mathbb{R}}$.

**Definition N.11** (Epigraph)**.** The *epigraph* of a function $f : \mathcal{C} \to \overline{\mathbb{R}}$, where $\mathcal{C} \subseteq \mathbb{R}^n$, is the subset of $\mathbb{R}^{n+1}$ given by

$$\operatorname{epi}(f) \doteq \{\, (x, w) \ : \ x \in \mathcal{C}, w \in \mathbb{R}, f(x) \leq w \,\}.$$

**Definition N.12** (Domain)**.** The (effective) *domain* of an extended real-valued function $f : \mathcal{C} \to \overline{\mathbb{R}}$, where $\mathcal{C} \subseteq \mathbb{R}^n$, is given by

$$\operatorname{dom}(f) \doteq \{\, x \in \mathcal{C} \ : \ f(x) < +\infty \,\}.$$

**Definition N.13** (Proper)**.** Let $\mathcal{C} \subseteq \mathbb{R}^n$ be a convex set. An extended real-valued function $f : \mathcal{C} \to \overline{\mathbb{R}}$ is said to be *proper* if $f(x) < +\infty$ for at least one $x \in \mathcal{C}$ and $f(x) > -\infty$ for all $x \in \mathcal{C}$.

**Definition N.14** (Convex function)**.** Let $\mathcal{C} \subseteq \mathbb{R}^n$ be a convex set. A function $f : \mathcal{C} \to \mathbb{R}$ is *convex* if

$$f(\alpha x + (1 - \alpha)y) \leq \alpha f(x) + (1 - \alpha)f(y), \quad \forall x, y \in \mathcal{C}, \ \forall \alpha \in [0, 1].$$

The function is called *strictly convex* is the above expression holds with strict inequality for all $x, y \in \mathcal{C}$ with $x \neq y$ and all $\alpha \in (0, 1)$.

If, in addition to convex, function $f$ is differentiable, then it satisfies

$$f(x) \geq f(y) + \langle \nabla f(y), x - y \rangle, \ \forall x, y \in \mathcal{C}.$$

If $f : \mathbb{R}^n \to \overline{\mathbb{R}}$ is a proper extended real-valued function, then it is convex if $\operatorname{epi}(f)$ is a convex subset of $\mathbb{R}^{n+1}$. In this case, the above inequality and Definition N.14 hold for all $x, y \in \operatorname{dom}(f)$.

**Definition N.15** (Concave function)**.** A function $f$ is concave if $-f$ is convex.

**Definition N.16** (Closed function)**.** Let $\mathcal{C} \subseteq \mathbb{R}^n$ be a convex set. A function $f : \mathcal{C} \to \overline{\mathbb{R}}$ is *closed* if its epigraph is a closed set.

**Notation N.17** (Space of proper closed convex functions)**.** We denote the space of proper closed convex functions $f : \mathbb{R}^n \to \overline{\mathbb{R}}$ by $\Gamma(\mathbb{R}^n)$.

**Definition N.18** (Subdifferential)**.** Let $f \in \Gamma(\mathbb{R}^n)$. A vector $g$ is said to be a *subgradient* of $f$ at $x \in \mathrm{dom}(f)$ if

$$f(y) \geq f(x) + \langle g, y - x \rangle, \quad \forall y \in \mathbb{R}^n.$$

The set of all subgradients of $f$ at $x$ is called the *subdifferential of $f$ at $x$*, and is denoted by $\partial f(x)$. If $f$ is differentiable, then its gradient $\nabla f(x)$ at any point $x \in \mathrm{int}(\mathrm{dom}(f))$ is the unique subgradient of $f$ and $x$, i.e., $\partial f(x) = \{\nabla f(x)\}$.

**Definition N.19** (Strong convexity)**.** A function $f : \mathbb{R}^n \to \overline{\mathbb{R}}$ is $\mu$-strongly convex, for a given $\mu \in \mathbb{R}_{>0}$, over the convex set $\mathcal{C} \subseteq \mathbb{R}^n$ if and only if for all pairs $(x, y) \in \mathcal{C} \times \mathcal{C}$

$$f(x) \geq f(y) + \langle g, x - y \rangle + \frac{\mu}{2}\|x - y\|_2^2, \; \forall g \in \partial f(y),$$

or equivalently

$$\langle g_x - g_y, x - y \rangle \geq \mu\|x - y\|_2^2, \; \forall g_x \in \partial f(x), \; \forall g_y \in \partial f(y).$$

If the set $\mathcal{C}$ is not specified, then it is implicitly assumed that $\mathcal{C} = \mathrm{int}(\mathrm{dom}(f))$.

**Definition N.20** (Smoothness)**.** A function $f : \mathbb{R}^n \to \mathbb{R}$ is said to be $L$-smooth (otherwise referred to as Lipschitz) over $\mathcal{C} \subseteq \mathbb{R}^n$, for a given $L \in \mathbb{R}_{>0}$, if it is differentiable over $\mathcal{C}$ and satisfies

$$\|\nabla f(x) - \nabla f(y)\|_2 \leq L\|x - y\|_2, \; \forall x, y \in \mathcal{C}.$$

If the set $\mathcal{C}$ is not specified, then it is implicitly assumed that $\mathcal{C} = \mathbb{R}^n$.

**Descent lemma.** A well known result on $L$-smooth functions is the so called *descent lemma*, which states that, an $L$-smooth function $f : \mathbb{R}^n \to \mathbb{R}$ over $\mathcal{C} \subseteq \mathbb{R}^n$ satisfies,

$$f(y) \leq f(x) + \langle \nabla f(x), y - x \rangle + \frac{L}{2}\|x - y\|_2^2, \; \forall x, y \in \mathcal{C}.$$

**Definition N.21** ($\mathcal{K}_\infty$-class functions)**.** A function $\alpha : \mathbb{R}_{\geq 0} \to \mathbb{R}$ is said to be a $\mathcal{K}_\infty$-class function (denoted by $\alpha \in \mathcal{K}_\infty$), if it is continuous, strictly increasing, unbounded above and $\alpha(0) = 0$.

### Operators

**Definition N.22** (Indicator function)**.** The *indicator function* of a non-empty set $\mathcal{C} \subseteq \mathcal{D}$ is the function $\mathcal{I}_\mathcal{C} : \mathcal{D} \to \{0, +\infty\}$ where, given $x \in \mathcal{D}$, $\mathcal{I}_\mathcal{C}(x) = 0$ if $x \in \mathcal{C}$ and $\mathcal{I}_\mathcal{C}(x) = +\infty$ if $x \notin \mathcal{C}$.

**Definition N.23** (Projection onto a set)**.** Given a non-empty set $\mathcal{C} \subseteq \mathbb{R}^n$ and $R \in \mathbb{S}^n_{++}$, the operator $\pi^R_\mathcal{C} : \mathbb{R}^n \to \mathbb{R}^n$ given by

$$\pi^R_\mathcal{C}(x) = \arg\min_{z \in \mathcal{C}} \|z - x\|_R,$$

is the *R-weighted projection* of $x \in \mathbb{R}^n$ onto $\mathcal{C}$. If $R = I_n$, we simply write $\pi_\mathcal{C}(x)$, which is the standard Euclidean projection onto $\mathcal{C}$.

**Definition N.24** (Proximal operator)**.** Given a function $f : \mathbb{R}^n \to \overline{\mathbb{R}}$, its *proximal operator* is given by

$$\operatorname{prox}_f(x) = \arg\min_{z \in \mathbb{R}^n} f(z) + \frac{1}{2}\|z - x\|_2^2, \text{ for any } x \in \mathbb{R}^n.$$

**Notation N.25** (Ceiling and floor)**.** Given $a \in \mathbb{R}$, $\lceil a \rceil$ denotes the smallest integer greater than of equal to $a$, and $\lfloor a \rfloor$ the largest integer smaller than or equal to $a$.

### Well-known inequalities and laws

**Triangle inequality.** For an $x, y \in \mathbb{R}^n$, $\|x + y\| \leq \|x\| + \|y\|$.

**Cauchy-Schwarz inequality.** For any $x, y \in \mathbb{R}^n$, we have $|\langle x, y \rangle| \leq \|x\|_2 \|y\|_2$.

**Jensen's inequality.** For a convex function $f : \mathbb{R}^n \to \mathbb{R}$, vectors $x_j \in \mathbb{R}^n$ and scalars $a_j \in \mathbb{Z}_{>0}$ for $j \in \mathbb{Z}_1^m$, we have

$$f\left(\frac{\sum_{j=1}^m a_j x_j}{\sum_{j=1}^m a_j}\right) \leq \frac{\sum_{j=1}^m a_j f(x_j)}{\sum_{j=1}^m a_j}.$$

Moreover, if $f$ is concave, the inequality is reversed (from $\leq$ to $\geq$).

**Parallelogram law.** For any vector norm $\|\cdot\|$ on $\mathbb{R}^n$ and $x, y \in \mathbb{R}^n$, we have

$$\|x\|_2^2 + \|y\|_2^2 = \frac{1}{2}\|x + y\|_2^2 + \frac{1}{2}\|x - y\|_2^2.$$

# Chapter 1

# Introduction

## 1.1 Motivation and objectives

The results presented in this dissertation originate from our objective of implementing model predictive control (MPC) in *programmable logic controllers* (PLC), which are the most widely used embedded system in the industry for the implementation of low-level control loops. These devices are small ruggedized digital computers that have been adapted to the harsh industrial environments and to their use for controlling manufacturing and chemical processes. However, they are typically characterized for having very small memory and computational resources; mainly due to the fact that they are mostly used for the implementation of PID controllers and simple automatons.

The possibility of routinely implementing MPC in PLCs, given their current prevalence, could lead to important economical and competitive advantages, as well as help bridge the gap between the great interest that the academic control community has shown towards MPC in recent years with the reticence that many industries have towards changing their current control schemes, in spite of the potential advantages that MPC offers.

There are many aspects to tackle and decisions to make to achieve this goal, as well as many different avenues to do so. Indeed, the implementation of MPC in embedded systems is a widely researched topic, as will become apparent throughout this manuscript. In this dissertation, we focus on the implementation of linear MPC in embedded systems. We center our attention on two of the main aspects that one must consider: the MPC formulation to be implemented, and the optimization algorithm employed to solve its optimization problem.

There are other aspects to consider: practical aspects such as the fact that the system state might not be measured, the use of additional ingredients to guarantee an offset-free control, or time certification of the optimization algorithm. Although we have studied some of these practical issues in the context of the implementation of MPC in PLCs, our focus in this dissertation is on the development of solvers and MPC formulations particularly adapted to the limitations of embedded systems. This has led to the development of an open-source Matlab toolbox for the automatic code-generation of library-free sparse solvers for MPC. The toolbox is called *SPCIES: Suite of Predictive Controllers for Industrial Embedded Systems* [4] and is available at `https://github.com/GepocUS/Spcies`.

The solvers we develop are based on *first order methods*, since we found that the structures that arise in certain MPC formulations could be exploited with them. Even though these optimization methods typically require more iterations to converge to a close proximity of the optimal solution than other optimization algorithms, their use may be preferable due to the fact that each iteration requires less operations than the ones of other alternative iterative methods, particularly when dealing with large optimization problems.

Even though each iteration of first order methods may be computationally cheep, relatively speaking, a reduction of the number of iterations can have a significant impact on the applicability of the method, including its use as the embedded solver for MPC. Therefore, improvements of the convergence of these methods, both from a practical and theoretical standpoint, can be beneficial in many different areas. There are many ways in which this objective may be tackled: proper selection of the method's parameters, numerical preconditioning of the optimization problem, warmstarting procedures, etc. We have focused on improving the convergence of a subclass of first order methods known as *accelerated first order methods*, which have a quicker theoretical convergence than the non-accelerated variants and often perform better in practice. However, they may suffer from an oscillating behaviour that increases the number of iterations required to obtain a solution. This phenomenon can be addressed with the use of *restart schemes*, which restart the algorithm when certain conditions are met. We propose various restart schemes with theoretical properties.

Finally, in line with our study of MPC formulations suitable for their implementation in embedded systems, we developed a novel formulation that may be particularly suitable for its implementation in embedded systems due to its good performance with low prediction horizons, both in terms of its performance as well as its domain of attraction. This formulation, which we label *harmonic MPC*, tackles the problem of implementing MPC in embedded systems from another perspective: that of developing formulations that by their nature require fewer resources for their implementation.

Even though the initial focus was on the implementation of MPC in industrial embedded systems, the results shown in this dissertation will not focus on PLCs. Instead, we describe the solvers in general terms, focusing on showing how the structures are exploited and on comparing them with other alternatives from the literature. Therefore, we do not focus on the other practical aspects that must be addressed for the implementation of MPC in a real industrial setting. Additionally, the restart schemes for accelerated first order methods are presented from the point of view of the broader field of *convex optimization*. That is, they are not particularized to their use for MPC, although we show numerical results of their application to this problem.

## 1.2 Outline of the dissertation

This dissertation presents contributions in two separate (although closely related) fields: restart schemes for accelerated first order methods, and the implementation of MPC in embedded systems. Thus, we divide the manuscript into two parts.

**Part I** presents novel restart schemes for accelerated first order methods. These methods tend to exhibit a degradation of their convergence when close to the optimal solution due to the appearance of an oscillating behavior. Restart schemes aim to suppress this oscillatory behavior by stopping and restarting (thus their name) the algorithm when certain conditions are met. Previous restart schemes in the literature have been shown to provide a practical performance improvement. However, some do not prove linear convergence of the restarted algorithm or they require the knowledge of hard-to-obtain parameters of the function to be minimized. We present three novel restart schemes that are applicable to a wide range of (not necessarily strongly convex) optimization problems without requiring prior knowledge of its characterizing parameters. We present the iteration complexity of the proposed schemes and present numerical results comparing them to other restart schemes from the literature.

Part I is divided into the following chapters:

**Chapter 2** presents the preliminary concepts and first order methods that will be used and discussed throughout the rest of this dissertation. In particular, Section 2.1 introduces the first order methods and Section 2.2 the accelerated first order methods. The results presented in this chapter are well known, although we introduce some propositions that will be used in following chapters.

**Chapter 3** discusses restart schemes for accelerated first order methods. We start by providing a general explanation of restart schemes and motivating their usefulness. We then present a brief overview of the literature, explaining some of the most well known restart schemes in Section 3.1. Our three novel restart schemes are presented in Section 3.2. The first two are particularized to a certain accelerated first order method known as FISTA (which is essentially the accelerated proximal gradient method), whereas the third one considers a wider range of accelerated first order methods. We conclude the chapter with some numerical results comparing the proposed restart schemes with some of the ones presented in the literature in Section 3.3.

**Part II** deals with the implementation of linear MPC in embedded systems. The on-line implementation of MPC in embedded systems is challenging due to the low computational and memory resources that typically characterize them. We present an efficient approach for solving various MPC formulations using some of the first order methods presented in Part I. This approach allows us to take advantage of the inner structures of the matrices of the MPC's optimization problem, leading to a sparse optimization algorithm that does not require storing the sparse matrices using the typical sparse matrix representations (such as *compressed sparse column/row* or *dictionary of keys*). Instead, we only need to store the repeating submatrices, not needing to store the information of their location. The developed solvers have been made available in the open-source Matlab toolbox *SPCIES* [4]. The toolbox takes the model of the system and parameters of the MPC controller, and automatically generates the library-free code of the solver for its implementation in the embedded system of choice. Additionally, we present a novel MPC formulation, which we call *harmonic MPC* and label by *HMPC*, which is particularly suitable for its use with small prediction horizons, since it exhibits a larger domain of attraction that reduces the feasibility and performance issues that some MPC formulations display when the prediction horizon is not large enough. Due to its good performance with small prediction horizons, this MPC formulation may be suitable for its use in embedded systems.

Part II is divided into the following chapters

**Chapter 4** introduces the topic of this part of the dissertation. We start by describing the problem formulation in Section 4.1 and briefly describing MPC in Section 4.2. Finally, in Section 4.3 we present a non-exhaustive review of the literature on the topic of the implementation of MPC in embedded systems.

**Chapter 5** presents our approach to sparsely solving MPC formulations using first order methods. We start by presenting some rather straightforward sparse algorithms for solving specific QP problems and systems of equations in Section 5.1. Section 5.2 and Section 5.3 show how we solve QP problems using two of the first order methods presented in Part I of this dissertation. This approach to solving QP problems, along with the algorithms presented in Section 5.1, will be used in Sections 5.4 to 5.6 to develop sparse solvers for the MPC formulations presented in each one of these sections. Section 5.7 presents three systems used as test benches in the numerical results presented in Section 5.8, where we compare our proposed solvers between each other and against other alternatives from the literature. Additionally, we show the results of incorporating the restart schemes of Part I to some of the proposed solvers.

**Chapter 6** presents the *harmonic* MPC formulation, starting with its description in Section 6.1. Sections 6.2 and 6.3 prove the recursive feasibility and asymptotic stability of the proposed formulation, respectively. A discussion related to the design of one of its main ingredients is presented in Section 6.4. Finally, we showcase its potential advantages in Section 6.5, where we compare its performance and domain of attraction against other MPC formulations with a particular focus on using small prediction horizons.

## 1.3 Publications

The results shown in this dissertation are supported by several journal and congress publications, some of which are currently under review.
The novel restart schemes and proofs presented in Part I can be found in:

[5] T. Alamo, P. Krupa, and D. Limon, "Restart FISTA with global linear convergence," in *Proceedings of the European Control Conference (ECC)*, pp. 1969–1974, IEEE, 2019. Extended version available at arXiv:1906.09126.

[6] T. Alamo, P. Krupa, and D. Limon, "Gradient based restart FISTA," in *Proceedings of the Conference on Decision and Control (CDC)*, pp. 3936–3941, IEEE, 2019.

[7] T. Alamo, P. Krupa, and D. Limon, "Restart of accelerated first order methods with linear convergence for non-strongly convex optimization," *arXiv preprint: 2102.12387*, 2021. Submitted to *Transactions on Automatic Control.*

The sparse solvers for the implementation of MPC in embedded systems presented in Chapter 5, as well as some additional numerical results to the ones presented in this dissertation, can be found in:

[8] P. Krupa, D. Limon, and T. Alamo, "Implementation of model predictive control in programmable logic controllers," *IEEE Transactions on Control Systems Technology*, vol. 29, no. 3, pp. 1117–1130, 2021.

[9] P. Krupa, D. Limon, and T. Alamo, "Implementation of model predictive controllers in programmable logic controllers using IEC 61131-3 standard," in *Proceedings of the European Control Conference (ECC)*, pp. 288–293, IEEE, 2018.

[10] P. Krupa, I. Alvarado, D. Limon, and T. Alamo, "Implementation of model predictive control for tracking in embedded systems using a sparse extended ADMM algorithm," *arXiv preprint: 2008.09071*, 2020. Submitted to *Transactions on Control Systems Technology.*

[11] P. Krupa, J. Camara, I. Alvarado, D. Limon, and T. Alamo, "Real-time implementation of MPC for tracking in embedded systems: Application to a two-wheeled inverted pendulum," in *Proceedings of the European Control Conference (ECC)*, 2021. Preprint available at arXiv:2103.14571.

[12] P. Krupa, R. Jaouani, D. Limon, and T. Alamo, "A sparse ADMM-based solver for linear MPC subject to terminal quadratic constraint," *arXiv preprint: 2105.08419*, 2021. Submitted to *Automatica.*

The *harmonic* MPC formulation presented in Chapter 6 can be found in:

[13] P. Krupa, D. Limon, and T. Alamo, "Harmonic based model predictive control for set-point tracking," *IEEE Transactions on Automatic Control*, 2020.

[14] P. Krupa, M. Pereira, D. Limon, and T. Alamo, "Single harmonic based model predictive control for tracking," in *Proceedings of the Conference on Decision and Control (CDC)*, pp. 151–156, IEEE, 2019.

Other articles related to this dissertation but not directly discussed in it are:

[15] P. Krupa, N. Saraf, D. Limon, and A. Bemporad, "PLC implementation of a real-time embedded MPC algorithm based on linear input/output models," *IFAC-PapersOnLine*, vol. 53, no. 2, pp. 6987–6992, 2020.

[16] I. Alvarado, P. Krupa, D. Limon, and T. Alamo, "Tractable robust MPC design based on nominal predictions," *arXiv preprint: 2104.06088*, 2021. Submitted to *Journal of Process Control*.

[17] P. Krupa, C. Danielson, C. Laughman, S. A. Bortoff, D. J. Burns, S. Di Cairano, and D. Limon, "Modelica implementation of centralized MPC controller for a multi-zone heat pump," in *Proceedings of the European Control Conference (ECC)*, 2019.

# Part I

# Restart schemes for accelerated first order methods

# Chapter 2

# Preliminaries: Accelerated first order methods

First order methods (FOM) are iterative numerical methods for solving optimization problems that only require knowledge of the gradient/subgradient (but not Hessian) of the objective function [18, 19]. The origin of these methods can be traced back to the works of Cauchy in 1847 [20]. Since then, many different FOMs have been proposed and developed [18, 19, 21, 22].

These methods have been used in a wide variety of fields, including machine learning and model predictive control [18, 23]. The rising number of fields and applications that deal with large optimization problems has increased the interest in these methods due to their lower computational cost per iteration and memory requirements when compared to other optimization methods, such as *interior point* or *active set* methods.

As a gentle introduction to the topic, the most most well know FOM is the *gradient descent* method. Consider the optimization problem

$$\min_{z \in \mathbb{R}^{n_z}} f(z), \tag{2.1}$$

where the convex function $f : \mathbb{R}^{n_z} \to \mathbb{R}$ is Lipschitz continuously differentiable. Then, the gradient descent method applies, for $k \geq 0$, the recursion

$$z_{k+1} = z_k - \rho_k \nabla f(z_k)$$

starting at an initial point $z_0 \in \mathbb{R}^{n_z}$ and where $\rho_k \in \mathbb{R}_{>0}$ is the step size. If $\rho_k$ is chosen appropriately, the above recursion converges to an optimal solution $z^*$ of (2.1) as $k \to \infty$ [18].

In 1983, Yurii Nestevor proposed a variation of the gradient descent method known as the *accelerated* gradient descent method [24]. This new method provided a $O(1/k^2)$ convergence rate (in terms of objective function value), which is the optimal convergence rate for the class of optimization problems (2.1); a significant improvement over the $O(1/k)$ rate of the gradient descent method. The acceleration is obtained by providing a certain "momentum" to the algorithm in the form of an overrelaxation step. In particular, the iteration of the accelerated gradient descent method is, for $k \geq 0$,

$$\begin{aligned} z_{k+1} &= y_k - \rho_k \nabla f(z_k), \\ t_{k+1} &= \frac{1}{2}\left(1 + \sqrt{1 + 4t_k^2}\right), \\ y_{k+1} &= z_{k+1} + \frac{t_k - 1}{t_{k+1}}(z_{k+1} - z_k), \end{aligned}$$

starting at an initial point $z_0, y_0 \in \mathbb{R}^{n_z}$ and taking $t_0 = 1$.

The idea behind this method has been applied to many other FOMs, ranging from methods for smooth convex optimization [24, 25], to methods for composite non-smooth convex optimization [22, 26, 27, 28]. These methods are referred to as *accelerated*, of *fast*, first order methods (AFOM).

This chapter introduces a few FOMs and AFOMs. In particular, we present two popular FOMs: the *proximal gradient method* and the *Alternating Direction Method of Multipliers* (ADMM); and the accelerated variant of the proximal gradient method, which we label as FISTA. Additionally, we present a variation of ADMM known as *extended* ADMM (EADMM), and a variant of FISTA known as *monotone* FISTA (MFISTA). We show the pseudocode of each method and some of the their properties (in particular, those that will be of use for future developments throughout this dissertation). The properties we show in this chapter are for the most part well established in the literature. However, we include our own proofs for some of them. We do this to show that the properties we cite are still applicable to our particular problem formulation and/or assumptions, which are not always the same as the ones that can be found in the references provided.

We refer the reader to [18, 19, 22, 29, 30] for a few references in the field of first order methods, and to [1, 23] for general references on convex optimization.

## 2.1 First order methods

This section describes two popular FOMs: the *proximal gradient method* and the *alternating direction method of multipliers* (ADMM). Additionally, we describe an extension of ADMM to composite convex optimization problems with three functions in the objective function.

### 2.1.1 Proximal gradient method

The *proximal gradient method* is a FOM for convex composite optimization problems that is based on performing a *proximal operator* at each iteration of the algorithm [18, 29]. We follow, however, an equivalent description of the method, in which we make use of the so called *composite gradient mapping* [27, §2], which we formerly define further ahead. For an in depth explanation of the proximal gradient method and its different variants, we refer the reader to [18, §10.2] and, as a more accessible reference, to [29].

Consider the convex composite optimization problem

$$f^* = \min_{z \in \mathcal{Z}} \{ f(z) \doteq \Psi(z) + h(z) \} \tag{2.2}$$

under the following assumption.

**Assumption 2.1.** We assume that:

*(i)* The function $\Psi \in \Gamma(\mathbb{R}^{n_z})$ may be non-smooth.

*(ii)* $h : \mathbb{R}^{n_z} \to \mathbb{R}$ is a smooth differentiable convex function. That is, there exists $R \in \mathbb{S}_{++}^{n_z}$ such that

$$h(z) \leq h(y) + \langle \nabla h(y), z - y \rangle + \frac{1}{2} \|z - y\|_R^2, \; \forall (z, y) \in \mathbb{R}^{n_z} \times \mathbb{R}^{n_z}. \tag{2.3}$$

*(iii)* $\mathcal{Z} \subseteq \mathbb{R}^{n_z}$ is a non-empty closed convex set.

*(iv)* Problem (2.2) is solvable. That is, there exists $z^* \in \mathcal{Z} \cap \operatorname{dom}(\Psi)$ such that

$$f^* = f(z^*) = \inf_{z \in \mathcal{Z}} f(z).$$

We denote by $\Omega_f \doteq \{ z \in \mathbb{R}^{n_z} \,:\, z \in \mathcal{Z}, \; f(z) = f^* \}$ the *optimal set* of problem (2.2), which, due to Assumtion 2.1.*(iv)* is non-empty. It is well known that this set is a singleton if $f$ is strictly convex.

**Remark 2.2.** *It is standard to write the smoothness condition* (2.3) *as the well known descent lemma*

$$h(z) \leq h(y) + \langle \nabla h(y), z - y \rangle + \frac{L}{2} \|z - y\|_S^2, \; \forall (z, y) \in \mathbb{R}^{n_z} \times \mathbb{R}^{n_z}, \tag{2.4}$$

*where* $L \in \mathbb{R}_{>0}$ *provides a bound on the Lipschitz constant of* $\nabla h$ *[3, §2.1] and* $S \in \mathbb{S}_{++}^{n_z}$ *is usually taken as the identity matrix. However, since*

$$\frac{L}{2} \|z - y\|_S^2 = \frac{1}{2} \|z - y\|_{LS}^2, \; \forall (z, y) \in \mathbb{R}^{n_z} \times \mathbb{R}^{n_z},$$

*we have that* (2.3) *implies* (2.4) *if we take* $R = LS$. *We use expression* (2.3) *because it simplifies the algebraic expressions and results we present throughout this dissertation, but analogous results can be obtained if the standard L-smoothness definition (Definition N.20) and descent lemma are used instead of* (2.3).

Given $y \in \mathbb{R}^{n_z}$, one could use the local information given by $\nabla h(y)$ to minimize the value of $f$ around $y$. Under Assumtion 2.1, this can be done obtaining the optimal solution of the strictly convex optimization problem

$$\min_{z \in \mathcal{Z}} \; \Psi(z) + \langle \nabla h(y), z - y \rangle + \frac{1}{2} \|z - y\|_R^2.$$

The solution of this optimization problem leads to the notion of the *composite gradient mapping* [27, §2], which constitutes a generalization of the gradient mapping (see [3, §2.2.3] for the case when $\Psi(\cdot) = 0$ and [26] for the case when $\mathcal{Z} = \mathbb{R}^{n_z}$).

**Definition 2.3** (Composite gradient mapping)**.** Let $f \in \Gamma(\mathbb{R}^{n_z})$ be the composition of two functions $f = \Psi + h$, and let $\mathcal{Z} \subseteq \mathbb{R}^{n_z}$. Let Assumtion 2.1 hold. We denote the following mappings $\mathcal{T}_R^{f,\mathcal{Z}} : \mathbb{R}^{n_z} \to \mathbb{R}^{n_z}$ and $\mathcal{G}_R^{f,\mathcal{Z}} : \mathbb{R}^{n_z} \to \mathbb{R}^{n_z}$:

$$\mathcal{T}_R^{f,\mathcal{Z}}(y) = \arg\min_{z\in\mathcal{Z}} \Psi(z) + \langle \nabla h(y), z - y\rangle + \frac{1}{2}\|z - y\|_R^2, \tag{2.5a}$$

$$\mathcal{G}_R^{f,\mathcal{Z}}(y) = R\left(y - \mathcal{T}_R^{f,\mathcal{Z}}(y)\right). \tag{2.5b}$$

When the identities of $f$, $R$ and $\mathcal{Z}$ are clear from the context, we will often omit the superscripts and subscript and simply write $\mathcal{T}(y)$ and $\mathcal{G}(y)$.

The following proposition states the uniqueness of the composite gradient mapping under Assumtion 2.1. This result is well known in the literature (see, for instance, [18, §6.1] for an analogous result). However, we include a proof of the proposition for our particular notation and description of problem (2.2).

**Proposition 2.4** (On the unique solution of the composite gradient mapping)**.** Let $f \in \Gamma(\mathbb{R}^{n_z})$ be the composition of two functions $f = \Psi + h$, and let $\mathcal{Z} \subseteq \mathbb{R}^{n_z}$. Let Assumtion 2.1 hold. Then, $\mathcal{T}_R^{f,\mathcal{Z}}(y)$ is a singleton for any $y \in \mathbb{R}^{n_z}$.

***Proof:*** Let $\Psi_{\mathcal{Z}} = \Psi + \mathcal{I}_{\mathcal{Z}}$. Then, after some simple algebraic manipulations and cancellation of constant terms, we have that problem (2.5) can be equivalently written as

$$\mathcal{T}_R^{f,\mathcal{Z}}(y) = \arg\min_{z\in\mathbb{R}^{n_z}} \Psi_{\mathcal{Z}} + \frac{1}{2}\|z - \left(y - R^{-1}\nabla h(y)\right)\|_R^2. \tag{2.6}$$

From Assumtion 2.1.*(iv)*, we have that $\mathrm{dom}(\Psi_{\mathcal{Z}}) = \mathrm{dom}(\Psi) \cap \mathcal{Z} \neq \emptyset$. Moreover, since both $\Psi \in \Gamma(\mathbb{R}^{n_z})$ (Assumtion 2.1.*(i)*) and $\mathcal{I}_{\mathcal{Z}} \in \Gamma(\mathbb{R}^{n_z})$ (since Assumtion 2.1.*(iii)* states that $\mathcal{Z}$ is non-empty), we conclude that $\Psi_{\mathcal{Z}} \in \Gamma(\mathbb{R}^{n_z})$. Therefore, (2.6) is a strongly convex function that is not everywhere infinite. As such, it is well known that it has a unique minimizer for any $y \in \mathbb{R}^{n_z}$. ■

Equation (2.6) shows that the composite gradient mapping is closely related to the proximal operator. Indeed, if $R = L\mathbf{I}_{n_z}$ for some $L \in \mathbb{R}_{>0}$, then by definition of the proximal operator (Definition N.24), we have that (2.6) is equivalent to $\mathrm{prox}_{L^{-1}\Psi_{\mathcal{Z}}}(y - R^{-1}\nabla h(y))$.

In the context of optimal gradient methods, it is assumed that the computation of $\mathcal{T}(y)$ is cheap. This is the case when $\mathcal{Z}$ is a simple set (a box, $\mathbb{R}^{n_z}$, etc.), $R \in \mathbb{D}_{++}^{n_z}$, and $\Psi(\cdot)$ is a separable function. For example, in the well known Lasso optimization problem, the computation of $\mathcal{T}(y)$ resorts to the computation of the shrinkage operator [26]. See [18, §6], [29, §6],[31], or [32, §28] for numerous examples in which the computation of the composite gradient mapping is simple.

The following proposition gathers some well-known properties of the composite gradient mapping and the dual norm $\|\mathcal{G}(\cdot)\|_{R^{-1}}$. We include the proof for completeness.

**Proposition 2.5.** Consider problem (2.2) and let Assumtion 2.1 hold. Then:

*(i)* For every $y \in \mathbb{R}^n$ and $z \in \mathcal{Z}$,

$$f(\mathcal{T}(y)) - f(z) \leq \langle \mathcal{G}(y), \mathcal{T}(y) - z \rangle + \frac{1}{2}\|\mathcal{G}(y)\|^2_{R^{-1}} \tag{2.7a}$$

$$= \langle \mathcal{G}(y), y - z \rangle - \frac{1}{2}\|\mathcal{G}(y)\|^2_{R^{-1}} \tag{2.7b}$$

$$= -\frac{1}{2}\|\mathcal{T}(y) - z\|^2_R + \frac{1}{2}\|y - z\|^2_R. \tag{2.7c}$$

*(ii)* For every $y \in \mathcal{Z}$, $\frac{1}{2}\|\mathcal{G}(y)\|^2_{R^{-1}} \leq f(y) - f(\mathcal{T}(y)) \leq f(y) - f^*$.

***Proof:*** Let $h_y : \mathbb{R}^{n_z} \to \mathbb{R}$ be the given by $h_y(z) = \langle \nabla h(y), z - y \rangle + \frac{1}{2}\|z - y\|^2_R$, $\Psi_{\mathcal{Z}} : \mathbb{R}^{n_z} \to \overline{\mathbb{R}}$ be given by $\Psi_{\mathcal{Z}}(z) = \Psi(z) + \mathcal{I}_{\mathcal{Z}}(z)$, and $F_y : \mathbb{R}^{n_z} \to \overline{\mathbb{R}}$ be given by $F_y(z) = \Psi_{\mathcal{Z}}(z) + h_y(z)$. From Proposition 2.4, we have that, for any given $y \in \mathbb{R}^{n_z}$, $\mathcal{T}(y)$ is the unique minimizer of problem (2.5a). Therefore,

$$\Psi(\mathcal{T}(y)) + h_y(\mathcal{T}(y)) \leq \Psi(z) + h_y(z),\ \forall z \in \mathcal{Z},$$

which, since $\mathcal{T}(y) \in \mathcal{Z}$, implies

$$\Psi_{\mathcal{Z}}(\mathcal{T}(y)) + h_y(\mathcal{T}(y)) \leq \Psi_{\mathcal{Z}}(z) + h_y(z),\ \forall z \in \mathbb{R}^{n_z}.$$

From this inequality we have that $F_y(\mathcal{T}(y)) \leq F_y(z)$, $\forall z \in \mathbb{R}^{n_z}$, which due to the definition of the subdifferential (Definition N.18) implies

$$\mathbf{0}_{n_z} \in \partial F_y(\mathcal{T}(y)). \tag{2.8}$$

From Assumtion 2.1.*(iv)*, we have that $\mathrm{dom}(\Psi_{\mathcal{Z}}) = \mathrm{dom}(\Psi) \cap \mathcal{Z} \neq \emptyset$. Moreover, since both $\Psi \in \Gamma(\mathbb{R}^{n_z})$ (due to Assumtion 2.1.*(i)*) and $\mathcal{I}_{\mathcal{Z}} \in \Gamma(\mathbb{R}^{n_z})$ (since Assumtion 2.1.*(iii)* states that $\mathcal{Z}$ is non-empty), we conclude that $\Psi_{\mathcal{Z}} \in \Gamma(\mathbb{R}^{n_z})$. Additionally, since $h_y$ is a continuous real-valued function in $\mathbb{R}^{n_z}$, we have that it is closed (see [1, Proposition 1.1.3]) and that its domain is $\mathbb{R}^{n_z}$. Therefore, $F_y$ is the sum of two closed convex functions satisfying $\mathrm{ri}(\mathrm{dom}(\Psi_{\mathcal{Z}})) \cap \mathrm{ri}(\mathrm{dom}(h_y)) \neq \emptyset$, which along [1, Proposition 5.4.6], lets us conclude that

$$\partial F_y(\mathcal{T}(y)) = \partial \Psi_{\mathcal{Z}}(\mathcal{T}(y)) + \partial h_y(\mathcal{T}(y)),$$

where the subdifferential $\partial h_y$ of the differentiable function $h_y$ is given by

$$\partial h_y(\mathcal{T}(y) = \nabla h_y(\mathcal{T}(y)) = \nabla h(y) + R(\mathcal{T}(y) - y).$$

Thus, we obtain from (2.8) that

$$0 \in \partial\Psi_{\mathcal{Z}}(\mathcal{T}(y)) + \partial h_y(\mathcal{T}(y)) = \partial\Psi_{\mathcal{Z}}(\mathcal{T}(y)) + \nabla h(y) + R(\mathcal{T}(y) - y),$$

which, using the definition of $\mathcal{G}(y) \doteq R(y - \mathcal{T}(y))$, can be rewritten as

$$\mathcal{G}(y) - \nabla h(y) \in \partial\Psi_{\mathcal{Z}}(\mathcal{T}(y)).$$

Therefore, by the definition of $\partial\Psi_{\mathcal{Z}}$ (Definition N.18), we have that

$$\Psi_{\mathcal{Z}}(z) \geq \Psi_{\mathcal{Z}}(\mathcal{T}(y)) + \langle \mathcal{G}(y) - \nabla h(y), z - \mathcal{T}(y)\rangle, \; \forall z \in \mathbb{R}^{n_z}.$$

Since, $\mathcal{Z} \subseteq \mathbb{R}^{n_z}$, $\mathcal{T}(y) \in \mathcal{Z}$ and $\Psi_{\mathcal{Z}} = \Psi$ for every $z \in \mathcal{Z}$, the previous inequality implies

$$\Psi(z) \geq \Psi(\mathcal{T}(y)) + \langle \mathcal{G}(y) - \nabla h(y), z - \mathcal{T}(y)\rangle, \; \forall z \in \mathcal{Z}. \tag{2.9}$$

From the convexity of $h$, we have that $h(z) \geq h(y) + \langle \nabla h(y), z - y\rangle$, $\forall z \in \mathcal{Z}$. Adding this to (2.9) yields

$$\begin{aligned} f(z) &= \Psi(z) + h(z) \\ &\geq \Psi(\mathcal{T}(y)) + \langle \mathcal{G}(y) - \nabla h(y), z - \mathcal{T}(y)\rangle + h(y) + \langle \nabla h(y), z - y\rangle \\ &= \Psi(\mathcal{T}(y)) + \langle \mathcal{G}(y), z - \mathcal{T}(y)\rangle + h(y) + \langle \nabla h(y), \mathcal{T}(y) - y\rangle, \; \forall z \in \mathcal{Z}. \end{aligned} \tag{2.10}$$

From Assumtion 2.1.*(ii)* we have

$$\begin{aligned} h(y) &\geq h(\mathcal{T}(y)) - \langle \nabla h(y), \mathcal{T}(y) - y\rangle - \frac{1}{2}\|\mathcal{T}(y) - y\|_R^2 \\ &= h(\mathcal{T}(y)) - \langle \nabla h(y), \mathcal{T}(y) - y\rangle - \frac{1}{2}\|R^{-1}\mathcal{G}(y)\|_R^2 \\ &= h(\mathcal{T}(y)) - \langle \nabla h(y), \mathcal{T}(y) - y\rangle - \frac{1}{2}\|\mathcal{G}(y)\|_{R^{-1}}^2. \end{aligned}$$

Adding this inequality to (2.10) yields

$$\begin{aligned} f(z) &\geq \Psi(\mathcal{T}(y)) + h(\mathcal{T}(y)) + \langle \mathcal{G}(y), z - \mathcal{T}(y)\rangle - \frac{1}{2}\|\mathcal{G}(y)\|_{R^{-1}}^2 \\ &= f(\mathcal{T}(y)) + \langle \mathcal{G}(y), z - \mathcal{T}(y)\rangle - \frac{1}{2}\|\mathcal{G}(y)\|_{R^{-1}}^2, \; \forall z \in \mathcal{Z}, \end{aligned}$$

which, if rearranged, proves (2.7a). We now prove (2.7b) and (2.7c) by means of simple algebraic manipulations.

$$\begin{aligned} f(\mathcal{T}(y)) - f(z) &\leq \langle \mathcal{G}(y), \mathcal{T}(y) - z\rangle + \frac{1}{2}\|\mathcal{G}(y)\|_{R^{-1}}^2 \\ &= \langle \mathcal{G}(y), y - z + \mathcal{T}(y) - y\rangle + \frac{1}{2}\|\mathcal{G}(y)\|_{R^{-1}}^2 \\ &= \langle \mathcal{G}(y), y - z\rangle + \langle \mathcal{G}(y), \mathcal{T}(y) - y\rangle + \frac{1}{2}\|\mathcal{G}(y)\|_{R^{-1}}^2 \\ &= \langle \mathcal{G}(y), y - z\rangle - \langle \mathcal{G}(y), R^{-1}\mathcal{G}(y)\rangle + \frac{1}{2}\|\mathcal{G}(y)\|_{R^{-1}}^2 \\ &= \langle \mathcal{G}(y), y - z\rangle - \|\mathcal{G}(y)\|_{R^{-1}}^2 + \frac{1}{2}\|\mathcal{G}(y)\|_{R^{-1}}^2 \\ &= \langle \mathcal{G}(y), y - z\rangle - \frac{1}{2}\|\mathcal{G}(y)\|_{R^{-1}}^2, \; \forall z \in \mathcal{Z}. \end{aligned} \tag{2.11}$$

---

**Algorithm 1:** Proximal gradient method

---

**Require:** $z_0 \in \mathbb{R}^n$, $\epsilon \in \mathbb{R}_{>0}$

$k \leftarrow 0$

**repeat**

  $k \leftarrow k + 1$

  $z_k \leftarrow \mathcal{T}_R^{f,\mathcal{Z}}(z_{k-1})$

**until** $\|\mathcal{G}_R^{f,\mathcal{Z}}(z_k)\|_{R^{-1}} \leq \epsilon$

**Output:** $\tilde{z}^* \leftarrow z_k$

---

This proves (2.7b). Next, from the definition of $\mathcal{G}(y)$, we have

$$\begin{aligned} f(\mathcal{T}(y)) - f(z) &\leq \langle R(y - \mathcal{T}(y)), y - z\rangle - \frac{1}{2}\|R(y - \mathcal{T}(y))\|^2_{R^{-1}} \\ &= -\langle R(y - \mathcal{T}(y)), z - y\rangle - \frac{1}{2}\|y - \mathcal{T}(y)\|^2_R \\ &= -\frac{1}{2}\|y - \mathcal{T}(y) + z - y\|^2_R + \frac{1}{2}\|z - y\|^2_R \\ &= -\frac{1}{2}\|\mathcal{T}(y) - z\|^2_R + \frac{1}{2}\|y - z\|^2_R, \ \forall z \in \mathcal{Z}, \end{aligned}$$

which proves (2.7c) and concludes the proof of claim *(i)*. Assume now that $y \in \mathcal{Z}$. Particularizing (2.11) to $z = y$ yields

$$\frac{1}{2}\|\mathcal{G}(y)\|^2_{R^{-1}} \leq f(y) - f(\mathcal{T}(y)), \ \forall y \in \mathcal{Z},$$

which, along $f^* \leq f(\mathcal{T}(y))$, proves claim *(ii)*. ■

The proximal gradient method solves problem (2.2) by iteratively applying the operator $\mathcal{T}(\cdot)$ at each iteration. Algorithm 1 shows the pseudocode of this method for a given exit tolerance $\epsilon \in \mathbb{R}_{>0}$ and an initial condition $z_0 \in \mathbb{R}^{n_z}$. It returns a suboptimal solution $\tilde{z}^*$ of problem (2.2), where the suboptimality is determined by $\epsilon$. The fact that the exit condition $\|\mathcal{G}(z_k)\|_{R^{-1}} \leq \epsilon$ serves as a measure of optimality of $z_k$ is provided by the following proposition. Once again, this proposition is well known in the literature (see, for instance, [18, Corollary 10.8]), but we include a proof for completeness.

**Proposition 2.6.** Consider problem (2.2). Let Assumtion 2.1 hold and $\mathcal{G}_R^{f,\mathcal{Z}}$ be given by Definition 2.3. Then $z \in \mathbb{R}^{n_z}$ belongs to the optimal set

$$\Omega_f \doteq \{\, z \in \mathbb{R}^{n_z} \ : \ z \in \mathcal{Z}, \ f(z) = f^* \,\},$$

if and only if $\mathcal{G}_R^{f,\mathcal{Z}}(z) = \mathbf{0}_{n_z}$.

***Proof:*** We first show the implication $\mathcal{G}(z) = \mathbf{0}_{n_z} \implies z \in \Omega_f$. Since $\mathbb{R} \in \mathbb{S}^{n_z}_{++}$, we infer from $\mathcal{G}(z) \doteq R(z - \mathcal{T}(z))$ that $\mathcal{G}(z) = \mathbf{0}_{n_z} \implies z = \mathcal{T}(z)$. Let $z^* \in \Omega_f \subseteq \mathcal{Z}$.

**Algorithm 2:** Alternating direction method of multipliers

**Require:** $v_0 \in \mathbb{R}^{n_v}$, $\lambda_0 \in \mathbb{R}^p$, $\rho \in \mathbb{R}_{>0}$, $\epsilon_p \in \mathbb{R}_{>0}$, $\epsilon_d \in \mathbb{R}_{>0}$

$k \leftarrow 0$

**repeat**

  $z_{k+1} \leftarrow \arg\min\limits_{z} \mathcal{L}_\rho(z, v_k, \lambda_k)$

  $v_{k+1} \leftarrow \arg\min\limits_{v} \mathcal{L}_\rho(z_{k+1}, v, \lambda_k)$

  $\lambda_{k+1} \leftarrow \lambda_k + \rho(Cz_{k+1} + Dv_{k+1} - d)$

  $k \leftarrow k + 1$

**until** $r_p \leq \epsilon_p$ **and** $r_d \leq \epsilon_d$

**Output:** $\tilde{z}^* \leftarrow z_k$, $\tilde{v}^* \leftarrow v_k$, $\tilde{\lambda}^* \leftarrow \lambda_k$

Then, from the implication $\mathcal{G}(z) = \mathbf{0}_{n_z} \implies z = \mathcal{T}(z)$ and Proposition 2.5.*(i)* we have that

$$f^* = f(z^*) \geq f(\mathcal{T}(z)) - \langle \mathcal{G}(z), \mathcal{T}(z) - z^* \rangle - \frac{1}{2}\|\mathcal{G}(z)\|^2_{R^{-1}} = f(\mathcal{T}(z)) = f(z).$$

Therefore, $z \in \Omega_f$. Next, we show the implication $z \in \Omega_f \implies \mathcal{G}(z) = \mathbf{0}_{n_z}$. In this case, we have that $f(z) = f^*$, which along Proposition 2.5.*(ii)* leads to $\frac{1}{2}\|\mathcal{G}(z)\|^2_{R^{-1}} \leq f(z) - f^* = 0$. Thus, $\mathcal{G}(z) = \mathbf{0}_{n_z}$. ■

**Remark 2.7.** *We note that it is common to substitute the exit condition of Algorithm 1 (step 5) for* $\|\mathcal{G}(z_{k-1})\|_{R^{-1}} \leq \epsilon$*, since* $\mathcal{G}(z_{k-1})$ *requires* $\mathcal{T}(z_{k-1})$*, which has already been computed in step 4.*

### 2.1.2 Alternating direction method of multipliers

The *alternating direction method of multipliers* (ADMM) [21] is a FOM for solving convex composite optimization problems of the form

$$\min_{z,v} \; f(z) + g(v) \tag{2.12a}$$

$$s.t. \; Cz + Dv = d, \tag{2.12b}$$

where $z \in \mathbb{R}^{n_z}$, $v \in \mathbb{R}^{n_v}$, $f \in \Gamma(\mathbb{R}^{n_z})$, $g \in \Gamma(\mathbb{R}^{n_v})$, $C \in \mathbb{R}^{p \times n_z}$, $D \in \mathbb{R}^{p \times n_v}$ and $d \in \mathbb{R}^p$. Let the *augmented* Lagrangian function $\mathcal{L}_\rho : \mathbb{R}^{n_z} \times \mathbb{R}^{n_v} \times \mathbb{R}^p \to \overline{\mathbb{R}}$ be given by

$$\mathcal{L}_\rho(z, v, \lambda) = f(z) + g(v) + \frac{\rho}{2}\|Cz + Dv - d + \frac{1}{\rho}\lambda\|^2_2,$$

where $\lambda \in \mathbb{R}^p$ are the dual variables and $\rho \in \mathbb{R}_{>0}$ is the *penalty parameter*. We denote a solution point of (2.12) by $(z^*, v^*, \lambda^*)$, assuming that one exists.

Algorithm 2 shows the ADMM algorithm applied to problem (2.12) for the given exit tolerances $(\epsilon_p, \epsilon_d) \in \mathbb{R}_{>0} \times \mathbb{R}_{>0}$ and initial point $(v_0, \lambda_0) \in \mathbb{R}^{n_v} \times \mathbb{R}^p$.

The exit condition of the algorithm is determined by the primal ($r_p \in \mathbb{R}_{\geq 0}$) and dual ($r_d \in \mathbb{R}_{\geq 0}$) residuals given by

$$r_p = \|Cz_k + Dv_k - d\|_\infty,$$
$$r_d = \|v_k - v_{k-1}\|_\infty.$$

The use of these residuals, as well as the justification for using $r_d$ as a residual for dual optimality, can be found in [21, §3.3]. The algorithm returns a suboptimal solution point $(\tilde{z}^*, \tilde{v}^*, \tilde{\lambda}^*)$ of (2.12), where the subopimality is determined by the values of $\epsilon_p$ and $\epsilon_d$.

There are numerous results on the convergence of ADMM applied to (2.12) depending on the assumptions made on its ingredients. For instance, we refer the reader to [21], [18, Theorem 15.4], [33] or [34]. Additionally, it has been shown that the use of different values of $\rho$ for each constraint can improve the convergence of the algorithm in practical settings [35, §5.2].

### 2.1.3 The extended alternating direction method of multipliers

This section introduces the *extended* ADMM (EADMM) algorithm [36], which, as its name suggests, is an extension of the ADMM algorithm [21] (see Section 2.1.2 and Algorithm 2) to optimization problems with more than two separable functions in the objective function. In particular, we focus on the case in which the objective function has three separable functions.

Let $\theta_i : \mathbb{R}^{n_i} \to \mathbb{R}$ for $i \in \mathbb{Z}_1^3$ be convex functions, $\mathcal{Z}_i \subseteq \mathbb{R}^{n_i}$ for $i \in \mathbb{Z}_1^3$ be closed convex sets, $C_i \in \mathbb{R}^{m_z \times n_i}$ for $i \in \mathbb{Z}_1^3$ and $d \in \mathbb{R}^{m_z}$. Consider the problem

$$\min_{z_1, z_2, z_3} \sum_{i=1}^{3} \theta_i(z_i) \tag{2.13a}$$
$$s.t. \sum_{i=1}^{3} C_i z_i = d \tag{2.13b}$$
$$z_i \in \mathcal{Z}_i, \ i \in \mathbb{Z}_1^3. \tag{2.13c}$$

and let its augmented Lagrangian $\mathcal{L}_\rho : \mathbb{R}^{n_1} \times \mathbb{R}^{n_2} \times \mathbb{R}^{n_3} \times \mathbb{R}^{m_z} \to \mathbb{R}$ be given by

$$\mathcal{L}_\rho(z_1, z_2, z_3, \lambda) = \sum_{i=1}^{3} \theta_i(z_i) + \langle \lambda, \sum_{i=1}^{3} C_i z_i - d \rangle + \frac{\rho}{2} \left\| \sum_{i=1}^{3} C_i z_i - d \right\|_2^2, \tag{2.14}$$

where $\lambda \in \mathbb{R}^{m_z}$ are the dual variables and $\rho \in \mathbb{R}_{>0}$ is the penalty parameter. We denote a solution point of (2.13) by $(z_1^*, z_2^*, z_3^*, \lambda^*)$, assuming that one exists.

Algorithm 3 shows the implementation of the EADMM algorithm for a given exit tolerance $\epsilon > 0$ and initial points $(z_2^0, z_3^0, \lambda^0)$. Algorithm 3 returns a suboptimal solution $(\tilde{z}_1^*, \tilde{z}_2^*, \tilde{z}_3^*, \tilde{\lambda}^*)$ of problem (2.13), where the suboptimality is determined by primal and dual exit tolerances $\epsilon_p \in \mathbb{R}_{>0}$ and $\epsilon_d \in \mathbb{R}_{>0}$, respectively.

**Algorithm 3:** Extended ADMM

**Require:** $z_2^0 \in \mathbb{R}^{n_2}$, $z_3^0 \in \mathbb{R}^{n_3}$, $\lambda^0 \in \mathbb{R}^{m_z}$, $\rho \in \mathbb{R}_{>0}$, $\epsilon \in \mathbb{R}_{>0}$

1 $k \leftarrow 0$

2 **repeat**

3 $\quad z_1^{k+1} \leftarrow \arg\min_{z_1} \left\{ \mathcal{L}_\rho(z_1, z_2^k, z_3^k, \lambda^k) \mid z_1 \in \mathcal{Z}_1 \right\}$

4 $\quad z_2^{k+1} \leftarrow \arg\min_{z_2} \left\{ \mathcal{L}_\rho(z_1^{k+1}, z_2, z_3^k, \lambda^k) \mid z_2 \in \mathcal{Z}_2 \right\}$

5 $\quad z_3^{k+1} \leftarrow \arg\min_{z_3} \left\{ \mathcal{L}_\rho(z_1^{k+1}, z_2^{k+1}, z_3, \lambda^k) \mid z_3 \in \mathcal{Z}_3 \right\}$

6 $\quad r \leftarrow \sum_{i=1}^{3} C_i z_i^{k+1} - d$

7 $\quad \lambda^{k+1} \leftarrow \lambda^k + \rho r$

8 $\quad k \leftarrow k + 1$

9 **until** $\|r\|_\infty \leq \epsilon_p$ **and** $\|z_2 - z_2^{k-1}\|_\infty \leq \epsilon_d$ **and** $\|z_3^k - z_3^{k-1}\|_\infty \leq \epsilon_d$

**Output:** $\tilde{z}_1^* \leftarrow z_1^k$, $\tilde{z}_2^* \leftarrow z_2^k$, $\tilde{z}_3^* \leftarrow z_3^k$, $\tilde{\lambda}^* \leftarrow \lambda^k$

The fact that the exit condition shown in step 9 of Algorithm 3 serves as an indicator of the (sub-)optimality of the current iterate can be found in [36, §5].

As shown in [37], the EADMM algorithm is not necessarily convergent under the typical assumptions of the classical ADMM algorithm. However, multiple results have shown its convergence under additional assumptions [36, 38, 39] or by adding additional steps [40, 41]. In particular, [36] proved its convergence under the following assumption, as stated in the following theorem.

**Assumption 2.8** ([36], Assumption 3.1)**.** The functions $\theta_1$ and $\theta_2$ are convex, function $\theta_3$ is $\mu_3$-strongly convex, and $C_1$ and $C_2$ are full column rank.

**Theorem 2.9** (Convergence of EADMM; [36], Theorem 3.1)**.** Suppose that Assumtion 2.8 holds and that the penalty parameter $\rho \in \left(0, \frac{6\mu_3}{17\|C_3^\top C_3\|}\right)$. Then, the sequence of points $(z_1^k, z_2^k, z_3^k, \lambda^k)$ generated by Algorithm 3 converges to a point in the optimal set of problem (2.13) as $k \to \infty$.

## 2.2 Accelerated first order methods

This section describes two accelerated variants of the proximal gradient method: a non-monotone variant, and a monotone variant.

### 2.2.1 Fast proximal gradient method (FISTA)

The fast proximal gradient method is the accelerated variant of the proximal gradient method (see Section 2.1.1), which constitutes an extension of the fast gradient method [24] to composite convex optimization problems given by (2.2)

**Algorithm 4:** FISTA

**Require:** $z \in \mathbb{R}^n$, $\epsilon \in \mathbb{R}_{>0}$

$y_0 \leftarrow \mathcal{T}(z)$, $z_0 \leftarrow \mathcal{T}(z)$, $t_0 \leftarrow 1$, $k \leftarrow 0$

**repeat**

  $k \leftarrow k + 1$

  $z_k \leftarrow \mathcal{T}(y_{k-1})$

  $t_k \leftarrow \frac{1}{2}\left(1 + \sqrt{1 + 4t_{k-1}^2}\right)$

  $y_k \leftarrow z_k + \frac{t_{k-1} - 1}{t_k}(z_k - z_{k-1})$

**until** $\|\mathcal{G}(z_k)\|_{R^{-1}} \leq \epsilon$

**Output:** $\tilde{z}^* \leftarrow z_k$

under Assumtion 2.1. Its interest stems from its improved convergence rate with respect to the non-accelerated variant, achieving a $O(1/k^2)$ rate of convergence in terms of objective function values, whereas the non-accelerated variant achieves a rate of $O(1/k)$ [26].

We will refer to the fast proximal gradient method using the acronym FISTA (Fast Iterative Shrinking-Threshold Algorithm), which was the term coined by Beck and Teboulle in 2009 [26].

Algorithm 4 shows the FISTA algorithm applied to (2.2) under Assumtion 2.1, starting from an initial point $z \in \mathbb{R}^{n_z}$ and for a given exit tolerance $\epsilon \in \mathbb{R}_{>0}$. It returns a suboptimal solution $\tilde{z}^*$ of problem (2.2), where the suboptimality is determined by $\epsilon$. Note that Algorithm 4 has a similar complexity to the proximal gradient method (Algorithm 1), since they both apply the operator $\mathcal{T}_R^{f,\mathcal{Z}}$ at each iteration and the others steps of FISTA are simple operations.

**Remark 2.10.** *We note that Remark 2.7 also applies to the exit condition of Algorithm 4, where, in this case, the exit condition $\|\mathcal{G}(z_k)\|_{R^{-1}} \leq \epsilon$ can be substituted by $\|\mathcal{G}(y_{k-1})\|_{R^{-1}} \leq \epsilon$.*

The following proposition shows two well known properties satisfied by the iterates of FISTA (see [24], [26, §10.7]) that will be of importance for future developments and results presented in this dissertation. The proof makes use of the following lemma, which shows two important facts about the sequence $\{t_k\}_{\geq 0}$ generated by Algorithm 4. We include the proofs of the lemma and proposition for completeness.

**Lemma 2.11.** Let the sequence $\{t_k\}_{\geq 0}$ be defined by $t_0 = 1$ and

$$t_k = \frac{1}{2}\left(1 + \sqrt{1 + 4t_{k-1}^2}\right), \ k \in \mathbb{Z}_{>0}.$$

Then, *(i)* $t_{k-1}^2 = t_k^2 - t_k$, $\forall k \geq 1$, and *(ii)* $t_k \geq \frac{k+2}{2}$, $\forall k \geq 0$.

***Proof:*** Claim *(i)* is trivial, since $t_k$ is defined as the roots of $t_k^2 - t_k - t_{k-1}^2 = 0$ for all $k \geq 1$. Claim *(ii)* is proven by induction on $k$. The claim is trivially satisfied for $k = 0$. Suppose that the claim is satisfied for $k - 1$, i.e., $t_{k-1} \geq \frac{k+2}{2}$. Then,

$$t_k = \frac{1}{2}\left(1 + \sqrt{1 + 4t_{k-1}^2}\right) \geq \frac{1}{2}\left(1 + \sqrt{4t_{k-1}^2}\right) = \frac{1}{2} + t_{k-1} \geq \frac{1}{2} + \frac{k+1}{2} = \frac{k+2}{2},$$

which shows that it is also satisfied for $k$. ■

**Proposition 2.12.** Consider problem (2.2) and let Assumtion 2.1 hold. Let $\{z_k\}_{\geq 0}$ and $\{y_k\}_{\geq 0}$ be the sequences generated by FISTA (Algorithm 4) for problem (2.2) starting at $z_0 = y_0 = \mathcal{T}_R^{f,\mathcal{Z}}(z)$, where $z \in \mathbb{R}^{n_z}$ is given. Then:

*(i)* $f(z_k) - f^* \leq \dfrac{2}{(k+1)^2}\|z_0 - \pi_{\Omega_f}^R(z_0)\|_R^2$, $\forall k \geq 1$.

*(ii)* $\|\mathcal{G}_R^{f,\mathcal{Z}}(y_k)\|_{R^{-1}} \leq \dfrac{4}{k+2}\|z_0 - \pi_{\Omega_f}^R(z_0)\|_R$, $\forall k \geq 0$.

***Proof:*** For the sake of simplicity, we use the shorthand notations $\bar{z}_0 \doteq \pi_{\Omega_f}^R(z_0)$, $\mathcal{G}_k \doteq \mathcal{G}_R^{f,\mathcal{Z}}(y_k)$, $y_k^+ \doteq \mathcal{T}_R^{f,\mathcal{Z}}(y_k)$, $\delta f_k \doteq f(z_k) - f^*$, $\delta z_k \doteq z_k - \bar{z}_0$ and $\delta y_k \doteq y_k - \bar{z}_0$, where the subscripts may change, e.g., $y_{k+1}^+ \doteq \mathcal{T}_R^{f,\mathcal{Z}}(y_{k+1})$.

We start by proving claim *(i)*. Particularizing inequality (2.7c) from Proposition 2.5.*(i)* to $y_0 \in \mathcal{Z} \subseteq \mathbb{R}^{n_z}$ and $\bar{z}_0 \in \Omega_f \subseteq \mathcal{Z}$, we obtain

$$f(y_0^+) - f(\bar{z}_0) \leq -\frac{1}{2}\|y_0^+ - \bar{z}_0\|_R^2 + \frac{1}{2}\|y_0 - \bar{z}_0\|_R^2.$$

From steps 1 and 4 of the algorithm we have $z_0 = y_0$ and $z_1 = y_0^+$. Furthermore, by definition of $\bar{z}_0$ we have $f(\bar{z}_0) = f^*$. Therefore we can rewrite the previous inequality as

$$f(z_1) - f^* \leq -\frac{1}{2}\|z_1 - \bar{z}_0\|_R^2 + \frac{1}{2}\|z_0 - \bar{z}_0\|_R^2 \leq \frac{1}{2}\|z_0 - \bar{z}_0\|_R^2, \tag{2.15}$$

which proves claim *(i)* for $k = 1$. We now proceed to prove the claim for $k \geq 2$. From step 4 of Algorithm 4 we have

$$z_k = y_{k-1}^+, \ \forall k \geq 1. \tag{2.16}$$

Therefore, inequality (2.7b) from Proposition 2.5.*(i)* leads to

$$f(z) \geq f(z_{k+1}) + \frac{1}{2}\|\mathcal{G}_k\|_{R^{-1}}^2 - \langle \mathcal{G}_k, y_k - z\rangle, \ \forall z \in \mathcal{Z}, \ \forall k \geq 1.$$

We notice that, by construction, $z_k \in \mathcal{Z}$, $\forall k \geq 1$. Thus, particularizing the previous inequality to $z = z_k$ and to $z = \bar{z}_0 \in \Omega_f \subseteq \mathcal{Z}$, we obtain the inequalities:

$$f(z_k) \geq f(z_{k+1}) + \frac{1}{2}\|\mathcal{G}_k\|_{R^{-1}}^2 - \langle \mathcal{G}_k, y_k - z_k\rangle, \ \forall k \geq 1,$$

$$f(\bar{z}_0) \geq f(z_{k+1}) + \frac{1}{2}\|\mathcal{G}_k\|_{R^{-1}}^2 - \langle \mathcal{G}_k, y_k - \bar{z}_0\rangle, \ \forall k \geq 1,$$

which can be rewritten as

$$\delta f_k - \delta f_{k+1} \geq \frac{1}{2}\|\mathcal{G}_k\|^2_{R^{-1}} - \langle \mathcal{G}_k, \delta y_k - \delta z_k \rangle, \ \forall k \geq 1, \tag{2.17a}$$

$$-\delta f_{k+1} \geq \frac{1}{2}\|\mathcal{G}_k\|^2_{R^{-1}} - \langle \mathcal{G}_k, \delta y_k \rangle, \ \forall k \geq 1. \tag{2.17b}$$

We introduce now the auxiliary nomenclature $\Gamma_k$ defined as

$$\Gamma_k \doteq t^2_{k-1}\delta f_k - t^2_k \delta f_{k+1}, \ \forall k \geq 1,$$

which, making use of Lemma 2.11.*(i)* allows us to write

$$\Gamma_k = (t^2_k - t_k)\delta f_k - t^2_k \delta f_{k+1} = (t^2_k - t_k)(\delta f_k - \delta f_{k+1}) - t_k \delta f_{k+1}, \ \forall k \geq 1. \tag{2.18}$$

In view of Lemma 2.11.*(ii)*, $t_k \geq 1$, $\forall k \geq 0$. This implies that we can replace the terms $\delta f_k - \delta f_{k+1}$ and $-\delta f_{k+1}$ in inequality (2.18) by the lower bounds given in inequalities (2.17). Doing so, we obtain:

$$\begin{aligned}\Gamma_k &\geq (t^2_k - t_k)\left(\frac{1}{2}\|\mathcal{G}_k\|^2_{R^{-1}} - \langle \mathcal{G}_k, \delta y_k - \delta z_k \rangle\right) + t_k \left(\frac{1}{2}\|\mathcal{G}_k\|^2_{R^{-1}} - \langle \mathcal{G}_k, \delta y_k \rangle\right) \\ &= \frac{t^2_k}{2}\|\mathcal{G}_k\|^2_{R^{-1}} - \langle \mathcal{G}_k, t^2_k(\delta y_k - \delta z_k) + t_k \delta z_k \rangle, \ \forall k \geq 1. \end{aligned} \tag{2.19}$$

Step 6 of the algorithms can be rewritten as

$$\delta y_k - \delta z_k = \frac{t_{k-1} - 1}{t_k}(\delta z_k - \delta z_{k-1}), \ \forall k \geq 1,$$

which, introducing the notation $s_k \doteq \delta z_{k-1} + t_{k-1}(\delta z_k - \delta z_{k-1})$, leads to

$$\begin{aligned} s_k - \delta z_k &= \delta z_{k-1} + t_{k-1}(\delta z_k - \delta z_{k-1}) - \delta z_k = (t_{k-1} - 1)(\delta z_k - \delta z_{k-1}) \\ &= t_k(\delta y_k - \delta z_k), \ \forall k \geq 1. \end{aligned} \tag{2.20}$$

Using the definition of $s_k$ along with (2.20), we show that $\mathcal{G}_k$ can be written in terms of $s_k$ and $s_{k+1}$:

$$\begin{aligned} t_k \mathcal{G}_k &\overset{(*)}{=} t_k R(y_k - z_{k+1}) = t_k R(\delta y_k - \delta z_{k+1}) \\ &= t_k R(\delta y_k - \delta z_k + \delta z_k - \delta z_{k+1}) \\ &= R(s_k - \delta z_k + t_k(\delta z_k - \delta z_{k+1})) \\ &= R(s_k - s_{k+1}), \ \forall k \geq 1, \end{aligned} \tag{2.21}$$

where $(*)$ is due to (2.16). We can now keep developing (2.19) as follows

$$\begin{aligned} \Gamma_k &\geq \frac{1}{2}\|t_k \mathcal{G}_k\|^2_{R^{-1}} - \langle \mathcal{G}_k, t_k(s_k - \delta z_k) + t_k \delta z_k \rangle = \frac{1}{2}\|t_k \mathcal{G}_k\|^2_{R^{-1}} - \langle t_k \mathcal{G}_k, s_k \rangle \\ &\overset{(2.21)}{\geq} \frac{1}{2}\|R(s_k - s_{k+1})\|^2_{R^{-1}} - \langle R(s_k - s_{k+1}), s_k \rangle \\ &= \frac{1}{2}\|s_{k+1} - s_k\|^2_R + \langle R(s_{k+1} - s_k), s_k \rangle \\ &= \frac{1}{2}\|(s_{k+1} - s_k) + s_k\|^2_R - \frac{1}{2}\|s_k\|^2_R = \frac{1}{2}\|s_{k+1}\|^2_R - \frac{1}{2}\|s_k\|^2_R, \ \forall k \geq 1. \end{aligned}$$

Thus, for every $k \geq 1$,

$$\Gamma_k = t_{k-1}^2 \delta f_k - t_k^2 \delta f_{k+1} \geq \frac{1}{2}\|s_{k+1}\|_R^2 - \frac{1}{2}\|s_k\|_R^2.$$

Equivalently,

$$t_k^2 \delta f_{k+1} + \frac{1}{2}\|s_{k+1}\|_R^2 \leq t_{k-1}^2 \delta f_k + \frac{1}{2}\|s_k\|_R^2, \ \forall k \geq 1.$$

Since this inequality holds $\forall k \geq 1$, we can apply it in a recursive way to obtain

$$\begin{aligned} t_k^2 \delta f_{k+1} + \frac{1}{2}\|s_{k+1}\|_R^2 &\leq t_0^2 \delta f_1 + \frac{1}{2}\|s_1\|_R^2 \\ &= \delta f_1 + \frac{1}{2}\|\delta z_0 + t_0(\delta z_1 - \delta z_0)\|_R^2 \\ &= \delta f_1 + \frac{1}{2}\|z_1 - \bar{z}_0\|_R^2 \overset{(2.15)}{\leq} \frac{1}{2}\|z_0 - \bar{z}_0\|_R^2, \ \forall k \geq 1. \end{aligned}$$

Therefore,

$$t_k^2 (f(z_{k+1}) - f^*) + \frac{1}{2}\|s_{k+1}\|_R^2 \leq \frac{1}{2}\|z_0 - \bar{z}_0\|_R^2, \ \forall k \geq 1,$$

from where, making use of Lemma 2.11.*(ii)*, we finally conclude that

$$f(z_{k+1}) - f^* \leq \frac{\|z_0 - \bar{z}_0\|_R^2}{2t_k^2} \leq \frac{2\|z_0 - \bar{z}_0\|_R^2}{(k+2)^2}, \ \forall k \geq 1.$$

Next, we prove claim *(ii)*. We start by proving the claim for $k = 0$. From (2.15) we derive

$$\|z_1 - \bar{z}_0\|_R \leq \|z_0 - \bar{z}_0\|_R. \tag{2.22}$$

Thus

$$\begin{aligned} \|\mathcal{G}(y_0)\|_{R^{-1}} &= \|R(y_0 - y_0^+)\|_{R^{-1}} = \|y_0 - y_0^+\|_R = \|z_0 - z_1\|_R \\ &= \|z_0 - \bar{z}_0 + \bar{z}_0 - z_1\|_R \leq \|z_0 - \bar{z}_0\|_R + \|z_1 - \bar{z}_0\|_R \\ &\leq 2\|z_0 - \bar{z}_0\|_R. \end{aligned}$$

We now prove the claim for $k > 0$. From $s_1 = \delta z_0 + t_0(\delta z_1 - \delta z_0) = z_1 - \bar{z}_0$, along with (2.22), we derive $\|s_1\|_R = \|z_1 - \bar{z}_0\|_R \leq \|z_0 - \bar{z}_0\|_R$, which in addition to the following inequality that can be deduced from (2.15):

$$\|s_{k+1}\|_R \leq \|z_0 - \bar{z}_0\|_R, \ \forall k \geq 1,$$

leads to

$$\|s_k\|_R \leq \|z_0 - \bar{z}_0\|_R, \ \forall k \geq 1.$$

From here we derive, for every $k \geq 1$,

$$\|s_{k+1} - s_k\|_R \leq \|s_{k+1}\|_R + \|s_k\|_R \leq \|z_0 - \bar{z}_0\|_R + \|z_0 - \bar{z}_0\|_R = 2\|z_0 - \bar{z}_0\|_R. \tag{2.23}$$

---

**Algorithm 5:** MFISTA

---

**Require:** $z \in \mathrm{dom}(f)$, $\epsilon \in \mathbb{R}_{>0}$

$y_0 \leftarrow z$, $z_0 \leftarrow z$, $t_0 \leftarrow 1$, $k \leftarrow 0$

**repeat**

$\quad k \leftarrow k + 1$

$\quad v_k \leftarrow \mathcal{T}_R^{f,\mathcal{Z}}(y_{k-1})$

$\quad t_k \leftarrow \dfrac{1}{2}\left(1 + \sqrt{1 + 4t_{k-1}^2}\right)$

$\quad z_k \leftarrow \begin{cases} v_k & \text{if } f(v_k) \leq f(z_{k-1}) \\ z_{k-1} & \text{otherwise} \end{cases}$

$\quad y_k \leftarrow z_k + \dfrac{t_{k-1}}{t_k}(v_k - z_k) + \dfrac{t_{k-1} - 1}{t_k}(z_k - z_{k-1})$

**until** $\|\mathcal{G}(z_k)\|_{R^{-1}} \leq \epsilon$

**Output:** $\tilde{z}^* \leftarrow z_k$

---

From (2.21) we have

$$\mathcal{G}_k = \frac{1}{t_k} R(s_k - s_{k+1}), \ \forall k \geq 1.$$

Therefore, taking into account Lemma 2.11.*(ii)*,

$$\|\mathcal{G}_k\|_{R^{-1}} = \frac{1}{t_k}\|s_k - s_{k+1}\|_R \overset{(2.23)}{\leq} \frac{2}{t_k}\|z_0 - \bar{z}_0\|_R \leq \frac{4}{k+2}\|z_0 - \bar{z}_0\|_R, \ \forall k \geq 1.$$

■

### 2.2.2 A monotone variant of FISTA

This section presents a *monotone* version of FISTA, which we label by MFISTA, and which is presented in [42, §V.A]. FISTA is not a monotone algorithm, in that the function values $\{f(z_k)\}_{\geq 0}$ generated by Algorithm 4 are not necessarily non-increasing. However, a small adjustment can be made to it to provide a monotone behaviour.

Algorithm 5 shows the MFISTA algorithm applied to problem (2.2) under Assumption 2.1 for an initial point $z \in \mathrm{dom}(f)$ and an exit tolerance $\epsilon \in \mathbb{R}_{>0}$. It returns a suboptimal solution $\tilde{z}^*$ of problem (2.2), where the suboptimality is determined by $\epsilon$. We note that Remark 2.10 also applies to Algorithm 5.

Note that MFISTA and FISTA both require a single computation of the composite gradient mapping $\mathcal{T}_R^{f,\mathcal{Z}}(y_{k-1})$ and a few vector-vector operations (with MFISTA requiring a few more than FISTA). The main difference between the two, computationally speaking, is that MFISTA requires the two function evaluations $f(v_k)$ and $f(z_{k-1})$ at each iteration, while FISTA requires none.

The interesting aspect of MFISTA is that it shares FISTA's convergence rate result given in Proposition 2.12.*(i)*, as stated in the following proposition.

**Proposition 2.13.** Consider problem (2.2) and let Assumtion 2.1 hold. Let $\{z_k\}_{\geq 0}$ be the sequence generated by MFISTA (Algorithm 5) for problem (2.2) starting at $z \in \text{dom}(f)$. Then,

$$f(z_k) - f^* \leq \frac{2}{(k+1)^2} \|z_0 - \pi^R_{\Omega_f}(z_0)\|^2_R, \; \forall k \geq 1.$$

***Proof:*** The proof of this proposition requires modifications to the proof given for FISTA (see Proposition 2.12). We refer the reader to [42] for its proof.

## 2.3 Conclusions

In this chapter we have presented some well known FOMs and AFOMs. As we discussed previously, the popularity of these methods has increased in recent years do to the rising number of applications that deal with large optimization problems. In this case, even though FOMs typically require more iterations to solve than methods that make use of the Hessian of the optimization problem, their lower computational cost per iteration may make them a more suitable choice.

However, these methods are not without issue. In addition to the fact that they may require a large number of iterations to converge, there are other aspects that have drawn the attention of the scientific community. One of the main challenges is the study and improvement of their theoretical and practical convergence properties. As an example, numerous results address the problem of selecting the penalty parameter of ADMM [22, 43, 44].

In the following chapter we will address one of these challenges: improvements of the convergence of AFOMs related to the undesirable oscillating behaviour that they display when close to the optimal solution.

# Chapter 3

# Restart schemes for accelerated first order methods

One of the drawbacks of accelerated first order methods is that they may suffer from undesirable oscillatory behavior which slows down their convergence [45]. The intuitive reason behind this is that AFMOs can be thought of as *momentum* driven, in the sense that each iterations depends on the previous ones with an added momentum that is increased over time. This momentum, coupled with the fact that most AFOMs are non-monotone in the objective function value, leads to the appearance, in many applications, of a periodic oscillation when close to the optimal solution. To illustrate this issue, let us take a look at a simple example.

**Example 3.1.** Consider the quadratic programming (QP) problem

$$\min_{z \in \mathbb{R}^2} \{ f(z) \doteq \frac{1}{2} z^\top H z + q^\top z \}, \tag{3.1}$$

where $H = \text{diag}(0.5, 1)$ and $q = -(0.1, 1)$. It is well known that the solution of this problem is given by $z^* = -H^{-1} q = (0.2, 1)$. We solve problem (3.1) using FISTA algorithm (Algorithm 4) with the initial condition $z = -(2, 5)$, and taking $R = 100\mathbf{I}_2$ and $\epsilon = 10^{-6}$.

**Remark 3.2.** *We note that our selection of $R$ in Example 3.1, while it does satisfy Assumtion 2.1.(ii) for the given smalll-sized optimization problem, could have easily been selected to provide a much better performance of FISTA. For larger optimization problems, however, it is common for $R$ to be overestimated, since the computation of a better value is not worth the computational burden. For the purposes of this example, it was selected this way to enhance the oscillatory behavior of the algorithm, thus providing visually obvious oscillations. We remark that the oscillatory behavior would be present, in a greater or lesser extent, when solving most optimization problems if the exit tolerance is sufficiently small. In Section 3.3 we show examples with larger-sized optimization problems.*

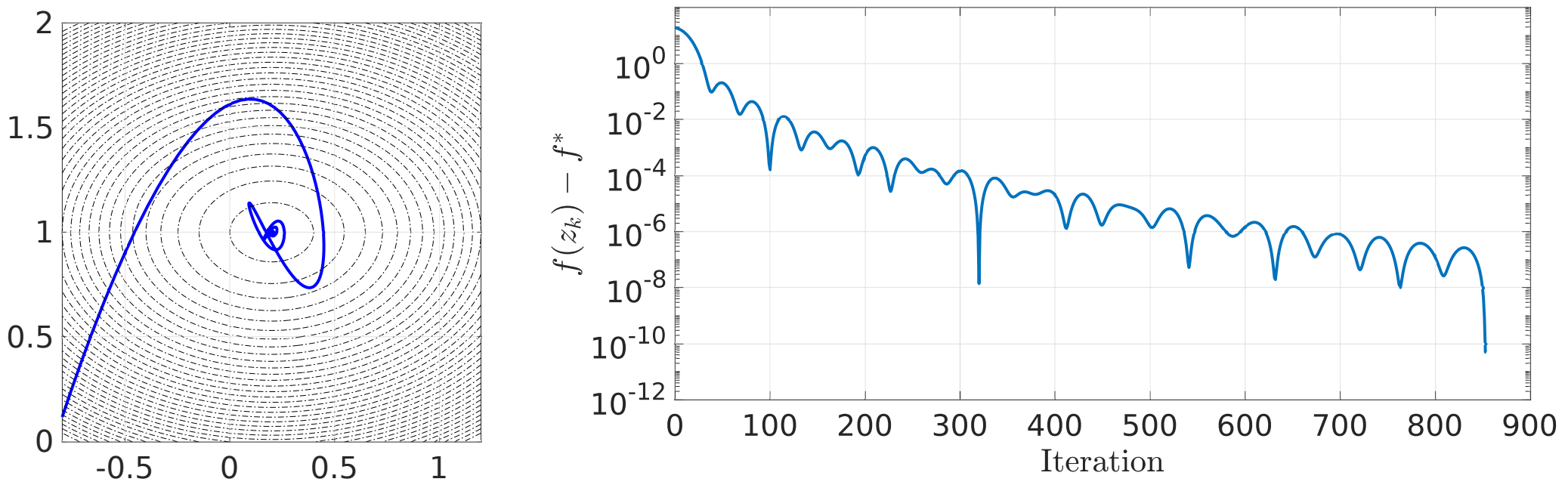

(a) Trajectory of $z_k$.

(b) Convergence to optimal objective function value.

Figure 3.1: Example 3.1. FISTA without a restart scheme.

Figure 3.1 shows the result of solving Example 3.1 using FISTA (Algorithm 4). Figure 3.1a shows the (incomplete) trajectory of the iterates $z_k$ produced by FISTA, whereas Figure 3.1b shows the difference $f(z_k) - f^*$ in a logarithmic scale. Notice that during the initial iterations $f(z_k)$ is monotonically converging towards the optimal solution, as can be seen in Figure 3.1b. However, it soon becomes non-monotone with the appearance of the oscillations. The general trend still converges towards the optimal solution, but at a slower rate. The algorithm converged after 853 iterations to the given tolerance.

This issue can be addressed by means of a *restart scheme*. Intuitively, a restart scheme is a procedure in which an iterative method is stopped when a certain condition, which we label the *restart condition*, is met, and is then restarted using the last value of the iterates as the new initial condition. This process continues until a *terminal* exit condition is satisfied.

In the context of AFOMs, and as a more formal description of a restart scheme, we consider the following. Let $\mathcal{A}$ be a non-specific AFOM that starts at a given initial point $z_0$, and generates a sequence $\{z_k\}_{\geq 0}$ that converges towards the optimal solution $z^*$ of the optimization problem as $k \to \infty$. The algorithm terminates at a finite iterate $k_{\text{out}}$ when a certain condition $E_c$ is met. In a non-restarted paradigm, this restart condition will be a certificate for the degree of suboptimality of $z_{\text{out}} \doteq z_{k_{\text{out}}}$. For instance, the exit condition of FISTA (Algorithm 4) is given by

$$E_c(z_k) = \begin{cases} \text{true} & \text{if } \|\mathcal{G}(z_k)\|_{R^{-1}} \leq \epsilon, \\ \text{false} & \text{otherwise}, \end{cases}$$

where $\epsilon \in \mathbb{R}_{>0}$ is the given exit tolerance. The exit condition $E_c$ is therefore a map to $\{\text{true}, \text{false}\}$ in relation to the current state of the algorithm. If the exit condition $E_c$ is evaluated to *true*, then the algorithm terminates, returning $z_{\text{out}} \leftarrow z_k$ and $k_{\text{out}} \leftarrow k$. We will write this as $[z_{\text{out}}, k_{\text{out}}] \leftarrow \mathcal{A}(z, E_c)$.

A restarted paradigm follows the procedure shown in Algorithm 6. The procedure uses an AFOM algorithm $\mathcal{A}$; an initial condition $r_0 \in \mathcal{C} \subseteq \mathbb{R}^{n_z}$, where $\mathcal{C}$ will depend on the algorithm $\mathcal{A}$; a *restart condition* $E_c$; and a *terminal* exit

**Algorithm 6:** General restart procedure for AFOMs.

**Require:** $r_0 \in \mathcal{C}$, $E_c$, $E_t$

$k \leftarrow 0$, $j \leftarrow 0$

**repeat**

  $j \leftarrow j + 1$

  $[r_j, k_j] \leftarrow \mathcal{A}(r_{j-1}, E_c)$

  $k \leftarrow k + k_j$

  Evaluate $E_t$;

**until** $E_t = true$

**Output:** $r_{\text{out}} \leftarrow r_j$, $j_{\text{out}} \leftarrow j$, $k_{\text{out}} \leftarrow k$

condition $E_t$. The algorithm makes successive calls to $\mathcal{A}$, generating a sequence $\{r_j\}_{\geq 0}$ of *restart points*, where $j$ is the *restart counter.* At each iteration $j$ (with the exception of the first iteration) algorithm $\mathcal{A}$ is given the value it returned in its previous call and is executed until the restart condition $E_c$ is satisfied. The *terminal* exit condition $E_t$ should be a measure of the optimality of $r_j$, whereas the *restart condition* $E_c$ can be any condition; although, obviously, it will be chosen to improve the (practical) convergence of this scheme when compared with the non-restarted variant. The algorithm returns, in $r_{\text{out}}$, the last $r_j$ provided by $\mathcal{A}$, the total number of restarts $j_{\text{out}}$, and the total number $k_{\text{out}}$ of iterations performed by $\mathcal{A}$ among all its calls. In the following, we will use the letters $r$ and $j$ to refer to the restart points and restart counter, respectively, and the letters $z$ and $k$ to refer to the *inner* iterates of $\mathcal{A}$.

## 3.1 A brief review of the literature

Various restart schemes for AFOMs have been presented in the literature [27, 45, 46, 47, 48, 49, 50, 51], each considering a particular class of optimization problem and with different assumptions on the knowledge (or lack thereof) of the parameters that characterize it.

Additionally, each one derives convergence properties of their restart scheme under certain assumptions. These convergence properties are often based on the notion of the *iteration complexity* of the algorithm, which is an upper bound on the number of iterations required to find an $\epsilon$-accurate optimal solution. Typically, and unless we explicitly state otherwise, this is expressed as the number of iterations $k$ required to reach a suboptimal solution $f(z_k)$ satisfying

$$f(z_k) - f^* \leq \epsilon,$$

but other measures of suboptimality may be used.

Complexity bounds have been derived for many first order methods and their accelerated variants. For instance, when applied to an $L$-smooth convex optimization problem, the projected gradient method attains a complexity $O(L/\epsilon)$, whereas the fast projected gradient method has the much better complexity $O(\sqrt{L/\epsilon})$ [3, 26, 52]. If, in addition, the problem is $\mu$-strongly convex, then the iteration complexities are $O\left(L/\mu \log(1/\epsilon)\right)$ for the projected gradient method and

$$O\left(\sqrt{\frac{L}{\mu}}\log\frac{1}{\epsilon}\right) \tag{3.2}$$

for the fast projected gradient method [3, 27, 49]. The $\log(1/\epsilon)$ term is particularly relevant here, since it is indicative of a *linear convergence rate*, i.e., when the iterates satisfy an expression of the form

$$f(z_k) - f^* \leq C\rho^k,$$

for some $C \in \mathbb{Z}_{>0}$ and $\rho \in (0, 1)$ [53, §3.1]. Furthermore, the iteration complexity (3.2) is the optimal one for this class of optimization problems using a first order black box oracle [3], meaning that a better iteration complexity cannot be attained (up to a constant, that is) with a first order method that does not take advantage of particularities of the problem structure.

One of the main questions when dealing with restart schemes is whether or not iteration complexities can be derived; in particular, which of the above and under what assumptions. Since restart schemes will be applied to AFOMs, we are particularly interested in determining under what conditions (if any), does a restart scheme attain the optimal iteration complexity (3.2).

### 3.1.1 Fixed-rate restart schemes

The most simple class of restart schemes are what we refer to as *fixed-rate* restart schemes, in which the AFOM is restarted every time it performs a predetermined fixed number of iterations. That is, the restart condition is given by

$$E_{k_m}(k) = \begin{cases} \text{true} & \text{if } k \geq k_m \\ \text{false} & \text{otherwise,} \end{cases}$$

where $k_m \in \mathbb{Z}_{>0}$ is fixed.

In [27, §5.1], a fixed-rate restart scheme is presented which considers the class of optimization problems $\min_{z \in \mathcal{Z}} f(z)$, where $\mathcal{Z} \subseteq \mathbb{R}^{n_z}$ is a closed convex set and $f : \mathbb{R}^{n_z} \to \mathbb{R}$ is a $\mu$-strongly convex and $L$-smooth function. By choosing the restart rate appropriately, this restart scheme recovers the optimal iteration complexity (3.2) for this class of optimization problem. Similar (and oftentimes equivalent) fixed-rate restart schemes are presented in other publications, such as in [48, §3.2], [45, §3.1].

Another particularly relevant fixed-rate restart scheme was presented in [46, §5.2.2]. The interest of this scheme is that it proves linear convergence of the restarted scheme for non-necessarily strongly convex optimization problems. In particular, the paper considers the class of optimization problems

$$f^* = \min_{z \in \mathcal{Z}} f(z), \tag{3.3}$$

where $\mathcal{Z} \subseteq \mathbb{R}^{n_z}$ is a closed convex set and $f : \mathbb{R}^{n_z} \to \mathbb{R}$ is a $L$-smooth convex function satisfying the following condition, known as the *quadratic functional growth condition* (see [46, Definition 4]), for some $\mu \in \mathbb{R}_{>0}$.

**Definition 3.3** (Quadratic functional growth)**.** Consider problem (3.3), where $\mathcal{Z} \subseteq \mathbb{R}^{n_z}$ is a closed convex set and $f : \mathbb{R}^{n_z} \to \mathbb{R}$, and assume that its optimal set $\Omega_f$ is nonempty and that $f^*$ is finite. We say that $f$ has a *quadratic functional growth* on $\mathcal{Z}$ if there exists a constant $\mu \in \mathbb{R}_{>0}$ such that

$$f(z) - f^* \geq \frac{\mu}{2}\|z - \pi_{\Omega_f}(z)\|_2^2, \; z \in \mathcal{Z}. \tag{3.4}$$

**Remark 3.4.** *Inequality* (3.4) *is satisfied, at least locally, for a large class of not necessarily strongly convex nor strictly convex functions [46, 54]. We note that it is always satisfied if $f$ is strongly convex.*

The scheme presented in [46, §5.2.2] can be implemented in two different ways. The first is by performing the fixed-rate restart given by

$$E_{k^*}(k) = \begin{cases} \text{true} & \text{if } k \geq k^*, \\ \text{false} & \text{otherwise,} \end{cases}$$

where

$$k^* = \left\lceil 2e\sqrt{\frac{L}{\mu}} \right\rceil.$$

The second alternative is to restart $\mathcal{A}$ whenever the following restart condition $E_f^* : \mathbb{R}^{n_z} \to \{\text{true}, \text{false}\}$ is satisfied:

$$E_f^*(z_k) = \begin{cases} \text{true} & \text{if } f(z_k) - f^* \leq \dfrac{f(z_0) - f^*}{e^2}, \\ \text{false} & \text{otherwise.} \end{cases} \tag{3.5}$$

We note that this second approach, while equivalent to the first in the sense that the same convergence results are obtained, is not (strictly speaking) a fixed-rate restart scheme.

There are two main problems with fixed-rate restart schemes. First, they require knowledge of parameters that are not typically known beforehand or that are often computationally expensive to compute, such as the parameters that characterize the smoothness, the strong convexity, the quadratic functional growth or the optimum value $f^*$. Second, the restart rate is derived from global parameters of the optimization problem, and may thus result in an inappropriate restart rate in better conditioned regions [45, §3.1].

**Example 3.5.** As an example of the use of a fixed-rate restart scheme, Figure 3.2 shows the use of the restart scheme from [46, §5.2.2] using the exit condition (3.5) applied to FISTA to solve Example 3.1, along with the result of the non-restarted

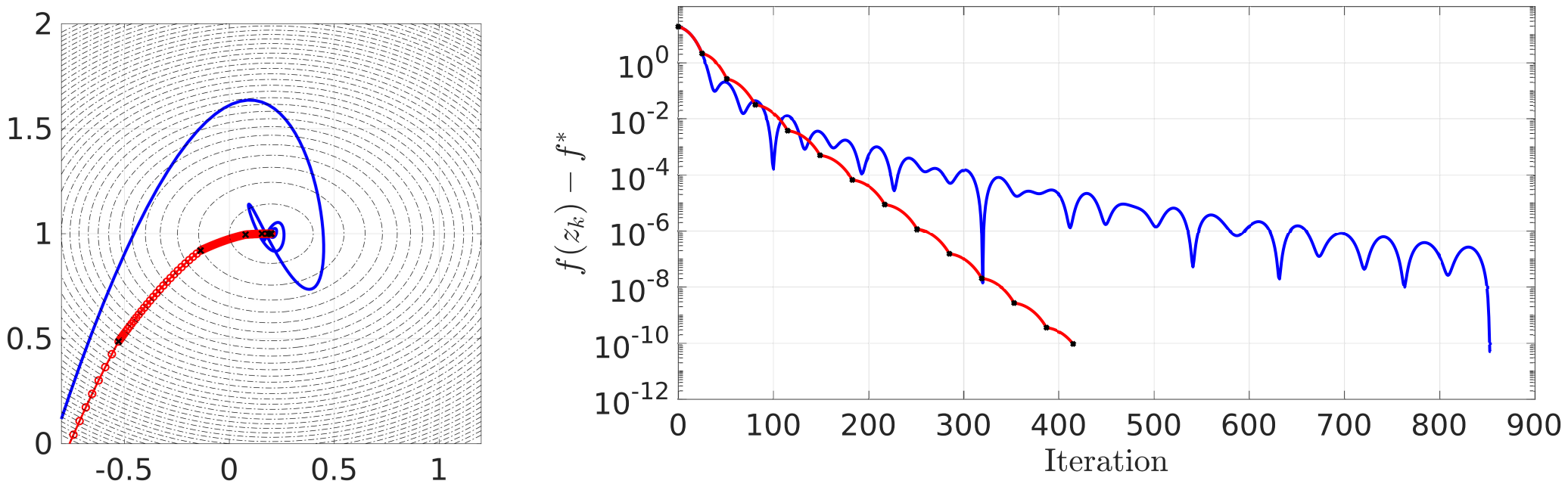

(a) Trajectory of $z_k$.

(b) Convergence to optimal objective function value.

Figure 3.2: Example 3.1. FISTA with the fixed-rate restart scheme from [46, §5.2.2]. The restarted version is depicted in red. The non-restarted variant is depicted in blue. The red circles in (a) depict each of the iterates $z_k$. The black crosses depict each of the restart points $r_j$.

variant of FISTA. The terminal exit condition was taken as the exit condition of (non-restarted) FISTA (step 7 of Algorithm 4) with the exit tolerance $\epsilon$ given in Example 3.1. The total number of FISTA iterations is $k_{\text{out}} = 415$ and the number of restarts is $j_{\text{out}} = 13$.

### 3.1.2 Adaptive restart schemes

The second main class of restart schemes are commonly referred to as *adaptive* restart schemes, because they either adapt the restart rate online or directly use a restart condition that serves to evaluate the performance of the AFOM.

In [45], two simple and well-known adaptive restart schemes are presented; which we label as the *objective function value* scheme, and the *gradient alignment* scheme. It is shown that they both preserve the optimal convergence rate of the fast gradient method for smooth strongly convex optimization problems, provided that the objective function is quadratic.

The *objective function value restart scheme* performs a restart of the AFOM each time the current iterate does not decrease with respect to the one of the previous iterate. That is, the restart condition $E_f : \mathbb{R}^{n_z} \times \mathbb{R}^{n_z} \to \{\text{true}, \text{false}\}$ of this scheme is given by

$$E_f(z_k, z_{k-1}) = \begin{cases} \text{true} & \text{if } f(z_k) - f(z_{k-1}) > 0 \\ \text{false} & \text{otherwise.} \end{cases} \tag{3.6}$$

The idea behind this restart scheme is very intuitive: the oscillations that appear in AFOMs can be viewed as the algorithm "overshooting" the optimal value. Thus, one of the simplest indicators of the presence of the oscillations is an increase of the objective function value.

**Example 3.6.** Figure 3.3 shows the application of the *objective function value restart scheme* using FISTA, alongside the non-restarted variant, to solve Exam-

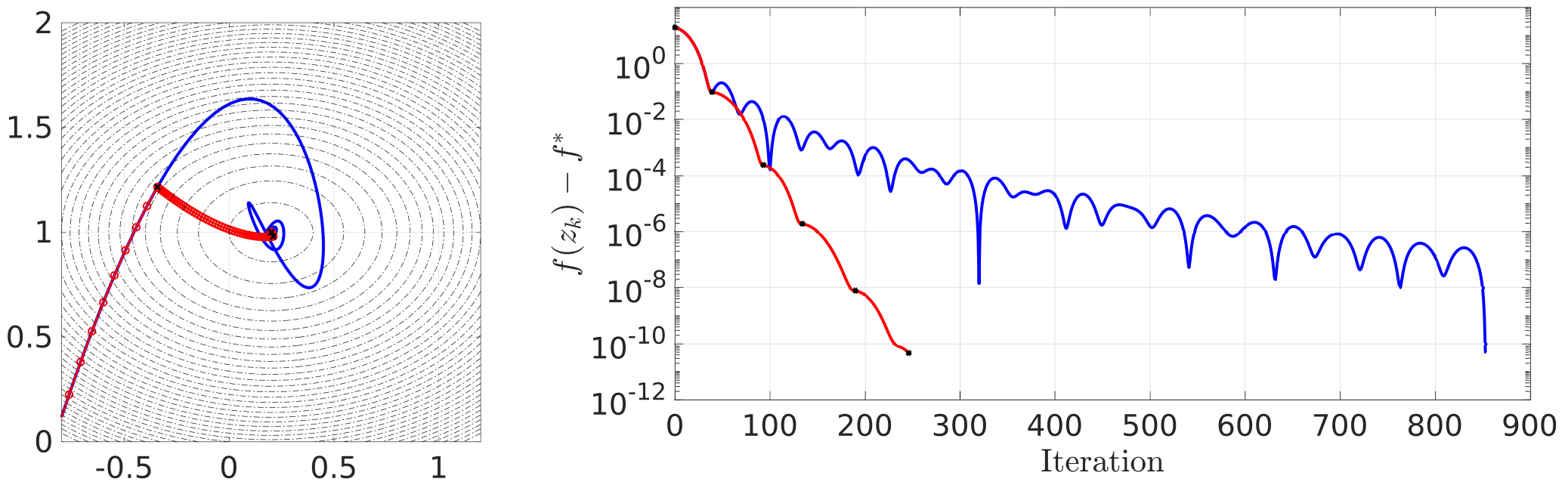

(a) Trajectory of $z_k$.

(b) Convergence to optimal objective function value.

Figure 3.3: Example 3.1. FISTA with the objective function value restart scheme. The restarted version is depicted in red. The non-restarted variant is depicted in blue. The red circles in (a) depict each of the iterates $z_k$. The black crosses depict each of the restart points $r_j$.

ple 3.1. The terminal exit condition was taken as the exit condition of (non-restarted) FISTA (step 7 of Algorithm 4) with the exit tolerance $\epsilon$ given in Example 3.1. The total number of FISTA iterations is $k_{\text{out}} = 246$ and the number of restarts is $j_{\text{out}} = 5$.

In its simplest form, the *gradient alignment restart scheme* considers the class of non-constrained optimization problems $\min_{z \in \mathbb{R}^{n_z}} f(z)$ with a continuously differentiable objective function $f : \mathbb{R}^{n_z} \to \mathbb{R}$. The scheme performs a restart of the AFOM each time the direction in which the iterates are "moving", i.e., $z_k - z_{k-1}$, is not aligned with the negative of the gradient $\nabla f(z_{k-1})$ of the objective function. The idea is to restart the AFOM whenever the iterates are moving in a "bad" direction, as measured by the negative of the gradient.

The extension of this scheme to constrained optimization problems and/or to problems with non-differentiable objective functions is to substitute the gradient $\nabla f(\cdot)$ with whichever gradient operator is used by the AFOM (the projected gradient, the projected subgradient, the proximal operator, the composite gradient mapping, etc.). Thus, the restart condition $E_g : \mathbb{R}^{n_z} \times \mathbb{R}^{n_z} \times \mathbb{R}^{n_z} \to \{\text{true}, \text{false}\}$ of this scheme can be expressed as

$$E_g(z_k, z_{k-1}, y; g) = \begin{cases} \text{true} & \text{if } \langle g(y), z_k - z_{k-1} \rangle > 0 \\ \text{false} & \text{otherwise,} \end{cases} \tag{3.7}$$

where $g : \mathbb{R}^{n_z} \to \mathbb{R}^{n_z}$ is the gradient operator used by the AFOM, and $y$ is the point where the gradient operator is evaluated. By default, we consider $g \equiv \mathcal{G}_R^{f,\mathcal{Z}}$, where $\mathcal{G}_R^{f,\mathcal{Z}}$ is the composite gradient mapping operator given by Definition 2.3. In this case, we simply label $E_c(z_k, z_{k-1}, y; \mathcal{G}_R^{f,\mathcal{Z}})$ by $E_c(z_k, z_{k-1}, y)$.

**Remark 3.7.** *The most straightforward choice of $y$ in the restart condition* (3.7) *is to take $y = z_{k-1}$, but this may not always be the case. In fact, to avoid*

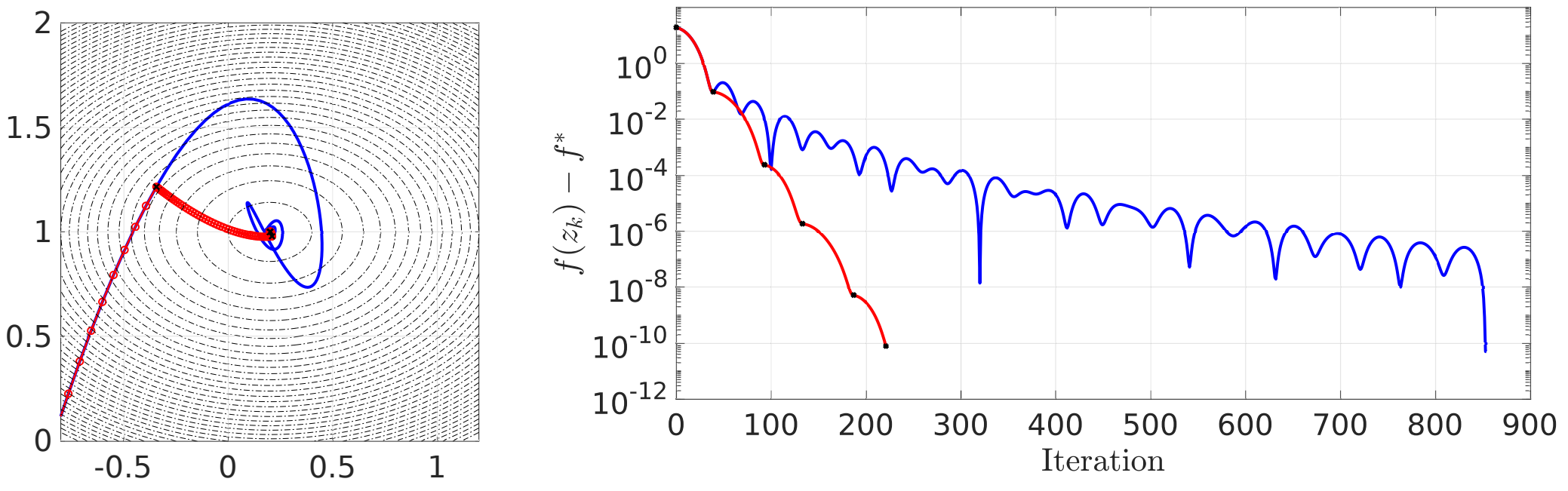

(a) Trajectory of $z_k$. (b) Convergence to optimal objective function value.

Figure 3.4: Example 3.1. FISTA with the gradient alignment restart scheme. The restarted version is depicted in red. The non-restarted variant is depicted in blue. The red circles in (a) depict each of the iterates $z_k$. The black crosses depict each of the restart points $r_j$.

*additional computations, the best choice is to choose $y$ as whichever point is used in the AFOM to evaluate the gradient operator $g(\cdot)$ at. For instance, in FISTA (Algorithm 4), where $g(y) = \mathcal{G}(y) = R(y - \mathcal{T}(y))$, the best choice is to take $y = y_{k-1}$, since $\mathcal{T}(y_{k-1})$ is evaluated in step 4.*

As stated in [45, §3.2], this scheme has two advantages when compared to the *objective function value* restart scheme. First, it is more numerically stable, since the gradient will tend to a fixed point as we approach the optimum, whereas $f(z_k) - f(z_{k-1})$ will tend to present numerical issues due to cancellation errors as $f(z_k)$ approaches $f^*$. Second, the evaluation of its restart condition requires no additional computations, since $g(y)$ is evaluated in AFOMs (see Remark 3.7).

**Example 3.8.** Figure 3.4 shows the result of applying the *gradient alignment restart scheme* to Example 3.1 using FISTA, along with the result of the non-restarted variant of FISTA. The terminal exit condition was taken as the exit condition of (non-restarted) FISTA (step 7 of Algorithm 4) with the exit tolerance $\epsilon$ given in Example 3.1. We apply the exit condition $E_g(z_k, z_{k-1}, y_{k-1})$, where $y_{k-1}$ is chosen as the point in which the composite gradient mapping is evaluated because $\mathcal{T}(y_{k-1})$ is computed in step 4 of Algorithm 4, thus avoiding additional computations (see Remark 3.7). In this case, since the QP problem is unconstrained, the composite gradient mapping $\mathcal{G}$ of FISTA is simply the gradient $\nabla f$. The total number of FISTA iterations is $k_{\rm out} = 221$ and the number of restarts is $j_{\rm out} = 5$.

In [47], the restart schemes $E_f$ and $E_g$ from [45] are applied to the class of composite optimization problems

$$\min_{z \in \mathbb{R}^{n_z}} \Psi(z) + h(z), \tag{3.8}$$

where $\Psi \in \Gamma(\mathbb{R}^{n_z})$ and $h : \mathbb{R}^{n_z} \to \mathbb{R}$ is a smooth convex function. Additionally, they present a novel variation of the restart schemes proposed in [45]. Finally, they present some numerical results that indicate that these restart methods provide good practical convergence for this class of composite optimization problems, even though, for the most part, they lack theoretical linear convergence rates.

In [51], the adaptive restart schemes of [45] are applied to the optimized gradient method [55], as well as a similar, but novel, adaptive restart scheme for this FOM. However, no convergence rates nor iteration complexities are derived.

In [49], a restart scheme for the accelerated proximal gradient method is presented. The authors extend the results from [27, §5.3] by presenting a scheme with the same iteration complexity but that requires less computations of the gradient operator. The scheme considers the class composite optimization problems (3.8) where $\Psi \in \Gamma(\mathbb{R}^{n_z})$ and $h : \mathbb{R}^{n_z} \to \mathbb{R}$ is $L$-smooth and $\mu$-strongly convex. The proposed scheme starts with an estimate $\hat{\mu}$ of $\mu$, which is hopefully a lower bound, i.e., $\hat{\mu} \leq \mu$. It then runs the accelerated proximal gradient method, where $L$ is determined by line search, using $\hat{\mu}$ and performs checks to determine if it is indeed a lower bound of $\mu$. It it is determined that $\hat{\mu}$ does not satisfy $\hat{\mu} \leq \mu$, then it is decreased and the algorithm starts from the previous restart point, i.e., the iterates obtained from the last call to the accelerated proximal gradient method are discarded. Even though this may result in the execution of "useless" iterations, linear convergence is still derived. Additionally, they adapt the scheme to the case of high-dimensional Lasso problems, which are not generally strongly convex over the entire domain of $h$.

In [50], a restart scheme that improves upon the results of [27, §5.3] and [49] is presented, wherein an (improved) iteration complexity (3.2) is derived for the same class of optimization problems but where $h$ only needs to satisfy the quadratic functional growth condition with parameter $\mu$ (Definition 3.3), instead of it needing to be strongly convex. The proposed approach performs calls to the AFOM with a fixed number of iterations that is determined by estimating the parameter $\mu$. After each call to the AFOM, it is determined whether the estimate was correct, and adjusted if it is not. An additional improvement of this scheme, when compared with the ones in [27, §5.3] and [49], is that no iterates of the AFOM are "wasted" if the estimation of $\mu$ is changed. Finally, the authors in [50] show that a fixed-rate restart scheme for AFOMs applied to the above class of optimization problem is linearly convergent (in the sense described in [50, Corollary 3]), for any positive choice of the restart rate. However, the optimal iteration complexity (3.2) is only recovered if the problem parameters are known.

## 3.2 Implementable restart schemes with linear convergence

This section presents three novel restart schemes for AFOMs. Section 3.2.1 presents a restart scheme for FISTA that is based on monitoring the evolution of the objective function values, whereas Section 3.2.2 presents a restart scheme for FISTA that is based on monitoring the composite gradient mapping. Both schemes consider the class of optimization problems (2.2) under Assumtion 2.1 and under the satisfaction of a weighted quadratic functional growth. Under these assumptions, linear convergence of the restarted FISTA algorithm is recovered. The restart scheme from Section 3.2.1 is then extended to a wider class of AFOMs in Section 3.2.3.

### 3.2.1 A restart scheme for FISTA

This section presents a restart scheme, originally presented in [5], for FISTA (Algorithm 4) that exhibits linear convergence for the class of optimization problems

$$f^* = \min_{z \in \mathcal{Z}} \{ f(z) \doteq \Psi(z) + h(z) \} \tag{3.9}$$

under Assumtion 2.1 and the following assumption:

**Assumption 3.9** (Weighted quadratic functional growth)**.** Problem (3.9) satisfies the following weighted quadratic functional growth condition (see Definition 3.3) with parameter $\mu \in \mathbb{R}_{>0}$ for the $R$-weighted Euclidean norm $\| \cdot \|_R$:

$$f(z) - f^* \geq \frac{\mu}{2} \| z - \pi^R_{\Omega_f}(z) \|^2_R, \; \forall z \in \mathcal{Z}. \tag{3.10}$$

**Remark 3.10.** *We note that Assumtion 3.9 is equivalent to the quadratic functional growth condition (Definition 3.3) in the event that $h$ is $L$-smooth in the sense of Definition N.20 instead of being smooth in the sense of Assumtion 2.1.(ii).*

The proposed restart scheme is based on the following proposition, which presents some novel results that further characterize the convergence properties of FISTA under the above quadratic functional growth assumption.

**Proposition 3.11.** Consider problem (3.9) under Assumtions 2.1 and 3.9. Then, the iterates of FISTA (Algorithm 4) satisfy:

*(i)* $f(z_k) - f^* \leq \dfrac{4}{\mu(k+1)^2}(f(z_0) - f^*), \; \forall k \geq 1.$

*(ii)* $f(z_k) \leq f(z_0), \; \forall k \geq \left\lfloor \frac{2}{\sqrt{\mu}} \right\rfloor.$

*(iii)* $f(z_k) - f^* \leq \dfrac{f(z_0) - f(z_k)}{e}, \; \forall k \geq \left\lfloor \frac{2\sqrt{e+1}}{\sqrt{\mu}} \right\rfloor.$

***Proof:*** For convenience, we denote $\bar{z}_0 \doteq \pi^R_{\Omega_f}(z_0)$. We recall that, from Proposition 2.12.*(i)*, we have

$$f(z_k) - f^* \leq \frac{2}{(k+1)^2}\|z_0 - \bar{z}_0\|^2_R, \ \forall k \geq 1,$$

which along (3.10) leads to

$$f(z_k) - f^* \leq \frac{4}{\mu(k+1)^2}(f(z_0) - f^*), \ \forall k \geq 1,$$

where the fact that $z_k \in \mathcal{Z}$ follows from step 4 of Algorithm 4 and the definition of $\mathcal{T}$ (Definition 2.3). This proves claim *(i)*. Denote

$$\alpha_k \doteq \frac{4}{\mu(k+1)^2}, \ \forall k \geq 1,$$

and suppose that $k \geq \left\lfloor \frac{2}{\sqrt{\mu}} \right\rfloor$. Then,

$$\alpha_k = \frac{4}{\mu(k+1)^2} \leq \frac{4}{\mu\left(\left\lfloor \frac{2}{\sqrt{\mu}} \right\rfloor + 1\right)^2} < \frac{4}{\mu\left(\frac{2}{\sqrt{\mu}}\right)^2} = 1.$$

Therefore, $\alpha_k \in (0,1)$, $\forall k \geq \left\lfloor \frac{2}{\sqrt{\mu}} \right\rfloor$, which along with claim *(i)*, leads to

$$f(z_k) - f^* \leq f(z_0) - f^*, \ \forall k \geq \left\lfloor \frac{2}{\sqrt{\mu}} \right\rfloor,$$

thus proving claim *(ii)*. Next, in view of claim *(i)*, we have

$$\begin{aligned} f(z_k) - f^* &\leq \alpha_k(f(z_0) - f^*) = \alpha_k(f(z_0) - f(z_k) + f(z_k) - f^*) \\ &= \alpha_k(f(z_0) - f(z_k)) + \alpha_k(f(z_k) - f^*). \end{aligned}$$

Therefore

$$(1 - \alpha_k)(f(z_k) - f^*) \leq \alpha_k(f(z_0) - f(z_k)). \tag{3.11}$$

Suppose now that $k \geq \left\lfloor \frac{2\sqrt{e+1}}{\sqrt{\mu}} \right\rfloor$. Note that this implies that $k \geq \left\lfloor \frac{2}{\sqrt{\mu}} \right\rfloor$, which in turn implies that $1 - \alpha_k > 0$. Therefore, we can divide both terms of inequality (3.11) to obtain

$$\begin{aligned} f(z_k) - f^* &\leq \frac{\alpha_k}{1 - \alpha_k}(f(z_0) - f(z_k)) = \frac{\frac{4}{\mu(k+1)^2}}{1 - \frac{4}{\mu(k+1)^2}}(f(z_0) - f(z_k)) \\ &= \frac{4(f(z_0) - f(z_k))}{\mu(k+1)^2 - 4} \leq \frac{4(f(z_0) - f(z_k))}{\mu(\left\lfloor \frac{2\sqrt{e+1}}{\sqrt{\mu}} \right\rfloor + 1)^2 - 4} \\ &\leq \frac{4(f(z_0) - f(z_k))}{\mu(\frac{2\sqrt{e+1}}{\sqrt{\mu}})^2 - 4} = \frac{4(f(z_0) - f(z_k))}{4(e+1) - 4} = \frac{f(z_0) - f(z_k)}{e}, \end{aligned}$$

which proves claim *(iii)*. ■

**Algorithm 7:** Restart FISTA based on objective function values

**Require:** $r_0 \in \mathcal{Z}$, $\epsilon > 0$

$n_0 \leftarrow 0$, $j \leftarrow 1$ $k \leftarrow 0$

$[r_1, n_1] \leftarrow \text{FISTA}\left(r_0, E_f^e(n_0)\right)$

**repeat**

  $j \leftarrow j + 1$

  $[r_j, n_j] \leftarrow \text{FISTA}\left(r_{j-1}, E_f^e(n_{j-1})\right)$

  $k \leftarrow k + n_j$

  **if** $f(r_{j-1}) - f(r_j) > \dfrac{f(r_{j-2}) - f(r_{j-1})}{e}$ **then**

    $n_j \leftarrow 2n_{j-1}$

  **end if**

**until** $\|\mathcal{G}(r_j)\|_{R^{-1}} \leq \epsilon$

**Output:** $r_{\text{out}} \leftarrow r_j$, $j_{\text{out}} \leftarrow j$, $k_{\text{out}} \leftarrow k$

Proposition 3.11.*(iii)* tells us that there exists a finite number of iterations after which the distance to the optimal value $f^*$ is smaller than a fraction of the distance covered from $z_0$ to the current iterate. Furthermore, notice that $\left\lfloor \frac{2\sqrt{e+1}}{\sqrt{\mu}} \right\rfloor > \left\lfloor \frac{2}{\sqrt{\mu}} \right\rfloor$. Therefore, if $k \geq \left\lfloor \frac{2\sqrt{e+1}}{\sqrt{\mu}} \right\rfloor$, not only is Proposition 3.11.*(iii)* satisfied, but also, from Proposition 3.11.*(ii)*, we have that $f(z_0) - f(z_k) \geq 0$. Therefore, a fixed-rate restart scheme with the rate set to $k_m = \left\lfloor \frac{2\sqrt{e+1}}{\sqrt{\mu}} \right\rfloor > \left\lfloor \frac{2}{\sqrt{\mu}} \right\rfloor$ would provide linear convergence. However, this restart scheme would suffer from the same issues of other fixed-rate restart schemes: the value of $\mu$ would be required for its implementation, and the use of the global $\mu$ might lead to a worse practical convergence than if local information was used instead (see Section 3.1.1).

**Remark 3.12.** *We note that claims (ii) and (iii) of Proposition 3.11 may be satisfied for values of $k$ smaller that the given bounds. That is, the proposition states that the satisfaction of both inequalities is guaranteed if $k$ is larger or equal to the given bounds, but they may also be satisfied for smaller values of $k$.*

We propose the restart scheme shown in Algorithm 7, which uses the restart condition $E_f^e : \mathbb{R}^{n_z} \times \mathbb{R}^{n_z} \times \mathbb{Z}_{>0} \to \{\text{true}, \text{false}\}$ given by

$$
E_f^e(z_k, z_m, k; n) = \text{true} \iff \begin{cases} f(z_m) - f(z_k) \leq \dfrac{f(z_0) - f(z_m)}{e} & \text{(3.12a)} \\ f(z_k) \leq f(z_0) & \text{(3.12b)} \\ k \geq n, & \text{(3.12c)} \end{cases}
$$

where $m = \left\lfloor \frac{k}{2} \right\rfloor + 1$. For convenience, we drop the $(z_k, z_m, k)$ notation, and instead simply write $E_f^e(n)$. Starting from an initial point $r_0 \in \mathcal{Z}$, and for a given exit tolerance $\epsilon \in \mathbb{R}_{>0}$, Algorithm 7 makes successive calls to FISTA using the above exit condition. The algorithm returns an $\epsilon$-accurate solution $r_{\text{out}}$ of problem (3.9), in terms of $\|\mathcal{G}(r_{\text{out}})\|_{R^{-1}} \leq \epsilon$ (which is a valid condition for suboptimality of $r_{\text{out}}$, as shown in Proposition 2.6); the total number of restart iterations $j_{\text{out}}$; and the total number of iterations $k_{\text{out}}$ performed by FISTA.

The idea behind this restart scheme is to use an approximation of claims *(ii)* and *(iii)* of Proposition 3.11, where for claim *(iii)* $z_k$ is used as our best estimation of $z^*$ and $z_m$ is the middle point of the sequence of iterates generated by the current call to FISTA. Therefore, the restart condition is only true if *(i)* the current iterate $z_k$ of FISTA is no worse, in terms of its objective function value, than the current initial condition $z_0$, and *(ii)* if the amount that the objective function value has decreased from $z_0$ to $z_m$ is $e$ times smaller than the amount decreased from $z_m$ to $z_k$. That is, (3.12a) serves as a way of detecting a degradation of the performance of FISTA, whereupon the condition is satisfied if the second half of the iterates are providing significantly less benefits that the first half, measured in terms of the decrease obtained in the objective function value; whereas (3.12b) is guaranteeing that, at the very least, the new restart point is no worse than the previous one.

**Remark 3.13.** *In fact, the imposition of* (3.12b) *guarantees that* $f(r_{j-1}) > f(r_j)$ *provided that* $r_{j-1} \neq z^*$*. This follows from noting that at every restart iteration* $j$ *of Algorithm 7, step 1 of Algorithm 4 performs the assignment* $z_0 \leftarrow \mathcal{T}(r_{j-1})$*. Thus, in view of Proposition 2.5.(ii) and the use of the exit condition* (3.12b)*, we have that* $f(r_{j-1}) > f(z_0) \geq f(z_{out}) = f(r_j)$ *if* $r_{j-1} \neq z^*$*. Therefore, in general, a net improvement is attained at every iteration* $j$*.*

The problem is that, due to the non-monotone nature of FISTA, the amounts $f(z_m) - f(z_k)$ and $f(z_0) - f(z_m)$ are not necessarily positive. However, we know from Proposition 3.11 that this will not the case if $k$ is sufficiently large. For this reason, the restart condition also includes a minimum number of iterations $n$ that is adapted between calls to FISTA. In particular, at each iterate $j$ of Algorithm 7, we set the minimum number of iterations of FISTA as $n_j$, which is generally taken as the number of iterations performed by FISTA in iteration $j-1$. However, this minimum number of iterations is doubled if the condition

$$f(r_{j-1}) - f(r_j) > \frac{f(r_{j-2}) - f(r_{j-1})}{e} \tag{3.13}$$

is satisfied, which, once again, is an approximation of Proposition 3.11.*(iii)* but this time between restart points. The idea is that $n_j$ will eventually converge to a number of iterations in which claims *(ii)* and *(iii)* of Proposition 3.11 are being satisfied by the outputs of FISTA obtained at each iteration $j$. If these conditions are being satisfied for the current value of $n_j$, then it is not increased. However, if we detect that they are not satisfied, either because FISTA exits with a number of iterations $n_j$ larger than $n_{j-1}$ or because condition (3.13) is satisfied, then $n$ is increased. The use of $n_j$ can be viewed as a way of finding the optimal restart rate, in the sense of the rate used by the fixed-rate restart schemes (see

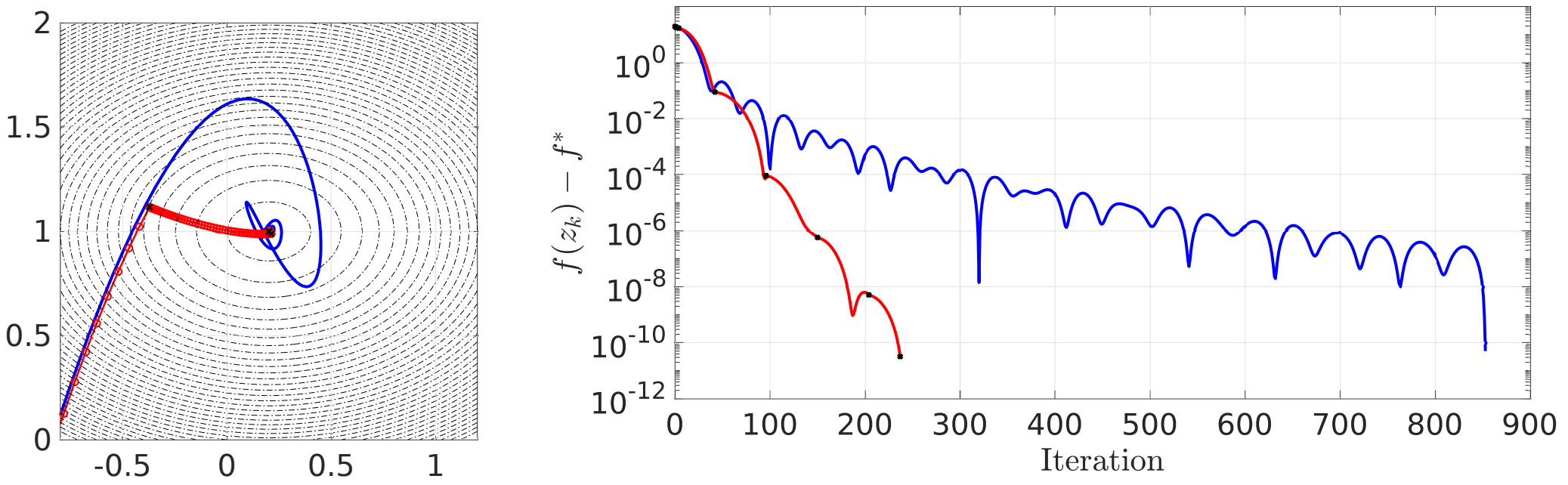

(a) Trajectory of $z_k$.

(b) Convergence to optimal objective function value.

Figure 3.5: In red, Algorithm 7 applied to Example 3.1. The non-restarted FISTA is depicted in blue. The red circles in (a) depict each of the iterates $z_k$. The black crosses depict each of the restart points $r_j$.

Section 3.1.1). However, we do not use global information. Instead, if the iterates a providing a good performance, then we do not increase $n_j$ artificially.

The main drawbacks of this restart scheme are that *(i)* it requires evaluating the objective function value at every iteration of FISTA, which is not necessary in the non-restarted FISTA algorithm, and *(ii)* the evaluation of (3.12a) requires storing the historic of the past objective function values for the current call to FISTA. This latter issue is not a significant problem in practical terms, since the use of a maximum desired number of iterations in a real implementation of iterative algorithms is commonplace. Thus, memory can be allocated for storing the past values of the objective function values even if the value of $\mu$ is unknown.

**Example 3.14.** Figure 3.5 shows the use of Algorithm 7 to solve Example 3.1 along with the result of the non-restarted FISTA algorithm. The total number of FISTA iterations is $k_{\text{out}} = 237$ and the number of restarts is $j_{\text{out}} = 8$.

The main properties of Algorithm 7 are summarized in the following theorem. Claim *(i)* simply provides an upper bound of $\|\mathcal{G}(r_{j-1})\|_{R^{-1}}$ in terms of the difference $f(r_{j-1}) - f(r_j)$, which due to (3.12b) is non-negative; claim *(ii)* provides an upper bound to the value of $n_j$; and claim *(iii)* provides the iteration complexity of the algorithm to reach an $\epsilon$-accurate solution of problem (3.9) in terms of $\|\mathcal{G}(r_{\text{out}})\|_{R^{-1}} \leq \epsilon$. Note that the iteration complexity is given by

$$O\left(\frac{1}{\sqrt{\mu}} \log \frac{1}{\epsilon}\right), \tag{3.14}$$

which we note is optimal in the sense discussed in Section 3.1 and shown in (3.2).

**Remark 3.15.** *The reader may have noted that the Lipschitz constant $L$ does not seem to be present in the iteration complexity* (3.14)*, as is the case of the iteration complexities provided in Section 3.1 for the AFOMs. However, this is not the case, since parameter $R$ effectively accounts for $L$ (see Remark 2.2). If we were to rewrite the characterization of the smoothness of $h$ using the $L$-smoothness* *condition given in Definition N.20 instead of using* (2.3)*, then* $\sqrt{L}$ *would appear in* (3.14)*.*

**Theorem 3.16** (Convergence properties of Algorithm 7)**.** Consider problem (3.9) under Assumtions 2.1 and 3.9. Then, the sequences $\{r_j\}$ and $\{n_j\}$ generated by Algorithm 7 and its output $k_{\text{out}}$ satisfy:

*(i)* $\frac{1}{2}\|\mathcal{G}_R^{f,\mathcal{Z}}(r_{j-1})\|_{R^{-1}}^2 \leq f(r_{j-1}) - f(r_j), \ \forall j \geq 1.$

*(ii)* $n_j \leq \dfrac{4\sqrt{e+1}}{\sqrt{\mu}}, \ \forall j \geq 0.$

*(iii)* The total number of iterations of FISTA required to attain the exit condition $\|\mathcal{G}_R^{f,\mathcal{Z}}(r_{\text{out}})\|_{R^{-1}} \leq \epsilon$ is upper bounded by

$$k_{\text{out}} \leq \frac{16}{\sqrt{\mu}} \left\lceil \log\left(1 + \frac{2(f(r_0) - f^*)}{\epsilon^2}\right) \right\rceil.$$

***Proof:*** By construction, $r_{j-1} \in \mathcal{Z}$, for all $j \geq 1$. Therefore, we have from Proposition 2.5.*(ii)*, that

$$\frac{1}{2}\|\mathcal{G}(r_{j-1})\|_{R^{-1}}^2 \leq f(r_{j-1}) - f(\mathcal{T}(r_{j-1})), \ \forall j \geq 1. \tag{3.15}$$

We also notice that $r_j$ is computed invoking FISTA algorithm using $r_{j-1}$ as initial condition ($z = r_{j-1}$). That is, $[r_j, n_j] \leftarrow FISTA(r_{j-1}, E_f^e(n_{j-1}))$. Since the output value $f(r_j)$ is forced to be no larger than the one corresponding to $z_0 = \mathcal{T}(z) = \mathcal{T}(r_{j-1})$, we have $f(r_j) \leq f(\mathcal{T}(r_{j-1}))$. Therefore, we obtain from inequality (3.15) that

$$\frac{1}{2}\|\mathcal{G}(r_{j-1})\|_{R^{-1}}^2 \leq f(r_{j-1}) - f(\mathcal{T}(r_{j-1})) \leq f(r_{j-1}) - f(r_j), \ \forall j \geq 1,$$

which proves claim *(i)*.

Let us now prove that if $n_{j-1} \leq \frac{4\sqrt{e+1}}{\sqrt{\mu}}$, then the value $n_j$ obtained from $[r_j, n_j] \leftarrow FISTA(r_{j-1}, E_f^e(n_{j-1}))$, also satisfies

$$n_j \leq \frac{4\sqrt{e+1}}{\sqrt{\mu}}. \tag{3.16}$$

Let us denote $\bar{m} \doteq \left\lfloor \frac{2\sqrt{e+1}}{\sqrt{\mu}} \right\rfloor$. Since $\bar{m} \geq \left\lfloor \frac{2\sqrt{e+1}}{\sqrt{\mu}} \right\rfloor$, we infer, from Proposition 3.11.*(iii)* that

$$f(z_{\bar{m}}) - f^* \leq \frac{f(z_0) - f(z_{\bar{m}})}{e}.$$

From this inequality, we obtain

$$f(z_{\bar{m}}) - f(z_k) \leq f(z_{\bar{m}}) - f^* \leq \frac{f(z_0) - f(z_{\bar{m}})}{e}.$$

Therefore, the restart condition (3.12a) is satisfied for $m = \bar{m}$. Since $m = \lfloor \frac{k}{2} \rfloor + 1$ we have $m \geq \frac{k}{2}$. This means that for $m = \bar{m}$, the corresponding value for $k$ is no larger than

$$k \leq 2\bar{m} = 2\left\lfloor \frac{2\sqrt{e+1}}{\sqrt{\mu}} \right\rfloor \leq \frac{4(\sqrt{e+1})}{\sqrt{\mu}}.$$

We also notice that, in view of Proposition 3.11.*(ii)*, the restart condition (3.12b) is satisfied for every $k \geq \left\lfloor \frac{2}{\sqrt{\mu}} \right\rfloor$. Therefore, $n_{j-1} \leq \frac{4\sqrt{e+1}}{\sqrt{\mu}}$ implies that $n_j$, obtained from $[r_j, n_j] \leftarrow FISTA(r_{j-1}, E_f^e(n_{j-1}))$, also satisfies (3.16).

Making use of the previous discussion, we now prove claim *(ii)* by reduction to the absurd. Suppose that

$$n_j > \frac{4\sqrt{e+1}}{\sqrt{\mu}}. \tag{3.17}$$

Because of the previous discussion, an taking into account that $n_0 = 0$, inequality (3.17) could only be attained by the doubling step $n_j = 2n_{j-1}$ of the algorithm (see step 8). That is, inequality 3.17 is possible only if there is $s \in \mathbb{Z}_{>0} > 1$ such that $n_{s-1} > \frac{2\sqrt{e+1}}{\sqrt{\mu}}$ and

$$f(r_{s-1}) - f(r_s) > \frac{f(r_{s-2}) - f(r_{s-1})}{e}.$$

Since $[r_{s-1}, n_{s-1}] \leftarrow FISTA(r_{s-2}, E_f^e(n_{s-2}))$, we have that $r_{s-1}$ is obtained from $r_{s-2}$ applying $n_{s-1} > \frac{2\sqrt{e+1}}{\sqrt{\mu}}$ iterations of FISTA algorithm. However, we have from Proposition 3.11.*(iii)* that this number of iterations implies

$$f(r_{s-1}) - f(r_s) \leq f(r_{s-1}) - f^* \leq \frac{f(\mathcal{T}(r_{s-2})) - f(r_{s-1})}{e}.$$

From Proposition 2.5.*(ii)* we also have $f(\mathcal{T}(r_{s-2})) \leq f(r_{s-2})$, which when combined with the above inequality leads to

$$f(r_{s-1}) - f(r_s) \leq \frac{f(r_{s-2}) - f(r_{s-1})}{e}.$$

That is, there is no doubling step if $n_{s-1} \geq \frac{2\sqrt{e+1}}{\sqrt{\mu}}$, which proves claim *(ii)*.

We now show that there is a doubling step (i.e., step 8 is executed) at least every

$$T \doteq \left\lceil \log\left(1 + \frac{2(f(r_0) - f^*)}{\epsilon^2}\right) \right\rceil$$

iterations of the algorithm. Suppose that there is no doubling step from iteration $j = s+1$ to $j = s+T$, where $s \in \mathbb{Z}_{>0} \geq 1$. That is,

$$f(r_{j-1}) - f(r_j) \leq \frac{f(r_{j-2}) - f(r_{j-1})}{e}, \ \forall j \in [s+1, s+T].$$

From this, and claim *(i)* of the theorem, we obtain the following sequence of inequalities:

$$\begin{aligned}
\frac{1}{2}\|\mathcal{G}(r_{s+T-1})\|^2_{R^{-1}} &\le f(r_{s+T-1}) - f(r_{s+T}) \le \frac{f(r_{s+T-2}) - f(r_{s+T-1})}{e} \\
&\le \left(\frac{1}{e}\right)^T (f(r_{s-1}) - f(r_s)) \le \left(\frac{1}{e}\right)^T (f(r_{s-1}) - f^*) \\
&\le \left(\frac{1}{e}\right)^T (f(r_0) - f^*) = \left(\frac{1}{e}\right)^{\left\lceil \log\left(1+\frac{2(f(r_0)-f^*)}{\epsilon^2}\right)\right\rceil} (f(r_0) - f^*) \\
&\le \left(\frac{1}{e}\right)^{\log\left(1+\frac{2(f(r_0)-f^*)}{\epsilon^2}\right)} (f(r_0) - f^*) \\
&= \left(\frac{1}{1+\frac{2(f(r_0)-f^*)}{\epsilon^2}}\right)(f(r_0) - f^*) \le \frac{\epsilon^2}{2}.
\end{aligned}$$

We conclude that $T$ consecutive iterations without doubling step implies that the exit condition $\|\mathcal{G}(r_{s+T-1})\|_{R^{-1}} \le \epsilon$ is satisfied. Therefore, there must be at least one doubling step every $T$ iterations. This implies that there exist $j \in [s+1, s+T]$ such that

$$f(r_{j-1}) - f(r_j) > \frac{f(r_{j-2}) - f(r_{j-1})}{e},$$

which, in view of step 7, implies that $n_j = 2n_{j-1}$. Moreover, since $\{n_j\}$ is a non-decreasing sequence, we get $n_{s+T} \ge n_j = 2n_{j-1} \ge 2n_s$, $\forall s \ge 1$. That is,

$$n_s \le \frac{n_{s+T}}{2},\ \forall s \ge 1. \tag{3.18}$$

Let us rewrite $j$ as $j = m + nT$, where $0 \le m < T$ and $n \ge 0$. From the non decreasing nature of $\{n_j\}$, we have that

$$\begin{aligned}
\sum_{i=0}^{j} n_i = \sum_{i=0}^{m+nT} n_i = \sum_{i=0}^{m} n_i + \sum_{\ell=0}^{n-1}\sum_{i=1}^{T} n_{m+i+\ell T} &\le T n_m + T\sum_{\ell=1}^{n} n_{m+\ell T} \\
= T\sum_{\ell=0}^{n} n_{m+\ell T} = T\sum_{\ell=0}^{n} n_{j-\ell T}.
\end{aligned} \tag{3.19}$$

Also, from inequality (3.18), we have $n_{j-T} \le \frac{n_j}{2}$. Using this inequality in a recursive manner we obtain

$$n_{j-\ell T} \le \left(\frac{1}{2}\right)^{\ell} n_j,\ \ell = 0, \dots, n,$$

which along with (3.19) leads to

$$\sum_{i=0}^{j} n_i \le T\sum_{\ell=0}^{n}\left(\frac{1}{2}\right)^{\ell} n_j \le T\sum_{\ell=0}^{\infty}\left(\frac{1}{2}\right)^{\ell} n_j = 2Tn_j.$$

**Algorithm 8:** Gradient Based Restart FISTA

**Require:** $r_0 \in \mathbb{R}^n$, $\epsilon \in \mathbb{R}_{>0}$

1 $y_0 \leftarrow \mathcal{T}(r_0)$, $z_0 \leftarrow \mathcal{T}(r_0)$, $t_0 \leftarrow 1$, $k \leftarrow 0$ $j \leftarrow 0$, $\rho_0 \leftarrow \|\mathcal{G}(r_0)\|_{R^{-1}}$
2 **repeat**
3     $k \leftarrow k+1$
4     $z_k \leftarrow \mathcal{T}(y_{k-1})$
5     $t_k \leftarrow \dfrac{1}{2}\left(1+\sqrt{1+4t_{k-1}^2}\right)$
6     $y_k \leftarrow z_k + \dfrac{t_{k-1}-1}{t_k}(z_k - z_{k-1})$
7     **if** $\|\mathcal{G}(y_k)\|_{R^{-1}} \leq \dfrac{\rho_j}{e}$ **then**
8         $j \leftarrow j+1$
9         $r_j \leftarrow y_k$
10         $\rho_j \leftarrow \|\mathcal{G}(r_j)\|_{R^{-1}}$
11         $y_k \leftarrow \mathcal{T}(r_j)$, $z_k \leftarrow \mathcal{T}(r_j)$, $t_k \leftarrow 1$
12     **end if**
13 **until** $\rho_j \leq \epsilon$

**Output:** $z_{\text{out}} \leftarrow z_k$, $j_{\text{out}} \leftarrow j$, $k_{\text{out}} \leftarrow k$

By now making the summation up to $j = j_{\text{out}}$, claim *(iii)* now directly follows by noting that claim *(ii)* states that $n_j \leq \frac{4\sqrt{e+1}}{\sqrt{\mu}}$, $\forall j \geq 0$, and by noting that $k_{\text{out}} = \sum_{j=0}^{j_{\text{out}}} n_j$ if no doubling steps (step 8 of the algorithm) are performed and $k_{\text{out}} \leq \sum_{j=0}^{j_{\text{out}}} n_j$ otherwise. ■

### 3.2.2 A gradient based restart scheme for FISTA

This section presents a computationally cheap and simple to implement restart scheme for FISTA (Algorithm 4) that was originally presented in [6] and that exhibits linear convergence when applied to the class of optimization problems (3.9) under Assumtions 2.1 and 3.9. The scheme is based on monitoring the evolution of the composite gradient mapping $\mathcal{G}$ during the iterates of FISTA through the use of the exit condition $E_g^e : \mathbb{R}^{n_z} \to \{\text{true}, \text{false}\}$ given by

$$E_g^e(y_k) = \begin{cases} \text{true} & \text{if } \|\mathcal{G}(y_k)\|_{R^{-1}} \leq \dfrac{1}{e}\|\mathcal{G}(z)\|_{R^{-1}} \\ \text{false} & \text{otherwise,} \end{cases} \tag{3.20}$$

where $y_k$ is the iterate of FISTA computed at step 6 of Algorithm 4. That is, FISTA is restarted every time the value of the composite gradient mapping is reduced by $e$ with respect to the one corresponding to the initial condition.

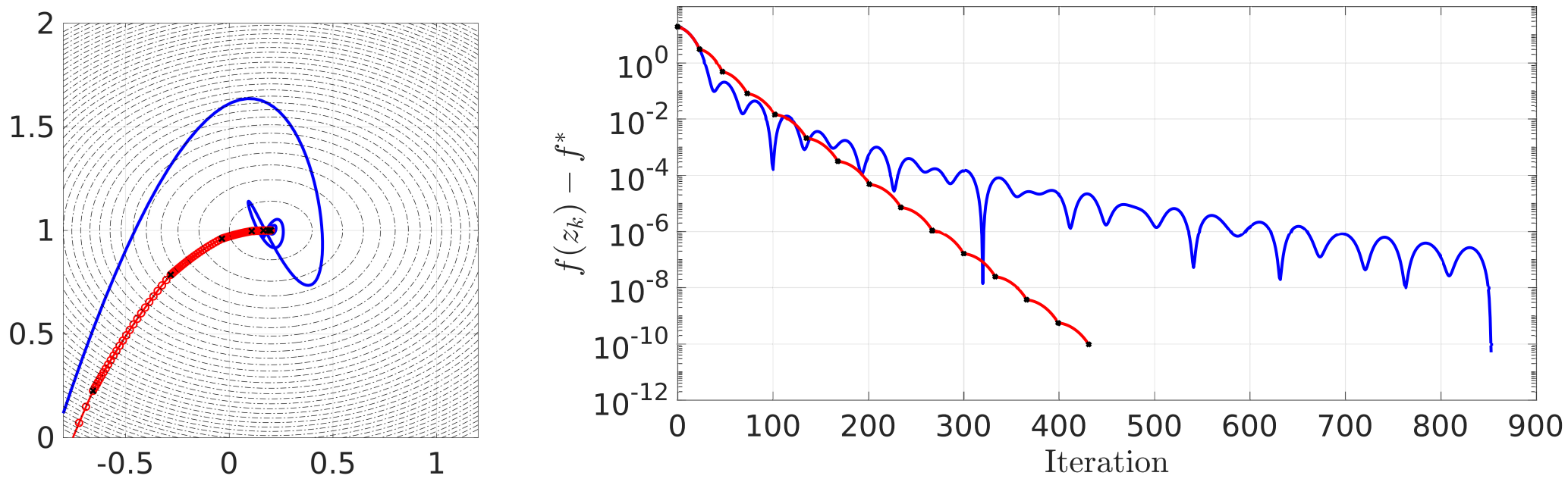

(a) Trajectory of $z_k$.

(b) Convergence to optimal objective function value.

Figure 3.6: In red, Algorithm 8 applied to Example 3.1. The non-restarted FISTA is depicted in blue. The red circles in (a) depict each of the iterates $z_k$. The black crosses depict each of the restart points $r_j$.

However, there is a minor difference with the implementation of this restart scheme that does not fit within the general restart scheme shown in Algorithm 6, which is that FISTA is not restarted using the last point $z_k$ that it generated, but instead the last $y_k$ is used. Due to this small difference, we now present the complete pseudocode of this restart scheme.

For an initial starting point $r_0 \in \mathbb{R}^{n_z}$, and an exit tolerance $\epsilon \in \mathbb{R}_{>0}$, Algorithm 8 implements FISTA with the restart scheme (3.20) for solving the class of optimization problems (3.9) under Assumtions 2.1 and 3.9. Steps 3 to 6 are the steps of the non-restarted FISTA algorithm (see Algorithm 4), whereas steps 8 to 11 implement the restart procedure and step 7 checks the restart condition (3.20). Notice that in step 9 the new restart point $r_j$ is not assigned $z_k$, as would be the case of the general restart scheme described by Algorithm 6. Instead, the current iterate $y_k$ is used, both as the new initial condition (see step 9) and as the new "initial" value of the gradient mapping (see step 10). The algorithm generates a sequence $\rho_j = \|\mathcal{G}(r_j)\|_{R^{-1}}$ and exits when $\rho_j \leq \epsilon$, i.e., when an $\epsilon$-accurate solution of problem (3.9) has been found. This solution is returned in $z_{\text{out}}$, alongside the total number of restart iterations $j_{\text{out}}$ and total number of FISTA iterations $k_{\text{out}}$.

**Example 3.17.** Figure 3.6 shows the use of Algorithm 8 to solve Example 3.1 along with the result of the non-restarted FISTA algorithm. The total number of FISTA iterations is $k_{\text{out}} = 431$ and the number of restarts is $j_{\text{out}} = 14$.

The convergence properties of this restart scheme are derived from the following proposition, which shows that, under the assumption of the quadratic functional growth condition, $\|\mathcal{G}(y)\|_{R^{-1}}$ can be used to bound the distance of $\mathcal{T}(y)$ to the optimal set $\Omega_f$ of problem (3.9). A similar result can be found in [46, Theorem 7]. However, our proposition can be applied to any non-necessarily feasible point $y \in \mathbb{R}^{n_z}$, i.e., $y$ need not belong to $\mathcal{Z}$. Analogous results can also be found in other publications, such as [56, 57].

**Proposition 3.18.** Consider problem (3.9) under Assumtions 2.1 and 3.9 and let $\Omega_f$ be the optimal set of problem (3.9). Then,

$$\|\mathcal{T}_R^{f,\mathcal{Z}}(y) - \pi_{\Omega_f}^R(\mathcal{T}_R^{f,\mathcal{Z}}(y))\|_R \leq \frac{1+\sqrt{1+\mu}}{\mu}\|\mathcal{G}_R^{f,\mathcal{Z}}(y)\|_{R^{-1}}, \ \ \forall y \in \mathbb{R}^{n_z}.$$

***Proof:*** To simplify the notation, let us define $y^+ \doteq \mathcal{T}_R^{f,\mathcal{Z}}(y)$, $\bar{y} \doteq \pi_{\Omega_f}^R(y)$ and $\bar{y}^+ \doteq \pi_{\Omega_f}^R(\mathcal{T}_R^{f,\mathcal{Z}}(y))$. From equation (2.7a) in Proposition 2.5 we have that

$$f(y^+) - f(z) \leq \langle \mathcal{G}(y), y^+ - z\rangle + \frac{1}{2}\|\mathcal{G}(y)\|_{R^{-1}}^2, \ \forall y \in \mathbb{R}^{n_z}, \ \forall z \in \mathcal{Z}.$$

Particularizing this inequality to $z = \bar{y}^+$ we obtain

$$f(y^+) - f(\bar{y}^+) \leq \langle \mathcal{G}(y), y^+ - \bar{y}^+\rangle + \frac{1}{2}\|\mathcal{G}(y)\|_{R^{-1}}^2.$$

Since $\bar{y}^+ \in \Omega_f$ we have $f(\bar{y}^+) = f^*$. Therefore, we infer from Assumtion 3.9 that

$$\begin{aligned}
\frac{\mu}{2}\|y^+ - \bar{y}^+\|_R^2 &\leq f(y^+) - f^* = f(y^+) - f(\bar{y}^+) \\
&\leq \langle \mathcal{G}(y), y^+ - \bar{y}^+\rangle + \frac{1}{2}\|\mathcal{G}(y)\|_{R^{-1}}^2 \\
&= \langle R^{-\frac{1}{2}}\mathcal{G}(y), R^{\frac{1}{2}}(y^+ - \bar{y}^+)\rangle + \frac{1}{2}\|\mathcal{G}(y)\|_{R^{-1}}^2 \\
&\overset{(*)}{\leq} \|\mathcal{G}(y)\|_{R^{-1}}\|y^+ - \bar{y}^+\|_R + \frac{1}{2}\|\mathcal{G}(y)\|_{R^{-1}}^2,
\end{aligned}$$

where in $(*)$ we are making use of the Cauchy-Schwarz inequality. Then, adding $\frac{1}{2}\|y^+ - \bar{y}^+\|_R^2$ to both terms of the last inequality, we obtain

$$\frac{1}{2}(\mu+1)\|y^+ - \bar{y}^+\|_R^2 \leq \frac{1}{2}\left(\|\mathcal{G}(y)\|_{R^{-1}} + \|\bar{y}^+ - \bar{y}\|_R\right)^2.$$

That is,

$$\begin{aligned}
\sqrt{\mu+1}\|y^+ - \bar{y}^+\|_R &\leq \|\mathcal{G}(y)\|_{R^{-1}} + \|\bar{y}^+ - \bar{y}\|_R, \\
(\sqrt{\mu+1} - 1)\|y^+ - \bar{y}^+\|_R &\leq \|\mathcal{G}(y)\|_{R^{-1}},
\end{aligned}$$

from where we conclude that

$$\|y^+ - \bar{y}^+\|_R \leq \frac{1}{\sqrt{1+\mu}-1}\|\mathcal{G}(y)\|_{R^{-1}} = \frac{1+\sqrt{\mu+1}}{\mu}\|\mathcal{G}(y)\|_{R^{-1}}.$$

■

The following theorem gathers the main convergence properties of the restart scheme presented in Algorithm 8. Claim *(i)* states the sublinear convergence of $\{\|\mathcal{G}(y_k)\|_{R^{-1}}\}_{\geq 0}$ under the quadratic functional growth assumption without taking into account the restart procedure; claim *(ii)* provides an upper bound to the number of iterations required to satisfy the restart condition (3.20) (step 7 of Algorithm 8); and claim *(iii)* states the linear convergence of the composite gradient mapping under the restart procedure. Finally, claim *(iii)* provides the iteration complexity.

**Remark 3.19.** *We note that we use $\|\mathcal{G}(r_j)\|_{R^{-1}} \leq \epsilon$ as the exit condition of Algorithm 8 to simplify the proofs and expressions of Theorem 3.20. In a practical setting, however, the objective is to find an $\epsilon$-accurate solution in as fewer iterations as possible. Therefore, the exit condition $\|\mathcal{G}(y_{k-1})\|_{R^{-1}} \leq \epsilon$ may be used instead (see Remark 2.10).*

**Theorem 3.20** (Convergence properties of Algorithm 8)**.** Consider problem (3.9) under Assumtions 2.1 and 3.9. Then, the sequences generated by Algorithm 8 satisfy:

*(i)* $\|\mathcal{G}^{f,\mathcal{Z}}_R(y_k)\|_{R^{-1}} \leq \dfrac{4(1+\sqrt{\mu+1})}{\mu k}\|\mathcal{G}^{f,\mathcal{Z}}_R(r_0)\|_{R^{-1}},\ \forall k \geq 1.$

*(ii)* $\|\mathcal{G}^{f,\mathcal{Z}}_R(y_k)\|_{R^{-1}} \leq \dfrac{1}{e}\|\mathcal{G}^{f,\mathcal{Z}}_R(r_0)\|_{R^{-1}},\ \forall k \geq \dfrac{4e}{\mu}(1+\sqrt{\mu+1}).$

*(iii)* Let $a(\mu) \doteq \max\left\{1,\ \dfrac{4e(1+\sqrt{\mu+1})}{\mu}\right\}$. Then,

$$\min_{i=0,\dots,k}\{\|\mathcal{G}(y_i)\|_{R^{-1}},\ \|\mathcal{G}(r_0)\|_{R^{-1}}\} \leq e\left(\frac{1}{e}\right)^{\left\lfloor \frac{k}{a(\mu)}\right\rfloor}\|\mathcal{G}(r_0)\|_{R^{-1}}.$$

*(iv)* The total number of iterations $k_{\text{out}}$ required to exit the algorithm is upper bounded by

$$k_{\text{out}} \leq \max\left\{1,\ \frac{4e(1+\sqrt{\mu+1})}{\mu}\right\}\left\lceil \log\left(\frac{\|\mathcal{G}(r_0)\|_{R^{-1}}}{\epsilon}\right)\right\rceil.$$

***Proof:*** Notice that, from step 1 of Algorithm 8, we have that $z_0 = y_0 = \mathcal{T}(r_0)$, which along with Proposition 3.18 implies that

$$\|z_0 - \pi^R_{\Omega_f}(z_0)\|_R \leq \frac{1+\sqrt{\mu+1}}{\mu}\|\mathcal{G}(r_0)\|_{R^{-1}}.$$

From Proposition 2.12.*(ii)* we also have that

$$\|\mathcal{G}(y_k)\|_{R^{-1}} \leq \frac{4}{k+2}\|z_0 - \pi^R_{\Omega_f}(z_0)\|_R,\ \forall k \geq 0.$$

Therefore, we obtain

$$\begin{aligned}\|\mathcal{G}(y_k)\|_{R^{-1}} &\leq \frac{4}{k+2}\|z_0 - \pi^R_{\Omega_f}(z_0)\|_R \leq \frac{4(1+\sqrt{\mu+1})}{(k+2)\mu}\|\mathcal{G}(r_0))\|_{R^{-1}} \\ &\leq \frac{4(1+\sqrt{\mu+1})}{k\mu}\|\mathcal{G}(r_0)\|_{R^{-1}},\ \forall k \geq 1,\end{aligned}$$

which proves claim *(i)*. Let $m \in \mathbb{R}$ be the scalar that satisfies

$$\frac{4(1+\sqrt{\mu+1})}{\mu m} = \frac{1}{e},$$

from where we obtain that $m = \frac{4e}{\mu}(1+\sqrt{\mu+1})$. This means that

$$\frac{4(1+\sqrt{\mu+1})}{\mu k} \leq \frac{1}{e},\ \forall k \geq m.$$

Therefore, a sufficient condition for $\|\mathcal{G}(y_k)\|_{R^{-1}} \leq \frac{1}{e}\|\mathcal{G}(r_0)\|_{R^{-1}}$ is

$$k \geq \frac{4e}{\mu}(1+\sqrt{\mu+1}),$$

which proves claim *(ii)*. Notice that the restart condition $\|\mathcal{G}(y_k)\|_{R^{-1}} \leq \frac{\rho_j}{e}$ (see step 7 of the algorithm) implies that

$$\|\mathcal{G}(r_j)\|_{R^{-1}} \leq \left(\frac{1}{e}\right)^j \|\mathcal{G}(r_0)\|_{R^{-1}}. \tag{3.21}$$

Additionally, in view of claim *(ii)*, we have that each restart occurs in a number of iterations no larger than $a(\mu)$. Thus, the number of restarts $j$ is no smaller than $\left\lfloor \frac{k}{a(\mu)} \right\rfloor$. Therefore, we obtain from (3.21) that

$$\min_{i=0,\dots,k}\{\|\mathcal{G}(y_i)\|_{R^{-1}}, \|\mathcal{G}(r_0)\|_{R^{-1}}\} \leq \left(\frac{1}{e}\right)^{\left\lfloor \frac{k}{a(\mu)} \right\rfloor} \|\mathcal{G}(r_0)\|_{R^{-1}} \leq e\left(\frac{1}{e}\right)^{\frac{k}{a(\mu)}} \|\mathcal{G}(r_0)\|_{R^{-1}},$$

which proves claim *(iii)*. Clearly, if $\|\mathcal{G}(r_0)\|_{R^{-1}} \leq \epsilon$ the algorithm exits in one iteration. Otherwise, we obtain from (3.21) that the number of restarts $j$ required to attain the desired accuracy $\epsilon$ is no larger than

$$\left\lceil \log\left(\frac{\|\mathcal{G}(r_0)\|_{R^{-1}}}{\epsilon}\right)\right\rceil.$$

Since, in view of claim *(ii)*, a restart occurs in a number of iterations no larger than $a(\mu)$, we conclude that the total number of iterations is upper bounded by

$$k_{\text{out}} \leq a(\mu)\left\lceil \log\left(\frac{\|\mathcal{G}(r_0)\|_{R^{-1}}}{\epsilon}\right)\right\rceil,$$

thus proving claim *(iv)*. ■

### 3.2.3 Restart scheme for accelerated first order methods

We remind the reader that $\|\cdot\|$ represents any vector norm and $\|\cdot\|_*$ represents its dual norm (see Definition N.1).

This section presents a restart scheme for AFOMs that can be viewed as an extension of the scheme presented in Section 3.2.1, in that it is based on a restart condition that checks the evolution of the objective function value and on the use of a minimum number of iterations that is adapted between restarts. However, the scheme presented here is applicable to a wider class of AFOMs and relies on a less restrictive interpretation of the quadratic functional growth condition (see Definition 3.3). The restart scheme presented in this section is currently under review, but a preprint of the article can be found in [7].

Let us consider a convex optimization problem

$$f^* = \min_{z \in \mathbb{R}^{n_z}} f(z) \tag{3.22}$$

that we assume is solvable and is subject to the following assumption, which is a relaxation of the quadratic functional growth condition (Definition 3.3).

**Assumption 3.21.** The function $f \in \Gamma(\mathbb{R}^{n_z})$ satisfies, for every $\rho \in \mathbb{R}_{\geq 0}$ a quadratic functional growth condition of the form

$$f(z) - f^* \geq \frac{\mu_\rho}{2} \|z - \pi_{\Omega_f}(z)\|^2, \ \forall z \in V_f(\rho),$$

for some $\mu_\rho \in \mathbb{R}_{>0}$, where $V_f(\rho) = \{\, z \in \mathbb{R}^{n_z} \ : \ f(z) - f^* \leq \rho \,\}$ is a level set of problem (3.22) with respect to $f^*$.

Additionally, let us consider a fixed point algorithm $\mathcal{A}$ applied to problem (3.22) starting from an initial condition $z_0 \in \text{dom}(f)$, i.e., given an initial condition $z_0$, algorithm $\mathcal{A}$ generates a sequence $\{z_k\}_{\geq 0}$ such that $\lim_{k\to\infty} f(z_k) = f^*$. We use the following notation to refer to the iterates provided by algorithm $\mathcal{A}$.

**Definition 3.22.** Let $\mathcal{A}$ be an iterative algorithm applied to solve problem (3.22) using as initial condition $z_0 \in \mathbb{R}^{n_z}$. Given $k \in \mathbb{Z}_{\geq 0}$, we denote by $\mathcal{A}_k(z_0) \in \mathbb{R}^{n_z}$ the vector corresponding to iteration $k$ of the algorithm.

**Remark 3.23.** *We note that the notation given in Definition 3.22 bares a close resemblance to the one presented in Chapter 3 for the output of a restart scheme* $[z_{out}, k_{out}] \leftarrow \mathcal{A}(z_0, E_c)$. *However, they are not equivalent, since* $\mathcal{A}(z_0, E_c)$ *is a notation to the call of the algorithm with a certain exit condition, whereas* $\mathcal{A}_k(z_0)$ *is the vector in* $\mathbb{R}^{n_z}$ *corresponding to iterate* $k$ *of the algorithm. Even so, the two notations are strongly related. In fact,*

$$[\mathcal{A}_{k_m}(z_0), k_m] \leftarrow \mathcal{A}(z_0, E_{k_m}),$$

*where, for* $k_m \in \mathbb{R}_{>0}$, $E_{k_m}$ *is the exit condition given by* $E_{k_m}(k) = true \iff k = k_m$.

The following assumption characterizes the class of fixed point algorithms we consider in this section. We note that $g$ is to be taken as the gradient operator used by the AFOM $\mathcal{A}$, e.g., the composite gradient mapping (see Definition 2.3), the proximal operator, the projected gradient, the gradient, the subgradient, etc.

**Assumption 3.24.** The fixed point iterative algorithm $\mathcal{A}$, applied to solve (3.22) under Assumtion 3.21 satisfies, for every $z_0 \in \text{dom}(f)$:

*(i)* $f(\mathcal{A}_k(z_0)) - f^* \leq \dfrac{a_f}{(k+1)^2}\|z_0 - \pi_{\Omega_f}(z_0)\|^2, \ \forall k \geq 1,$

*(ii)* $f(\mathcal{A}_1(z_0)) \leq f(z_0) - \dfrac{1}{2L_f}\|g(z_0)\|_*^2,$

for some $a_f \in \mathbb{R}_{>0}$, $L_f \in \mathbb{R}_{>0}$ and where $g : \mathbb{R}^{n_z} \to \mathbb{R}^{n_z}$ is a gradient operator satisfying $g(z) = 0 \iff z \in \Omega_f$, where $\Omega_f$ is the optimal set of (3.22).

The parameters $a_f$ and $L_f$ shown in Assumtion 3.24 will take different values depending on the AFOM being used and on the characteristics of problem (3.22). For instance, if $f = h + \Psi$ is the composition of a smooth function $h$ and a non-smooth function $\Psi$, then $L_f$ is the Lipschitz constant of $h$ (see Sections 3.2.1 and 3.2.2). The value of $a_f$ is equal to $2L_f$ in the FISTA and MFISTA algorithms, but it will take other values in different AFOMs or if certain strategies, such as backtracking strategies, are used (see [18, §10] or [27]).

The two conditions listed in Assumtion 3.24 are satisfied by most AFOMs. The sublinear convergence stated in Assumtion 3.24.*(i)* is one of the main properties of AFOMs. For instance, we show the satisfaction of this property for FISTA (see Proposition 2.12.*(i)*) and MFISTA (see Proposition 2.13). For other AFOMs satisfying this condition, we refer the reader to [18, 19, 58], to [30] for an accelerated variant of ADMM, and to [59] for an accelerated version of the *alternating minimization algorithm* with a particular focus on its certification for model predictive control.

The condition stated in Assumtion 3.24.*(ii)*, on the other hand, is not necessarily satisfied by most AFOMs shown in the literature. However, any AFOM can be easily modified to satisfy this condition. For instance, in view of Proposition 2.5.*(ii)*, all that an AFOM based on the composite gradient mapping needs to include to satisfy Assumtion 3.24.*(ii)* is to start by obtaining $z_1 = \mathcal{T}(z_0)$, which is a simple modification. We note that the FISTA and MFISTA algorithms that we present in Algorithms 4 and 5, respectively, already include this initial step.

In short, by considering problem (3.22) under to Assumtion 3.21 and a fixed point algorithm $\mathcal{A}$ under Assumtion 3.24, we are encompassing a wide range of optimization problems and AFOMs.

We now present a proposition regarding the iterates of $\mathcal{A}$ under Assumtions 3.21 and 3.24 that will serve as a basis for the development of the convergence analysis of the restart scheme presented in this section. An equivalent result can be found in [46, §5.2.2]. The proposition makes use of a scalar $n_\rho$, which, due to

**Algorithm 9:** Delayed exit condition on $\mathcal{A}$

**Prototype:** $[z_m, m] \leftarrow \mathcal{A}_d(z_0, n)$
**Require:** $z_0 \in \mathrm{dom}(f)$, $n \in \mathbb{R}$
1 $k \leftarrow 0$
2 **repeat**
3 $\quad k \leftarrow k + 1$
4 $\quad z_k \leftarrow \begin{cases} \mathcal{A}_k(z_0) & \text{if } f(\mathcal{A}_k(z_0)) \leq f(z_{k-1}) \\ z_{k-1} & \text{otherwise} \end{cases}$
5 $\quad \ell \leftarrow \left\lfloor \frac{k}{2} \right\rfloor$
6 **until** $k \geq n$ *and* $f(z_\ell) - f(z_k) \leq \dfrac{1}{3}(f(z_0) - f(z_\ell))$
**Output:** $z_m \leftarrow z_k$, $m \leftarrow k$

its prevalence and importance in future developments, we characterize separately in the following definition.

**Definition 3.25.** We define $n_\rho \in \mathbb{R}_{>0}$ as the scalar satisfying

$$n_\rho \doteq \max\left\{\frac{1}{2}, \sqrt{\frac{2a_f}{\mu_\rho}}\right\},$$

where $\mu_\rho$ and $a_f$ are given in Assumtions 3.21 and 3.24, respectively.

**Proposition 3.26.** Let Assumtion 3.24 hold. Then, for every $z_0 \in V_f(\rho)$,

$$f(\mathcal{A}_k(z_0)) - f^* \leq \left(\frac{n_\rho}{k+1}\right)^2 (f(z_0) - f^*), \; \forall k \geq 1.$$

***Proof:*** Denote $f_0 \doteq f(z_0)$, $f_k \doteq f(\mathcal{A}_k(z_0))$, $\forall k > 1$. Then,

$$f_k - f^* \leq \frac{a_f}{(k+1)^2}\|z_0 - \pi_{\Omega_f}(z_0)\|^2 \leq \frac{2a_f}{\mu_\rho(k+1)^2}(f_0 - f^*) \leq \frac{n_\rho^2}{(k+1)^2}(f_0 - f^*).$$

■

Before we present the restart scheme, we must introduce Algorithm 9, which implements a delayed exit condition on algorithm $\mathcal{A}$. For a given initial condition $z_0 \in \mathrm{dom}(f)$ and a scalar $n \in \mathbb{R}_{\geq 0}$, Algorithm 9 uses a given AFOM $\mathcal{A}$ to generate a sequence $\{z_k\}_{\geq 0}$ that satisfies (see step 4)

$$f(z_k) = \min\{f(z_{k-1}), f(\mathcal{A}_k(z_0))\}, \; \forall k \geq 1.$$

Therefore,

$$f(z_k) = \min_{i=0,\dots,k} f(\mathcal{A}_i(z_0)). \tag{3.23}$$

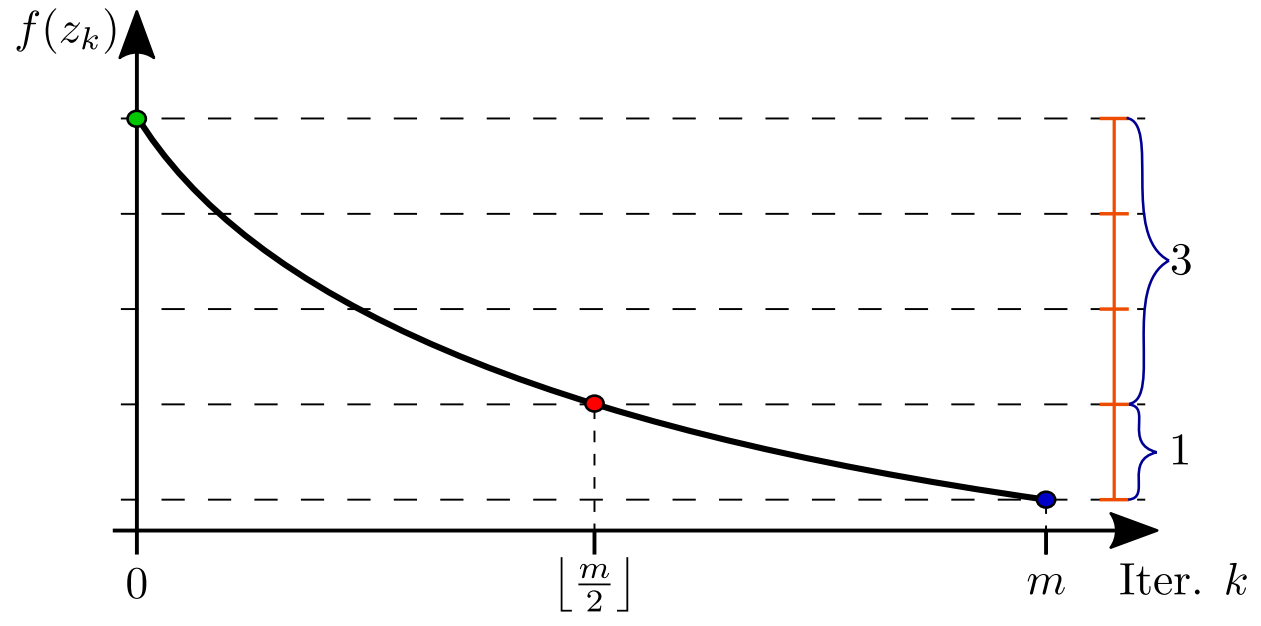

Figure 3.7: Satisfaction of the delayed exit condition (3.24).

The algorithm terminates after $k \geq n$ iterations if the following exit condition is satisfied (see step 6):

$$f(z_\ell) - f(z_k) \leq \frac{1}{3}\left(f(z_0) - f(z_\ell)\right), \tag{3.24}$$

where $\ell = \left\lfloor \frac{k}{2} \right\rfloor$. The outputs of the algorithm are the current iterate $z_m = z_k$, and the total number of iterations $m \geq n$ required to satisfy the exit condition (3.24).

Intuitively, as illustrated in Figure 3.7, exit condition (3.24) detects a degradation in the performance of the iterations of $\mathcal{A}$. Notice that at iteration $m$, the reduction corresponding to the last half of the iterations (from $\left\lfloor \frac{m}{2} \right\rfloor$ to $m$) is no larger than one third of the reduction achieved in the first half of the iterations (from 0 to $\left\lfloor \frac{m}{2} \right\rfloor$).

**Remark 3.27.** *We note that algorithm $\mathcal{A}$ does not need to be started from $z_0$ every time it is called in step 4 of Algorithm 9. Instead, the algorithm is started at $z_0$ the first time it is called, and then one more iteration is performed at each iteration $k$ of Algorithm 9. Then, either the current value of $\mathcal{A}_k(z_0)$ or the previous value of $z_{k-1}$ is stored in $z_k$.*

**Remark 3.28.** *Step 4 of Algorithm 9 generates a sequence $\{z_k\}_{\geq 0}$ whose associated objective function values are non-increasing. Therefore, if $\mathcal{A}$ is a monotone AFOM, e.g., MFISTA (see Section 2.2.2), then the assignment $z_k \leftarrow z_{k-1}$ will never be performed. In this case, the same sequence $\{z_k\}_{\geq 0}$ can be obtained by simply calling the monotone algorithm $\mathcal{A}$ with the following exit condition $E_f^3 : \mathbb{R}^{n_z} \times \mathbb{R}^{n_z} \times \mathbb{Z}_{>0} \to \{true, false\}$:*

$$E_f^3(z_k, z_l, k; z_0, n) = true \iff \begin{cases} f(z_l) - f(z_k) \leq \dfrac{1}{3}(f(z_0) - f(z_l)) \\ k \geq n. \end{cases}$$

The following property provides three conditions satisfied by the inputs and outputs of Algorithm 9. This result is instrumental to prove the convergence results of the restart scheme we present further ahead. Claim *(i)* is similar to Assumtion 3.24.*(ii)*, in that it states that the output of the algorithm provides a net improvement with respect to the input (unless $z_0 = z^*$, in which case, obviously, there is no net gain). Claim *(ii)* states the convergence rate of the algorithm. Finally, claim *(iii)* provides an upper bound to the output $m$, provided that $n$ is no larger than $\lceil 4n_\rho \rceil$.

**Proposition 3.29.** Let Assumtion 3.24 hold. Then, given an initial condition $z_0 \in V_f(\rho) \subseteq \text{dom}(f)$ and a scalar $n \in \mathbb{R}$, the output $[z_m, m]$ of Algorithm 9 satisfies

*(i)* $f(z_m) \leq f(z_0) - \dfrac{1}{2L_f}\|g(z_0)\|_*^2$,

*(ii)* $f(z_m) - f^* \leq \left(\dfrac{n_\rho}{m+1}\right)^2 (f(z_0) - f^*)$,

*(iii)* $n \in (0, \lceil 4n_\rho \rceil] \implies m \in [n, \lceil 4n_\rho \rceil]$.

***Proof:*** From (3.23) and Assumtion 3.24.*(ii)* we have

$$f(z_m) = \min_{i=0,\dots,m} f(\mathcal{A}_i(z_0)) \leq f(\mathcal{A}_1(z_0)) \leq f(z_0) - \frac{1}{2L_f}\|g(z_0)\|_*^2,$$

which proves claim *(i)*. In view of (3.23) and Proposition 3.26 we have that

$$f(z_k) - f^* \overset{(3.23)}{\leq} f(\mathcal{A}_k(z_0)) - f^* \leq \left(\frac{n_\rho}{k+1}\right)^2 (f(z_0) - f^*), \ \forall k \in \mathbb{Z}_1^m. \tag{3.25}$$

Claim *(ii)* immediately follows by taking $k = m$. The inequality $n \leq m$ trivially follows from step 6 of Algorithm 9. Therefore, in order to prove claim *(iii)*, it suffices to show that inequality (3.24) is satisfied for $\hat{k} = \lceil 4n_\rho \rceil$ and $\hat{\ell} = \left\lfloor \frac{\lceil 4n_\rho \rceil}{2} \right\rfloor \geq \lfloor 2n_\rho \rfloor \geq 1$ (where the last inequality follows from Definition 3.25 which states that $n_\rho \geq 1/2$). Indeed,

$$f(z_{\hat{\ell}}) - f^* \overset{(3.25)}{\leq} \left(\frac{n_\rho}{\hat{\ell}+1}\right)^2 (f(z_0) - f^*) \leq \left(\frac{n_\rho}{2n_\rho}\right)^2 (f(z_0) - f^*) = \frac{1}{4}(f(z_0) - f^*),$$

which implies that

$$f(z_{\hat{\ell}}) \leq \frac{1}{4} f(z_0) + \frac{3}{4} f^* \leq \frac{1}{4} f(z_0) + \frac{3}{4} f(z_{\hat{k}}).$$

Thus,

$$f(z_{\hat{\ell}}) - f(z_{\hat{k}}) \leq \frac{1}{4}(f(z_0) - f(z_{\hat{k}})) = \frac{1}{4}(f(z_0) - f(z_{\hat{\ell}})) + \frac{1}{4}(f(z_{\hat{\ell}}) - f(z_{\hat{k}})),$$

from where we conclude that $f(z_{\hat{\ell}}) - f(z_{\hat{k}}) \leq \dfrac{1}{3}(f(z_0) - f(z_{\hat{\ell}}))$. ■

**Algorithm 10:** Restart scheme for AFOMs based on $\mathcal{A}_d$

**Require:** $r_0 \in \text{dom}(f)$, $\epsilon \in \mathbb{R}_{>0}$

1 $m_0 \leftarrow 1$, $m_{-1} \leftarrow 1$, $j \leftarrow -1$, $k \leftarrow 0$

2 **repeat**

3 $\quad j \leftarrow j + 1$

4 $\quad s_j \leftarrow \begin{cases} \sqrt{\dfrac{f(r_{j-1}) - f(r_j)}{f(r_{j-2}) - f(r_j)}} & \text{if } j \geq 2 \\ 0 & \text{otherwise} \end{cases}$

5 $\quad n_j \leftarrow \max\{m_j, 4 s_j m_{j-1}\}$

6 $\quad [r_{j+1}, m_{j+1}] \leftarrow \mathcal{A}_d(r_j, n_j)$

7 $\quad k \leftarrow k + m_{j+1}$

8 **until** $f(r_j) - f(r_{j+1}) \leq \epsilon$

**Output:** $r_{out} \leftarrow r_{j+1}$, $j_{out} \leftarrow j$, $k_{\text{out}} \leftarrow k$

Algorithm 10 shows the proposed restart scheme. Starting from an initial condition $r_0 \in \text{dom}(f)$ and given a AFOM $\mathcal{A}$ satisfying Assumtion 3.24, Algorithm 10 makes successive calls to $\mathcal{A}_d$ (Algorithm 9), with a minimum number of iterations $n_j$ that is adapted at each restart iteration $j$ to take into account the evolution of the previous objective function values. The algorithm generates a sequence $\{r_j\}_{\geq 0}$ that converges to an optimal solution $z^*$ of problem (3.22) under Assumtion 3.21 as $j \to +\infty$. The algorithm returns an $\epsilon$-accurate solution $r_{\text{out}}$ in terms of the exit condition shown in step 8 (see Remark 3.32 for some discussion on this exit condition); the total number of restart iterations $j_{\text{out}}$; and the total number of iterations $k_{\text{out}}$ performed by algorithm $\mathcal{A}$.

The idea behind the restart scheme proposed in Algorithm 10 bares some similarities to the one presented in Section 3.2.1 (see Algorithm 7), in that it makes use of a minimum number of iterations $n_j$, whose value at each iteration $j$ is determined by the evolution of the past objective function values. The similarity is fairly obvious if step 4 of Algorithm 10 is compared with step 7 of Algorithm 7. In both algorithms, the minimum number of iterations is increased if the fraction

$$\frac{f(r_{j-1}) - f(r_j)}{f(r_{j-2}) - f(r_{j-1})} \tag{3.26}$$

is larger that a certain amount. In this case, however, the reason behind using (3.26) as a measure for determining the (possible) increase of $n_j$ does not have such an intuitive explanation as the one discussed in Section 3.2.1.

**Example 3.30.** Figure 3.8 shows the use of Algorithm 10 to solve Example 3.1 along with the result of the non-restarted FISTA algorithm. We take $\mathcal{A}$ as MFISTA, which can be implemented as described in Remark 3.28. The total number of MFISTA iterations is $k_{\text{out}} = 239$ and the number of restarts is $j_{\text{out}} = 5$.

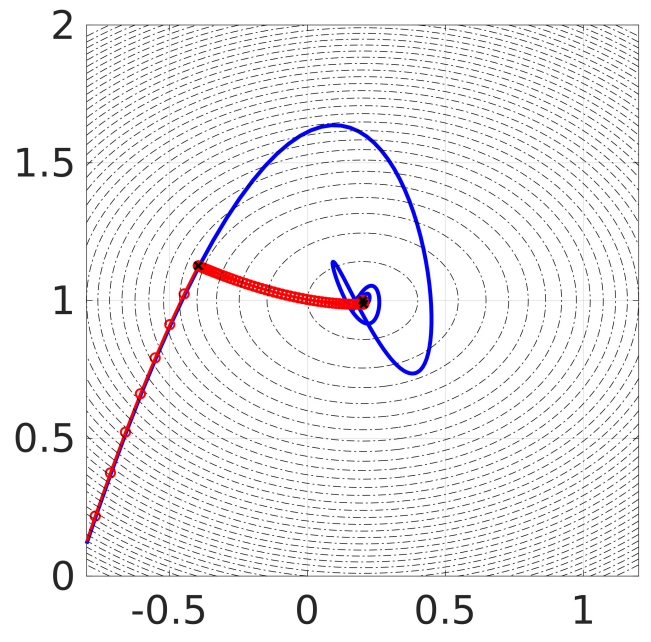

(a) Trajectory of $z_k$.

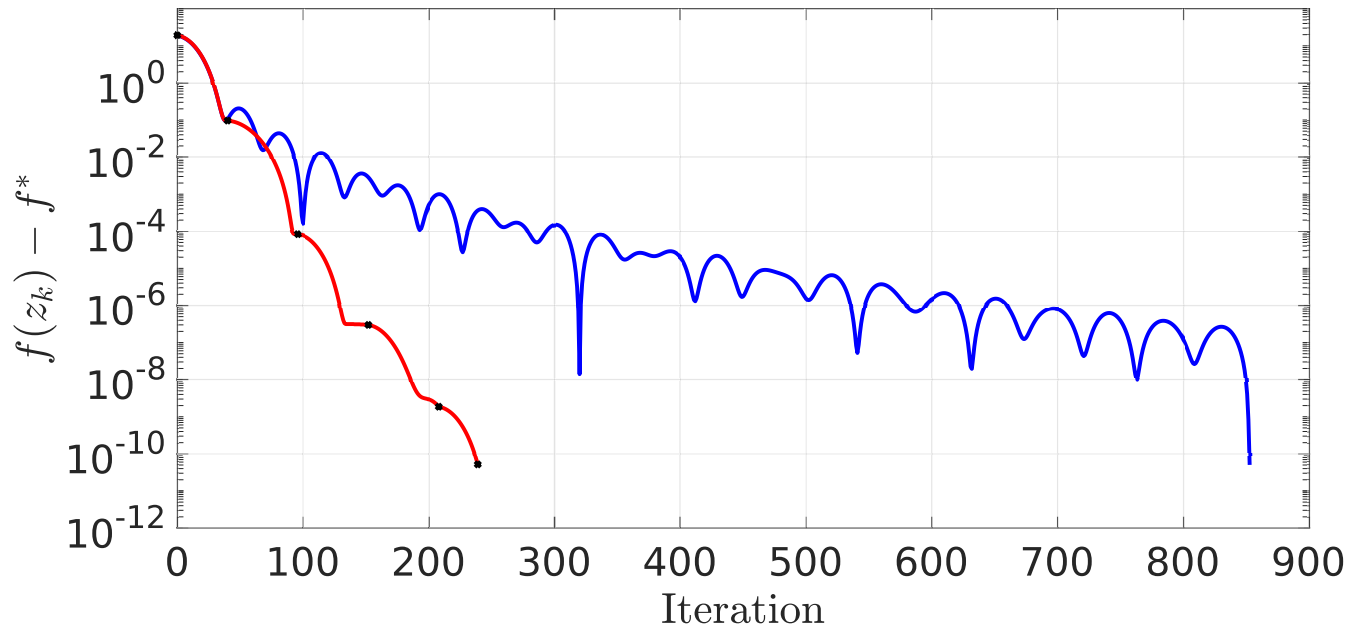

(b) Convergence to optimal objective function value.

Figure 3.8: In red, Algorithm 10 applied to Example 3.1 using MFISTA. The non-restarted FISTA is depicted in blue. The red circles in (a) depict each of the iterates $z_k$. The black crosses depict each of the restart points $r_j$.

The following proposition provides some key results on the iterates of Algorithm 10.

**Proposition 3.31.** Let Assumtion 3.24 hold and consider Algorithm 10 for a given initial condition $r_0 \in V_f(\rho) \subseteq \text{dom}(f)$ and an exit tolerance $\epsilon \in \mathbb{R}_{>0}$. Then:

*(i)* Proposition 3.29 can be applied to the iterates of Algorithm 10, i.e., taking $z_0 \equiv r_j$, $n \equiv n_j$, $z_m \equiv r_{j+1}$ and $m \equiv m_{j+1}$.

*(ii)* The sequence $\{m_j\}_{\geq 0}$ is non-decreasing. In particular,

$$m_j \leq n_j \leq m_{j+1}, \; \forall j \in \mathbb{Z}_0^{j_{out}}. \tag{3.27}$$

*(iii)* The sequence $\{s_j\}_{\geq 0}$ satisfies $s_j \in (0, 1]$, $\forall j \in \mathbb{Z}_2^{j_{out}}$.

***Proof:*** Since $r_0 \in V_f(\rho)$ for some $\rho \in \mathbb{R}_{>0}$, and each $r_j$, for $j \geq 1$, is obtained from a call to Algorithm 9 (see step 6), we have in view of Proposition 3.29.*(i)* that the iterates $z_k$ also satisfy $z_j \in V_f(\rho)$, $\forall j \in \mathbb{Z}_0^{j_{\text{out}}}$, for the same value of $\rho$. Therefore, Proposition 3.29 can be applied to each call to $\mathcal{A}_d$, thus proving claim *(i)*. That is, for every $j \in \mathbb{Z}_{\geq 0}$, the iterates of Algorithm 10 satisfy

$$f(r_{j+1}) \leq f(r_j) - \frac{1}{2L_f}\|g(r_j)\|_*^2, \tag{3.28a}$$

$$f(r_{j+1}) - f^* \leq \left(\frac{n_\rho}{m_{j+1}+1}\right)^2 (f(r_j) - f^*), \tag{3.28b}$$

$$n_j \in (0, \lceil 4n_\rho \rceil] \Rightarrow m_{j+1} \in [n_j, \lceil 4n_\rho \rceil]. \tag{3.28c}$$

Next, due to step 5 we have $m_j \leq n_j$, $j \in \mathbb{Z}_0^{j_{\text{out}}}$. Moreover, from (3.28c), we have that $n_j \leq m_{j+1}$, $\forall j \in \mathbb{Z}_0^{j_{out}}$, which proves claim *(ii)*.

Finally, from the exit condition of the algorithm (step 8), we have

$$f(r_{j-1}) - f(r_j) > \epsilon, \ \forall j \in \mathbb{Z}_1^{j_{out}}. \tag{3.29}$$

Additionally, from (3.28a) we have $f(r_{j-2}) \geq f(r_{j-1})$, $\forall j \in \mathbb{Z}_2^{j_{out}}$. Thus,

$$f(r_{j-2}) - f(r_j) \geq f(r_{j-1}) - f(r_j) \overset{(3.29)}{>} \epsilon > 0, \ \forall j \in \mathbb{Z}_2^{j_{out}}.$$

Therefore, from step 4, and taking $j \geq 2$, we have

$$0 < s_j = \sqrt{\frac{f(r_{j-1}) - f(r_j)}{f(r_{j-2}) - f(r_j)}} \leq 1, \ \forall j \in \mathbb{Z}_2^{j_{out}},$$

which proves claim *(iii)*. ■

**Remark 3.32.** *From Proposition 3.31.(i), we have that we can rearrange Proposition 3.29.(i) to read as*

$$\|g(r_j)\|_*^2 \leq 2L_f(f(r_j) - f(r_{j+1})).$$

*Therefore, the exit condition $f(r_j) - f(r_{j+1}) \leq \epsilon$ implies $\|g(r_j)\|_*^2 \leq 2L_f\epsilon$. Since, as per Assumtion 3.24, $g(r_j)$ serves to characterize the optimality of $r_j$, we conclude that the exit condition of Algorithm 10 also serves to characterize the optimality of $r_{j+1}$. This means that the exit condition of Algorithm 10 could be replaced by $\|g(r_j)\|_* \leq \tilde{\epsilon}$, where $\tilde{\epsilon} \in \mathbb{R}_{>0}$.*

The following theorem presents the main convergence properties of Algorithm 10. In particular, claim *(iii)* provides its iteration complexity, which we note is of the same order than the optimal iteration complexity (3.2) that can be obtained using AFOMs applied to the class of optimization problems (3.22) under Assumtion 3.21.

**Theorem 3.33.** Let Assumtion 3.24 hold and consider Algorithm 10 for a given initial condition $r_0 \in V_f(\rho) \subseteq \text{dom}(f)$ and an exit tolerance $\epsilon \in \mathbb{R}_{>0}$. Then:

*(i)* The number of calls to $\mathcal{A}_d$ (step 6) is bounded. That is, $j_{out}$ is finite.

*(ii)* The number of iterations of $\mathcal{A}$ at each call of $\mathcal{A}_d$ (step 6) is upper bounded by $\lceil 4n_\rho \rceil$. That is,

$$m_{j+1} \leq \lceil 4n_\rho \rceil, \ \forall j \in \mathbb{Z}_0^{j_{out}}. \tag{3.30}$$

*(iii)* The total number of iterations of $\mathcal{A}$ required to attain the exit condition $f(r_{j_{\text{out}}}) - f(r_{j_{\text{out}}+1}) \leq \epsilon$ (see step 8), is upper bounded by

$$k_{\text{out}} = \sum_{j=0}^{j_{\text{out}}} m_{j+1} \leq \frac{e\lceil 4n_\rho \rceil}{2} \left\lceil 5 + \frac{1}{\log 15} \log\left(1 + \frac{f(r_0) - f^*}{\epsilon}\right) \right\rceil.$$

***Proof:*** In the following, we will make use of two technical lemmas: Lemmas 3.35 and 3.36, which we have included, for convenience, immediately after this proof.

Let $T \in \mathbb{Z}$ be such that

$$f(r_j) - f(r_{j+1}) > \epsilon, \ \forall j \in \mathbb{Z}_0^T, \tag{3.31}$$

is satisfied. Then, defining $d_j \doteq f(r_j) - f(r_{j+1})$, we have

$$f(r_0) - f(r_{T+1}) = \sum_{j=0}^{T} d_j \geq (T+1)\left(\min_{j=0,\dots,T} d_j\right) > (T+1)\epsilon.$$

Thus,

$$T + 1 < \frac{f(r_0) - f(r_{T+1})}{\epsilon} \leq \frac{f(r_0) - f^*}{\epsilon} \leq \frac{\rho}{\epsilon},$$

from where we infer that the largest integer $T$ satisfying (3.31) is bounded. Consequently, the exit condition of Algorithm 10 (step 8) is satisfied within a finite number of iterations, thus proving claim *(i)*.

To prove claim *(ii)*, we start by noting that both $m_1$ and $m_2$ are no larger than $\lceil 4n_\rho \rceil$. Indeed, from step 4 we have that $s_0 = s_1 = 0$, which, in virtue of step 5, implies that $n_0 = m_0 = 1$ and $n_1 = m_1$. Since $n_0 = 1$ is no larger than $\lceil 4n_\rho \rceil$ we have from (3.28c) that $m_1$ is also upper-bounded by $\lceil 4n_\rho \rceil$. Moreover, since $n_1 = m_1 \leq \lceil 4n_\rho \rceil$, we obtain by the same reasoning that $m_2 \leq \lceil 4n_\rho \rceil$. We now prove that if $j \geq 2$ and $m_j \leq \lceil 4n_\rho \rceil$, then $m_{j+1} \leq \lceil 4n_\rho \rceil$. From step 4 we have

$$\begin{aligned} s_j^2 &= \frac{f(r_{j-1}) - f(r_j)}{f(r_{j-2}) - f(r_j)} = 1 - \frac{f(r_{j-2}) - f(r_{j-1})}{f(r_{j-2}) - f(r_j)} \\ &\leq 1 - \frac{f(r_{j-2}) - f(r_{j-1})}{f(r_{j-2}) - f^*} = \frac{f(r_{j-1}) - f^*}{f(r_{j-2}) - f^*} \overset{(3.28b)}{\leq} \left(\frac{n_\rho}{m_{j-1}+1}\right)^2. \end{aligned}$$

Thus, we have $s_j m_{j-1} \leq n_\rho$. Therefore,

$$n_j = \max\{m_j, 4 s_j m_{j-1}\} \leq \max\{\lceil 4n_\rho \rceil, 4n_\rho\} = \lceil 4n_\rho \rceil,$$

which, along with (3.28c), leads to $m_{j+1} \leq \lceil 4n_\rho \rceil$, thus proving the claim.

Finally, to prove claim *(iii)*, we start by noting that the computation of each $r_{j+1}$ is obtained from $m_{j+1}$ iterations of $\mathcal{A}$ (see step 6). Thus,

$$k_{\text{out}} = \sum_{j=0}^{j_{\text{out}}} m_{j+1} \overset{(3.30)}{\leq} (1 + j_{out}) \lceil 4n_\rho \rceil. \tag{3.32}$$

Let us denote

$$D \doteq \left\lceil 5 + \frac{1}{\log 15} \log\left(1 + \frac{f(r_0) - f^*}{\epsilon}\right)\right\rceil.$$

Consider first the case $j_{\text{out}} < D$. Since both $j_{\text{out}}$ and $D$ are integers we infer from this inequality that $1 + j_{\text{out}} \leq D$. This, along with (3.32), implies that $k_{\text{out}} \leq \lceil 4n_\rho \rceil D$.

Suppose now that $j_{\text{out}} \geq D$. We first recall that Property 3.31.*(ii)* states that the sequence $\{m_{j+1}\}_{\geq 0}$ is non-decreasing. We now rewrite $j_{\text{out}}$ as $j_{\text{out}} = d + tD$, where $d \in \mathbb{Z}_{0,D-1}$ and $t \in \mathbb{Z}_{\geq 0}$. Thus,

$$k_{\text{out}} = \sum_{j=0}^{d} m_{j+1} + \sum_{j=1}^{tD} m_{d+j+1} \leq D m_{d+1} + D \sum_{i=1}^{t} m_{d+1+iD} = D \sum_{i=0}^{t} m_{d+1+iD}.$$

From Lemma 3.36.*(v)* we have

$$m_{d+1+iD} \leq \frac{m_{d+1+(i+1)D}}{\sqrt{15}}, \ \forall i \in \mathbb{Z}_0^{t-1}.$$

Thus,

$$k_{\text{out}} \leq D \sum_{i=0}^{t} m_{d+1+tD} \left( \frac{1}{\sqrt{15}} \right)^{t-i}.$$

Using now $m_{1+d+tD} \leq \lceil 4 n_\rho \rceil$, see (3.30), we obtain

$$\frac{N_{\mathcal{A}}}{D \lceil 4 n_\rho \rceil} \leq \sum_{i=0}^{t} \left( \frac{1}{\sqrt{15}} \right)^{t-i} = \sum_{j=0}^{t} \left( \frac{1}{\sqrt{15}} \right)^{j} \leq \sum_{j=0}^{\infty} \left( \frac{1}{\sqrt{15}} \right)^{j} = \frac{\sqrt{15}}{\sqrt{15}-1} \leq \frac{e}{2}.$$

Thus, $N_{\mathcal{A}} \leq \frac{e}{2} \lceil 4 n_\rho \rceil D.$ ■

**Remark 3.34.** *The iteration complexity provided in Theorem 3.33.(iii) considers the exit condition shown in step 8 of Algorithm 10, whereas the iteration complexities shown in previous sections consider exit conditions of the form $\|g(r_{out})\|_* \leq \tilde{\epsilon}$. However, in view of Remark 3.32, this exit condition can also be used in Algorithm 10, in which case the iteration complexity shown in Theorem 3.33.(iii) would be the same but replacing $\epsilon$ with $\tilde{\epsilon}/(2L_f)$.*

The following lemma is used exclusively in the proof of Lemma 3.36.

**Lemma 3.35.** The function $\varphi : \mathbb{R} \to \mathbb{R}$ defined as

$$\varphi(s) \doteq \left( \frac{1}{s^2} - 1 \right) \max \left\{ 1, (4s)^4 \right\},$$

satisfies $\varphi(s) \geq 15$, $\forall s \in (0, \frac{\sqrt{15}}{4}]$.

***Proof:*** We have that

$$\varphi(s) = \begin{cases} 4^4(s^2 - s^4) & \text{if } s > \frac{1}{4}, \\ \dfrac{1}{s^2} - 1 & \text{if } s \leq \frac{1}{4}. \end{cases}$$

It is clear that $\varphi(\cdot)$ is monotonically decreasing in $(0, \frac{1}{4}]$. Thus,

$$\min_{s\in(0,\frac{\sqrt{15}}{4}]} \varphi(s) = \min_{s\in[\frac{1}{4},\frac{\sqrt{15}}{4}]} \varphi(s) = \min_{s\in[\frac{1}{4},\frac{\sqrt{15}}{4}]} 4^4(s^2 - s^4).$$

We notice that the derivative of $s^2 - s^4$ is $2s(1-2s^2)$, which vanishes only once in the interval of interest (at $s = \frac{1}{\sqrt{2}}$). From here we infer that $s^2 - s^4$ is increasing in $[\frac{1}{4}, \frac{1}{\sqrt{2}})$ and decreasing in $(\frac{1}{\sqrt{2}}, \frac{\sqrt{15}}{4}]$. Thus, the minimum is attained at the extremes of the interval $[\frac{1}{4}, \frac{\sqrt{15}}{4}]$. That is, we conclude that

$$\min_{s\in(0,\frac{\sqrt{15}}{4}]} \varphi(s) = \min\{\varphi(\frac{1}{4}), \varphi(\frac{\sqrt{15}}{4})\} = \min\{15, 15\} = 15.$$

■

**Lemma 3.36** (A few technical results on the iterates of Algorithm 10)**.** Let Assumtion 3.24 hold and consider Algorithm 10 for a given initial condition $r_0 \in V_f(\rho) \subseteq \mathrm{dom}(f)$ and an exit tolerance $\epsilon \in \mathbb{R}_{>0}$. Assume that $j_{\mathrm{out}} \geq 2$ and that there exist $T \in \mathbb{Z}_2^{j_{\mathrm{out}}}$ and $\ell \in \mathbb{Z}_0^{j_{\mathrm{out}}-T}$ such that

$$m_{\ell+1} > \frac{1}{\sqrt{15}} m_{\ell+1+T}.$$

Then,

*(i)* $s_j \in \left(0, \frac{\sqrt{15}}{4}\right], \forall j \in \mathbb{Z}_{\ell+2}^{\ell+T}$,

*(ii)* $\sum_{j=\ell+2}^{\ell+T} \log\left(\max\{1, (4s_j)^4\}\right) < 4\log 15$,

*(iii)* $\sum_{j=\ell+2}^{\ell+T} \log\left(\frac{1}{s_j^2} - 1\right) \leq \log\left(1 + \frac{f(r_0) - f^*}{\epsilon}\right)$,

*(iv)* $T < 5 + \frac{1}{\log 15}\log\left(1 + \frac{f(r_0) - f^*}{\epsilon}\right)$.

Additionally, let

$$D \doteq \left\lceil 5 + \frac{1}{\log 15}\log\left(1 + \frac{f(r_0) - f^*}{\epsilon}\right)\right\rceil.$$

Then,

*(v)* $m_{\ell+1} \leq \frac{1}{\sqrt{15}} m_{\ell+1+D}, \ \forall \ell \in \mathbb{Z}_0^{j_{\mathrm{out}}-D}$.

***Proof:*** From step 4 of Algorithm 10 we have

$$s_j^2 = \frac{f(r_{j-1}) - f(r_j)}{f(r_{j-2}) - f(r_j)},\ j \in \mathbb{Z}_2^{j_{\rm out}}.$$

The inequality $s_j > 0$, $\forall j \in \mathbb{Z}_{\ell+2}^{\ell+T}$ follows from Property 3.31.*(iii)*. In order to prove the first claim it remains to prove the inequality $s_j \leq \frac{\sqrt{15}}{4}$, $\forall j \in \mathbb{Z}_{\ell+2}^{\ell+T}$. We proceed by reductio ad absurdum. Assume that there is $j \in \mathbb{Z}_{\ell+2}^{\ell+T}$ such that $s_j > \frac{\sqrt{15}}{4}$. In this case,

$$m_{j+1} \overset{(3.27)}{\geq} n_j = \max\{m_j, 4s_j m_{j-1}\} \geq 4s_j m_{j-1} > \sqrt{15} m_{j-1}.$$

From this and the non-decreasing nature of the sequence $\{m_j\}_{\geq 0}$ stated in Property 3.31.*(ii)* we obtain

$$m_{\ell+1+T} \geq m_{j+1} > \sqrt{15} m_{j-1} \geq \sqrt{15} m_{\ell+1},$$

which contradicts the initial assumptions, thus proving claim *(i)*.

From the non-decreasing nature of the sequence $\{m_j\}_{\geq 0}$ stated in Property 3.31.*(ii)* we have

$$m_{j+1} \overset{(3.27)}{\geq} n_j = \max\{m_j, 4s_j m_{j-1}\} \geq m_{j-1} \max\{1, 4s_j\},\ \forall j \in \mathbb{Z}_{\ell+2}^{\ell+T}.$$

Equivalently,

$$\log(\max\{1, 4s_j\}) \leq \ln \frac{m_{j+1}}{m_{j-1}},\ \forall j \in \mathbb{Z}_{\ell+2}^{\ell+T},$$

which implies

$$\begin{aligned}
\sum_{j=\ell+2}^{\ell+T} \log(\max\{1, 4s_j\}) &\leq \sum_{j=\ell+2}^{\ell+T} \log \frac{m_{j+1}}{m_{j-1}} = \log \frac{m_{\ell+T} m_{\ell+1+T}}{m_{\ell+1} m_{\ell+2}} \\
&\leq \log \frac{m_{\ell+1+T}^2}{m_{\ell+1}^2} = 2\log \frac{m_{\ell+1+T}}{m_{\ell+1}} \\
&< 2\log \sqrt{15} = \log 15.
\end{aligned}$$

Claim *(ii)* immediately follows from multiplying this inequality by 4.

To prove the third claim, we start by noticing that

$$\prod_{j=\ell+2}^{\ell+T} (\frac{1}{s_j^2} - 1) = \prod_{j=\ell+2}^{\ell+T} \frac{f(r_{j-2}) - f(r_{j-1})}{f(r_{j-1}) - f(r_j)} = \frac{f(r_\ell) - f(r_{\ell+1})}{f(r_{\ell+T-1}) - f(r_{\ell+T})}.$$

Since $\ell + T \leq j_{\rm out}$ we have $f(r_{\ell+T-1}) - f(r_{\ell+T}) > \epsilon > 0$, which leads to

$$\prod_{j=\ell+2}^{\ell+T} (\frac{1}{s_j^2} - 1) < \frac{f(r_\ell) - f(_{\ell+1})}{\epsilon} \overset{(3.28a)}{\leq} \frac{f(r_0) - f(r_{\ell+1})}{\epsilon} \leq \frac{f(r_0) - f^*}{\epsilon},$$

from where claim *(iii)* directly follows.

Next, we add the inequalities given in claims *(ii)* and *(iii)* to obtain

$$\sum_{j=\ell+2}^{\ell+T} \log\left(\left(\frac{1}{s_j^2}-1\right)\max\left\{1,(4s_j)^4\right\}\right) < \log\left(1+\frac{f(r_0)-f^*}{\epsilon}\right)+4\log 15. \quad (3.33)$$

From claim *(i)* we have $s_j \in \left(0, \frac{\sqrt{15}}{4}\right]$, $\forall j \in \mathbb{Z}_{\ell+2}^{\ell+T}$. Thus, making use of Lemma 3.35, the left term of (3.33) can be lower bounded by means of the following inequality

$$15 \le \left(\frac{1}{s^2}-1\right)\max\left\{1,(4s)^4\right\},\ \forall s \in \left(0, \frac{\sqrt{15}}{4}\right].$$

That is,

$$\sum_{j=\ell+2}^{\ell+T} \log 15 < \log\left(1+\frac{f(r_0)-f^*}{\epsilon}\right)+4\log 15.$$

Equivalently,

$$(T-1)\log 15 < \log\left(1+\frac{f(r_0)-f^*}{\epsilon}\right)+4\log 15,$$

thus proving claim *(iv)*.

Finally, we prove claim *(v)* by reductio ad absurdum. If there exist $\ell \in \mathbb{Z}_0^{j_{\text{out}}-D}$ such that $m_{\ell+1} > \frac{1}{\sqrt{15}} m_{\ell+1+D}$, then we obtain from claim *(iv)* that

$$D < 5 + \frac{1}{\ln 15}\ln\left(1+\frac{f(z_0)-f^*}{\epsilon}\right),$$

which contradicts the definition of $D$. ■

## 3.3 Numerical results

This section presents numerical results comparing the three restart schemes proposed in Section 3.2 with some of the restart schemes of the literature described in Section 3.1. In particular, we compare the restart schemes presented in Algorithms 7, 8 and 10 with the *objective function value restart scheme* [45], whose restart condition $E_f$ is given by (3.6); the *gradient alignment restart scheme* [45], whose restart condition $E_g$ is given by (3.7); and the optimal fixed-rate restart scheme from [46, §5.2.2], using the restart condition $E_f^*$ given by (3.5).

For convenience and space considerations, the tables and figures of this section will use the nomenclature of the exit conditions ($E_f$, $E_g$ and $E_f^*$) to refer to the restart schemes of the literature. Similarly, we will refer to the restart schemes we propose in Section 3.2 by their algorithms (Alg. 7, Alg. 8 and Alg. 10).

The results shown here use FISTA (Algorithm 4), with the exception of Algorithm 10, which uses MFISTA (Algorithm 5). We chose MFISTA because, due to its monotone behavior, step 4 of Algorithm 9 always evaluates to the upper expression (see Remark 3.28). Therefore, MFISTA can be directly used in place of $\mathcal{A}_d$ in step 6 of Algorithm 10 making use of the restart condition shown in Remark 3.28. The norm of the restart scheme presented in Algorithm 10 will use the same norm $\|\cdot\|_R$ as the other restart schemes. In order to provide a fair comparison between the different schemes, we exit them as soon as an iterate $z_k$ satisfying $\|\mathcal{G}(z_k)\|_{R^{-1}} \leq \epsilon$ is attained for the selected value of $\epsilon \in \mathbb{R}_{>0}$. This way, we can compare how quickly each restart scheme finds an $\epsilon$-accurate solution of the problem at hand. Restart scheme $E_f^*$ requires knowledge of $f^*$, which we obtain by solving the optimization problem using Algorithm 7 with and exit tolerance of $\epsilon = 10^{-8}$.

As additional numerical results, Section 5.8.7 shows the application of the above mentioned restart schemes to solve the optimization problems of different MPC formulations using FISTA.

### 3.3.1 Application to Lasso problems

This section presents the result of applying the restart schemes to weighted Lasso problems of the form

$$\min_{z \in \mathbb{R}^{n_z}} \frac{1}{2N} \|Az - b\|_2^2 + \|Wz\|_1, \tag{3.34}$$

where $z \in \mathbb{R}^{n_z}$, $A \in \mathbb{R}^{N \times n_z}$ is sparse with an average of 90% of its entries being zero (sparsity is generated by setting a 0.9 probability for each element of the matrix to be 0), $n_z > N$, and $b \in \mathbb{R}^N$. Each nonzero element in $A$ and $b$ is obtained from a Gaussian distribution with zero mean and variance 1. $W \in \mathbb{D}_{++}^{n_z}$ is a diagonal matrix with elements obtained from a uniform distribution on the interval $(0, \alpha]$, with $\alpha \in \mathbb{R}_{>0}$. We note that problems (3.34) can be reformulated in such a way that they satisfy the quadratic growth condition [46, §6.3].

We take $R$ as the diagonal matrix constructed as

$$R_{(i,i)} = \sum_{j=1}^{n} |H_{(i,j)}|, \tag{3.35}$$

where $H = \frac{1}{N} A^\top A$, which due to the Gershgorin Circle Theorem [60, §7.2] (see also [27, §6]) satisfies the smoothness condition given in Assumtion 2.1.*(ii)*.

Table 3.1 shows the results of solving 100 randomly generated problems (3.34) that share the values of $N = 600$, $n_z = 800$, and $\alpha = 0.003$. We take $\epsilon = 10^{-7}$. Figure 3.9 and Figure 3.10 show the evolution of $\|\mathcal{G}(z_k)\|_{R^{-1}}$ and $f(z_k) - f^*$, respectively, of each one of the restart schemes for one of the Lasso problems (3.34) used to obtain the results of Table 3.1. Additionally, they show the result of applying FISTA without a restart scheme.

Table 3.2 and Figures 3.11 and 3.12 show analogous results to Table 3.1 and Figures 3.9 and 3.10, respectively, but taking $N = 100$, $n_z = 200$ and $\alpha = 0.3$.

| Restart scheme | Alg. 7 | Alg. 8 | Alg. 10 | $E_f$ | $E_g$ | $E_f^*$ |
|---|---|---|---|---|---|---|
| Avg. Iter. | 914.69 | 1431.2 | 897.46 | 946.62 | 894.21 | 1350.3 |
| Med. Iter. | 907 | 1410.5 | 880.5 | 945.5 | 880.5 | 1331 |
| Max. Iter. | 1155 | 1800 | 1131 | 1349 | 1199 | 1696 |
| Min. Iter. | 704 | 1154 | 737 | 692 | 702 | 1084 |

Table 3.1: Comparison between restart schemes applied to FISTA to solve 100 problems (3.34) with $N = 600$, $n_z = 800$, $\alpha = 0.003$ and $\epsilon = 10^{-7}$.

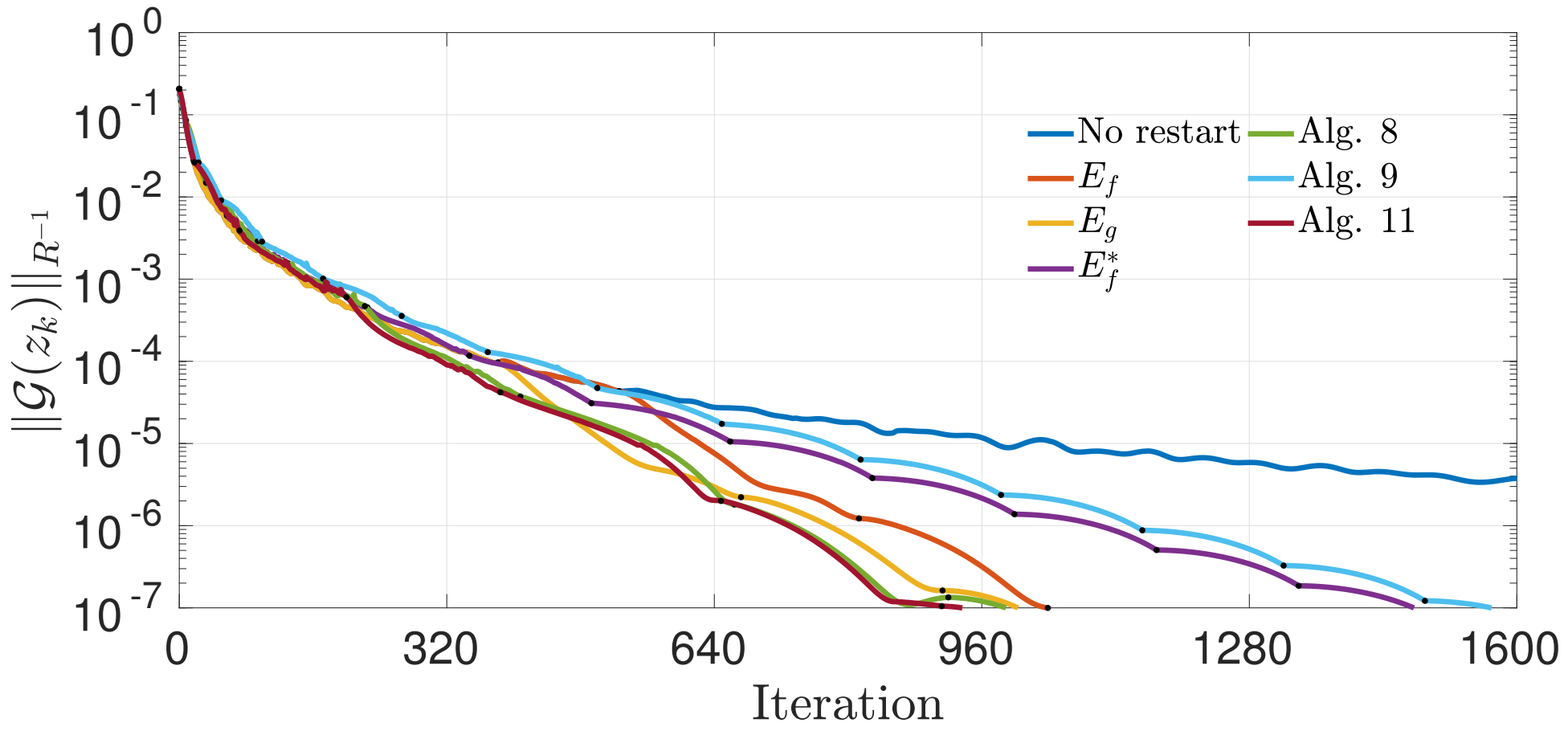

Figure 3.9: Evolution of the dual norm of the composite gradient mapping for different restart schemes applied to a randomly generated problem (3.34) with $N = 600$, $n_z = 800$, $\alpha = 0.003$ and $\epsilon = 10^{-7}$. Black dots represent the iterations in which a restart occurred.

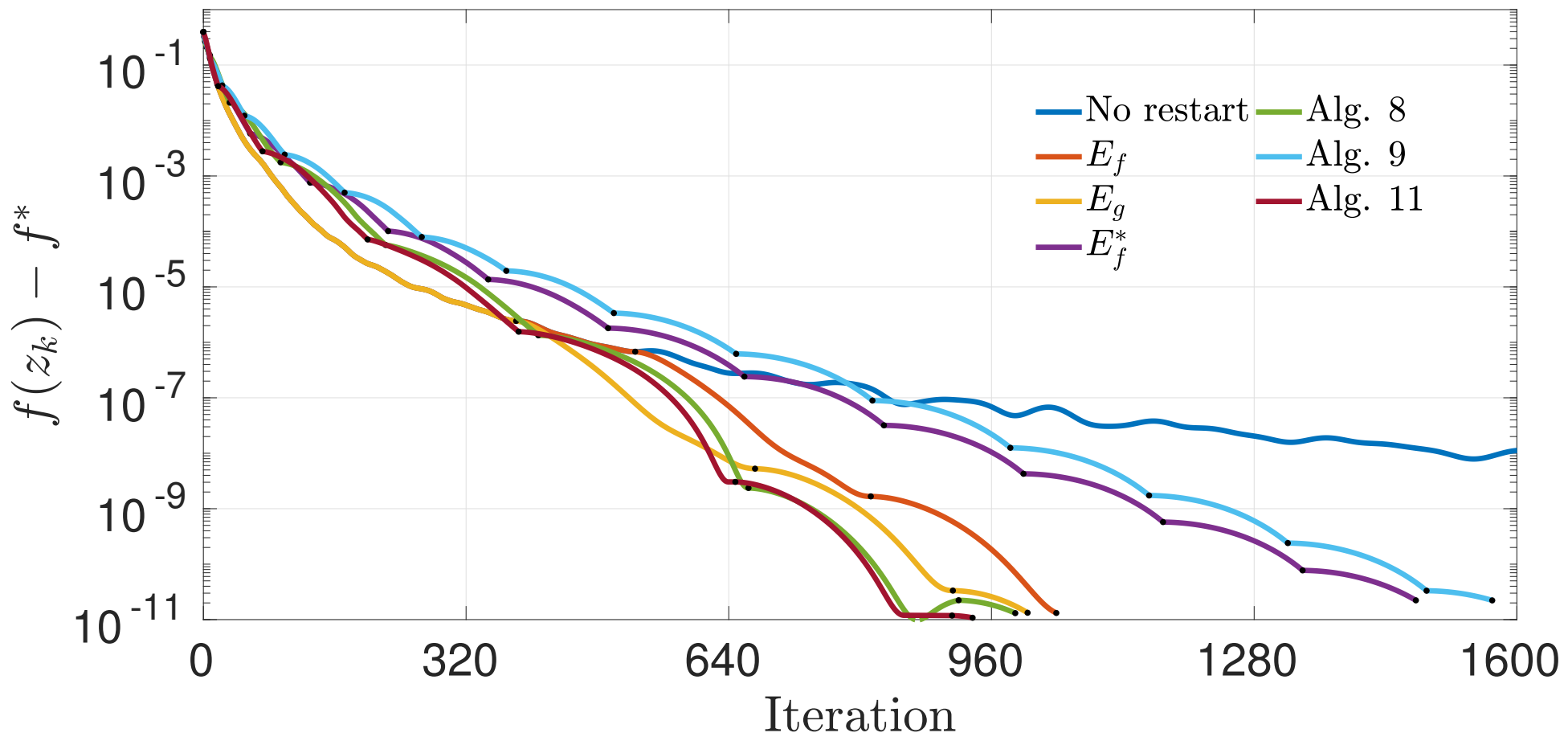

Figure 3.10: Evolution of the distance, in terms of function value, to the optimal solution for different restart schemes applied to a randomly generated problem (3.34) with $N = 600$, $n_z = 800$, $\alpha = 0.003$ and $\epsilon = 10^{-7}$. Black dots represent the iterations in which a restart occurred.

| Restart scheme | Alg. 7 | Alg. 8 | Alg. 10 | $E_f$ | $E_g$ | $E_f^*$ |
|---|---|---|---|---|---|---|
| Avg. Iter. | 91.33 | 146.5 | 83.04 | 87.53 | 82.96 | 136.42 |
| Med. Iter. | 91.5 | 141 | 80 | 86.5 | 81 | 131 |
| Max. Iter. | 195 | 372 | 158 | 164 | 146 | 350 |
| Min. Iter. | 69 | 104 | 61 | 58 | 58 | 98 |

Table 3.2: Comparison between restart schemes applied to FISTA to solve 100 problems (3.34) with $N = 100$, $n_z = 200$, $\alpha = 0.3$ and $\epsilon = 10^{-7}$.

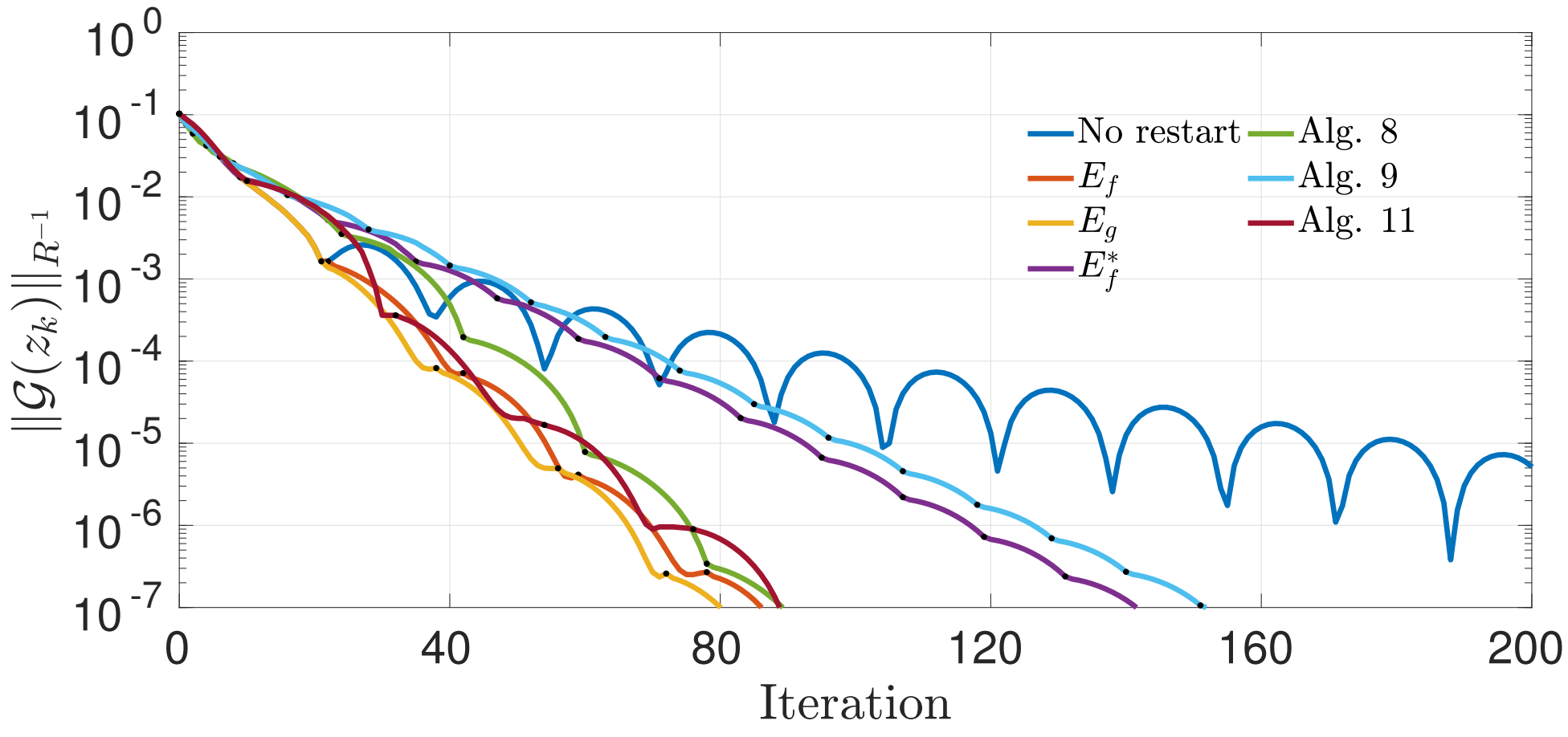

Figure 3.11: Evolution of the dual norm of the composite gradient mapping for different restart schemes applied to a randomly generated problem (3.34) with $N = 100$, $n_z = 200$, $\alpha = 0.3$ and $\epsilon = 10^{-7}$. Black dots represent the iterations in which a restart occurred.

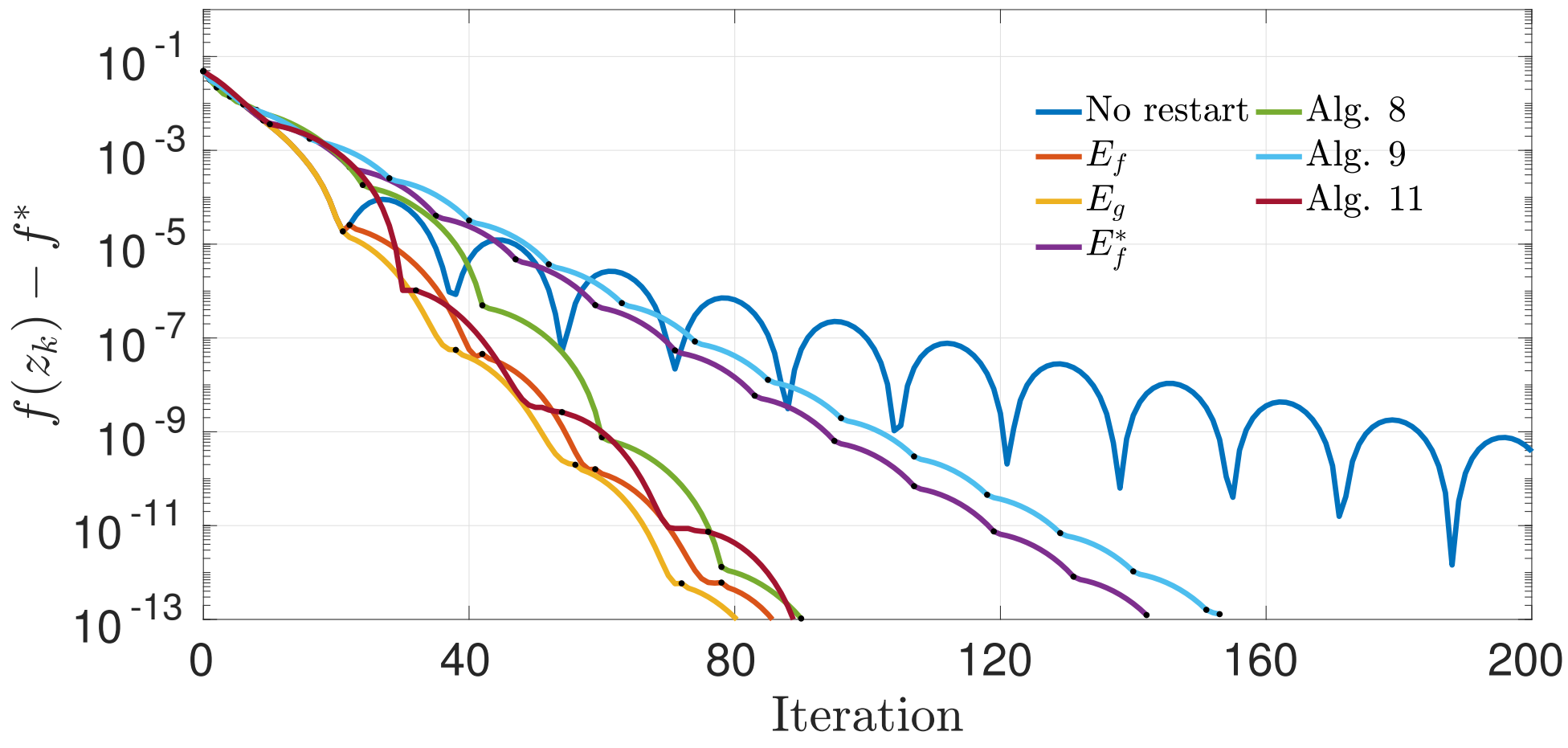

Figure 3.12: Evolution of the distance, in terms of function value, to the optimal solution for different restart schemes applied to a randomly generated problem (3.34) with $N = 100$, $n_z = 200$, $\alpha = 0.3$ and $\epsilon = 10^{-7}$. Black dots represent the iterations in which a restart occurred.

### 3.3.2 Application to QP problems

This section presents the results of applying the restart schemes to unconstrained QP problems of the form

$$\min_{z \in \mathbb{R}^{n_z}} \frac{1}{2} z^\top H z + q^\top z, \tag{3.36}$$

where $z \in \mathbb{R}^{n_z}$, $H \in \mathbb{S}_{++}^{n_z}$ and $q \in \mathbb{R}^{n_z}$. Matrix $H$ is obtained by first computing a matrix $M$ whose elements are taken from a uniform distribution on the interval $(0, 1]$ and then taking

$$H = \frac{1}{2} M^\top M + \alpha \mathbf{I}_{n_z}$$

for some $\alpha \in \mathbb{R}_{>0}$ whose value affects the condition number of $H$ (smaller values of $\alpha$ tend to result in higher condition numbers). The elements of vector $q$ are taken from a uniform distribution on the interval $(0, \beta]$ for some $\beta \in \mathbb{R}_{>0}$. Note that problems (3.36) have no non-smooth term $\Psi$. As in Section 3.3.1, we take $R$ using (3.35).

Table 3.3 shows the results of solving 100 randomly generated problems (3.36) that share the values of $n_z = 200$, $\alpha = 10$ and $\beta = 20$, taking $\epsilon = 10^{-5}$. It also shows information about the condition numbers of the matrices $H$. Figure 3.13 and Figure 3.14 show the evolution of $\|\mathcal{G}(z_k)\|_{R^{-1}}$ and $f(z_k) - f^*$, respectively, of each one of the restart schemes for one of the QP problems (3.36) used to obtain the results of Table 3.3. Additionally, they show the result of applying FISTA without a restart scheme.

Table 3.4 and Figures 3.15 and 3.16 show analogous results to Table 3.3 and Figures 3.13 and 3.14, respectively, but taking $\alpha = 0.1$ instead of $\alpha = 10$, resulting in worse conditioned QP problems.

| Restart scheme | Alg. 7 | Alg. 8 | Alg. 10 | $E_f$ | $E_g$ | $E_f^*$ | cond$(H)$ |
|---|---|---|---|---|---|---|---|
| Avg. Iter. | 177.39 | 324.76 | 179.5 | 189.23 | 184.82 | 304.95 | 40.12 |
| Med. Iter. | 178 | 325 | 180 | 188 | 184.5 | 305 | 40.03 |
| Max. Iter. | 187 | 334 | 187 | 218 | 206 | 317 | 43.25 |
| Min. Iter. | 171 | 315 | 173 | 167 | 168 | 296 | 37.93 |

Table 3.3: Comparison between restart schemes applied to FISTA to solve 100 problems (3.36) with $n_z = 200$, $\alpha = 10$, $\beta = 20$ and $\epsilon = 10^{-5}$. The last column shows details about the condition numbers of the matrices $H$.

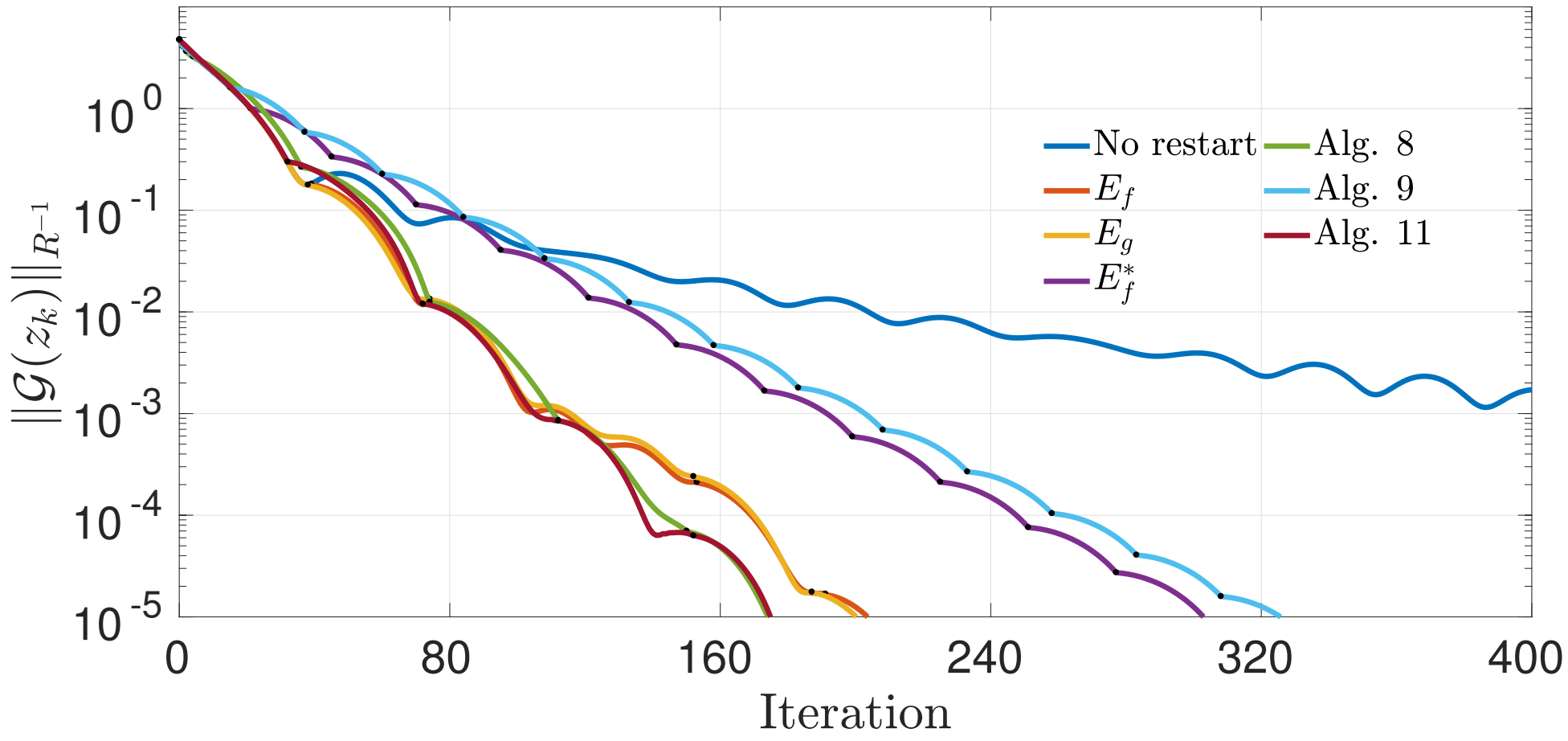

Figure 3.13: Evolution of the dual norm of the composite gradient mapping for different restart schemes applied to a randomly generated problem (3.36) with $n_z = 200$, $\alpha = 10$, $\beta = 20$ and $\epsilon = 10^{-5}$. Black dots represent the iterations in which a restart occurred.

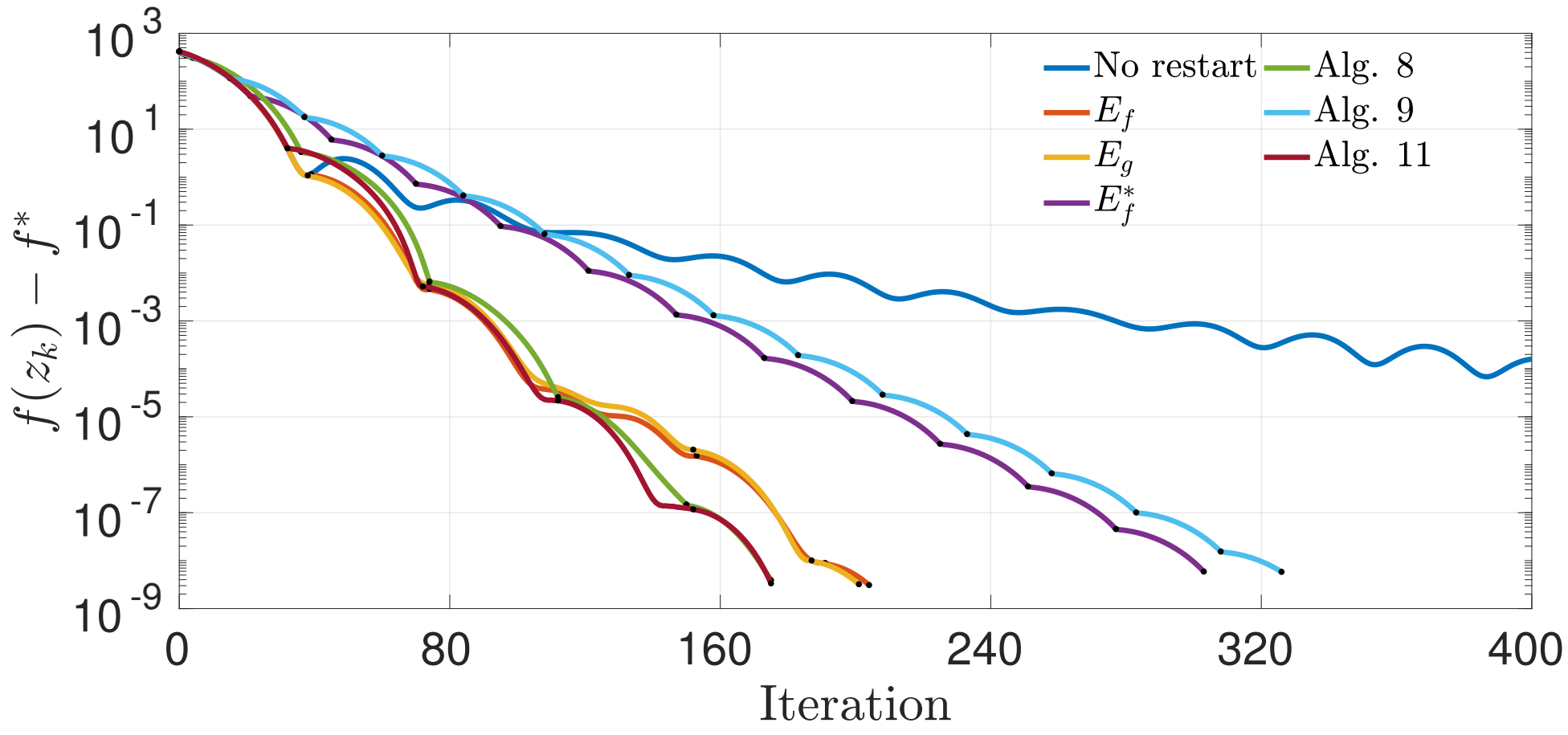

Figure 3.14: Evolution of the distance, in terms of function value, to the optimal solution for different restart schemes applied to a randomly generated problem (3.36) with $n_z = 200$, $\alpha = 10$, $\beta = 20$ and $\epsilon = 10^{-5}$. Black dots represent the iterations in which a restart occurred.

| Restart scheme | Alg. 7 | Alg. 8 | Alg. 10 | $E_f$ | $E_g$ | $E_f^*$ | cond$(H)$ |
|---|---|---|---|---|---|---|---|
| Avg. Iter. | 1816.4 | 3115.8 | 1784.8 | 1725.3 | 1724.1 | 3037.2 | 3864.9 |
| Med. Iter. | 1816.5 | 3126.5 | 1796 | 1727.5 | 1721.5 | 3053.5 | 3860.7 |
| Max. Iter. | 2099 | 3319 | 2041 | 2228 | 2228 | 3269 | 4200.7 |
| Min. Iter. | 1447 | 2719 | 1505 | 1396 | 1413 | 2621 | 3355.18 |

Table 3.4: Comparison between restart schemes applied to FISTA to solve 100 problems (3.36) with $n_z = 200$, $\alpha = 0.1$, $\beta = 20$ and $\epsilon = 10^{-5}$. The last column shows details about the condition numbers of the matrices $H$.

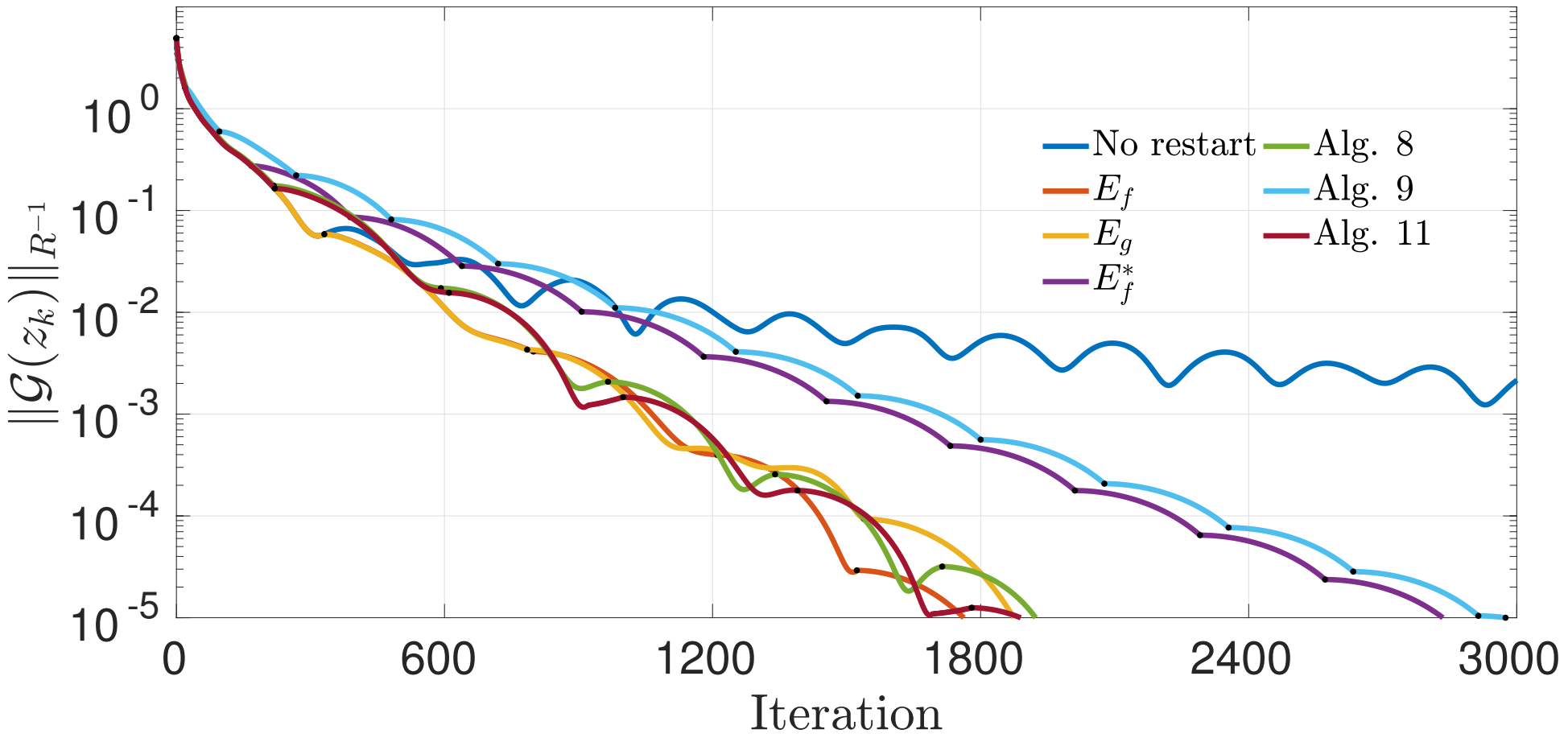

Figure 3.15: Evolution of the dual norm of the composite gradient mapping for different restart schemes applied to a randomly generated problem (3.36) with $n_z = 200$, $\alpha = 0.1$, $\beta = 20$ and $\epsilon = 10^{-5}$. Black dots represent the iterations in which a restart occurred.

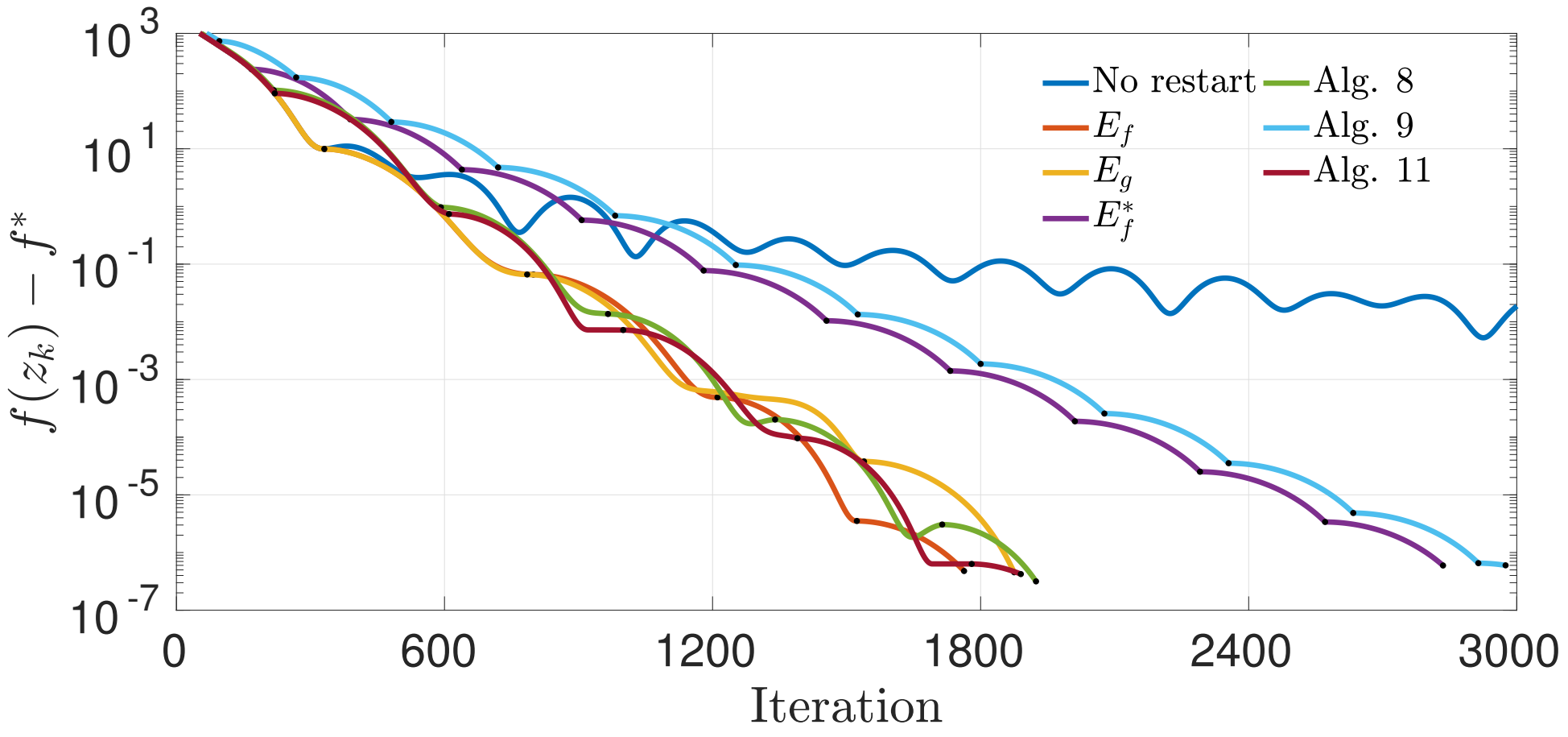

Figure 3.16: Evolution of the distance, in terms of function value, to the optimal solution for different restart schemes applied to a randomly generated problem (3.36) with $n_z = 200$, $\alpha = 0.1$, $\beta = 20$ and $\epsilon = 10^{-5}$. Black dots represent the iterations in which a restart occurred.

## 3.4 Conclusions and future lines of work

The numerical results presented in Section 3.3 showcase the good practical performance of Algorithms 7 and 10, which solve optimization problems in a number of iterations similar to the ones obtained using the adaptive restart schemes whose exit conditions are given by $E_f$ (3.6) and $E_g$ (3.7). Algorithms 7 and 10, however, are linearly convergent for the class of non-strongly convex optimization problems described in Section 3.2.

Algorithm 8, on the other hand, performs similarly to the restart scheme whose restart condition is given by $E_f^*$ (3.5), both in terms of the number of iterations and in the fact that a restart is produced at intervals of similar length. We find that Algorithm 8 tends to perform slightly worse than $E_f^*$ and that they both tend to have a much more predictable performance than the other restart schemes. Restart scheme $E_f^*$, however, requires knowledge of $f^*$, which is not always available.

A comparison between the figures showing $\|\mathcal{G}(z_k)\|_{R^{-1}}$ with those showing $f(z_k) - f^*$ clearly highlights the relation between $\|\mathcal{G}(z_k)\|_{R^{-1}}$ and the optimality of $z_k$ highlighted in Proposition 2.6.

The practical usefulness of the restart schemes would have to be further analyzed by performing tests comparing the computation times of each, since the advantages obtained by the reduction in the number of iterations may be overshadowed by the additional computations required to evaluate the restart conditions. Even though their usefulness would have to be determined in a case-by-case basis, additional research can be made in this line. For instance, the recent preprint [61] presents a restart scheme that is based on the restart scheme presented in Section 3.2.1 but that required less evaluations of the objective function. The results they show indicate that their proposed restart scheme required more iterations that ours, but that the computation times may be smaller.

Another future line of work is the proposal of a gradient-based restart scheme similar to the one presented in Section 3.2.2 but applicable to a wider class of accelerated first order methods.

Finally, the results shown here have focused on FISTA and its monotone variant MFISTA, but the application of the restart schemes to other accelerated first order methods, such as the accelerated ADMM, in an interesting line of research.

# Part II

# Implementation of MPC in embedded systems

# Chapter 4

# Preliminaries

This chapter addresses some of the preliminary concepts and ideas related to model predictive control (MPC) and its implementation in embedded systems.

We start by introducing the problem formulation in Section 4.1, i.e., by presenting the class of system and system model under consideration as well as the control objective. We then provide a brief description of MPC in Section 4.2 so as to introduce the concepts and nomenclature used throughout the remainder of the dissertation. Finally, in Section 4.3 we present a brief (and non-exhaustive) review of the literature regarding the implementation of MPC in embedded systems.

Before we proceed, we provide for the sake of clarity our definition of an embedded system. Some examples of well known embedded systems are Programmable Logic Controllers (PLC), FPGAs, Arduino and Raspberry Pi.

**Definition 4.1** (Embedded system)**.** We understand by *embedded system* any computer system that combines both hardware and software and that has a dedicated function within a larger system. In particular, we circumscribe our definition to those electronic systems that are based either on microcontrollers or microprocessors and that are used to control the physical device they are embedded in. Integral to our definition is the fact that embedded systems are characterized by their low computational and memory resources when compared to more powerful computing devices, such as the average desktop PC.

## 4.1 Problem formulation

We consider a discrete-time, linear, time-invariant, state space model

$$x(t+1) = Ax(t) + Bu(t) \tag{4.1}$$

where $x(t) \in \mathbb{R}^n$ and $u(t) \in \mathbb{R}^m$ are the state and control input of the system at the discrete time instant $t \in \mathbb{Z}$, $A \in \mathbb{R}^{n \times n}$, and $B \in \mathbb{R}^{n \times m}$. Additionally, we assume that the system is subject to state and input constraints

$$(x(t), u(t)) \in \mathcal{Y} \subseteq \mathbb{R}^n \times \mathbb{R}^m,\ \forall t, \tag{4.2}$$

where the set $\mathcal{Y}$ will be given either by box constraints of the form

$$\mathcal{Y} = \{ (x, u) \in \mathbb{R}^n \times \mathbb{R}^m \; : \; \underline{x} \le x \le \overline{x}, \; \underline{u} \le u \le \overline{u} \}, \tag{4.3}$$

where $\underline{x} \in \mathbb{R}^n$, $\overline{x} \in \mathbb{R}^n$, $\underline{u} \in \mathbb{R}^m$, and $\overline{u} \in \mathbb{R}^m$; or by coupled input-state constraints

$$\mathcal{Y} = \{ (x, u) \in \mathbb{R}^n \times \mathbb{R}^m \; : \; \underline{y} \le Ex + Fu \le \overline{y} \}, \tag{4.4}$$

where $\underline{y} \in \mathbb{R}^p$, $\overline{y} \in \mathbb{R}^p$, $E \in \mathbb{R}^{p \times n}$ and $F \in \mathbb{R}^{p \times m}$. We make the following assumption on model (4.1) and its constraints.

**Assumption 4.2.** We assume that:

*(i)* Model (4.1) is controllable.

*(ii)* The constraint set $\mathcal{Y}$ has a non-empty interior.

**Remark 4.3.** *If the constraints are given by* (4.3)*, then Assumtion 4.2.(ii) implies that* $\underline{x} < \overline{x}$ *and* $\underline{u} < \overline{u}$*.*

The control objective is to steer the system (4.1) to the given reference $(x_r, u_r) \in \mathbb{R}^n \times \mathbb{R}^m$ whilst satisfying the system constraints (4.2), where we assume that the reference is a (piecewise) constant set-point, i.e., a constant reference whose value can change at any unpredictable time. Obviously, this objective is only possible if the reference is an admissible steady state (see the following definition) of system (4.1) under (4.2).

**Definition 4.4** (Admissible steady state)**.** Consider a system (4.1) subject to (4.2) under Assumtion 4.2. We say that a pair $(\hat{x}, \hat{u}) \in \mathbb{R}^n \times \mathbb{R}^m$ is an *admissible steady state* of the (4.1) under (4.2) if $\hat{x} = A\hat{x} + B\hat{u}$ (i.e., it is a steady state) and $(\hat{x}, \hat{u}) \in \mathcal{Y}$. We say that it is *strictly* admissible if $(\hat{x}, \hat{u}) \in \text{int}(\mathcal{Y})$.

If the reference is not an admissible steady state, then we wish to steer the system to the "closest" admissible steady state, for some criterion of closeness.

In this dissertation we are concerned with linear systems (4.1). However, for future purposes, we now introduce some notation for the *real* system from where the linear model is derived.

In general, model (4.1) will be a linear representation of some real system governed by a continuous-time, non-linear, ordinary differential equation

$$\frac{d\chi}{dt} = f(\chi(t), u(t)), \tag{4.5}$$

where $\chi(t) \in \mathbb{R}^n$ and $u(t) \in \mathbb{R}^m$ are the state and control input of the system at the continuous time instant $t \in \mathbb{R}$, respectively, and $f : \mathbb{R}^n \times \mathbb{R}^m \to \mathbb{R}^n$ is a continuous differentiable function. Additionally, the constraints (4.2) will have been obtained so as to satisfy some real system constraints

$$(\chi(t), u(t)) \in \mathscr{Y} \subseteq \mathbb{R}^n \times \mathbb{R}^m, \; \forall t. \tag{4.6}$$

In the case of box constraints, the set $\mathcal{Y}$ will be given by

$$\mathcal{Y} = \{ (\chi, u) \in \mathbb{R}^n \times \mathbb{R}^m \, : \, \underline{\chi} \leq \chi \leq \overline{\chi}, \; \underline{u} \leq u \leq \overline{u} \}.$$

where $\underline{\chi} \in \mathbb{R}^n$, $\overline{\chi} \in \mathbb{R}^n$, $\underline{u} \in \mathbb{R}^m$, and $\overline{u} \in \mathbb{R}^m$.

In particular, we consider that model (4.1) has been obtained by linearizing the real system model (4.5) around an *operating point* $(\chi^\circ, u^\circ)$, which we assume to be an admissible steady state of (4.5) under (4.6). Additionally we consider that $x$ and $u$ may have been subject to scaling using the *scaling matrices* $N_x \in \mathbb{D}^n$ and $N_u \in \mathbb{D}^m$ as follows:

$$x = N_x(\chi - \chi^\circ), \quad u = N_u(u - u^\circ).$$

There are many other ways in which a linear model can be derived from a non-linear one or from data of the system [62]. Our above consideration is merely the one we take in the numerical examples we show in this dissertation.

We refer to $x$, $u$ and their bounds as being in *incremental units* to emphasize that they are (scaled) increments with respect to the operating point of the real system model. Additionally, we refer to $\chi$, $u$ and their bounds as being in *engineering units*, and we will represent them using calligraphic font. We use this naming convention to emphasize the fact that they relate to the real system (even though, in many occasions, the state has no tangible physical interpretation).

## 4.2 A brief introduction to model predictive control

Model predictive control (MPC) is an advanced control strategy that is very prevalent in the current control literature due to its inherent ability to provide constraint satisfaction and asymptotic stability to the given reference [63, 64]. In MPC, the control action is obtained, at each sample time, from the solution of an optimization problem in which a *prediction model* is used to predict the future evolution of the system over a given *prediction horizon*.

The need for solving an optimization problem in real-time and the fact that a model of the system is required have historically hindered the use of MPC in many industries and applications. In spite of this, there are numerous publications on the successful use of MPC in different areas, including water distribution networks [65], micro grids [66], power converters [67], automotive applications [68, 69], temperature control of buildings [70] or data centers [71], and of the heat, ventilation and air conditioning systems [17, 72, 73].

Different MPC formulations and control architectures are used in the above references. In this dissertation we focus on linear MPC formulations with simple constraints. This section is not intended as a thorough description nor explanation of MPC. Instead, it is a brief introduction to the class of MPC problems that we focus on in this dissertation, as well as to the nomenclature that we will use.

We consider linear MPC formulations described by the following parametric optimization problem:

$$\min_{\mathbf{x},\mathbf{u}} \left\{ J \doteq \sum_{j=0}^{N-1} \ell_j(\mathbf{x},\mathbf{u};x_r,u_r) + V_t(\mathbf{x},\mathbf{u};x_r,u_r) \right\} \tag{4.7a}$$

$$s.t.\ x_{j+1} = Ax_j + Bu_j,\ j \in \mathbb{Z}_0^{N-1} \tag{4.7b}$$

$$x_0 = x(t) \tag{4.7c}$$

$$(x_j, u_j) \in \mathcal{Y}_j,\ j \in \mathbb{Z}_0^{N-1} \tag{4.7d}$$

$$x_N \in \mathcal{X}_t, \tag{4.7e}$$

where $N \in \mathbb{Z}_{>0}$ is the prediction horizon; (4.7b) is the prediction model (4.1), with $A \in \mathbb{R}^{n \times n}$ and $B \in \mathbb{R}^{n \times m}$; the sequences $\mathbf{x} = (x_0, \dots, x_N)$ and $\mathbf{u} = (u_0, \dots, u_{N-1})$ are the predicted states $x_j \in \mathbb{R}^n$ and inputs $u_j \in \mathbb{R}^m$ throughout the prediction horizon; $x_r \in \mathbb{R}^n$ and $u_r \in \mathbb{R}^m$ are the given state and input reference; $x(t) \in \mathbb{R}^n$ is the system state at the current discrete time instant $t$; $\mathcal{Y}_j \subseteq \mathbb{R}^n \times \mathbb{R}^m$ are the constraints for the prediction step $j$; $\ell_j(\cdot)$ and $V_t(\cdot)$ are the *stage costs* for each prediction step $j$ and the *terminal cost* functions, respectively, which in general are functions of the predicted states $\mathbf{x}$ and inputs $\mathbf{u}$ and are parametrized by the reference $(x_r, u_r)$; and $\mathcal{X}_t$ is the *terminal set*. Constraint (4.7e) is known as the *terminal constraint*.

The MPC formulations that we discuss and present in the following chapters loosely fall within the above general formulation, although some of them include additional decision variables and constraints. However, we always consider that the functions $\ell_j$ and $V_t$ of the *cost function* $J$ (4.7a) are real valued convex functions and that $\mathcal{X}_t \subseteq \mathbb{R}^n$ is a closed convex set. Therefore, problem (4.7) is a convex optimization problem.

We denote an optimal solution of problem (4.7) (assuming that one exists) by $\mathbf{x}^* = (x_0^*, x_1^*, \dots, x_N^*)$ and $\mathbf{u}^* = (u_0^*, u_1^*, \dots, u_{N-1}^*)$. The MPC control law, i.e., the control action $u(t)$ to be applied to the system each discrete time instant $t$, is given by $u(t) = u_0^*$. At the next sample time $t+1$, the MPC optimization problem is recalculated taking the new system state $x(t+1)$ in (4.7c), and so forth.

The following two definitions play an important role in the stability of MPC.

**Definition 4.5** (Feasibility region)**.** Consider an MPC formulation defined by some optimization problem parametrized by the current system state at time instant $t \in \mathbb{Z}$. Its *feasibility region* is the set of states $x(t) \in \mathbb{R}^n$ for which said optimization problem is feasible.

**Definition 4.6** (Domain of attraction)**.** Consider an MPC formulation parametrized by the current system state at time instant $t \in \mathbb{Z}$ and whose objective is to steer the system to a reference $x_r \in \mathbb{R}^n$. Its *domain of attraction* is the set of states $x(t) \in \mathbb{R}^n$ such that the MPC controller steers the closed-loop system to the reference $x_r$ while satisfying the constraints.

We will often consider the following standard assumption.

**Assumption 4.7.** Let Assumtion 4.2 hold and assume that

*(i)* $(x_r, u_r)$ is an admissible steady state of system (4.1) under (4.2).

*(ii)* $x(t)$ strictly belongs to the feasibility region (Definition 4.5) of the MPC formulation, i.e., there exist $\mathbf{x}$ and $\mathbf{u}$ satisfying (4.7b) and (4.7c) for which (4.7d) and (4.7e) are strictly satisfied.

In MPC, the terminal set $\mathcal{X}_t$ is typically taken as an admissible invariant set (see Definition N.8) of system (4.1) for constraints (4.2) under some *terminal control law*. That is, for some predefined control law, any state belonging to $\mathcal{X}_t$ will remain in $\mathcal{X}_t$. The use of a terminal admissible invariant set provides the MPC controller with stability guarantees if its other ingredients are suitably designed [64].

The control objective of the MPC controller is to steer the system to the reference $(x_r, u_r)$ whilst satisfying the system constraints (4.2), which will only be possible if $(x_r, u_r)$ is an admissible steady state of (4.1) under (4.2).

**Remark 4.8.** *We note that our control objective, as stated above, is to steer the system to an admissible steady state of the linear model* (4.1) *whilst satisfying* (4.2). *Since this model is a linearization of a (possibly) non-linear model, our control objective does not guarantee that the real system will converge to the desired set point, unless additional ingredients are added to the overall control architecture [8, 74]. During this dissertation we will not delve into this topic any further, since our main interest is in the development of MPC formulations and solvers suitable for their implementation in embedded systems, and not on the additional ingredients required for their use in a real environment.*

## 4.3 A brief state of the art

The implementation of linear MPC in embedded systems is an extensively researched topic in the field of control, with many publications tackling different embedded systems, MPC formulations and approaches. This section presents a (non-exhaustive) review of the state of the art. We focus on the implementation of linear MPC, instead of non-linear, because the solvers proposed in this dissertation all fall into this paradigm. There are also multiple results on the implementation of nonlinear MPC in embedded systems [75, 76, 77, 78], as well as multiple tools and solvers [79, 80]. Additionally, we focus on the implementation of *nominal* MPC, although results of the implementation of *robust* MPC are also available, including the *min-max* robust MPC approach [81, 82, 83]. Henceforth, we will refer to "linear MPC" simply as "MPC", unless specifically stated otherwise.

The implementation of MPC in embedded systems can be divided into two main paradigms: *explicit* MPC, and *online* MPC.

In explicit MPC [84, 85], the solution of the MPC control law is computed offline and stored in the embedded system. In [86] the authors show that the control law of standard MPC subject to state and input constraints is given by a piecewise affine and continuous function of the state. This function can be computed offline and stored as a look-up table. Therefore, no optimization problem has to be solved online, since the control action is determined by evaluating the look-up table. This approach has been used in various publications to implement MPC in embedded systems [87, 88, 89, 90, 91]. However, its drawback is that the size of the look-up table, and thus the memory requirements and computational time needed to evaluate it, becomes prohibitively large for medium to large-sized systems and/or for problems involving many constraints.

The second main paradigm is to solve the MPC's optimization problem at each sample time for the current system state. A large portion of the results in this field come from the recent development of solvers that are tailored to their use in embedded systems, which can be used to solve the MPC optimization problem online. Many of them generate library-free code that is tailored to the problem to be solved, thus resulting in a rather efficient implementation. There are a wide range of solvers, each particularized to a certain class of optimization problem and based on a particular optimization method. A few of the most noteworthy and well known solvers that have been used for the embedded implementation of MPC are the following:

- CVXGEN [92] is an web-based "code generator for convex optimization problems that can be reduced to quadratic programs" [92] and that is better suited for small scale problems. It generates library-free C code tailored to the specific QP problem, which is solved using a primal-dual interior point method.

- qpOASES [93, 94], is an open-source C++ QP solver based on a parametric active set method that, according to the authors, is "particularly suited for model predictive control applications" [94].

- OSQP [35] is open-source sparse operator splitting solver for QP problems that is based on the ADMM algorithm. The solver is written in C, but has interfaces to various other programming languages.

- ODYS [95] is a proprietary solver for QP problems with a specific focus on MPC. It has been used to implement MPC in various industries, such as, for instance, in automotive mass production [69].

- FiOrdOs [96] is a Matlab toolbox for automated C-code generation of first order methods for parametric convex programs that is well suited for its application to MPC. The optimization problem is solved using the gradient method or the fast gradient method.

- FORCES [97] is a proprietary solver that generates code for convex optimization problems that is based on interior point methods and is well suited for MPC.

- qpDUNES [98] is an C-based open-source QP solver based on a dual Newton strategy.

Some examples of these tools being used to implement MPC in industrial embedded systems include [99], where qpOASES is implemented in a PLC to control a MISO system; [100], where qpOASES, CVXGEN and FiOrdOs are implemented in a PLC and compared (see also [101] for another implementation using FiOrdOs); or [102] which compares the scalability of FiOrdOs and qpOASES in a PLC. A very useful tool for building optimizations problems and solving them with the above solvers (and many others) is the YALMIP toolbox for Matlab [103].

The above solvers are not specific to MPC, although some consider and/or provide tools for it. As such, even though they can be used to solve a wider class of optimization problems, their use for MPC has been reported in numerous publications. However, another approach followed by some authors is to develop solvers that are tailored or more focused on MPC. These solvers may attain a more efficient implementation due to their narrow focus. Some noteworthy examples include:

- $\mu$AO-MPC [104] is a code generation tool for MPC subject to box constraints on the inputs. It generates self-contained code tailored to the specific MPC problem and system. The optimization problem is solved using the augmented Lagrangian method along with the fast gradient method.

- PRESAS [105] is a primal active-set solver for block-sparse QPs with a particular focus on MPC (both linear and non-linear).

- HPMPC [106] is a C library for the implementation of MPC. Various routines are used, depending on if the MPC is unconstrained or constrained (Riccati recursion, interior point and ADMM). Its successor, HPIPM [107], which is currently under active development, focuses on interior point methods for QPs, with a particular focus on MPC.

- SPCIES [4] is the Matlab toolbox that contains the solvers that have been developed as part of this dissertation.

We now briefly describe other noteworthy publications in this field that don't quite fall into the above paradigms, either because they implement MPC in embedded systems but without presenting a fully-fledged solver, or because they present novel methods/approaches for its implementation.

In [108] presents a primal barrier interior-point method that exploits the structure of the MPC optimization problem, achieving small computation times. In [109], an online implementation of embedded MPC is presented, where the MPC optimization problem is solved using a combination of Nesterov's fast gradient method and method of multipliers. A parallel coordinate descent algorithm for network systems is presented in [110] and applied to MPC in a PLC. In [111], an accelerated dual gradient projection algorithm for embedded MPC with an iteration complexity that grows linearly with the prediction horizon is presented. The article [112] presents the implementation of MPC in FPGAs using the fast gradient method and ADMM. In [113], Nesterov's fast gradient method is used to solve MPC in a PLC. The application of an accelerated dual proximal algorithm is proposed in [114] to finding the solution of the infinite horizon constrained LQR problem. A tool for software-hardware code generation of operator splitting methods for MPC in FPGAs is presented in [115], In [116], a primal-dual iterative algorithm particularized to MPC is developed for its implementation in FPGAs. In [117], MPC is programmed on a FPGA using high-level synthesis tools along with $\mu$AO-MPC. A primal active-set method for bounded-variable least-squares problems [118] is implemented in a PLC to solve MPC in [15]. An ADMM algorithm for MPC which exploits the problem structure by decoupling the system dynamics is presented in [119]. Decoupling the system dynamics allows for a parallel implementation of the algorithm. In [120], another ADMM based solver for MPC that decouples the system dynamics is presented, in this case by exploiting the symmetry of the system.

As evidences from the above references, research in this field is very extensive, encompassing a wide range of embedded systems, optimization methods and approaches. For additional references we refer the reader to the surveys [121, 122, 123, 124].

# Chapter 5

# Sparse solvers for model predictive control

This chapter considers system (4.1) subject to the box constraints (4.3).

This chapter presents sparse solvers for various MPC formulations in which the choice of formulation and optimization algorithm leads to an efficient solver that exploit the structure of the problem. That is, we define the MPC formulations and the decision variables in such a way that the ingredients of the first order methods (i.e., the resulting matrices) are simple enough to be able to perform matrix-vector operations without having to use the typical sparse matrix representations (e.g., the *compressed sparse column/row* or *dictionary of keys* formats). The advantage of doing this is that we do not need to store all the non-zero elements of the matrices nor the arrays that determine their positions within them. Instead, we only need to store the repeating sub-matrices once, and the sparse matrix-vector operations are performed by direct identification of their inner structure. This results in solvers with low iteration complexity and memory footprint; although, in some aspects, at the expense of some restrictions.

The solvers presented in this chapter have been included in the SPCIES toolbox for Matlab[1] [4], which is available at `https://github.com/GepocUS/Spcies`. We do not provide the exact pseudocode of the solvers, since the specific ways in which we perform the sparse matrix-vector operations, as well as certain details such as what variables are stored, are subject to future changes and improvements. Instead, we describe the pseudocode in more general terms. We show how the sparse matrix-vector operations arise, the inner structure of the matrices, and providing enough information to understand which variables are required and which are not. The reader can find the exact code of the solvers in [4].

[1]For future reference, the version of the toolbox as of the date of presentation of this dissertation is `v0.3.2`. Future versions may have made improvements upon the results shown here.

We start this chapter with three sections that provide the foundation for the sparse MPC solvers. Section 5.1 presents three solvers that will be used repeatedly throughout this chapter: one for a particular class of systems of equations, one for a particular class of equality-constrained QP problems and one for a particular class of box-constrained separable QP problems. Sections 5.2 and 5.3 show how the FISTA and ADMM algorithm can be applied to solve QP problems.

The ideas presented in the aforementioned section will be employed in the development of sparse solvers for the four MPC formulations described in Sections 5.4 to 5.6.

## 5.1 Various structure-exploiting solvers

This section presents three rather evident results that will be used throughout this chapter: a sparse solver for a particular class of systems of equations, a way to solve a particular class of equality-constrained QP problems, and, for future reference, the algorithm for solving a diagonal strongly-convex QP problem subject to box constraints.

### 5.1.1 Solving systems of equations with banded decomposition

Let us consider a matrix $W$ satisfying the following assumption.

**Assumption 5.1.** Matrix $W$ satisfies:

*(i)* $W \in \mathbb{S}^{m_z}_{++}$.

*(ii)* The Cholesky decomposition of $W$, i.e., the upper triangular matrix $W_c$ satisfying $W = W_c^\top W_c$, has the following structure:

$$W_c = \begin{pmatrix} \beta^1 & \alpha^1 & \mathbf{0} & .. & .. & \mathbf{0} \\ \mathbf{0} & \beta^2 & \alpha^2 & \mathbf{0} & .. & .. \\ .. & .. & .. & .. & .. & \mathbf{0} \\ .. & .. & .. & \mathbf{0} & \beta^{N-1} & \alpha^{N-1} \\ \mathbf{0} & .. & .. & .. & \mathbf{0} & \beta^N \end{pmatrix}, \tag{5.1}$$

where $\beta^i \in \mathbb{R}^{n\times n}$, for $i \in \mathbb{Z}_1^N$, are upper triangular matrices, and $\alpha^i \in \mathbb{R}^{n\times n}$, for $i \in \mathbb{Z}_1^{N-1}$. Obviously, $m_z = Nn$.

Algorithm 11 shows how to sparsely solve a system of equations $Wz = w$ where $W$ satisfies Assumtion 5.1. It performs forward and backward substitution using the Cholesky factorization of $W$. That is, it sequentially solves the two following triangular systems of equations

$$W_c^\top \hat{z} = w, \quad W_c z = \hat{z},$$

using the auxiliary variable $\hat{z}$.

**Algorithm 11:** Sparse solver for $Wz = w$ under Assumtion 5.1

**Prototype:** $z \leftarrow \text{solve_W}(w)$

$z \leftarrow w$
`// Forward substitution:`
**for** $j \leftarrow 1$ **to** $n$ **do** `// Compute first` $n$ `elements`
  **for** $i \leftarrow 1$ **to** $j-1$ **do** $z_{(j)} \leftarrow z_{(j)} - \hat{\beta}^{1}_{(i,j)} z_{(i)}$ **end for**
  $z_{(j)} \leftarrow \hat{\beta}^{1}_{(j,j)} z_{(j)}$
**end for**
**for** $k \leftarrow 1$ **to** $N-1$ **do** `// Compute the rest of the elements`
  **for** $j \leftarrow 1$ **to** $n$ **do**
    **for** $i \leftarrow 1$ **to** $n$ **do**
      $z_{(j+kn)} \leftarrow z_{(j+kn)} - \alpha^{k}_{(i,j)} z_{(i+(k-1)n)}$
    **end for**
    **for** $i \leftarrow 1$ **to** $j-1$ **do**
      $z_{(j+kn)} \leftarrow z_{(j+kn)} - \hat{\beta}^{k+1}_{(i,j)} z_{(i+kn)}$
    **end for**
    $z_{(j+kn)} \leftarrow \hat{\beta}^{k+1}_{(j,j)} z_{(j+kn)}$
  **end for**
**end for**
`// Backwards substitution:`
**for** $j \leftarrow n$ **to** $1$ **by** $-1$ **do** `// Compute last` $n$ `elements`
  **for** $i \leftarrow n$ **to** $j+1$ **by** $-1$ **do**
    $z_{(j+(N-1)n)} \leftarrow z_{(j+(N-1)n)} - \hat{\beta}^{N}_{(j,i)} z_{(i+(N-1)n)}$
  **end for**
  $z_{(j+(N-1)n)} \leftarrow \hat{\beta}^{N}_{(j,j)} z_{(j+(N-1)n)}$
**end for**
**for** $k \leftarrow N-2$ **to** $0$ **by** $-1$ **do** `// Compute the rest of the elements`
  **for** $j \leftarrow n$ **to** $1$ **by** $-1$ **do**
    **for** $i \leftarrow n$ **to** $1$ **by** $-1$ **do**
      $z_{(j+kn)} \leftarrow z_{(j+kn)} - \alpha^{k+1}_{(j,i)} z_{(i+(k+1)n)}$
    **end for**
    **for** $i \leftarrow n$ **to** $j+1$ **by** $-1$ **do**
      $z_{(j+kn)} \leftarrow z_{(j+kn)} - \hat{\beta}^{k+1}_{(j,i)} z_{(i+kn)}$
    **end for**
    $z_{(j+kn)} \leftarrow \hat{\beta}^{k+1}_{(j,j)} z_{(j+kn)}$
  **end for**
**end for**

**Output:** $z$

Algorithm 11 has been developed to have a small memory and computational footprint by taking advantage of the structure of $W_c$. In particular, variable $z$ is used to store $w$, $\hat{z}$ and $z$ itself, thus only requiring the declaration of a single array of dimension $m_z$, instead of three. That is, in a real implementation, the vector $w$ would be passed by reference, and the function would return the solution $z$ in the same array. To avoid performing divisions, which are more computationally expensive than multiplications, we store and use the matrices $\hat{\beta}^i$ given by

$$\hat{\beta}^i_{(j,k)} = \begin{cases} \beta^i_{(j,k)} & \text{if } j \neq k, \\ \dfrac{1}{\beta^i_{(j,k)}} & \text{if } j = k, \end{cases}$$

instead of matrices $\beta^i$.

**Remark 5.2.** *Only the matrices $\hat{\beta}^i$, for $i \in \mathbb{Z}_1^N$, and $\alpha^i$, for $i \in \mathbb{Z}_1^{N-1}$, need to be stored in a real implementation of Algorithm 11 (as well as an array of dimension $m_z$). That is, the position of the non-zero elements of $W_c$ do not need to be stored, as would be the case if $W_c$ where to be stored using the typical sparse matrix representations, such as, for instance, the compressed sparse column format.*

### 5.1.2 Solving equality-constrained QPs with banded structure

We now present a procedure for sparsely solving a particular class of equality-constrained QP problems by making use of Algorithm 11.

Let us consider the equality-constrained QP problem

$$\min_{z \in \mathbb{R}^{n_z}} \; z^\top H z + q^\top z \tag{5.2a}$$

$$s.t. \; Gz = b, \tag{5.2b}$$

under the following assumption.

**Assumption 5.3.** Assume that $H \in \mathbb{S}^{n_z}_{++}$, $G \in \mathbb{R}^{m_z \times n_z}$ satisfies $\text{rank}(G) = m_z$, $m_z < n_z$, and that $W = GH^{-1}G^\top$ satisfies Assumtion 5.1.

**Remark 5.4.** *Note that matrix $W$ defined in Assumtion 5.3 satisfies $W \in \mathbb{S}^{m_z}_{++}$ due to the assumptions on $H$ and $G$.*

Problem (5.2) has an explicit solution given by the following proposition, taken from [23, §10.1.1].

**Algorithm 12:** Solver for equality-constrained QP

**Prototype:** $z^* \leftarrow \text{solve_eqQP}(q, b; H, G)$

**Input:** $q \in \mathbb{R}^{n_z}$, $b \in \mathbb{R}^{m_z}$

1 $\mu \leftarrow -GH^{-1}q - b$

2 $\mu \leftarrow \text{solve_W}(\mu)$

3 $z^* \leftarrow -H^{-1}(G^\top \mu + q)$

**Output:** $z^*$

**Proposition 5.5** (Explicit solution of equality-contrained QP problem)**.** Consider a QP problem (5.2), where $H \in \mathbb{S}^{n_z}_{++}$, $q \in \mathbb{R}^{n_z}$, $G \in \mathbb{R}^{m_z \times n_z}$ and $b \in \mathbb{R}^{m_z}$. A vector $z^* \in \mathbb{R}^{n_z}$ is an optimal solution of this optimization problem if and only if there exists a vector $\mu \in \mathbb{R}^{m_z}$ such that

$$
\begin{aligned}
&Gz^* = b, \\
&Hz^* + q + G^\top \mu = 0,
\end{aligned}
$$

which using simple algebra and defining $W \doteq GH^{-1}G^\top$, leads to

$$W\mu = -(GH^{-1}q + b) \tag{5.3a}$$

$$z^* = -H^{-1}(G^\top \mu + q). \tag{5.3b}$$

Therefore, a solution to (5.2) under Assumtion 5.3 can be obtained by first solving (5.3a), which can be sparsely solved using Algorithm 11, and then evaluating (5.3b). We formally state this in the following corollary, where we consider (5.2) to be parametrized by $q \in \mathbb{R}^{n_z}$ and $b \in \mathbb{R}^{m_z}$.

**Corollary 5.6.** Problem (5.2) parametrized by $q \in \mathbb{R}^{n_z}$ and $b \in \mathbb{R}^{m_z}$ can be solved under Assumtion 5.3 using Algorithm 12, where step 2 calls Algorithm 11 for the matrix $W = GH^{-1}G^\top$.

**Remark 5.7.** *Algorithm 12 is a sparse solver if the matrix-matrix and matrix-vector operations in steps 1 and 3 can be performed sparsely, which, at the very least, requires $H^{-1}$ and $G$ to be sparse.*

### 5.1.3 Solving box-constrained separable QPs

This section presents a straightforward result for the sole purpose of facilitating future developments: how to solve QP problems with a diagonal Hessian subject to box constraints. Let us consider the box-constrained QP problem

$$\min_{z \in \mathbb{R}^{n_z}} \; z^\top H z + q^\top z \tag{5.4a}$$

$$s.t.\; \underline{z} \leq z \leq \overline{z}, \tag{5.4b}$$

under the following assumption.

---

**Algorithm 13:** Solver for box-constrained QP

---

**Prototype:** $z^* \leftarrow \text{solve_boxQP}(q; H, \underline{z}, \overline{z})$
**Input:** $q \in \mathbb{R}^{n_z}$
1 **for** $j \leftarrow 1$ **to** $n_z$ **do**
2 | $z^*_{(j)} \leftarrow \max\left\{\min\left\{-H^{-1}_{(j,j)} q_{(j)},\, \overline{z}_{(j)}\right\},\, \underline{z}_{(j)}\right\}$
3 **end for**
**Output:** $z^*$

---

**Assumption 5.8.** Assume that $H \in \mathbb{D}^{n_z}_{++}$ and that $\underline{z} \in \mathbb{R}^{n_z}$, $\overline{z} \in \mathbb{R}^{n_z}$ satisfy $\underline{z} \leq \overline{z}$.

Then, the components of $z$ are decoupled. It is clear that the optimal solution $z^*$ of (5.4) is given by

$$z^*_{(j)} = \max\left\{\min\left\{-H^{-1}_{(j,j)} q_{(j)},\, \overline{z}_{(j)}\right\},\, \underline{z}_{(j)}\right\},\; j \in \mathbb{Z}_1^{n_z}.$$

Algorithm 13 solves (5.4) under Assumtion 5.8 using the above expression, where we consider (5.4) to be parametrized by $q \in \mathbb{R}^{n_z}$.

**Remark 5.9.** *We note that only the inverse of the diagonal elements of $H$ are required in Algorithm 13. Thus, they can be stored in an array of dimension $n_z$.*

**Remark 5.10.** *If $H$ is given by $H = \rho \mathbf{I}_{n_z}$ for some $\rho \in \mathbb{R}_{>0}$, then we will use the slight abuse of notation $z^* \leftarrow$ solve_boxQP$(q; \rho, \underline{z}, \overline{z})$ as the prototype of Algorithm 13. In this case, step 2 of the algorithm will be substituted by*

$$z^*_{(j)} \leftarrow \max\left\{\min\left\{-\rho^{-1} q_{(j)},\, \overline{z}_{(j)}\right\},\, \underline{z}_{(j)}\right\}.$$

## 5.2 Solving QPs with FISTA through duality

Let us consider a QP problem

$$\min_{z} \left\{ J(z) \doteq \frac{1}{2} z^\top H z + q^\top z \right\} \tag{5.5a}$$

$$s.t.\; Gz = b \tag{5.5b}$$

$$z \in \mathcal{Z}, \tag{5.5c}$$

where $z \in \mathbb{R}^{n_z}$ are the *decision variables*, $J : \mathbb{R}^{n_z} \to \mathbb{R}$ is the *cost function*, $q \in \mathbb{R}^{n_z}$, $G \in \mathbb{R}^{m_z \times n_z}$ and $b \in \mathbb{R}^{m_z}$; under the following assumption.

**Assumption 5.11.** We assume that:

*(i)* The Hessian $H \in \mathbb{S}^{n_z}_{++}$.

*(ii)* $\mathcal{Z} \subseteq \mathbb{R}^{n_z}$ is a non-empty closed convex set.

*(iii)* $G$ satisfies $\mathrm{rank}(G) = m_z$.

*(iv)* There exist $\hat{z} \in \mathrm{ri}(\mathcal{Z})$ such that $G\hat{z} = b$.

We note that problem (5.5) is solvable under Assumtion 5.11.*(iv)*. Furthermore, since $H$ is assumed to be positive definite, problem (5.5) under Assumtion 5.11 is a convex optimization problem with a strongly convex cost function. Therefore, it is well known that it has a unique optimal solution $z^*$.

In this section we show how to solve problem (5.5) under Assumtion 5.11 by applying FISTA (Algorithm 4) to its dual formulation, which we now introduce.

### 5.2.1 QP problem's dual formulation

Let us consider the following Lagrangian function $\mathcal{L} : \mathbb{R}^{n_z} \times \mathbb{R}^{m_z} \to \mathbb{R}$ of (5.5):

$$\mathcal{L}(z, \lambda) = J(z) + \langle \lambda, b - Gz \rangle,$$

where $\lambda \in \mathbb{R}^{m_z}$ are the dual variables for the equality constraints (5.5b). Then, the *dual function* $\psi : \mathbb{R}^{m_z} \to \mathbb{R}$ is given by

$$\psi(\lambda) = \min_{z \in \mathcal{Z}} \mathcal{L}(z, \lambda), \tag{5.6}$$

and the *dual problem* by

$$\max_{\lambda \in \mathbb{R}^{m_z}} \psi(\lambda). \tag{5.7}$$

Let the optimal solutions of problems (5.5) and (5.7) be given by $z^*$ (*primal optimal solution*) and $\lambda^*$ (*dual optimal solution*), respectively. Then the *primal optimal value* and *dual optimal value* are given by $J^* = J(z^*)$ and $\psi^* = \psi(\lambda^*)$.

The following proposition states the necessary and sufficient conditions for optimality of problems (5.5) and (5.7).

**Proposition 5.12.** Consider the primal and dual problems (5.5) and (5.7). Let Assumtion 5.11 hold. Then, $J^* = \psi^*$, and $(z^*, \lambda^*)$ are a primal and dual optimal solution pair if and only if $z^*$ is feasible (i.e., $z^* \in \mathcal{Z}$ and $Gz^* = b$) and

$$z^* = \arg\min_{z \in \mathcal{Z}} \mathcal{L}(z, \lambda^*). \tag{5.8}$$

Moreover, at least one dual optimal solution $\lambda^*$ exists.

***Proof:*** This proposition follows directly from [1, Proposition 5.3.3], where the fact that $J^*$ is finite and that (5.8) has a unique optimal solution follow from $J$ being strongly convex (note Assumtion 5.11.*(i)* states that $H \in \mathbb{S}^{n_z}_{++}$). ■

For a given $\lambda$, let $z_\lambda \in \mathbb{R}^{n_z}$ be defined as

$$z_\lambda = \arg\min_{z \in \mathcal{Z}} \mathcal{L}(z, \lambda) \tag{5.9}$$

It is obvious that $\phi(\lambda) = \mathcal{L}(z_\lambda, \lambda)$. Since $J$ is a strongly convex function, we have that, under Assumtion 5.11.*(iv)*, the dual function $\psi$ is a continuously differentiable smooth concave function [125, Lemma 3.1] whose gradient is given by [19, Example 2.2.4]

$$\nabla\psi(\lambda) = b - Gz_\lambda. \tag{5.10}$$

The following proposition characterizes the smoothness of $\psi$ in terms of the matrix $W \doteq GH^{-1}G^\top$ (c.f., Assumtion 2.1.*(ii)*). Its proof is based on the following lemma. Note that, since $H \in \mathbb{S}^n_{++}$ and $\text{rank}(G) = m_z$, we have that $W \in \mathbb{S}^{m_z}_{++}$.

**Lemma 5.13.** Consider the primal and dual problems (5.5) and (5.7), where the dual function $\psi$ is given by (5.6). Let Assumtion 5.11 hold and $z_\lambda$ be given by (5.9). Assume that $\Delta z \in \mathbb{R}^{n_z}$ is such that $z_\lambda + \Delta z \in \mathcal{Z}$. Then,

$$\langle \Delta z, Hz_\lambda + q - G^\top\lambda \rangle \geq 0.$$

***Proof:*** From the convexity of $\mathcal{Z}$, have that $(1-\mu)z_\lambda + \mu(z_\lambda + \Delta z) \in \mathcal{Z}$, $\forall \mu \in [0, 1]$, (see Definition N.2) which leads to $z_\lambda + \mu\Delta z \in \mathcal{Z}$, $\forall \mu \in [0, 1]$. From the optimality of $z_\lambda$, we have

$$\begin{aligned}
\psi(\lambda) &= \frac{1}{2} z_\lambda^\top H z_\lambda + q^\top z_\lambda + \langle \lambda, b - Gz_\lambda \rangle \\
&\leq \frac{1}{2}(z_\lambda + \mu\Delta z)^\top H(z_\lambda + \mu\Delta z) + \langle q, z_\lambda + \mu\Delta z \rangle + \langle \lambda, b - G(z_\lambda + \Delta z) \rangle \\
&= \psi(\lambda) + \mu\langle \Delta z, Hz_\lambda + q - G^\top\lambda \rangle + \frac{\mu^2}{2}\Delta z^\top H \Delta z, \ \forall \mu \in [0, 1].
\end{aligned}$$

Therefore,

$$\mu\langle \Delta z, Hz_\lambda + q - G^\top\lambda \rangle + \frac{\mu^2}{2}\Delta z^\top H\Delta z \geq 0, \ \forall \mu \in [0, 1],$$

which is trivially satisfied if $\mu = 0$. Otherwise, dividing by $\mu$, we have that

$$\langle \Delta z, Hz_\lambda + q - G^\top\lambda \rangle + \frac{\mu}{2}\Delta z^\top H\Delta z \geq 0, \ \forall \mu \in (0, 1].$$

For this to be true $\forall \mu \in (0, 1]$, it must be true that $\langle \Delta z, Hz_\lambda + q - G^\top\lambda \rangle \geq 0$. ■

**Proposition 5.14** (Smoothness of the QP's dual function)**.** Consider the primal and dual problems (5.5) and (5.7), where the dual function $\psi$ is given by (5.6). Let Assumtion 5.11 hold and $z_\lambda$ be given by (5.9). Denote $W \doteq GH^{-1}G^\top$. Then for any $(\lambda, \Delta\lambda) \in \mathbb{R}^{m_z} \times \mathbb{R}^{m_z}$,

$$\psi(\lambda + \Delta\lambda) \geq \psi(\lambda) + \langle \Delta\lambda, b - Gz_\lambda \rangle - \frac{1}{2}\|\Delta\lambda\|_W^2.$$

***Proof:*** Let us define $\Delta z$ as

$$z_\lambda + \Delta z = \arg\min_{z\in\mathcal{Z}} \frac{1}{2} z^\top H z + q^\top z + \langle \lambda + \Delta\lambda, b - Gz\rangle.$$

It is clear that $z_\lambda + \Delta z \in \mathcal{Z}$. Then, from the definition of $\psi$, we have

$$\begin{aligned}
\psi(\lambda{+}\Delta\lambda) &= \frac{1}{2}(z_\lambda{+}\Delta z)^\top H (z_\lambda{+}\Delta z) + \langle q, z_\lambda{+}\Delta z\rangle + \langle \lambda{+}\Delta\lambda, b - G(z_\lambda{+}\Delta z)\rangle \\
&= \psi(\lambda) + \frac{1}{2}\Delta z^\top H \Delta z + \langle \Delta z, H z_\lambda + q - G^\top \lambda\rangle + \langle \Delta\lambda, b - G(z_\lambda{+}\Delta z)\rangle \\
&\overset{(*)}{\geq} \psi(\lambda) + \frac{1}{2}\Delta z^\top H \Delta z + \langle \Delta\lambda, b - G(z_\lambda{+}\Delta z)\rangle \\
&= \psi(\lambda) + \langle \Delta\lambda, b - G z_\lambda\rangle + \frac{1}{2}\Delta z^\top H \Delta z + \langle \Delta\lambda, -G\Delta z\rangle \\
&\geq \psi(\lambda) + \langle \Delta\lambda, b - G z_\lambda\rangle + \min_{\Delta z}\left\{\frac{1}{2}\Delta z^\top H \Delta z + \langle \Delta\lambda, -G\Delta z\rangle\right\} \\
&\overset{(**)}{=} \psi(\lambda) + \langle \Delta\lambda, b - G z_\lambda\rangle - \frac{1}{2}\Delta\lambda^\top G H^{-1} G^\top \Delta\lambda,
\end{aligned}$$

where $(*)$ is due to Lemma 5.13 and $(**)$ follows from

$$\min_{\Delta z}\left\{\frac{1}{2}\Delta z^\top H \Delta z + \langle \Delta\lambda, -G\Delta z\rangle\right\} = -\frac{1}{2}\Delta\lambda^\top G H^{-1} G^\top \Delta\lambda,$$

which is a well known result in the field of convex optimization. ■

### 5.2.2 Solving the QP's dual problem with FISTA

We now show how to solve the dual problem (5.7) using FISTA (Algorithm 4). The results we now present are well known and can be found in several prior publications, but are included here for completeness.

Since FISTA is expressed in terms of the minimization of a convex function, instead of maximization of a concave one, we consider the convex problem

$$-\min_{\lambda\in\mathbb{R}^{m_z}} -\psi(\lambda). \tag{5.11}$$

We note that, from (5.10), we have

$$\nabla(-\psi(\lambda)) = -(b - G z_\lambda). \tag{5.12}$$

Then, under Assumtion 5.11, and taking into consideration Proposition 5.14, we have that Assumtion 2.1 is satisfied for (5.11). In particular, problem (5.11) has no non-smooth term $\Psi$, and the matrix $R$ that characterizes its smoothness (Assumtion 2.1.*(ii)*) is $W \doteq G H^{-1} G^\top$. Therefore, the iterates of FISTA applied to problem (5.11) will converge to $\lambda^*$.

**Algorithm 14:** FISTA for solving QP problem (5.5)

**Require:** $\lambda \in \mathbb{R}^{m_z}$, $\epsilon \in \mathbb{R}_{>0}$

1 $t_0 \leftarrow 1$, $k \leftarrow 0$

2 $q_k \leftarrow q - G^\top \lambda$

3 $z_k \leftarrow \arg\min_{z \in \mathcal{Z}} \frac{1}{2} z^\top H z + q_k^\top z$

4 $\Gamma_k \leftarrow -(Gz_k - b)$

5 $\Delta\lambda_k \leftarrow$ solution of $W\Delta\lambda = \Gamma_k$

6 $y_k \leftarrow \lambda + \Delta\lambda_k$

7 $\lambda_k \leftarrow \lambda + \Delta\lambda_k$

8 **repeat**

9 $\quad k \leftarrow k + 1$

10 $\quad q_k \leftarrow q - G^\top y_{k-1}$

11 $\quad z_k \leftarrow \arg\min_{z \in \mathcal{Z}} \frac{1}{2} z^\top H z + q_k^\top z$

12 $\quad \Gamma_k \leftarrow -(Gz_k - b)$

13 $\quad \Delta\lambda_k \leftarrow$ solution of $W\Delta\lambda = \Gamma_k$

14 $\quad \lambda_k \leftarrow \Delta\lambda_k + y_{k-1}$

15 $\quad t_k \leftarrow \frac{1}{2}\left(1 + \sqrt{1 + 4t_{k-1}^2}\right)$

16 $\quad y_k \leftarrow \lambda_k + \frac{t_{k-1} - 1}{t_k}(\lambda_k - \lambda_{k-1})$

17 **until** $\|\Gamma_k\|_\infty \leq \epsilon$

**Output:** $\tilde{z}^* \leftarrow z_k$, $\tilde{\lambda}^* \leftarrow y_{k-1}$

Let us now take a closer look at step 4 of Algorithm 4, which performs the assignment $\lambda_k \leftarrow \mathcal{T}_W^{-\psi, \mathbb{R}^{m_z}}(y_{k-1})$. That is,

$$\begin{aligned}\lambda_k &\leftarrow \arg\min_{\lambda \in \mathbb{R}^{m_z}} \langle \nabla(-\psi(y_{k-1})), \lambda - y_{k-1} \rangle + \frac{1}{2}\|\lambda - y_{k-1}\|_W^2 \\ &= \arg\min_{\lambda \in \mathbb{R}^{m_z}} \langle Gz_{y_{k-1}} - b, \lambda - y_{k-1} \rangle + \frac{1}{2}\|\lambda - y_{k-1}\|_W^2,\end{aligned}$$

which is a strongly convex unconstrained QP problem. Therefore, $\lambda_k$ is the solution of the system of equations

$$W(\lambda_k - y_{k-1}) = -(Gz_{y_{k-1}} - b).$$

By defining $\Delta\lambda_k = \lambda_k - y_{k-1}$, step 4 of Algorithm 4 reduces to first computing $z_{y_{k-1}}$, then solving the system of equations $W\Delta\lambda_k = -(Gz_{y_{k-1}} - b)$, and then performing the assignment $\lambda_k \leftarrow \Delta\lambda_k + y_{k-1}$.

FISTA algorithm applied to the dual problem (5.11) is shown in Algorithm 14. Note that its exit condition (step 17) does not correspond to the exit condition of Algorithm 4, which would be $\|y_{k-1} - \lambda_k\|_W \leq \epsilon$. We use the exit condition $\|\Gamma_k\|_\infty \leq \epsilon$ because it is computationally cheaper to evaluate, since $\Gamma_k \in \mathbb{R}^{m_z}$ is

**Algorithm 15:** ADMM for solving QP problem (5.13)

**Require:** $v_0 \in \mathbb{R}^{n_z}$, $\lambda_0 \in \mathbb{R}^{n_z}$, $\rho \in \mathbb{R}_{>0}$, $\epsilon_p \in \mathbb{R}_{>0}$, $\epsilon_d \in \mathbb{R}_{>0}$

$k \leftarrow 0$

**repeat**

  $q_k \leftarrow q + \lambda_k - \rho v_k$

  $z_{k+1} \leftarrow \min_{z} \frac{1}{2} z^\top H_\rho z + q_k^\top z, \ s.t. \ Gz = b$

  $\hat{q}_k \leftarrow -\rho z_{k+1} - \lambda_k$

  $v_{k+1} \leftarrow \arg\min_{v \in \mathcal{Z}} \frac{\rho}{2} v^\top v + \hat{q}_k^\top v$

  $\lambda_{k+1} \leftarrow \lambda_k + \rho(z_{k+1} - v_{k+1})$

  $k \leftarrow k + 1$

**until** $\|z_k - v_k\|_\infty \leq \epsilon_p$ **and** $\|v_k - v_{k-1}\|_\infty \leq \epsilon_d$

**Output:** $\tilde{z}^* \leftarrow z_k$, $\tilde{v}^* \leftarrow v_k$, $\tilde{\lambda}^* \leftarrow \lambda_k$

already computed in step 12. The fact that $\Gamma_k$ is a measure of optimality of the current iterate follows from $\Gamma_k = -(Gz_k - b) = -(Gz_{y_{k-1}} - b) \overset{(5.12)}{=} -\nabla(-\psi(y_{k-1}))$. Then, since (5.11) is an unconstrained smooth convex problem, we have that $\lambda^* = y_{k-1} \iff \nabla(-\psi(y_{k-1})) = 0$. That is, $\lambda^* = y_{k-1} \iff \Gamma_k = 0$. Additionally, from Proposition 5.12, we have that $z^* = z_{\lambda^*}$. Thus, the outputs of the algorithm are the suboptimal dual solution $\tilde{\lambda}^* = y_{k-1}$ and the suboptimal primal solution $\tilde{z}^* = z_{\tilde{\lambda}^*} = z_{y_{k-1}} = z_k$, where the suboptimality is determined by $\epsilon$.

**Remark 5.15.** *Note that steps 3 and 11 of Algorithm 14 require solving a constrained QP problem, which in general is not trivial. However, in the following sections we will consider assumptions under which this problem will have a simple and explicit solution. Additionally, steps 5 and 13 require solving a system of equations, which can be computationally demanding. However, we will consider systems of equations that satisfy Assumtion 5.1, and may therefore be sparsely solved using Algorithm 11.*

## 5.3 Solving QPs with ADMM

This section explains how problem (5.5) under Assumtion 5.11 can be solved using ADMM (Algorithm 2). The results we now present are well known and can be found in several prior publications, but are included here for completeness.

We start by rewriting (5.5) into a problem of class (2.12) by taking

$$\min_{z,v} \frac{1}{2} z^\top H z + q^\top z \tag{5.13a}$$

$$s.t. \ Gz = b \tag{5.13b}$$

$$v \in \mathcal{Z} \tag{5.13c}$$

$$z - v = \mathbf{0}_{n_z}. \tag{5.13d}$$

That is, the ingredients of (2.12) are given by: $C = \mathbf{I}_{n_z}$, $D = -\mathbf{I}_{n_z}$, $d = \mathbf{0}_{n_z}$,

$$
\begin{aligned}
f(z) &= \frac{1}{2} z^\top H z + q^\top z + \mathcal{I}_{Gz=b}(z), \\
g(v) &= \mathcal{I}_{\mathcal{Z}}(v),
\end{aligned}
$$

where $\mathcal{I}_{Gz=b}$ is the indicator function of the set $\{\, z \in \mathbb{R}^{n_z} \,:\, Gz = b \,\}$.

Let us revisit steps 3 and 4 of Algorithm 2 when applied to (5.13). Step 3 of Algorithm 2 now reads as

$$
\begin{aligned}
z_{k+1} \leftarrow \arg\min_{z} \;& \frac{1}{2} z^\top H z + q^\top z + \frac{\rho}{2} \| z - v_k + \frac{1}{\rho} \lambda_k \|_2^2 \\
& s.t. \; Gz = b \\
= \arg\min_{z} \;& \frac{1}{2} z^\top H_\rho z + q_k^\top z \\
& s.t. \; Gz = b,
\end{aligned} \tag{5.14}
$$

where $H_\rho = H + \rho \mathbf{I}_{n_z}$ and $q_k = q + \lambda_k - \rho v_k$. Step 4 of Algorithm 2 now reads as

$$
\begin{aligned}
v_{k+1} \leftarrow \arg\min_{v \in \mathcal{Z}} \;& \frac{\rho}{2} \| z_{k+1} - v + \frac{1}{\rho} \lambda_k \|_2^2 \\
= \arg\min_{v \in \mathcal{Z}} \;& \frac{\rho}{2} v^\top v + \hat{q}_k^\top v,
\end{aligned} \tag{5.15}
$$

where $\hat{q}_k = -\rho z_{k+1} - \lambda_k$.

Algorithm 15 shows the result of particularizing Algorithm 2 to problem (5.13) following the above discussion.

**Remark 5.16.** *As in Remark 5.15, steps 4 and 6 of Algorithm 15 require solving an equality-constrained QP and a constrained QP, respectively. In general, these problems are not necessarily simple to solve. However, in this chapter we will consider assumptions under which both steps can be solved using Algorithms 12 and 13, respectively.*

## 5.4 Simple standard MPC formulations

This section considers two simple standard MPC formulations, which are particularizations of the general MPC formulation (4.7). We describe the formulation, and present sparse solvers based on the algorithms shown in Sections 5.2 and 5.3. The solvers were originally presented in [8, 9], although with a bigger emphasis on their implementation in PLCs.

The first one is an MPC formulation with a terminal equality constraint:

$$\min_{\mathbf{x},\mathbf{u}} \sum_{j=0}^{N-1} \|x_j - x_r\|_Q^2 + \|u_j - u_r\|_R^2 \tag{5.16a}$$

$$s.t.\ x_{j+1} = Ax_j + Bu_j,\ j \in \mathbb{Z}_0^{N-1} \tag{5.16b}$$

$$x_0 = x(t) \tag{5.16c}$$

$$\underline{x} \leq x_j \leq \overline{x},\ j \in \mathbb{Z}_1^{N-1} \tag{5.16d}$$

$$\underline{u} \leq u_j \leq \overline{u},\ j \in \mathbb{Z}_0^{N-1} \tag{5.16e}$$

$$x_N = x_r, \tag{5.16f}$$

and the second one is an MPC formulation without a terminal constraint:

$$\min_{\mathbf{x},\mathbf{u}} \sum_{j=0}^{N-1} \left(\|x_j - x_r\|_Q^2 + \|u_j - u_r\|_R^2\right) + \|x_N - x_r\|_T^2 \tag{5.17a}$$

$$s.t.\ x_{j+1} = Ax_j + Bu_j,\ j \in \mathbb{Z}_0^{N-1} \tag{5.17b}$$

$$x_0 = x(t) \tag{5.17c}$$

$$\underline{x} \leq x_j \leq \overline{x},\ j \in \mathbb{Z}_1^{N} \tag{5.17d}$$

$$\underline{u} \leq u_j \leq \overline{u},\ j \in \mathbb{Z}_0^{N-1}. \tag{5.17e}$$

We consider Assumtion 4.7 to hold for both the above MPC formulations. Assumptions on the *cost function matrices* $Q$, $R$ and $T$ will be stated further ahead when applicable.

The advantage of the MPC formulation (5.16) is that it avoids the computation of an admissible invariant set $\mathcal{X}_t$ (Definition N.8), which can be very demanding even for average-sized systems [126], by using the simplest one: the singleton $\{x_r\}$. To see that under Assumtion 4.7, the singleton $\{x_r\}$ is, indeed, an admissible invariant set of the system, note that the admissible terminal control law $u(t) = u_r$ keeps the system at $x_r$. This advantage, however, comes at the expense of a reduction of the domain of attraction (Definition 4.6) of the MPC controller due to the use of the smallest admissible invariant set of the system.

This drawback can be avoided by simply eliminating the terminal constraint, leading to the MPC formulation (5.17). However, this comes at the expense of losing the stability guarantees that accompany the use of a suitable terminal constraint. In particular, the controller may not stabilize every feasible initial state; only those contained within a certain region which may be difficult to characterize [127]. Thus, the feasibility region and the domain of attraction of the MPC formulation (5.17) may not be the same.

**Remark 5.17.** *We note that the MPC formulation* (5.16) *does not include a terminal cost because it does not require one for stability purposes and we chose to focus on a simple MPC formulation in order to attain solvers with low iteration complexity.*

The MPC formulations (5.16) and (5.17) can be posed as QPs (5.5) as follows.

**QP problem of the MPC formulation** (5.16)**:**

Take the decision variables as:

$$z = (u_0, x_1, u_1, x_2, u_2, \dots, x_{N-1}, u_{N-1}). \tag{5.18a}$$

Then, the resulting ingredients of the QP problem (5.5) are given by

$$H = \text{diag}(R, Q, R, Q, R, \dots, Q, R), \tag{5.18b}$$

$$q = -(Ru_r, Qx_r, Ru_r, Qx_r, Ru_r, \dots, Qx_r, Ru_r), \tag{5.18c}$$

$$G = \begin{bmatrix} B & -\mathbf{I}_n & \mathbf{0} & \cdots & \cdots & \mathbf{0} & \mathbf{0} \\ \mathbf{0} & A & B & -\mathbf{I}_n & \cdots & \mathbf{0} & \mathbf{0} \\ \mathbf{0} & \mathbf{0} & \ddots & \ddots & \ddots & \mathbf{0} & \mathbf{0} \\ \mathbf{0} & \mathbf{0} & \cdots & A & B & -\mathbf{I}_n & \mathbf{0} \\ \mathbf{0} & \mathbf{0} & \cdots & \mathbf{0} & \mathbf{0} & A & B \end{bmatrix}, \tag{5.18d}$$

$$b = (-Ax(t), \mathbf{0}_n, \mathbf{0}_n, \dots, \mathbf{0}_n, x_r), \tag{5.18e}$$

$$\mathcal{Z} = \{\, z \,:\, \underline{z} \leq z \leq \overline{z} \,\}, \tag{5.18f}$$

where

$$\underline{z} = (\underline{u}, \underline{x}, \underline{u}, \underline{x}, \underline{u}, \dots, \underline{x}, \underline{u}), \tag{5.18g}$$

$$\overline{z} = (\overline{u}, \overline{x}, \overline{u}, \overline{x}, \overline{u}, \dots, \overline{x}, \overline{u}). \tag{5.18h}$$

The dimensions of the QP problem are $n_z = (N-1)(n+m) + m$ and $m_z = Nn$.

**QP problem of the MPC formulation** (5.17)**:**

Take the decision variables as:

$$z = (u_0, x_1, u_1, x_2, u_2, \dots, x_{N-1}, u_{N-1}, x_N). \tag{5.19a}$$

Then, the resulting ingredients of the QP problem (5.5) are given by

$$H = \text{diag}(R, Q, R, Q, R, \dots, Q, R, T), \tag{5.19b}$$

$$q = -(Ru_r, Qx_r, Ru_r, Qx_r, Ru_r, \dots, Qx_r, Ru_r, Tx_r), \tag{5.19c}$$

$$G = \begin{bmatrix} B & -\mathbf{I}_n & \mathbf{0} & \cdots & \cdots & \mathbf{0} \\ \mathbf{0} & A & B & -\mathbf{I}_n & \cdots & \mathbf{0} \\ \mathbf{0} & \mathbf{0} & \ddots & \ddots & \ddots & \mathbf{0} \\ \mathbf{0} & \mathbf{0} & \mathbf{0} & A & B & -\mathbf{I}_n \end{bmatrix}, \tag{5.19d}$$

$$b = (-Ax(t), \mathbf{0}_n, \mathbf{0}_n, \dots, \mathbf{0}_n), \tag{5.19e}$$

$$\mathcal{Z} = \{\, z \,:\, \underline{z} \leq z \leq \overline{z} \,\}, \tag{5.19f}$$

where

$$\underline{z} = (\underline{u}, \underline{x}, \underline{u}, \underline{x}, \underline{u}, \dots, \underline{x}, \underline{u}, \underline{x}), \tag{5.19g}$$

$$\overline{z} = (\overline{u}, \overline{x}, \overline{u}, \overline{x}, \overline{u}, \dots, \overline{x}, \overline{u}, \overline{x}). \tag{5.19h}$$

The dimensions of the QP problem are given by $n_z = N(n+m)$ and $m_z = Nn$.

### 5.4.1 FISTA-based solver for standard MPC

We apply the results of Section 5.2 to the QP problems of the MPC formulations (5.16) and (5.17), whose ingredients are given by (5.18) and (5.19), respectively, under the following assumption.

**Assumption 5.18.** Let Assumtion 4.7 hold and assume that $Q \in \mathbb{D}^n_{++}$, $R \in \mathbb{D}^m_{++}$ and $T \in \mathbb{D}^n_{++}$.

It is easy to see that the satisfaction of Assumtion 5.18 implies the satisfaction of Assumtion 5.11 for both QP problems, where, additionally, $H \in \mathbb{D}^{n_z}_{++}$. Therefore, Algorithm 14 can be applied to them.

Let us take a closer look at the computationally expensive steps of Algorithm 14. Steps 5 and 13 of Algorithm 14 require solving $W\Delta\lambda = \Gamma_k$. However, due to the banded structures of $G$ and $H$, matrix $W = GH^{-1}G^\top$ satisfies Assumtion 5.1. Therefore, they can be solved using Algorithm 11.

Steps 3 and 11 of Algorithm 14 perform the assignment

$$\begin{aligned} z_k \leftarrow \min_z \; & \frac{1}{2} z^\top H z + q_k^\top z \\ s.t. \; & \underline{z} \leq z \leq \overline{z}, \end{aligned}$$

where

$$q_k = \begin{cases} q - G^\top \lambda & \text{if } k = 0 \\ q - G^\top y_{k-1} & \text{otherwise,} \end{cases}$$

which satisfies Assumtion 5.8, and can therefore be solved using Algorithm 13.

Algorithm 16 shows the particularization of Algorithm 14 to the QP problems of the MPC formulations (5.16) and (5.17) under Assumtion 5.18. The control action $u(t)$ is taken as the first $m$ elements of $\tilde{z}^*$.

**Remark 5.19.** *One of the key aspects of Algorithm 16 is that the matrix-matrix and matrix-vector operations can be performed sparsely without needing to store all the non-zero elements of the matrices. Instead, it only requires the computation/storage of the repeating elements once, such as matrices $A$ and $B$, for $G$; $Q^{-1}$, $R^{-1}$ and $T^{-1}$ (if applicable) for $H^{-1}$; $-Qx_r$, $-Ru_r$ and $-Tx_r$ for $q$; $-Ax(t)$ for $b$; and $\underline{x}$, $\overline{x}$, $\underline{u}$, $\overline{u}$ for $\underline{z}$ and $\overline{z}$.*

**Remark 5.20.** *We note that steps 5 and 13 of Algorithm 16 can be solved using Algorithm 13 because the Hessian is diagonal and we are considering box constraints. However, these conditions can be relaxed without greatly compromising the efficiency of the algorithm. In particular, (i) if non-diagonal cost function matrices are considered, and/or (ii) we consider coupled input-state constraints* (4.4)*, then these steps would result in $N$ decoupled small-scale, inequality-constrained QP problems, which could be individually solved using interior point or active set methods.*

**Algorithm 16:** Sparse FISTA solver for standard MPC formulations

**Require:** $x(t) \in \mathbb{R}^n$, $x_r \in \mathbb{R}^n$, $u_r \in \mathbb{R}^m$, $\lambda \in \mathbb{R}^{m_z}$, $\epsilon \in \mathbb{R}_{>0}$

1 Compute $q$ with $x_r$ and $u_r$
2 Compute $b$ with $x(t)$
3 $t_0 \leftarrow 1$, $k \leftarrow 0$
4 $q_k \leftarrow q - G^\top \lambda$
5 $z_k \leftarrow \text{solve_boxQP}(q_k; H, \underline{z}, \overline{z})$
6 $\Gamma_k \leftarrow -(Gz_k - b)$
7 $\Delta\lambda_k \leftarrow \text{solve_W}(\Gamma_k)$
8 $y_k \leftarrow \lambda + \Delta\lambda_k$
9 $\lambda_k \leftarrow \lambda + \Delta\lambda_k$
10 **repeat**
11 $\quad k \leftarrow k + 1$
12 $\quad q_k \leftarrow q - G^\top y_{k-1}$
13 $\quad z_k \leftarrow \text{solve_boxQP}(q_k; H, \underline{z}, \overline{z})$
14 $\quad \Gamma_k \leftarrow -(Gz_k - b)$
15 $\quad \Delta\lambda_k \leftarrow \text{solve_W}(\Gamma_k)$
16 $\quad \lambda_k \leftarrow \Delta\lambda_k + y_{k-1}$
17 $\quad t_k \leftarrow \frac{1}{2}\left(1 + \sqrt{1 + 4t_{k-1}^2}\right)$
18 $\quad y_k \leftarrow \lambda_k + \frac{t_{k-1} - 1}{t_k}(\lambda_k - \lambda_{k-1})$
19 **until** $\|\Gamma_k\|_\infty \leq \epsilon$

**Output:** $\tilde{z}^* \leftarrow z_k$, $\tilde{\lambda}^* \leftarrow y_{k-1}$

**Remark 5.21.** *The typical choice for the terminal cost function matrix $T$ is to take it as the solution of the discrete Riccati equation*

$$A^\top TA - T - (A^\top TB)(R + B^\top TB)^{-1}(B^\top TA) + Q = \mathbf{0}_{n \times n}.$$

*This matrix, however, is generally non-diagonal and thus cannot be used in Algorithm 16 unless the last $n$ elements of $z_k$ are updated in step 13 by solving a small-scale inequality-constrained QP problem, as discussed inn Remark 5.20.*

### 5.4.2 ADMM-based solver for standard MPC

We apply the results of Section 5.3 to the QP problems of the MPC formulations (5.16) and (5.17), whose ingredients are given by (5.18) and (5.19), respectively, under the following assumption.

**Assumption 5.22.** Let Assumtion 4.7 hold and assume that $Q \in \mathbb{S}_+^n$, $R \in \mathbb{S}_+^m$ and $T \in \mathbb{S}_+^n$.

**Algorithm 17:** Sparse ADMM solver for standard MPC formulations

**Require:** $x(t) \in \mathbb{R}^n$, $x_r \in \mathbb{R}^n$, $u_r \in \mathbb{R}^m$, $v_0 \in \mathbb{R}^{n_z}$, $\lambda_0 \in \mathbb{R}^{n_z}$, $\rho \in \mathbb{R}_{>0}$, $\epsilon_p \in \mathbb{R}_{>0}$, $\epsilon_d \in \mathbb{R}_{>0}$

1 Compute $q$ with $x_r$ and $u_r$
2 Compute $b$ with $x(t)$
3 $k \leftarrow 0$
4 **repeat**
5 $\quad q_k \leftarrow q + \lambda_k - \rho v_k$
6 $\quad z_{k+1} \leftarrow \text{solve_eqQP}(q_k, b; H_\rho, G)$
7 $\quad \hat{q}_k \leftarrow -\rho z_{k+1} - \lambda_k$
8 $\quad v_{k+1} \leftarrow \text{solve_boxQP}(\hat{q}_k; \rho, \underline{z}, \overline{z})$
9 $\quad \lambda_{k+1} \leftarrow \lambda_k + \rho(z_{k+1} - v_{k+1})$
10 $\quad k \leftarrow k + 1$
11 **until** $\|z_k - v_k\|_\infty \leq \epsilon_p$ **and** $\|v_k - v_{k-1}\|_\infty \leq \epsilon_d$

**Output:** $\tilde{z}^* \leftarrow z_k$, $\tilde{v}^* \leftarrow v_k$, $\tilde{\lambda}^* \leftarrow \lambda_k$

It is easy to see that the satisfaction of Assumtion 5.22 implies the satisfaction of Assumtion 5.11 for both QP problems. Therefore, Algorithm 15 can be applied to them. Let us take a closer look at steps 4 and 6 of Algorithm 15.

Step 4 of Algorithm 15 requires solving the QP problem (5.14), which due to the block diagonal structure of $H_\rho$ under Assumtion 5.22 and to the banded structure of $G$, satisfies Assumtion 5.3. Therefore, it can be solved using Algorithm 12 as described in Corollary 5.6.

Step 6 requires solving the QP problem (5.15), which satisfies Assumtion 5.8 and can therefore be solved using Algorithm 13.

Algorithm 17 shows the particularization of Algorithm 15 to the QP problems of the MPC formulations (5.16) and (5.17) under Assumtion 5.22. The control action $u(t)$ is taken as the first $m$ elements of $\tilde{v}^*$.

**Remark 5.23.** *We note that the reason why Assumtion 5.22 considers non-diagonal cost function matrices, when Assumtion 5.18 does not, is that step 8 of Algorithm 17 can always be solved using Algorithm 13, whereas in the FISTA-based solver the use of non-diagonal cost function matrices prevented its use. However, Remark 5.20.(ii) still also holds here for the case in which coupled input-state constraints* (4.4) *are used in place of box constraints.*

**Remark 5.24.** *Remark 5.19 also applies to Algorithm 17. That is, the matrix-matrix and matrix-vector operations are all performed sparsely thanks to the simple structures of the matrices and only the absolutely necessary variables are stored and computed.*

## 5.5 MPC with terminal quadratic constraint

This section presents the sparse ADMM-based solver presented in [12] for the MPC formulation

$$\min_{\mathbf{x},\mathbf{u}} \sum_{j=0}^{N-1} \left( \|x_j - x_r\|_Q^2 + \|u_j - u_r\|_R^2 \right) + \|x_N - x_r\|_T^2 \tag{5.20a}$$

$$s.t.\ x_0 = x(t) \tag{5.20b}$$

$$x_{j+1} = Ax_j + Bu_j,\ j \in \mathbb{Z}_0^{N-1} \tag{5.20c}$$

$$\underline{x}_j \leq x_j \leq \overline{x}_j,\ j \in \mathbb{Z}_1^{N-1} \tag{5.20d}$$

$$\underline{u}_j \leq u_j \leq \overline{u}_j,\ j \in \mathbb{Z}_0^{N-1} \tag{5.20e}$$

$$x_N \in \mathcal{E}(P, c, r), \tag{5.20f}$$

where $\mathcal{E}(P, c, r)$ is the ellipsoid defined by the given $P \in \mathbb{S}_{++}^n$, $c \in \mathbb{R}^n$ and $r \in \mathbb{R}_{>0}$ (see Definition N.7); under the following assumption.

**Assumption 5.25.** Let Assumtion 4.7 hold and assume that $Q \in \mathbb{S}_+^n$, $R \in \mathbb{S}_+^m$, $T \in \mathbb{S}_+^n$, and that the bounds $(\underline{x}_j, \overline{x}_j) \in \mathbb{R}^n \times \mathbb{R}^n$, $(\underline{u}_j, \overline{u}_j) \in \mathbb{R}^m \times \mathbb{R}^m$ given in (5.20d) and (5.20e) satisfy $\underline{x} \leq \underline{x}_j < \overline{x}_j \leq \overline{x}$, $j \in \mathbb{Z}_1^{N-1}$ and $\underline{u} \leq \underline{u}_j < \overline{u}_j \leq \overline{u}$, $j \in \mathbb{Z}_0^{N-1}$, where $\underline{x}$, $\overline{x}$, $\underline{u}$ and $\overline{u}$ are the system bounds (4.3).

As discussed in Section 4.2, MPC formulations typically rely on a terminal constraint $x_N \in \mathcal{X}_t$, see (4.7e), to guarantee stability of the closed-loop system [128], where the terminal set $\mathcal{X}_t$ is taken as an admissible invariant set of the system for its state and input constraints (see Definition N.8). This is true for both nominal MPC formulations [64, §2], as well as for many robust MPC formulations, in which case $\mathcal{X}_t$ takes the form of a robust admissible invariant set of the system [129]. Additionally, in some robust MPC approaches the constraints are tightened throughout the prediction horizon [129, 130], thus our consideration of step-dependent constraints in (5.20d) and (5.20e).

MPC formulation (5.20) captures these paradigms for the particular case in which the (robust) admissible invariant set is taken as an ellipsoid $\mathcal{E}(P, c, r)$. Typically, an admissible invariant set in the form of a polyhedral set is used, i.e., a set of the form $\{ x \in \mathbb{R}^n : A_t x \leq b_t \}$ where $A_t \in \mathbb{R}^{n_t \times n}$ and $b_t \in \mathbb{R}^{n_t}$, either so that the resulting optimization problem is a QP, or because the maximal (robust) admissible invariant set is used, which for controllable linear systems is a polyhedral set. The computation of the polyhedral admissible invariant set for MPC is a well researched field [131]. However, it typically results in the addition of a large amount of inequality constraints, i.e., $n_t$ is very large [126, §5], even for moderately large systems. As such, even though the resulting optimization problem is a QP, for which many efficient solvers have been developed (see the QP solvers listed in Section 4.3), it may be computationally demanding to solve.

In many cases, the polyhedral (robust) admissible invariant set can be substituted by one in the form of an ellipsoid $\mathcal{E}(P, c, r)$ [16], [126, §4.1] [132], thus leading to the MPC formulation (5.20). This comes at the expense of the resulting optimization problem no longer being a QP problem, but with the advantage of it (typically) having significantly fewer constraints. In fact, in many occasions, the computation of a polyhedral (robust) admissible invariant set is prohibitively expensive, whereas the computation of an ellipsoidal one is attainable, since it can be posed, for instance, as a linear matrix inequality (LMI) problem[2] [133]. This further motivates the usefulness of a solver tailored to this MPC formulation.

Finally, we note that optimization problem (5.20) is a quadratically constrained quadratic programming (QCQP) problem, which is generally challenging to solve but can have closed-form solutions in some cases [134]. In particular, in [135] the authors present the FalcOpt solver [80], which is suitable for embedded systems and considers nonlinear MPC subject to a terminal quadratic constraint. However, it does not consider state constraints and the use of (a variation of) the sequential quadratic programming approach may make it less efficient when dealing with the linear MPC case than a tailored linear MPC solver.

### 5.5.1 ADMM solver for MPC with terminal quadratic constraint

Optimization problem (5.20) cannot be cast as a QP (5.13) due to the terminal quadratic constraint (5.20f), so the procedure used in Section 5.3 cannot be employed here. However, a very similar procedure can be used.

The MPC formulation (5.20) can be cast as an optimization problem (2.12) as follows. Let us define the auxiliary variables $(\tilde{x}_1, \dots, \tilde{x}_N)$ and $(\tilde{u}_0, \tilde{u}_1, \dots, \tilde{u}_{N-1})$, and take

$$
\begin{aligned}
z &= (u_0, x_1, u_1, x_2, u_2, \dots, x_{N-1}, u_{N-1}, x_N), \\
v &= (\tilde{u}_0, \tilde{x}_1, \tilde{u}_1, \tilde{x}_2, \tilde{u}_2, \dots, \tilde{x}_{N-1}, \tilde{u}_{N-1}, \tilde{x}_N).
\end{aligned}
$$

To facilitate readability, we divide $z$ and $v$ into two parts, given by $z = (z_\circ, z_f) \in \mathbb{R}^{n_\circ} \times \mathbb{R}^n$ and $v = (v_\circ, v_f) \in \mathbb{R}^{n_\circ} \times \mathbb{R}^n$, where

$$
\begin{aligned}
z_\circ &\doteq (u_0, x_1, u_1, x_2, u_2, \dots, x_{N-1}, u_{N-1}), \\
v_\circ &\doteq (\tilde{u}_0, \tilde{x}_1, \tilde{u}_1, \tilde{x}_2, \tilde{u}_2, \dots, \tilde{x}_{N-1}, \tilde{u}_{N-1}),
\end{aligned}
$$

$z_f \doteq x_N$ and $v_f \doteq \tilde{x}_N$. Therefore, $n_\circ = (N-1)(n+m)+m$. Then, problem (5.20) can be recast as (2.12) by taking

$$
\begin{aligned}
f(z) &= \frac{1}{2} z^\top H z + q^\top z + \mathcal{I}_{(Gz-b=0)}(z), \\
g(v) &= \mathcal{I}_{(\underline{v} \leq v_\circ \leq \overline{v})}(v_\circ) + \mathcal{I}_{\mathcal{E}(P,c,r)}(v_f),
\end{aligned}
$$

[2]See Section 5.5.2, where, for completeness, we provide a design procedure for an ellipsoidal admissible positive invariant set $\mathcal{E}(P, c, r)$.

where

$$H = \operatorname{diag}(R,\, Q,\, R,\, Q,\, R,\, \dots, Q,\, R,\, T), \tag{5.21a}$$
$$q = -(Ru_r,\, Qx_r,\, Ru_r,\, \dots, Qx_r,\, Ru_r,\, Tx_r), \tag{5.21b}$$
$$\underline{v}_\circ = (\underline{u}_0,\, \underline{x}_1,\, \underline{u}_1,\, \dots, \underline{x}_{N-1},\, \underline{u}_{N-1}), \tag{5.21c}$$
$$\overline{v}_\circ = (\overline{u}_0,\, \overline{x}_1,\, \overline{u}_1,\, \dots, \overline{x}_{N-1},\, \overline{u}_{N-1}), \tag{5.21d}$$
$$G = \begin{bmatrix} B & -\mathbf{I}_n & \mathbf{0} & \cdots & \cdots & \mathbf{0} \\ \mathbf{0} & A & B & -\mathbf{I}_n & \cdots & \mathbf{0} \\ \mathbf{0} & \cdots & \ddots & \ddots & \ddots & \mathbf{0} \\ \mathbf{0} & \cdots & \mathbf{0} & A & B & -\mathbf{I}_n \end{bmatrix}, \tag{5.21e}$$
$$b = (-Ax(t),\, \mathbf{0}_n,\, \mathbf{0}_n, \dots,\, \mathbf{0}_n), \tag{5.21f}$$

and by imposing

$$z_\circ - v_\circ = \mathbf{0}_{n_\circ}, \tag{5.22a}$$
$$P^{1/2}(z_f - v_f) = \mathbf{0}_n, \tag{5.22b}$$

where $P^{1/2} \in \mathbb{S}^n_{++}$ is the matrix that satisfies $P = P^{1/2}P^{1/2}$. Thus, the matrices $C$, $D$ and vector $d$ of (2.12) are given by

$$\begin{aligned} C &= \operatorname{diag}(\mathbf{I}_m,\, \mathbf{I}_n,\, \mathbf{I}_m,\, \dots, \mathbf{I}_n,\, \mathbf{I}_m,\, P^{1/2}), \\ D &= -\operatorname{diag}(\mathbf{I}_m,\, \mathbf{I}_n,\, \mathbf{I}_m,\, \dots, \mathbf{I}_n,\, \mathbf{I}_m,\, P^{1/2}), \\ d &= \mathbf{0}_{n_\circ + n}. \end{aligned}$$

The reason for imposing $P^{1/2}(z_f - v_f) = \mathbf{0}_n$, instead of $z_f - v_f = \mathbf{0}_n$, will be clear further ahead (see Remark 5.27).

Let us now also divide $\lambda$ into two parts: $\lambda = (\lambda_\circ, \lambda_f) \in \mathbb{R}^{n_\circ} \times \mathbb{R}^n$, where $\lambda_\circ$ are the dual variables associated to the constraints (5.22a) and $\lambda_f$ are the dual variables associated to the constraints (5.22b).

Then, step 3 of Algorithm 2 requires solving the optimization problem:

$$\min_z \frac{1}{2} z^\top \hat{H} z + \hat{q}_k^\top z \tag{5.23a}$$
$$s.t.\ Gz = b, \tag{5.23b}$$

where

$$\begin{aligned} \hat{H} &= H + \rho C^\top C = H + \rho \operatorname{diag}(\mathbf{I}_m, \mathbf{I}_n, \mathbf{I}_m, \dots, \mathbf{I}_n, \mathbf{I}_m, P), \\ \hat{q}_k &= q + \rho C^\top D v^k + C^\top \lambda^k = q + (\lambda^k_\circ - \rho v^k_\circ, P^{1/2}\lambda^k_f - \rho P v^k_f). \end{aligned}$$

Due to the block diagonal structure of $\hat{H}$ and the banded structure of $G$ shown in (5.21e), problem (5.23) satisfies Assumtion 5.3 and can therefore be solved using Algorithm 12.

Step 4 of Algorithm 2 has a separable structure that allows the problem to be divided in two parts, one for $v_\circ$ and one for $v_f$. The update of $v_\circ^{k+1}$ is the solution of the following optimization problem:

$$\begin{aligned} v_\circ^{k+1} = \arg\min_{v_\circ} \ & \frac{\rho}{2} v_\circ^\top v_\circ - (\rho z_\circ^{k+1} + \lambda_0^k)^\top v_\circ \\ s.t. \ & \underline{v}_\circ \leq v_\circ \leq \overline{v}_\circ, \end{aligned}$$

which satisfies Assumtion 5.8 and can therefore be solved using Algorithm 13. The update of $v_f^{k+1}$ is the solution of the optimization problem

$$\begin{aligned} v_f^{k+1} = \arg\min_{v_f} \ & \frac{\rho}{2} v_f^\top P v_f - (\rho P z_f^{k+1} + P^{1/2}\lambda_f^k)^\top v_f \\ s.t. \ & v_f \in \mathcal{E}(P, c, r). \end{aligned}$$

Dividing the objective function by $\rho$ and defining $P^{-1/2} \doteq P^{-1}P^{1/2}$, this problem can be recast as

$$v_f^{k+1} = \arg\min_{v_f} \frac{1}{2} \|v_f - (z_f^{k+1} + \rho^{-1}P^{-1/2}\lambda_f^k)\|_P^2 \tag{5.24a}$$

$$s.t. \ v_f \in \mathcal{E}(P, c, r), \tag{5.24b}$$

which has an explicit solution given by the following theorem.

**Theorem 5.26.** Let $a \in \mathbb{R}^n$, $P \in \mathbb{S}_{++}^n$, $c \in \mathbb{R}^n$, $r \in \mathbb{R}_{>0}$. Then, the solution $v^*$ of the convex optimization problem

$$\min_v \frac{1}{2} \|v - a\|_P^2 \tag{5.25a}$$

$$s.t. \ v \in \mathcal{E}(P, c, r), \tag{5.25b}$$

is given by

$$v^* = \begin{cases} a & \text{if } a \in \mathcal{E}(P, c, r) \\ \dfrac{r(a-c)}{\|a-c\|_P} + c & \text{otherwise.} \end{cases}$$

***Proof:*** Let $\tilde{a} \doteq a - c$, and let $\Pi_{(P,c,r)}(a)$ denote the argument that minimizes (5.25) for the given $P$, $c$ and $r$. Then, it is clear that

$$\Pi_{(P,c,r)}(a) = \Pi_{(P,\mathbf{0}_n,r)}(\tilde{a}) + c, \tag{5.26}$$

since this simply corresponds to shifting the origin to the center of the ellipsoid in (5.25b), and then undoing it. Therefore, it suffices to find a closed expression for the solution of $\Pi_{(P,\mathbf{0}_n,r)}(\tilde{a})$, i.e., to find an explicit solution to

$$\min_v \frac{1}{2} \|v - \tilde{a}\|_P^2 \tag{5.27a}$$

$$s.t. \ v^\top P v \leq r^2. \tag{5.27b}$$

The Lagrangian $\mathcal{L} : \mathbb{R}^n \times \mathbb{R} \to \mathbb{R}$ of (5.27) is given by

$$\mathcal{L}(v, y) = \frac{1}{2}\|v - \tilde{a}\|_P^2 + y(v^\top P v - r^2), \tag{5.28}$$

where $y \in \mathbb{R}$ is the dual variable for the constraint (5.27b). The dual function associated with (5.27) is

$$\Psi(y) = \inf_v \mathcal{L}(v, y). \tag{5.29}$$

Let us define $v(y) \in \mathbb{R}^n$ as

$$v(y) = \arg\min_v \mathcal{L}(v, y). \tag{5.30}$$

It is clear from the definition of $v(y)$ that the dual function (5.29) can be expressed as

$$\Psi(y) = \mathcal{L}(v(y), y).$$

The dual problem of (5.27) is to maximize the dual function $\Psi(y)$ subject to $y \geq 0$ [23, §5], i.e., to find the optimal solution of

$$\max_{y \geq 0} \mathcal{L}(v(y), y). \tag{5.31}$$

We start by rewriting optimization problem (5.30) as

$$v(y) = \arg\min_v \frac{1}{2} v^\top (1{+}2y) P v - \tilde{a}^\top P v.$$

whose optimal solution is the vector $v(y)$ for which the gradient of the objective function is equal to zero, since it is an unconstrained convex problem with a differentiable real-valued objective function [23, §4.2.3]. That is,

$$(1{+}2y)Pv(y) - P\tilde{a} = \mathbf{0}_n,$$

which leads to

$$v(y) = \frac{1}{1{+}2y}\tilde{a}. \tag{5.32}$$

Substituting expression (5.32) into (5.28) leads to

$$\begin{aligned} \mathcal{L}(v(y), y) &= \left[ \frac{1}{2} \left( \frac{2y}{1 + 2y} \right)^2 + \frac{y}{(1 + 2y)^2} \right] \tilde{a}^\top P \tilde{a} - r^2 y \\ &= \frac{y}{1{+}2y} \tilde{a}^\top P \tilde{a} - r^2 y, \end{aligned} \tag{5.33}$$

which, for $y > -1/2$, is a differentiable real-valued concave function[3]. Therefore, given that $y$ is a scalar, the optimal solution $y^*$ of (5.31) is given by

$$y^* = \max\{\hat{y}, 0\}, \tag{5.34}$$

[3] An easy way to see this is to note that

$$\frac{y}{1 + 2y} = \frac{1}{2}\left(1 - \frac{1}{1 + 2y}\right),$$

where $-1/(1{+}2y)$ is differentiable, real-valued and concave for $y > -1/2$.

where $\hat{y}$ is the scalar such that

$$\left.\frac{d\mathcal{L}(v(y), y)}{dy}\right|_{\hat{y}} = 0,$$

which, differentiating (5.33), leads to

$$\frac{\tilde{a}^\top P \tilde{a} - r^2(1+2\hat{y})^2}{(1+2\hat{y})^2} = 0,$$

$$\hat{y} = \frac{1}{2}\left(\frac{\sqrt{\tilde{a}^\top P \tilde{a}}}{r} - 1\right). \tag{5.35}$$

Given that strong duality holds, we have that the optimal solution of (5.27), is given by $v^* = v(y^*)$. Therefore, from (5.32), (5.34) and (5.35), we have that

$$v(y^*) = \begin{cases} \dfrac{r\tilde{a}}{\sqrt{\tilde{a}^\top P \tilde{a}}}, & \text{if } \tilde{a}^\top P \tilde{a} > r^2, \\ \tilde{a}, & \text{if } \tilde{a}^\top P \tilde{a} \leq r^2, \end{cases}$$

which, noting that $v(y^*) \equiv \Pi_{(P, \mathbf{0}_n, r)}(\tilde{a})$ and taking into account (5.26), proves the claim. ■

**Remark 5.27.** *We note that the reason for imposing* (5.22b) *is so that* (5.24) *can be solved using the simple explicit solution provided in Theorem 5.26. If, instead, we had taken the more simple constraint $z_f - v_f = \mathbf{0}_n$, problem* (5.24) *would have been a standard Euclidean projection onto the ellipsoid $\mathcal{E}(P, c, r)$, which does not have an explicit solution and would have thus required an iterative algorithm, such as [136, §2], to be solved.*

Algorithm 18 shows the particularization of Algorithm 2 to problem (5.20) under Assumtion 5.25 obtained from the above discussion. The matrices used in the algorithm are computed offline and stored in the embedded system. Therefore, the value of $\rho$ cannot be updated online. However, the values of $c$ and $r$ can change between sample times, as well as the values of $x(t)$, $x_r$, $u_r$ and the bounds in (5.20d) and (5.20e). The control action $u(t)$ to be applied to the system is taken as the first $m$ elements of $\tilde{v}^*$.

**Remark 5.28.** *Remark 5.19 also applies to Algorithm 18. That is, the matrix-matrix and matrix-vector operations are all performed sparsely thanks to the simple structures of the matrices and only the absolutely necessary variables are stored and computed.*

**Algorithm 18:** Sparse ADMM-based solver for (5.20)

**Require:** $x(t)$, $(x_r, u_r)$, $v^0$, $\lambda^0$, $\epsilon_p > 0$, $\epsilon_d > 0$

Compute $q$ with $x_r$ and $u_r$

Compute $b$ with $x(t)$

$k \leftarrow 0$

**repeat**

  $\hat{q}_k \leftarrow q + (\lambda_\circ^k - \rho v_\circ^k, P^{1/2}\lambda_f^k - \rho P v_f^k)$

  $z^{k+1} \leftarrow \text{solve_eqQP}(\hat{q}_k, b; \hat{H}, G)$

  $v_\circ^{k+1} \leftarrow \text{solve_boxQP}(-\rho z_\circ^{k+1} - \lambda_\circ^k; \rho, \underline{v}_\circ, \overline{v}_\circ)$

  $v_f^{k+1} \leftarrow z_f^{k+1} + \rho^{-1} P^{-1/2} \lambda_f^k$

  **if** $(v_f^{k+1} - c)^\top P (v_f^{k+1} - c) > r^2$ **then**

    $$v_f^{k+1} \leftarrow \frac{r(v_f^{k+1} - c)}{\sqrt{(v_f^{k+1} - c)^\top P (v_f^{k+1} - c)}} + c$$

  **end if**

  $\lambda_\circ^{k+1} \leftarrow \lambda_\circ^k + \rho(z_\circ^{k+1} - v_\circ^{k+1})$

  $\lambda_f^{k+1} \leftarrow \lambda_f^k + \rho P^{1/2}(z_f^{k+1} - v_f^{k+1})$

  $k \leftarrow k + 1$

**until** $r_p \leq \epsilon_p$ **and** $r_d \leq \epsilon_d$

**Output:** $\tilde{z}^* \leftarrow z^k$, $\tilde{v}^* \leftarrow v^k$, $\tilde{\lambda}^* \leftarrow \lambda^k$

### 5.5.2 Computation of admissible ellipsoidal invariant sets

Consider a system described by (4.1) subject to

$$Cx \leq c, \quad Du \leq d, \tag{5.36}$$

where $C \in \mathbb{R}^{p_x \times n}$, $c \in \mathbb{R}^{p_x}$, $D \in \mathbb{R}^{p_u \times m}$, $d \in \mathbb{R}^{p_u}$; and a steady state reference $(x_r, u_r)$ satisfying the above constraints that we wish the system to converge to.

We describe a procedure taken from various articles in the literature (see, for instance, [137] [138, §C.8.1], [133]) for computing a control gain $K$ and an admissible positive invariant ellipsoidal set $\mathcal{E}(P, c, r)$ (see Definition N.7 and Definition N.8) of system (4.1) subject to (5.36) for the control law $u(t){=}K(x(t){-}x_r){+}u_r$ such that the closed-loop system admissibly converges to $x_r$ for any initial state $x(t) \in \mathcal{E}(P, c, r)$. The procedure is based on solving an optimization problem subject to LMIs.

To simplify the procedure, we start by shifting the origin to the reference $(x_r, u_r)$. That is, let $\hat{x} \doteq x - x_r$ and $\hat{u} \doteq u - u_r$. Then, we compute a gain $K$ and an ellipsoid $\mathcal{E}(P, \mathbf{0}_n, r)$ (i.e., taking $c = x_r$), such that the closed loop system $\hat{x}(t+1) = A_K \hat{x}(t)$, where $A_K \doteq A + BK$, admissibly converges to the origin for the constraints

$$C\hat{x} \leq \hat{c}, \quad DK\hat{x} \leq \hat{d}, \tag{5.37}$$

where $\hat{c} \in \mathbb{R}^{p_x}$ and $\hat{d} \in \mathbb{R}^{p_u}$ are given by

$$\hat{c} = c - Cx_r, \quad \hat{d} = d - Du_r.$$

The implication for invariance given in Definition N.8 reads as

$$\hat{x}^\top P \hat{x} \le r^2 \implies \hat{x}^\top A_K^\top P A_K \hat{x} \le r^2,$$

which, applying the S-procedure [23, §B.2], is satisfied if there exists a scalar $\lambda \ge 0$ such that

$$\begin{bmatrix} \lambda P - A_K^\top P A_K & \mathbf{0}_n \\ \mathbf{0}_n^\top & r^2(1-\lambda) \end{bmatrix} \in \mathbb{S}_{++}^{n+1},$$

where we note that $\lambda$ must therefore satisfy $\lambda \in [0, 1)$. This expression can be rewritten as

$$\begin{bmatrix} \lambda P & \mathbf{0}_n \\ \mathbf{0}_n^\top & r^2(1-\lambda) \end{bmatrix} - \begin{bmatrix} A_K^\top \\ 0 \end{bmatrix} P \begin{bmatrix} A_K & 0 \end{bmatrix} \in \mathbb{S}_{++}^{n+1}.$$

Applying the Schur complement [23, §A.5.5] leads to:

$$\begin{bmatrix} \lambda P & \mathbf{0}_n & A_K^\top \\ \mathbf{0}_n^\top & r^2(1-\lambda) & \mathbf{0}_n^\top \\ A_K & \mathbf{0}_n & P^{-1} \end{bmatrix} \in \mathbb{S}_{++}^{2n+1},$$

Finally, pre- and post-multiplying by $\operatorname{diag}(P^{-1}, 1, I_n)$ and taking the transformation $W \doteq P^{-1}$, $Y \doteq KP^{-1}$, leads to the LMI:

$$\begin{bmatrix} \lambda W & \mathbf{0}_n & WA^\top + Y^\top B^\top \\ \mathbf{0}_n^\top & r^2(1-\lambda) & \mathbf{0}_n^\top \\ AW + BY & \mathbf{0}_n & W \end{bmatrix} \in \mathbb{S}_{++}^{2n+1}. \tag{5.38}$$

Next, (5.37) must be satisfied for all $\hat{x} \in \mathcal{E}(P, \mathbf{0}_n, r)$. It is well known that

$$\max_{x \in \mathcal{E}(P, \mathbf{0}_n, r)} v^\top x = r\sqrt{v^\top P^{-1} v},$$

where $v \in \mathbb{R}^n$ [133, §5.2.2]. Therefore this condition can be imposed by finding $P$ and $K$ such that

$$r^2 C_{(j)}^\top P^{-1} C_{(j)} \le \hat{c}_{(j)}^2,\ j \in \mathbb{Z}_1^{p_x},$$
$$r^2 (D_{(j)}K)^\top P^{-1} (D_{(j)}K) \le \hat{d}_{(j)}^2,\ j \in \mathbb{Z}_1^{p_u},$$

where the subindex $(j)$ indicates the $j$-th row of the matrix or component of the vector. Taking the above definitions of $W$ and $Y$, this can be posed as the LMIs:

$$r^2 C_{(j)}^\top W C_{(j)} \le \hat{c}_{(j)}^2,\ j \in \mathbb{Z}_1^{p_x}, \tag{5.39}$$

$$\begin{bmatrix} r^{-2}\hat{d}_{(j)}^2 & D_{(j)}Y \\ Y^\top D_{(j)}^\top & W \end{bmatrix} \in \mathbb{S}_{++}^{n+1},\ j \in \mathbb{Z}_1^{p_u}. \tag{5.40}$$

The procedure is to select a value of $r$ (typically, for convenience, we pick $r = 1$) and then to solve the following convex optimization problem subject to LMI constraints

$$\min_{W,Y} \; -\operatorname{tr}(W) \tag{5.41a}$$

$$s.t.\ (5.38), (5.39) \text{ and } (5.40), \tag{5.41b}$$

for increasing values of $\lambda$ in the range $0 \leq \lambda < 1$ until a feasible solution is found. Finally, $P$ and $K$ are recovered from the solution of (5.41). The minimization of $-\operatorname{tr}(W)$ is done to maximize the size of the resulting ellipsoid.

**Remark 5.29.** *Similar procedures can be used to compute the terminal ingredients of a robust MPC controller that follows formulation* (5.20) *[16, §V.A].*

## 5.6 MPC for tracking

The solvers presented in the previous sections all consider variations of the standard MPC formulations (4.7). This section, on the other hand, presents a sparse solver, based on the extended ADMM algorithm (Algorithm 3), for a non-standard MPC formulation known as *MPC for tracking* (MPCT) [139, 140]. This formulation differs from standard MPC formulations (4.7) in the inclusion of a pair of decision variables $(x_s, u_s) \in \mathbb{R}^n \times \mathbb{R}^m$ known as the *artificial reference*. Out of the different variants of the MPCT formulation, this dissertation is concerned with the following one, which uses a terminal equality constraint:

$$\min_{\substack{\mathbf{x},\mathbf{u},\\ x_s,u_s}} \sum_{j=0}^{N-1} \left(\|x_j - x_s\|_Q^2 + \|u_j - u_s\|_R^2\right) + \|x_s - x_r\|_T^2 + \|u_s - u_r\|_S^2 \tag{5.42a}$$

$$s.t.\ x_0 = x(t) \tag{5.42b}$$

$$x_{j+1} = Ax_j + Bu_j,\ j \in \mathbb{Z}_0^{N-1} \tag{5.42c}$$

$$\underline{x} \leq x_j \leq \overline{x},\ j \in \mathbb{Z}_1^{N-1} \tag{5.42d}$$

$$\underline{u} \leq u_j \leq \overline{u},\ j \in \mathbb{Z}_0^{N-1} \tag{5.42e}$$

$$x_s = Ax_s + Bu_s \tag{5.42f}$$

$$\underline{x} + \varepsilon_x \leq x_s \leq \overline{x} - \varepsilon_x \tag{5.42g}$$

$$\underline{u} + \varepsilon_u \leq u_s \leq \overline{u} - \varepsilon_u \tag{5.42h}$$

$$x_N = x_s, \tag{5.42i}$$

where $\varepsilon_x \in \mathbb{R}^n$ and $\varepsilon_u \in \mathbb{R}^m$ are vectors with arbitrarily small positive components which are added to avoid a possible loss of controllability when the constraints are active at the equilibrium point [140]. We consider the following assumption, which we note does not assume that $(x_r, u_r)$ is an admissible steady state of system (4.1) subject to (4.3).

**Assumption 5.30.** Let Assumtion 4.2 hold and assume that:

*(i)* $Q \in \mathbb{S}^n_+$, $R \in \mathbb{S}^m_+$, $T \in \mathbb{S}^n_+$ and $S \in \mathbb{S}^m_+$.

*(ii)* $x(t)$ strictly belongs to the feasibility region of (5.42). That is, there exist $\mathbf{x}$, $\mathbf{u}$, $x_s$ and $u_s$ such that (5.42b), (5.42c), (5.42f) and (5.42i) are satisfied and (5.42d), (5.42e), (5.42g) and (5.42h) are strictly satisfied.

**Remark 5.31.** *We note that Assumtion 5.30.(ii) implies that* $\underline{x} + \varepsilon_x < \overline{x} - \varepsilon_x$ *and* $\underline{u} + \varepsilon_u < \overline{u} - \varepsilon_u$*.*

The cost function (5.42a) penalizes, on one hand, the difference between the predicted states $x_j$ and control actions $u_j$ with the artificial reference $x_s$ and $u_s$, respectively, and on the other, the discrepancy between the artificial reference and the reference $(x_r, u_r)$ given by the user.

The inclusion of the artificial reference provides the MPCT formulation with a series of advantages with respect to other (standard) MPC formulations.

First, a common issue of standard MPC formulations with stability guarantees is that the domain of attraction of the controller can become insufficient if the prediction horizon is chosen too small. However, the use of small prediction horizons is desirable in order to help overcome the computational and memory limitations typically imposed by embedded systems. The MPCT formulation provides significantly larger domains of attraction than standard MPC formulations [140], especially for small prediction horizons.

Second, it intrinsically deals with references that are not attainable, i.e., that are not a steady state of the system and/or that violate the system constraints [140]. In this case, it will steer the closed-loop system to the admissible steady state $(x_a, u_a) \in \mathbb{R}^n \times \mathbb{R}^m$ that minimizes the cost $\|x_a - x_r\|^2_T + \|u_a - u_r\|^2_S$. Finally, it also guarantees recursive feasibility of the closed-loop system even in the event of a sudden reference change [139].

We now present a sparse solver, originally presented in [10], for the MPCT formulation (5.42) that is based on the EADMM algorithm (Algorithm 3). The MPCT formulation (5.42) can be expressed as a QP problem (5.13). However, the inclusion of the artificial reference $(x_s, u_s)$ leads to a Hessian matrix that is not block diagonal. Therefore, if we were to follow the same procedure presented in Section 5.3, we would not be able to solve step 4 of Algorithm 15 using Algorithm 12 because the resulting $W$ matrix would not satisfy Assumtion 5.1. Therefore, we propose to use instead use EADMM (Algorithm 3), which will allow us to recover the structures that emerged in the previous sections, thus resulting in a sparse solver with a very similar iteration complexity to the previous ones.

### 5.6.1 Recasting the MPCT formulation for EADMM

Let us recast problem (5.42) by defining the auxiliary variables $\tilde{x}_i \doteq x_i - x_s$ and $\tilde{u}_i \doteq u_i - u_s$ as follows:

$$\min_{\substack{\tilde{\mathbf{x}},\tilde{\mathbf{u}},\mathbf{x},\\ \mathbf{u},x_s,u_s}} \sum_{i=0}^{N} \left( \|\tilde{x}_i\|_Q^2 + \|\tilde{u}_i\|_R^2 \right) + \|x_s - x_r\|_T^2 + \|u_s - u_r\|_S^2 \tag{5.43a}$$

$$s.t.\ x_0 = x(t) \tag{5.43b}$$

$$\tilde{x}_{i+1} = A\tilde{x}_i + B\tilde{u}_i,\ i \in \mathbb{Z}_0^{N-1} \tag{5.43c}$$

$$\underline{x} \leq x_i \leq \overline{x},\ i \in \mathbb{Z}_1^{N-1} \tag{5.43d}$$

$$\underline{u} \leq u_i \leq \overline{u},\ i \in \mathbb{Z}_0^{N-1} \tag{5.43e}$$

$$\underline{x} + \varepsilon_x \leq x_N \leq \overline{x} - \varepsilon_x, \tag{5.43f}$$

$$\underline{u} + \varepsilon_u \leq u_N \leq \overline{u} - \varepsilon_u, \tag{5.43g}$$

$$x_s = Ax_s + Bu_s \tag{5.43h}$$

$$\tilde{x}_i + x_s - x_i = \mathbf{0}_n,\ i \in \mathbb{Z}_0^N \tag{5.43i}$$

$$\tilde{u}_i + u_s - u_i = \mathbf{0}_m,\ i \in \mathbb{Z}_0^N \tag{5.43j}$$

$$x_N = x_s, \tag{5.43k}$$

$$u_N = u_s, \tag{5.43l}$$

where the decision variables are $x_s$, $u_s$, $\tilde{\mathbf{x}} = (\tilde{x}_0, \dots, \tilde{x}_N)$, $\tilde{\mathbf{u}} = (\tilde{u}_0, \dots, \tilde{u}_N)$, $\mathbf{x} = (x_0, \dots, x_N)$, and $\mathbf{u} = (u_0, \dots, u_N)$. Equality constraints (5.43i) and (5.43j) impose the congruence of the decision variables with the original problem (5.42). Inequalities (5.42g) and (5.42h) are omitted because they are already imposed by (5.43f) and (5.43g) alongside the inclusion of (5.43k) and (5.43l).

**Remark 5.32.** *Note that the summation in the cost function* (5.43a) *now includes $i = N$. However, this does not change the solution of the optimization problem due to the inclusion of* (5.43k) *and* (5.43l)*.*

We can now obtain a problem of form (2.13) by taking

$$z_1 = (x_0, u_0, x_1, u_1, \dots, x_{N-1}, u_{N-1}, x_N, u_N), \tag{5.44a}$$

$$z_2 = (x_s, u_s), \tag{5.44b}$$

$$z_3 = (\tilde{x}_0, \tilde{u}_0, \tilde{x}_1, \tilde{u}_1, \dots, \tilde{x}_{N-1}, \tilde{u}_{N-1}, \tilde{x}_N, \tilde{u}_N), \tag{5.44c}$$

which leads to

$$\theta_1(z_1) = 0,$$
$$\theta_2(z_2) = \frac{1}{2} z_2^\top \operatorname{diag}(T, S) z_2 - (Tx_r, Su_r)^\top z_2,$$
$$\theta_3(z_3) = \frac{1}{2} z_3^\top \operatorname{diag}(Q, R, Q, R, \dots, Q, R) z_3,$$

$$C_1 = \begin{bmatrix} [\mathbf{I}_n \, \mathbf{0}_{n\times m}] & \mathbf{0} & \mathbf{0} \\ \hline -\mathbf{I}_{n+m} & \mathbf{0} & \mathbf{0} \\ \mathbf{0} & \ddots & \mathbf{0} \\ \mathbf{0} & \mathbf{0} & -\mathbf{I}_{n+m} \\ \hline \mathbf{0} & \mathbf{0} & -\mathbf{I}_{n+m} \end{bmatrix}, \qquad C_2 = \begin{bmatrix} \mathbf{0} \\ \hline \mathbf{I}_{n+m} \\ \vdots \\ \mathbf{I}_{n+m} \\ \hline \mathbf{I}_{n+m} \end{bmatrix},$$

$$C_3 = \begin{bmatrix} \mathbf{0} & \dots & \mathbf{0} \\ \hline \mathbf{I}_{n+m} & \mathbf{0} & \mathbf{0} \\ \mathbf{0} & \ddots & \mathbf{0} \\ \mathbf{0} & \mathbf{0} & \mathbf{I}_{n+m} \\ \hline \mathbf{0} & \dots & \mathbf{0} \end{bmatrix}, \qquad b = \begin{bmatrix} x(t) \\ \hline \mathbf{0} \\ \vdots \\ \mathbf{0} \\ \hline \mathbf{0} \end{bmatrix}.$$

Matrices $C_1$, $C_2$ and $C_3$ contain the equality constraints (5.43b), (5.43i), (5.43j), (5.43k) and (5.43l). Specifically, the first $n$ rows impose constraint (5.43b), the last $n+m$ rows impose the constraints (5.43k) and (5.43l), and the rest of the rows impose the constraints (5.43i) and (5.43j). Set $\mathcal{Z}_1$ is the set of vectors $z_1$ (5.44a) for which the box constraints (5.43d)-(5.43g) are satisfied; set $\mathcal{Z}_2$ is the set of vectors $z_2$ (5.44b) that satisfy the equality constraint (5.43h); and set $\mathcal{Z}_3$ is the set of vectors $z_3$ (5.44c) that satisfy the equality constraints (5.43c).

**Remark 5.33.** *Our selection of $z_i$ and $C_i$ for $i \in \mathbb{Z}_1^3$ results in an optimization problem that satisfies Assumption 2.8. Therefore, under a proper selection of $\rho$, the iterates of the EADMM algorithm will converge to the optimal solution of the MPCT controller. In practice, the parameter $\rho$ may be selected outside the range shown in Theorem 2.9 in order to improve the convergence rate of the algorithm [36]. In this case, the convergence will not be guaranteed and will have to be extensively checked with simulations.*

### 5.6.2 EADMM-based solver for the MPCT formulation

We now particularize Algorithm 3 to the optimization problem (5.43). Let us take a closer look at steps 3, 4 and 5 of Algorithm 3.

Step 3 of Algorithm 3 minimizes the Lagrangian (2.14) over $z_1$, resulting in the following box-constrained optimization problem:

$$\min_{z_1} \frac{1}{2} z_1^\top H_1 z_1 + q_1^\top z_1 \tag{5.45a}$$

$$s.t. \; \underline{z}_1 \leq z_1 \leq \overline{z}_1, \tag{5.45b}$$

where

$$\begin{aligned} H_1 &= \rho C_1^\top C_1, \\ q_1 &= \rho C_1^\top C_2 z_2^k + \rho C_1^\top C_3 z_3^k + C_1^\top \lambda^k - \rho C_1^\top b, \\ \underline{z}_1 &= (-M_n, \underline{u}, \underline{x}, \dots, \underline{u}, \underline{x} + \varepsilon_x, \underline{u} + \varepsilon_u), \\ \overline{z}_1 &= (M_n, \overline{u}, \overline{x}, \dots, \overline{u}, \overline{x} - \varepsilon_x, \overline{u} - \varepsilon_u), \end{aligned}$$

and $M_n \in \mathbb{R}^n_{>0}$ has arbitrarily large components. Due to the structure of $C_1$, we have that $H_1 \in \mathbb{D}^{(N+1)(n+m)}_{++}$. Therefore, problem (5.45) can be solved using Algorithm 13, since it satisfies Assumtion 5.8.

**Remark 5.34.** *Remark 5.20.(ii) also applies to problem* (5.45).

Step 4 of Algorithm 3 minimizes the Lagrangian (2.14) over $z_2 = (x_s, u_s)$, resulting in the following equality-constrained QP problem:

$$\min_{z_2} \frac{1}{2} z_2^\top H_2 z_2 + q_2^\top z_2 \tag{5.46a}$$

$$s.t.\ G_2 z_2 = b_2, \tag{5.46b}$$

where

$$\begin{aligned} H_2 &= \text{diag}(T, S) + \rho C_2^\top C_2, \\ q_2 &= -(T x_r, S u_r) + \rho C_2^\top C_1 z_1^{k+1} + \rho C_2^\top C_3 z_3^k + C_2^\top \lambda^k - \rho C_2^\top b, \\ G_2 &= [(A - \mathbf{I}_n) \ \ B], \\ b_2 &= \mathbf{0}_n. \end{aligned}$$

The solution of this problem is given by Proposition 5.5. Note that, since the prediction model (5.42c) is assumed to be controllable, $G_2$ is full rank. Moreover, $H_2 \in \mathbb{S}^{n+m}_{++}$. Therefore, the optimal solution of (5.46) can be obtained by substituting (5.3a) into (5.3b), which leads to the expression

$$z_2^* = M_2 q_2,$$

where $M_2 = H_2^{-1} G_2^\top (G_2 H_2^{-1} G_2^\top)^{-1} G_2 H_2^{-1} - H_2^{-1} \in \mathbb{R}^{(n+m)\times(n+m)}$ is computed offline and $b_2$ does not appear because it is equal to zero.

Step 5 of Algorithm 3 minimizes the Lagrangian (2.14) over $z_3$, resulting in the following equality-constrained QP problem:

$$\min_{z_3} \frac{1}{2} z_3^\top H_3 z_3 + q_3^\top z_3$$

$$s.t.\ G_3 z_3 = b_3,$$

where

$$\begin{aligned} H_3 &= \text{diag}(Q, R, Q, R, \dots, Q, R) + \rho C_3^\top C_3, \\ q_3 &= \rho C_3^\top C_1 z_1^{k+1} + \rho C_3^\top C_2 z_2^{k+1} + C_3^\top \lambda^k - \rho C_3^\top b, \\ G_3 &= \begin{bmatrix} A & B & -\mathbf{I}_n & \mathbf{0} & \cdots & \cdots & \mathbf{0} & \mathbf{0} \\ \mathbf{0} & \mathbf{0} & A & B & -\mathbf{I}_n & \cdots & \mathbf{0} & \mathbf{0} \\ & & & & & & & \\ \mathbf{0} & \mathbf{0} & \cdots & \ddots & \ddots & \ddots & \mathbf{0} & \mathbf{0} \\ \mathbf{0} & \mathbf{0} & \cdots & \mathbf{0} & A & B & -\mathbf{I}_n & \mathbf{0} \end{bmatrix}, \\ b &= \mathbf{0}_{Nn}, \end{aligned}$$

**Algorithm 19:** Extended ADMM for MPCT

**Require:** $z_2^0$, $z_3^0$, $\lambda^0$, $\rho > 0$, $\epsilon > 0$

Update $b$ with $x(t)$

Update $q_2$ with $x_r$ and $u_r$

$k \leftarrow 0$

**repeat**

  $q_1 \leftarrow \rho C_1^\top C_2 z_2^k + \rho C_1^\top C_3 z_3^k + C_1^\top \lambda^k - \rho C_1^\top b$

  $z_1^{k+1} \leftarrow \text{solve_boxQP}(q_1; H_1, \underline{z}_1, \overline{z}_1)$

  $q_2 \leftarrow -(Tx_r, Su_r) + \rho C_2^\top C_1 z_1^{k+1} + \rho C_2^\top C_3 z_3^k + C_2^\top \lambda^k - \rho C_2^\top b$

  $z_2^{k+1} \leftarrow M_2 q_2$

  $q_3 \leftarrow \rho C_3^\top C_1 z_1^{k+1} + \rho C_3^\top C_2 z_2^{k+1} + C_3^\top \lambda^k - \rho C_3^\top b$

  $z_3^{k+1} \leftarrow \text{solve_eqQP}(q_3, \mathbf{0}; H_3, G_3)$

  $\Gamma \leftarrow \sum\limits_{i=1}^{3} C_i z_i^{k+1} - b$

  $\lambda^{k+1} \leftarrow \lambda^k + \rho\Gamma$

  $k \leftarrow k + 1$

**until** $\|\Gamma\|_\infty \leq \epsilon$, $\|z_2^k - z_2^{k-1}\|_\infty \leq \epsilon$, $\|z_3^k - z_3^{k-1}\|_\infty \leq \epsilon$

**Output:** $\tilde{z}_1^* \leftarrow z_1^k$, $\tilde{z}_2^* \leftarrow z_2^k$, $\tilde{z}_3^* \leftarrow z_3^k$, $\tilde{\lambda}^* \leftarrow \lambda^k$

which, due to the block diagonal structure of $H_3$ and the banded structure of $G_3$ satisfies Assumtion 5.3 and can therefore be solved using Algorithm 12.

Algorithm 19 shows the particularization of Algorithm 3 to the MPCT formulation (5.43) obtained from the above discussion. The control action $u(t)$ is obtained from the elements $u_0$ of the output $\tilde{z}^*$ (5.44a) returned by Algorithm 19.

**Remark 5.35.** *It has been shown that the performance of ADMM can be significantly improved by having different values of $\rho$ for different constraints [141, §5.2], i.e., by considering $\rho$ as a diagonal positive definite matrix. In particular, we find that the convergence of Algorithm 19 improves significantly if the equality constraints* (5.43b), (5.43k), (5.43l), (5.43j) *for $i = N$, and* (5.43i) *for $i = 0$ and $i = N$, are penalized more than the others.*

**Remark 5.36.** *We note the operations in Algorithm 19 with matrices $C_i$, $i \in \mathbb{Z}_1^3$, can be performed exclusively using vector-vector operations. This is due to the fact that the matrices $C_i^\top C_j$, $i \in \mathbb{Z}_1^3$, $j \in \mathbb{Z}_1^3$, are diagonal, with some having a very small number of off-diagonal elements that are easily identified. This is also true for the case in which $\rho$ is taken as a diagonal positive definite matrix (Remark 5.35). We maintain the expressions using matrices $C_i$ in Algorithm 19 for simplicity.*

**Remark 5.37.** *The theoretical upper bound for $\rho$ provided in Theorem 2.9 is easily computable in this case. Indeed, we have that $C_3^\top C_3$ is the identity matrix,*

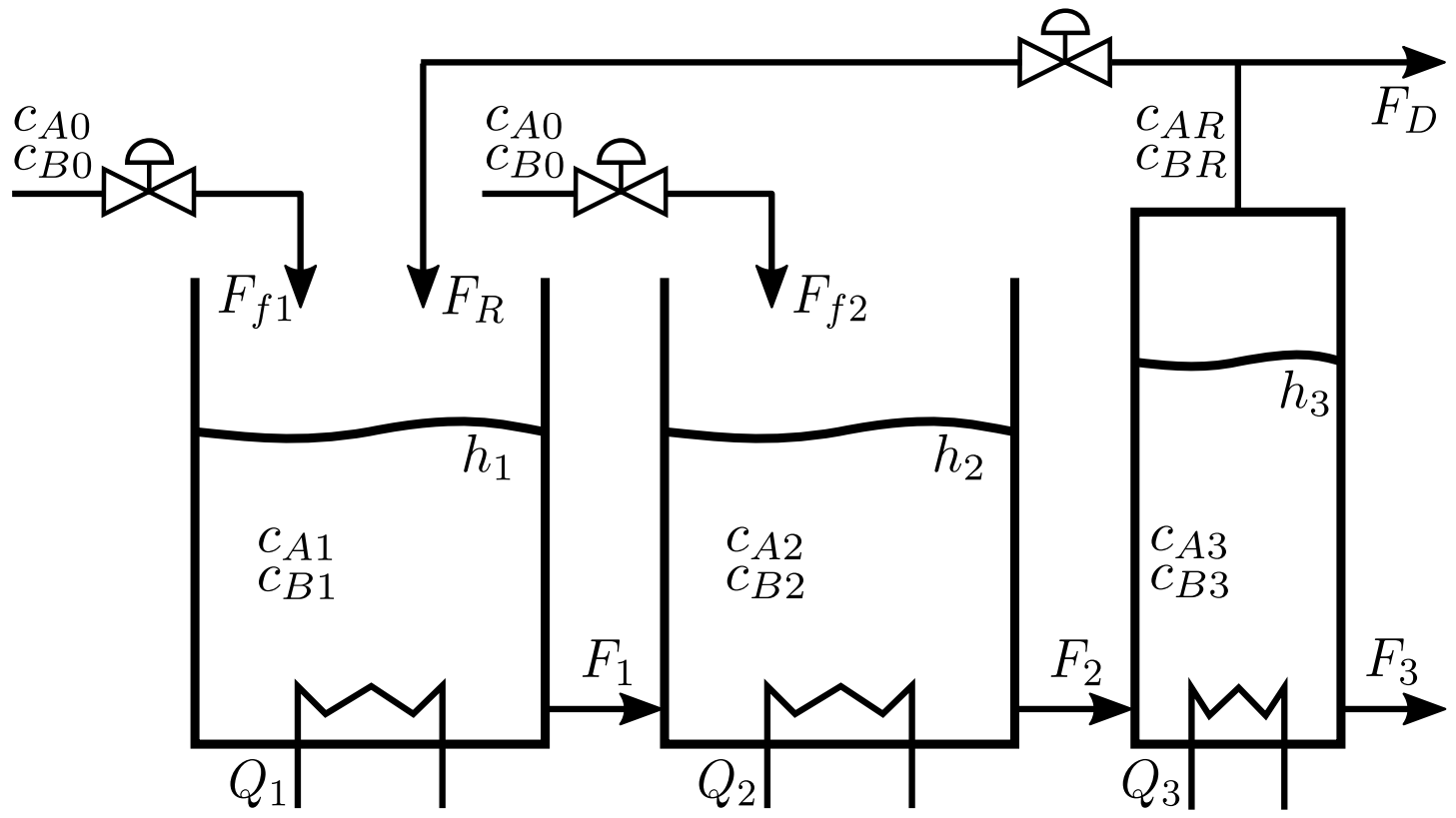

Figure 5.1: Double reactor and separator system.

*and therefore its spectral norm is* $\|C_3^\top C_3\| = 1$. *Furthermore,* $\mu_3$ *is the minimum eigenvalue of* $\text{diag}(Q, R)$*, which is simple to compute. However, in practice, we find that better results are obtained following Remark 5.35 for larger values of* $\rho$ *than the one described in Theorem 2.9. Therefore, the convergence of the algorithm will have to be extensively checked with simulations.*

## 5.7 Test Benches

This section presents various systems that will be used in the numerical results shown in the following sections. In particular, we present the following systems: a multivariable chemical plant consisting of two reactors and one separator, a ball and plate system, and a series of masses connected by springs.

### 5.7.1 Chemical plant: double reactor and separator

The chemical plant system, depicted in Figure 5.1 and inspired from [142], is a chemical plant in which the two following first-order reaction take place between the reactants A, B and C:

$$\text{A} \rightarrow \text{B}, \tag{5.47a}$$

$$\text{B} \rightarrow \text{C}. \tag{5.47b}$$

The reactions occur in two consecutive cylindrical reactors (labeled 1 and 2), which are fed by flows $F_{f1}$ and $F_{f2}$, respectively, of reactants A and B with concentrations $c_{A0}$ and $c_{B0}$ The reactants then flow into a cylindrical separator (labeled 3), where part of the reactants are redirected either to reactor 1 through flow $F_R$ or discarded through flow $F_D$. The concentration of the reactants in each cylinder is denoted by $c_{Xi}$, where $X \in \{\text{A}, \text{B}, \text{C}\}$ and $i \in \mathbb{Z}_1^3$. Each cylinder $i \in \mathbb{Z}_1^3$ has a bottom with area $A_i$ and a height of liquid $h_i$ that will be determined by the input and output flows of the cylinder, where the output flows $F_i$ are connected as shown in Figure 5.1. The kinetics of the reactions (5.47) taking place in reactors $i \in \mathbb{Z}_1^2$ are given by the Arrhenius equations

$$k_{Ai} = k_A e^{-\frac{E_A}{RT_i}}, \quad k_{Bi} = k_B e^{-\frac{E_B}{RT_i}},$$

respectively, where $e$ is Euler's number, $k_A$ and $k_B$ are the pre-exponential factor, $E_A$ and $E_B$ the activation energies, $R$ the universal gas constant, and $T_i$ the temperature of each reactor. The cylinders $i \in \mathbb{Z}_1^3$ have heating/cooling systems which transfer heats $Q_i$, thus affecting the temperatures and, as a result, the speed of the reactions.

| Parameter | Value | Units | Parameter | Value | Units |
|---|---|---|---|---|---|
| $A_1$ | 1 | $m^2$ | $T_0$ | 313 | K |
| $A_2$ | 1 | $m^2$ | $k_A$ | $10^{-5}$ | 1/s |
| $A_3$ | 1 | $m^2$ | $k_B$ | $5 \cdot 10^{-6}$ | 1/s |
| $\rho$ | 1100 | $kg/m^3$ | $E_A/R$ | -2840 | K |
| $C_p$ | 4 | kJ/kg K | $E_B/R$ | -2077 | K |
| $k_{v1}$ | 50 | kg/m s | $\Delta H_A$ | -100 | kJ/kg |
| $k_{v2}$ | 50 | kg/m s | $\Delta H_B$ | -39 | kJ/kg |
| $k_{v3}$ | 30 | kg/m s | $\alpha_A$ | 3.5 | - |
| $\alpha_C$ | 0.5 | - | $\alpha_B$ | 1.1 | - |
| $\alpha_D$ | 0.001 | - | $c_{A0}$ | 1 | wt(%) |
| $c_{B0}$ | 0 | wt(%) | | | |

Table 5.1: Parameters of the double reactor and separator system

| Var. | Val. [m] | Var. | Val. [wt(%)] | Var. | Val. [wt(%)] | Var. | Val. [K] |
|---|---|---|---|---|---|---|---|
| $h_1^\circ$ | 0.7 | $c_{A1}^\circ$ | 0.4155 | $c_{B1}^\circ$ | 0.5480 | $T_1^\circ$ | 329 |
| $h_2^\circ$ | 0.9 | $c_{A2}^\circ$ | 0.2581 | $c_{B2}^\circ$ | 0.6755 | $T_2^\circ$ | 333 |
| $h_3^\circ$ | 1.33 | $c_{A3}^\circ$ | 0.2282 | $c_{B3}^\circ$ | 0.7 | $T_3^\circ$ | 323 |

Table 5.2: Operating point of the double reactor and separator system

The non-linear model (4.5) of the plant is given by:

$$\frac{dh_1}{dt} = \frac{1}{\rho A_1}(F_{f1} + F_R - F_1), \tag{5.48a}$$

$$\frac{dc_{A1}}{dt} = \frac{1}{\rho A_1 h_1}\left(F_{f1}(c_{A0} - c_{A1}) + F_R(c_{AR} - c_{A1})\right) - k_{A1}c_{A1}, \tag{5.48b}$$

$$\frac{dc_{B1}}{dt} = \frac{1}{\rho A_1 h_1}\left(F_{f1}(c_{B0} - c_{B1}) + F_R(c_{BR} - c_{B1})\right) \tag{5.48c}$$

$$- k_{B1}c_{B1} + k_{A1}c_{A1}, \tag{5.48d}$$

$$\frac{dT_1}{dt} = \frac{1}{\rho A_1 h_1}\left(F_{f1}(T_0 - T_1) + F_R(T_R - T_1)\right) + \frac{Q_1}{\rho A_1 H_1 C_p} \tag{5.48e}$$

| Variable | Upper | Lower | Variable | Upper | Lower |
|---|---|---|---|---|---|
| $h_{1,2,3}$ | 2 | 0 | $c_i$ | 1 | 0 |
| $T_{1,2}$ | 348 | 320 | $T_3$ | 338 | 320 |
| $Q_{1,2,3}$ | 5000 | -5000 | $F_{f1,f2,R}$ | 50 | 0 |

Table 5.3: Upper and lower bounds for the double reactor and separator system

$$- \frac{1}{C_p}\left(k_{A1}c_{A1}\Delta H_A + k_{B1}c_{B1}\Delta H_B\right), \tag{5.48f}$$

$$\frac{dh_2}{dt} = \frac{1}{\rho A_2}(F_{f2} + F_1 - F_2), \tag{5.48g}$$

$$\frac{dc_{A2}}{dt} = \frac{1}{\rho A_2 h_2}\left(F_{f2}(c_{A0} - c_{A2}) + F_1(c_{A1} - c_{A2})\right) - k_{A2}c_{A2}, \tag{5.48h}$$

$$\frac{dc_{B2}}{dt} = \frac{1}{\rho A_2 h_2}\left(F_{f2}(c_{B0} - c_{B2}) + F_1(c_{B1} - c_{B2})\right) \tag{5.48i}$$

$$- k_{B2}c_{B2} + k_{A2}c_{A2}, \tag{5.48j}$$

$$\frac{dT_2}{dt} = \frac{1}{\rho A_2 h_2}\left(F_{f2}(T_0 - T_2) + F_1(T_1 - T_2)\right) + \frac{Q_2}{\rho A_2 H_2 C_p} \tag{5.48k}$$

$$- \frac{1}{C_p}\left(k_{A2}c_{A2}\Delta H_A + k_{B2}c_{B2}\Delta H_B\right), \tag{5.48l}$$

$$\frac{dh_3}{dt} = \frac{1}{\rho A_3}(F_2 - F_D - F_R - F_3), \tag{5.48m}$$

$$\frac{dc_{A3}}{dt} = \frac{1}{\rho A_3 h_3}\left(F_2(c_{A2} - c_{A3}) - (F_D + F_R)(c_{AR} - c_{A3})\right), \tag{5.48n}$$

$$\frac{dc_{B3}}{dt} = \frac{1}{\rho A_3 h_3}\left(F_2(c_{B2} - c_{B3}) - (F_D + F_R)(c_{BR} - c_{B3})\right), \tag{5.48o}$$

$$\frac{dT_3}{dt} = \frac{1}{\rho A_3 h_3}F_2(T_2 - T_3) + \frac{Q_3}{\rho A_3 H_3 C_p}, \tag{5.48p}$$

where $\Delta H_A$ and $\Delta H_B$ are the enthalpy of the reactions (5.47), respectively, $C_p$ is the specific heat of the reactants, and $\rho$ their density (we assume $C_p$ and $\rho$ to be the same for the three reactants), and the following relations hold for $i \in \mathbb{Z}_1^3$:

$$F_i = k_{vi}h_i, \qquad F_D = \alpha_D F_R,$$

$$c_{AR} = \frac{\alpha_A c_{A3}}{\bar{c}_3}, \qquad c_{BR} = \frac{\alpha_B c_{B3}}{\bar{c}_3},$$

$$\bar{c}_3 = \alpha_A c_{A3} + \alpha_B c_{B3} + \alpha_C c_{C3}, \qquad c_{C3} = 1 - c_{A3} - c_{B3},$$

where $k_{vi} \in \mathbb{R}_{>0}$ determine the relation between the flows $F_i$ and heights $h_i$; $\alpha_A$, $\alpha_D$ and $\alpha_C$ determine the relative amount of each reactant in flows $F_R$ and $F_D$; and $\alpha_D$ determines the amount of discarded material.

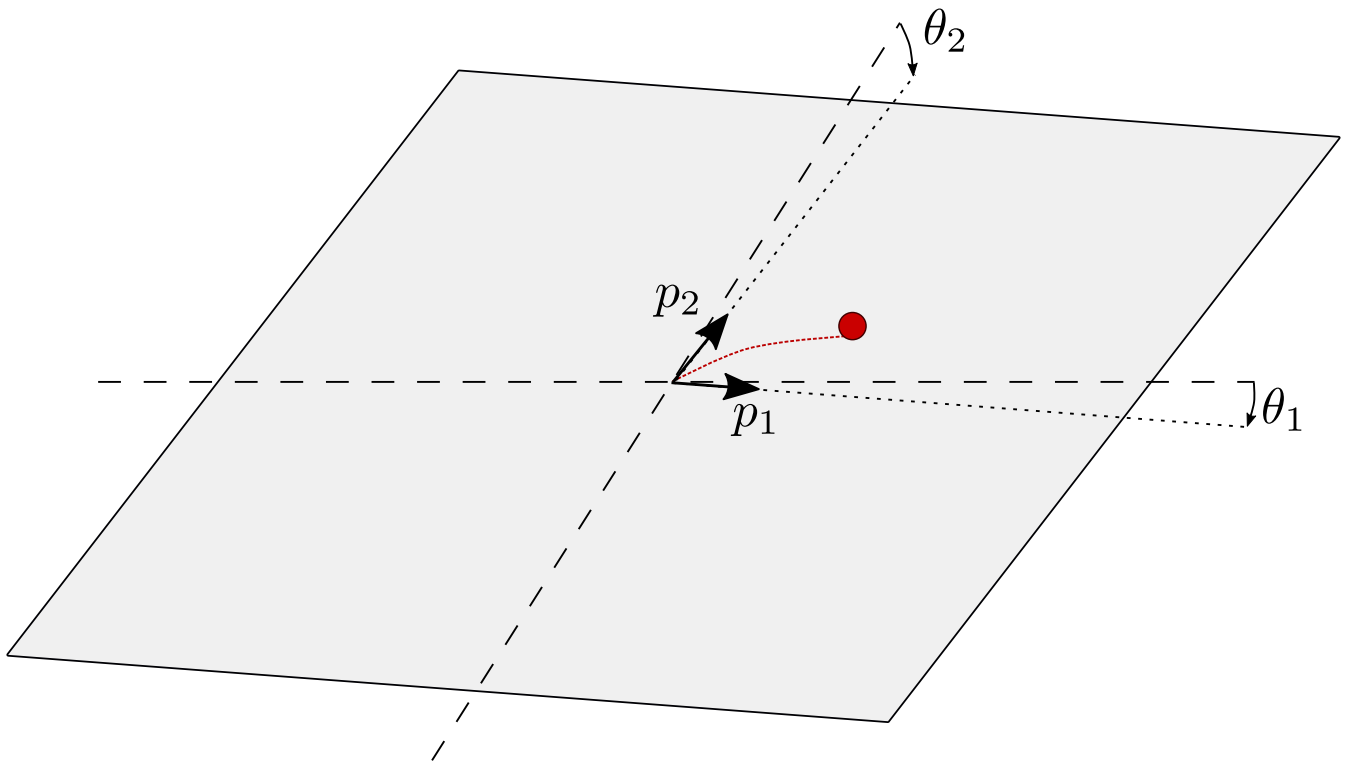

Figure 5.2: Ball and plate system.

The state and control input of the system are given by

$$\chi = (h_1, c_{A1}, c_{B1}, T_1, h_2, c_{A2}, c_{B2}, T_2, h_3, c_{A3}, c_{B3}, T_3),$$
$$u = (Q_1, Q_2, Q_3, F_{f1}, F_{f2}, F_R).$$

The values of the parameters of the system are given in Table 5.1

We obtain a linear model (4.1) of the system by linearizing (5.48) around the operating point $(\chi^\circ, u^\circ)$ described in Table 5.2 with a sample time of 3s and then scaling the resulting model with the scaling matrices

$$N_x = \text{diag}(1, 1, 1, 0.1, 1, 1, 1, 0.1, 1, 1, 1, 0.1),$$
$$N_u = \text{diag}(0.001, 0.001, 0.001, 0.1, 0.1, 0.1).$$

We consider the box constraints (4.3) on state and control inputs given in Table 5.3, where $c_i$ stands for all the concentrations $c_{A1}$, $c_{A2}$, $c_{A3}$, $c_{B1}$, $c_{B2}$ and $c_{B3}$.

### 5.7.2 Ball and plate

The ball and plate system, which is depicted in Figure 5.2, consists of a plate that pivots around its center point such that its slope can be manipulated by changing the angle of its two perpendicular axes. The objective is to control the position of a solid ball that rests on the plate. We assume that the ball is always in contact with the plate and that it does not slip when moving. The non-linear equations of the system are [143],

$$\ddot{p}_1 = \frac{m}{m + I_b/r^2}\left(p_1\dot{\theta}_1^2 + p_2\dot{\theta}_1\dot{\theta}_2 + g\sin\theta_1\right) \tag{5.49a}$$
$$\ddot{p}_2 = \frac{m}{m + I_b/r^2}\left(p_2\dot{\theta}_2^2 + p_1\dot{\theta}_1\dot{\theta}_2 + g\sin\theta_2\right), \tag{5.49b}$$

where $m = 0.05$Kg, $r = 0.01$m and $I_b = (2/5)mr^2 = 2 \cdot 10^{-6}$Kg·m$^2$ are the mass, radius, mass moment of inertia of a solid ball, respectively; $g = 9,81$m/s$^2$ is the gravitational acceleration; $p_1$ and $p_2$ are the position of the ball on the two axes of the plate relative to its center point; $\dot{p}_1$, $\dot{p}_2$, $\ddot{p}_1$ and $\ddot{p}_2$ their corresponding velocities and accelerations; $\theta_1$ and $\theta_2$ are the angle of the plate on each of its axes; and $\dot{\theta}_1$ and $\dot{\theta}_2$ their corresponding angular velocities.

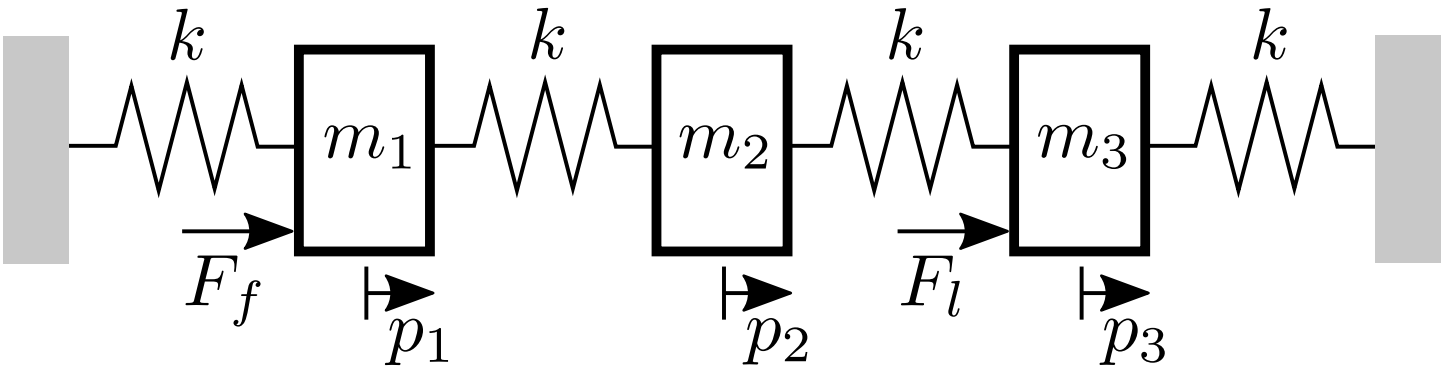

Figure 5.3: Chain of three masses connected by springs.

The state of the system is given by

$$x = (p_1, \dot{p}_1, \theta_1, \dot{\theta}_1, p_2, \dot{p}_2, \theta_2, \dot{\theta}_2),$$

and the control input $u = (\ddot{\theta}_1, \ddot{\theta}_2)$ is the angle acceleration of the plate in each one of its axes. We consider the following box constraints (4.3) on the state and control inputs:

$$|\dot{p}_i| \leq 0.5\,\mathrm{m/s}^2, \; |\theta_i| \leq \frac{\pi}{4}\,\mathrm{rad}, \; |\ddot{\theta}_i| \leq 0.4\,\mathrm{rad/s}^2, \; i \in \mathbb{Z}_1^2.$$

We obtain a linear model (4.1) by linearizing (5.49) around the operating point $x^\circ = \mathbf{0}_6$, $u^\circ = \mathbf{0}_2$ with a sample time of 0.2s and then scaling the resulting model with the scaling matrices

$$N_x = \mathrm{diag}(0.1, 1, 1, 1, 0.1, 1, 1, 1), \quad N_u = \mathbf{I}_2.$$

### 5.7.3 Oscillating masses

The oscillating masses system, which is inspired by the case study from [144], consists of three objects connected by springs as illustrated in Figure 5.3.

We take the mass of the outer objects as $m_1 = m_3 = 1$ and the mass of the central object as $m_2 = 0.5$. The spring constants are all taken as $k = 2$. There are two external forces acting on the system: a force $\mathcal{F}_f$ acting on the first object, and a force $\mathcal{F}_l$ acting on the last object, as illustrated in the figure. The state and control inputs of the system are given by

$$x = (p_1, p_2, p_3, v_1, v_2, v_3), \quad u = (\mathcal{F}_f, \mathcal{F}_l),$$

where $p_i$ and $v_i$, $i \in \mathbb{Z}_1^3$ are the position and velocity of each mass, respectively. The continuous-time (linear) dynamics of the system are given by

$$\begin{aligned} m_1 \dot{v}_1 &= k(p_2 - 2p_1) + \mathcal{F}_f \\ m_2 \dot{v}_2 &= k(p_1 + p_3 - 2p_2) \\ m_3 \dot{v}_3 &= k(p_2 - 2p_3) + \mathcal{F}_l. \end{aligned}$$

We compute a model (4.1) by taking a 0.2s sampling time and then scaling the resulting model using the scaling matrices

$$N_x = \text{diag}(10, 10, 10, 1, 1, 1), \quad N_u = I_2.$$

We consider the following box constraints (4.3) on the states and control inputs:

$$|p_i| \leq 0.3, \ i \in \mathbb{Z}_1^3, \quad |\mathcal{F}_j| \leq 0.8, \ j \in \{f, l\}.$$

The velocities $v_i$, $i \in \mathbb{Z}_1^3$, are not constrained.

## 5.8 Numerical results

This section presents numerical results using the proposed solvers to control the systems described in Section 5.7.

We present novel results using the latest version of the SPCIES toolbox [4], which at the time of writing of this dissertation is `v0.3.2`. Other solvers used for comparison are: OSQP (version `0.6.0`) [35], and qpOASES (version `3.2.0`) [93]. We use the default options of both solvers, with the exception of the OSQP solver, where the exit tolerances are set to $10^{-4}$ and its warmstart procedure is disabled to provide a better comparison between its underlying ADMM algorithm with the ones proposed in this manuscript. Furthermore, printing information to the console was disabled in both solvers due to the significant effect that it has on computation times.

Additional results showcasing the performance of the proposed solvers can be found in the papers [8, 9, 10, 11, 12]. In particular, [8] shows hardware-in-the-loop results of the standard MPC formulations (5.16) and (5.17) controlling the chemical plant described in Section 5.7.1 using a real industrial PLC; [9] shows similar hardware-in-the-loop results but for a smaller sized system; [10] shows some preliminary results of the solver for the MPCT formulation (5.42); [11] shows the implementation of the solver for the MPCT formulation in a Raspberry Pi to control an inverted pendulum robot, although the computation times are larger due to the use of the less refined C-code available at the time; and [12] shows the implementation of the MPC formulation with terminal quadratic constraint (5.20) in a Rapsberry Pi to control the oscillating masses system described in Section 5.7.3, where we compare it with other alternatives.

The results shown in this section have been obtained using a PC running the Ubuntu 18.04.5 LTS operating system on a Intel Core i5-8250U CPU operating at its base frequency of 1.60 GHz. For results on computation times in embedded systems we refer the reader to the above references.

The systems are simulated using the same linear models that are used as the prediction models of the MPC formulations. The plots and references, however, are shown in *engineering units* because they are more intuitive to visualize than the *incremental units* of the linear models.

| | $Q$ | $R$ | $N$ |
|---|---|---|---|
| Chemical plant | $5\mathbf{I}_{12}$ | $0.5\mathbf{I}_6$ | 20 |
| Ball and plate | $\text{diag}(10, 0.05, 0.05, 0.05, 10, 0.05, 0.05, 0.05)$ | $0.5\mathbf{I}_2$ | 30 |
| Oscillating masses | $\text{diag}(15, 15, 15, 1, 1, 1)$ | $0.1\mathbf{I}_2$ | 10 |

Table 5.4: Parameters of the MPC formulations used for each test bench systems.

| | $x_r$ | $u_r$ |
|---|---|---|
| Chemical plant | $(0.7, 0.419, 0.545, 329.571,$ $0.9, 0.261, 0.673, 333.435,$ $1.333, 0.231, 0.698, 337.602)$ | $(0, 0, 750, 30, 10, 5)$ |
| Ball and plate | $(1.8, 0, 0, 0, 1.4, 0, 0, 0)$ | $(0, 0)$ |
| Oscillating masses | $(0.25, 0.25, 0.25, 0, 0, 0)$ | $(0.5, 0.5)$ |

Table 5.5: References used for each test bench systems (in engineering units).

For convenience, throughout this section we will refer to the standard MPC formulations (5.16) and (5.17) as *equMPC* and *laxMPC*, respectively, the MPC formulation subject to a terminal ellipsoidal constraint (5.20) as *ellipMPC* and the MPC for tracking formulation (5.42) as *MPCT*.

### 5.8.1 Comparison between the proposed ADMM-based solvers

We start by comparing the ADMM-based solvers of each MPC formulation described in Sections 5.4 to 5.6 by performing tests on the three systems described in Section 5.7. In every test we use the matrices $Q$, $R$ and prediction horizon $N$ shown in Table 5.4 for the four MPC formulations. The ingredients $P$, $c$ and $r$ defining the terminal constraint of the *ellipMPC* formulation are computed using the procedure described in Section 5.5.2, where the LMI optimization problem is constructed using YALMIP [103] and solved using the SDPT3 solver [145]. In all cases, a feasible solution of problem (5.41) is found for $r = 1$ and $\lambda = 0.95$, from where we obtain $P$ and its associated control gain $K$. Using this gain $K$, matrix $T$ of the MPCT formulation is then computed as the solution of the Lyapunov equation

$$(A + BK)^\top T(A + BK) - T = -Q - K^\top RK. \tag{5.50}$$

The *laxMPC* and *equMPC* formulations are solved using Algorithm 17, *ellipMPC* using Algorithm 18 and *MPCT* using Algorithm 19. The penalty parameters are listed in Table 5.6. The exit tolerances are all set to $10^{-4}$.

In each test, the system is initialized in its operating point and the reference is set to the values provided in Table 5.5. We show the trajectory of one of the system states and one of the system outputs obtained with each one of the solvers, as well as the number of iterations and computation times at each sample time.

Figure 5.4 shows the comparison between the solvers applied to the chemical plant described in Section 5.7.1. Figure 5.4a shows the trajectory of the temper-ature $T_3$, Figure 5.4b the trajectory of the input $Q_3$, Figure 5.4c the computation times, and Figure 5.4d the number of iterations of each solver. Table 5.7 shows an analysis of the number of iterations and computation times of each solver.

| | laxMPC | equMPC | ellipMPC | MPCT |
|---|---|---|---|---|
| Chemical plant | 15 | 15 | 15 | $\rho_1 = 2,\ \rho_2 = 40$ |
| Ball and plate | 15 | 15 | 15 | $\rho_1 = 10,\ \rho_2 = 200$ |
| Oscillating masses | 15 | 15 | 15 | $\rho_1 = 2,\ \rho_2 = 40$ |

Table 5.6: Value of the penalty parameter $\rho$ used for each one of the MPC formulations and test benches. In the case of the MPCT formulation, $\rho_2$ indicates the value of $\rho$ for the constraints listed in Remark 5.35 and $\rho_1$ for the rest.

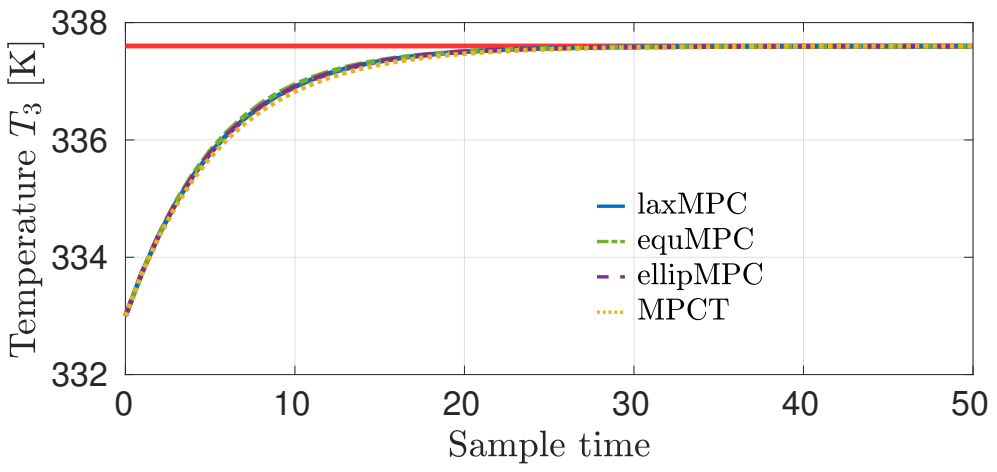

(a) State trajectory.

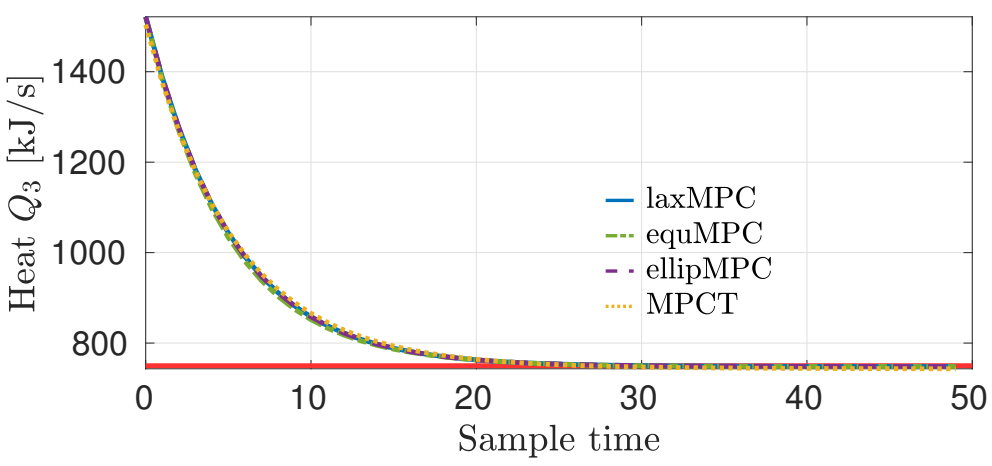

(b) Control action trajectory.

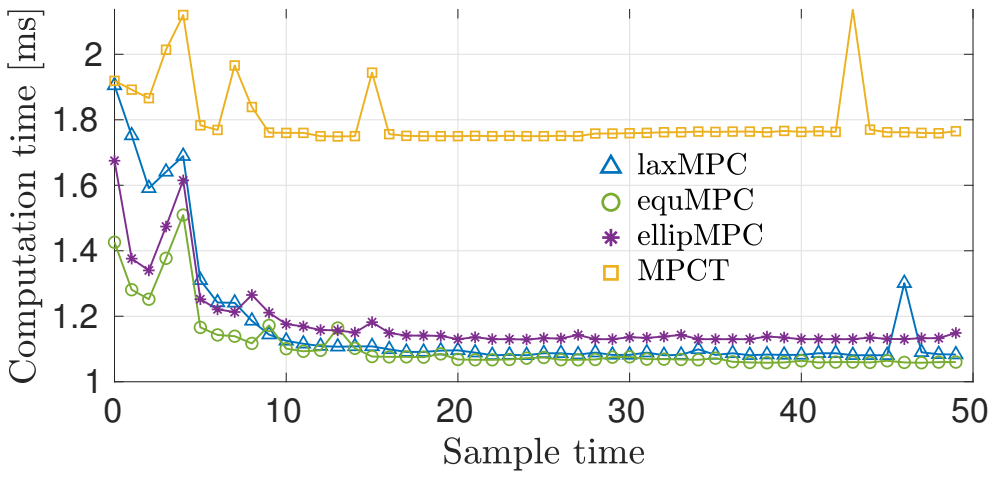

(c) Computation times.

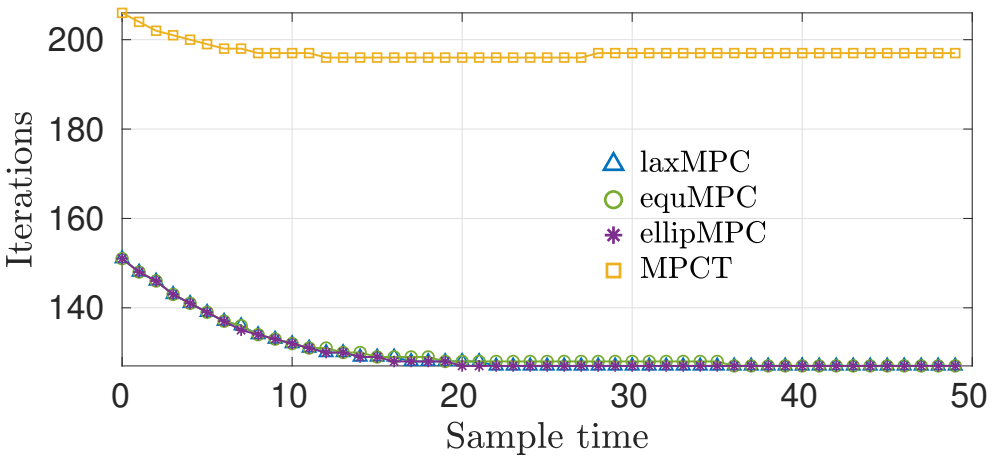

(d) Number of iterations.

Figure 5.4: Closed-loop simulation of the chemical plant: Comparison between the ADMM-based solvers.

Figure 5.5 and Table 5.8 show analogous results to Figure 5.4 and Table 5.7, respectively, but for the ball and plate system described in Section 5.7.2. Figure 5.5a shows the trajectory of the state $\dot{p}_1$ and Figure 5.5b the trajectory of the control input $\ddot{\theta}_1$.

Figure 5.6 and Table 5.9 show analogous results to Figure 5.4 and Table 5.7, respectively, but for oscillating masses system described in Section 5.7.3. Figure 5.6a shows the trajectory of the state $p_2$ and Figure 5.6b the trajectory of the control input $\mathcal{F}_l$.

The results indicate that there are no major differences, in terms of number of iterations and computation times, between the four solvers. The EADMM algorithm tends to present larger differences between the maximum and minimum number of iterations, whereas the number of iterations of the ADMM solvers show less variation.

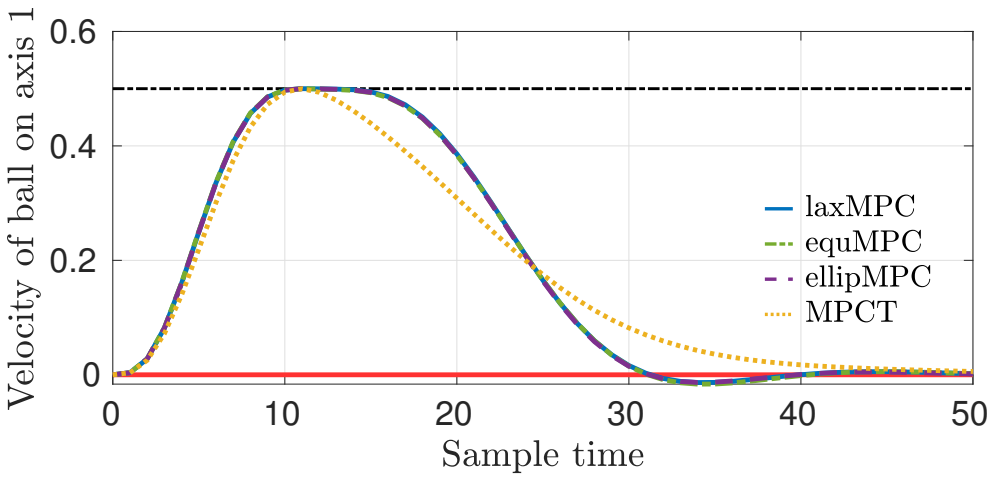

(a) State trajectory.

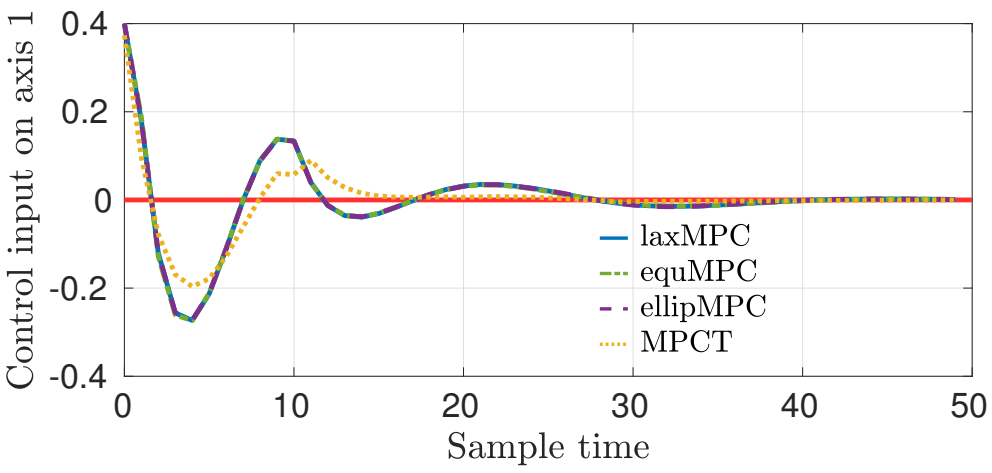

(b) Control action trajectory.

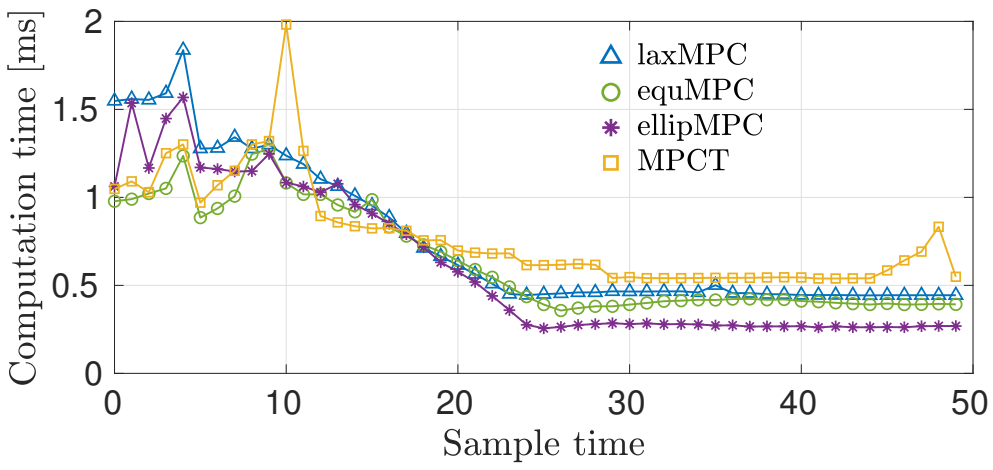

(c) Computation times.

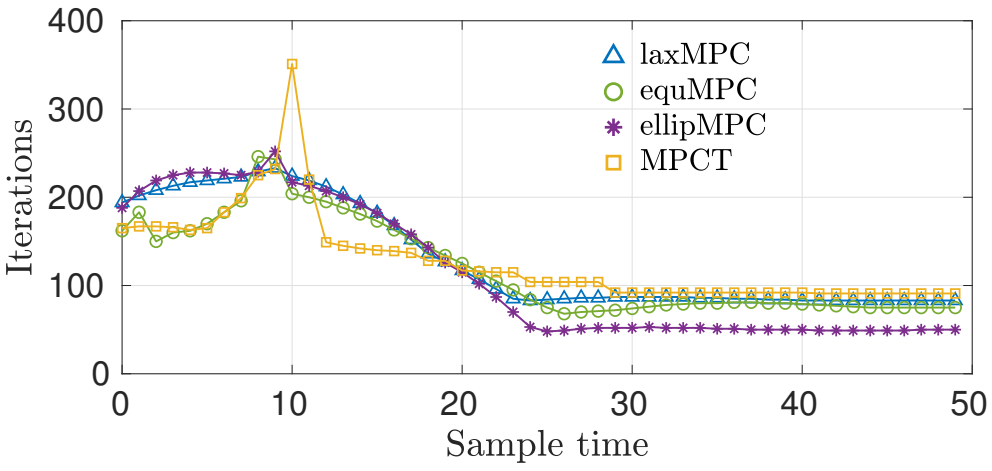

(d) Number of iterations.

Figure 5.5: Closed-loop simulation of the ball and plate: Comparison between the ADMM-based solvers. The black dashed/dotted line in figure (a) represents the upper bound on the velocity of the ball.

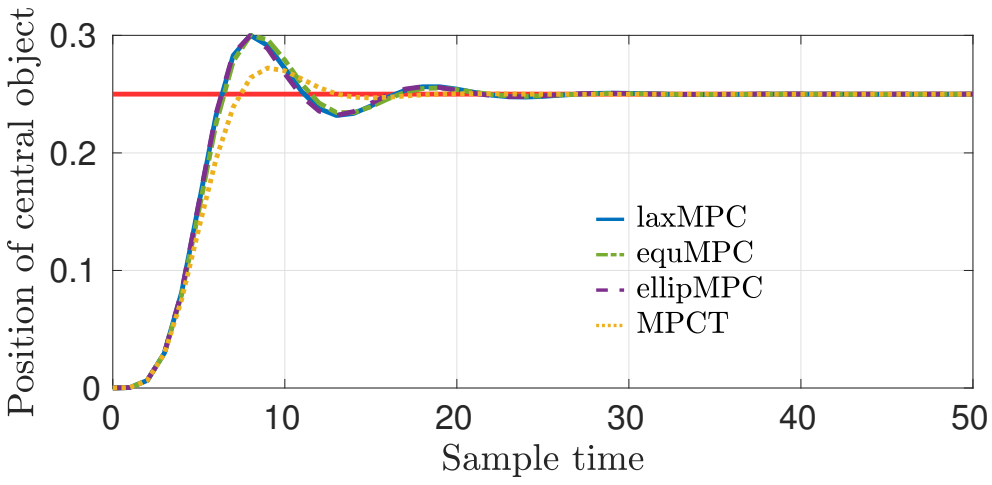

(a) State trajectory.

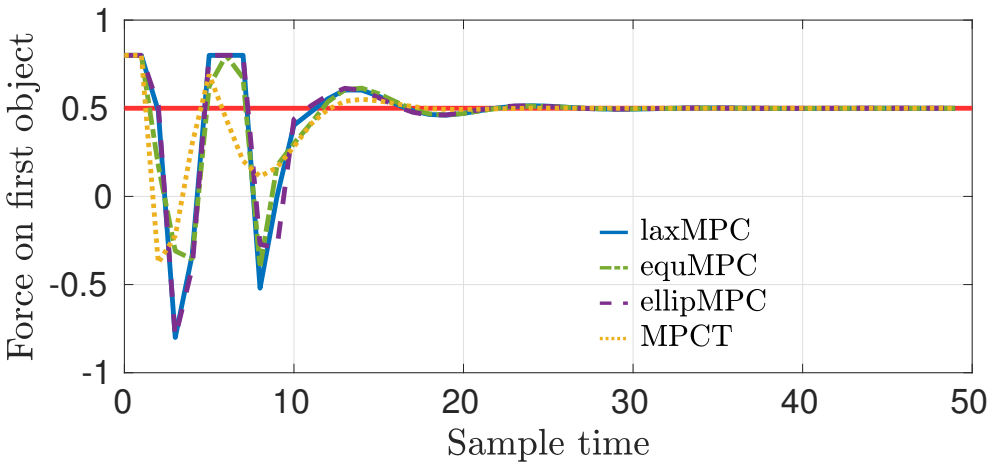

(b) Control action trajectory.

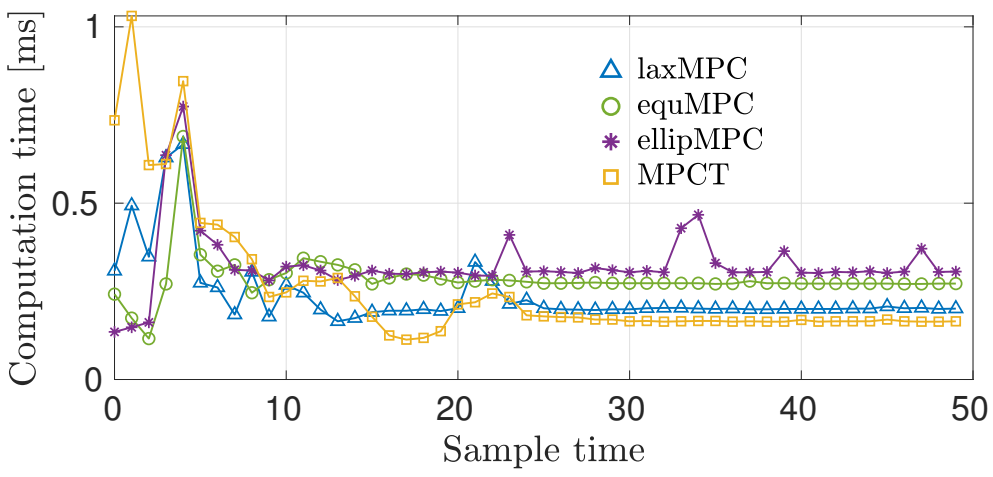

(c) Computation times.

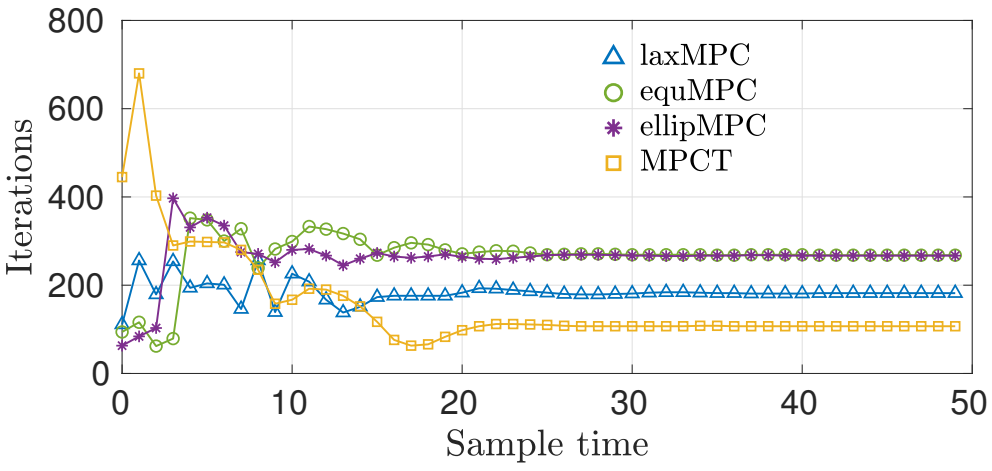

(d) Number of iterations.

Figure 5.6: Closed-loop simulation of the oscillating masses: Comparison between the ADMM-based solvers.

| Formulation | | laxMPC | equMPC | ellipMPC | MPCT |
|---|---|---|---|---|---|
| Iterations | Average | 130.28 | 130.64 | 130.2 | 197.32 |
| | Median | 127 | 128 | 127 | 197 |
| | Maximum | 151 | 151 | 151 | 206 |
| | Minimum | 127 | 127 | 127 | 196 |
| Comp. time | Average | 1.170 | 1.109 | 1.184 | 1.796 |
| | Median | 1.089 | 1.071 | 1.137 | 1.762 |
| | Maximum | 1.905 | 1.509 | 1.675 | 2.137 |
| | Minimum | 1.081 | 1.058 | 1.129 | 1.749 |

Table 5.7: Comparison between the ADMM-based solvers for the chemical plant: number of iterations and computation times [ms].

| Formulation | | laxMPC | equMPC | ellipMPC | MPCT |
|---|---|---|---|---|---|
| Iterations | Average | 131.58 | 120.36 | 114.44 | 128.34 |
| | Median | 87 | 82.5 | 53 | 104 |
| | Maximum | 233 | 246 | 252 | 351 |
| | Minimum | 83 | 68 | 48 | 91 |
| Comp. time | Average | 0.866 | 0.674 | 0.667 | 0.854 |
| | Median | 0.55 | 0.571 | 0.479 | 0.715 |
| | Maximum | 2.39 | 1.248 | 1.574 | 2.438 |
| | Minimum | 0.442 | 0.357 | 0.26 | 0.539 |

Table 5.8: Comparison between the ADMM-based solvers for the ball and plate: number of iterations and computation times [ms].

| Formulation | | laxMPC | equMPC | ellipMPC | MPCT |
|---|---|---|---|---|---|
| Iterations | Average | 183.84 | 265.9 | 262.52 | 157.76 |
| | Median | 182 | 269 | 267 | 107 |
| | Maximum | 256 | 352 | 397 | 680 |
| | Minimum | 111 | 62 | 63 | 63 |
| Comp. time | Average | 0.239 | 0.286 | 0.326 | 0.259 |
| | Median | 0.201 | 0.273 | 0.306 | 0.17 |
| | Maximum | 0.667 | 0.689 | 0.774 | 1.031 |
| | Minimum | 0.165 | 0.116 | 0.135 | 0.113 |

Table 5.9: Comparison between the ADMM-based solvers for the oscillating masses: number of iterations and computation times [ms].

We note that the algorithms still require a significant number of iterations even when the optimal solution does not have any active constraints or when close to the reference. These iterations could be reduced by incorporating a warmstart procedure, such as the ones presented in [146] or [147]. For instance, in [10, §7], we particularize the warmstart procedure from [146, §II] to the EADMM algorithm for solving the MPCT formulation (Algorithm 19). Its particularization to the other solvers would follow similarly.

Finally, we note that the results for the *ellipMPC* formulation are very similar to the ones obtained for the standard MPC formulations (*laxMPC* and *equMPC*), in spite of the addition of the terminal quadratic constraint. This is even the case during the first 7 sample times in Figure 5.6d, during which the terminal quadratic constraint was active. This indicates that the QCQP problem (5.20) solved using ADMM is comparable, in terms of number of iterations and computation times, to solving standard MPC formulations.

### 5.8.2 Comparison between the ADMM- and FISTA-based solvers

We now compare the sparse solvers for standard MPC formulations presented in Section 5.4. That is, we compare the ADMM-based solve given in Algorithm 17 with the FISTA-based solve given in Algorithm 16 applied to the standard MPC formulations (5.16) (labeled *equMPC*) and (5.17) (labeled *laxMPC*).

We present the results of applying the four solvers to the three systems described in Section 5.7. The values of $Q$, $R$ and $N$ are shown in Table 5.4. We obtain matrix $T$ by first computing the matrix $\hat{T}$ that solves the Riccati equation

$$A^\top \hat{T} A - \hat{T} - (A^\top \hat{T} B)(R + B^\top \hat{T} B)^{-1}(B^\top \hat{T} A) + Q = \mathbf{0}_{6\times 6},$$

and then taking $T$ as the diagonal matrix satisfying

$$T_{(i,i)} = \sum_{j=0}^{n} \hat{T}_{(i,j)}.$$

We note that the reason why matrix $\hat{T}$ is not used, as is typically done in MPC, is because the FISTA-based solvers require the cost function matrices to be diagonal (see Assumtion 5.18). We take the above $T$ because it is an easy choice, although in a real setting we should check that it satisfies the Lyapunov equation (5.50).

We start the systems at the operating point and set the reference to the admissible steady states provided in Table 5.5. The exit tolerances of the solvers (ADMM and FISTA) are all set to $10^{-4}$ and the penalty parameter of the ADMM algorithm is selected as $\rho = 15$ in all cases. We show the number of iterations and computation times at each sample time. The state and input trajectories are nearly indistinguishable to the naked eye when compared to the ones shown in Section 5.8.1.

Figure 5.7 shows the comparison between the solvers applied to the chemical plant described in Section 5.7.1. Figure 5.7a shows the computation times and

| Formulation | | laxMPC | | equMPC | |
|---|---|---|---|---|---|
| FOM | | ADMM | FISTA | ADMM | FISTA |
| Iterations | Average | 129.46 | 1 | 130.64 | 1 |
| | Median | 126 | 1 | 128 | 1 |
| | Maximum | 151 | 1 | 151 | 1 |
| | Minimum | 126 | 1 | 127 | 1 |
| Comp. time | Average | 1.193 | 0.048 | 1.222 | 0.048 |
| | Median | 1.083 | 0.031 | 1.15 | 0.033 |
| | Maximum | 2.12 | 0.382 | 1.903 | 0.417 |
| | Minimum | 1.063 | 0.029 | 1.055 | 0.032 |

Table 5.10: Comparison between the ADMM- and FISTA-based solvers for standard MPC applied to the chemical plant: number of iterations and computation times [ms].

| Formulation | | laxMPC | | equMPC | |
|---|---|---|---|---|---|
| FOM | | ADMM | FISTA | ADMM | FISTA |
| Iterations | Average | 113.78 | 138.86 | 120.36 | 130.1 |
| | Median | 73 | 1 | 82.5 | 1 |
| | Maximum | 243 | 1203 | 246 | 1089 |
| | Minimum | 59 | 1 | 68 | 1 |
| Comp. time | Average | 0.682 | 0.789 | 0.653 | 0.728 |
| | Median | 0.453 | 0.0245 | 0.433 | 0.025 |
| | Maximum | 1.55 | 6.288 | 1.587 | 5.59 |
| | Minimum | 0.33 | 0.023 | 0.357 | 0.023 |

Table 5.11: Comparison between the ADMM- and FISTA-based solvers for standard MPC applied to the ball and plate: number of iterations and computation times [ms].

| Formulation | | laxMPC | | equMPC | |
|---|---|---|---|---|---|
| FOM | | ADMM | FISTA | ADMM | FISTA |
| Iterations | Average | 193.26 | 24.24 | 265.9 | 26.96 |
| | Median | 186 | 1 | 269 | 1 |
| | Maximum | 307 | 360 | 352 | 279 |
| | Minimum | 102 | 1 | 62 | 1 |
| Comp. time | Average | 0.249 | 0.058 | 0.289 | 0.06 |
| | Median | 0.204 | 0.017 | 0.275 | 0.017 |
| | Maximum | 0.849 | 0.517 | 0.762 | 0.5 |
| | Minimum | 0.186 | 0.016 | 0.119 | 0.016 |

Table 5.12: Comparison between the ADMM- and FISTA-based solvers for standard MPC applied to the oscillating masses: number of iterations and computation times [ms].

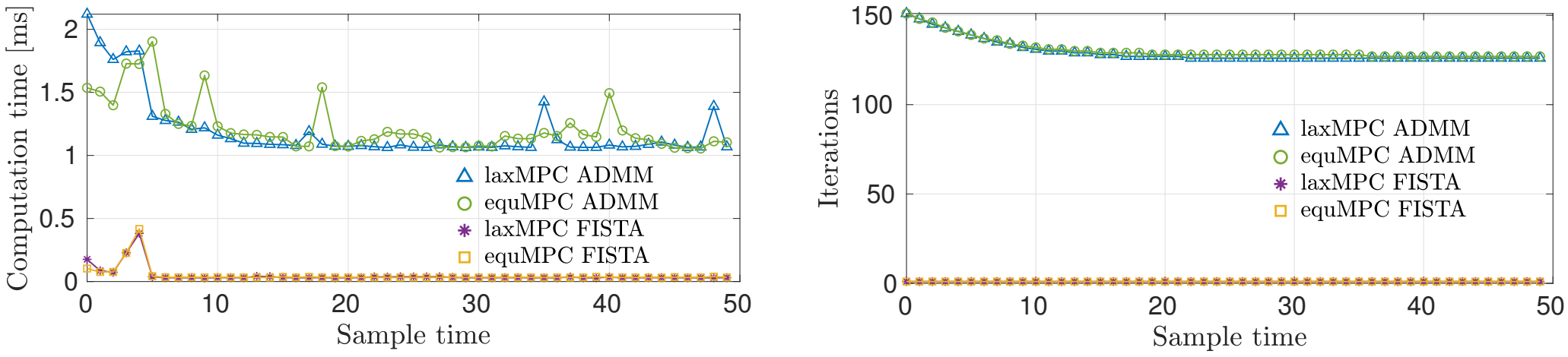

(a) Computation times. (b) Number of iterations.

Figure 5.7: Closed-loop simulation of the chemical plant: Comparison between the ADMM- and FISTA-based solvers for the standard MPC formulations.

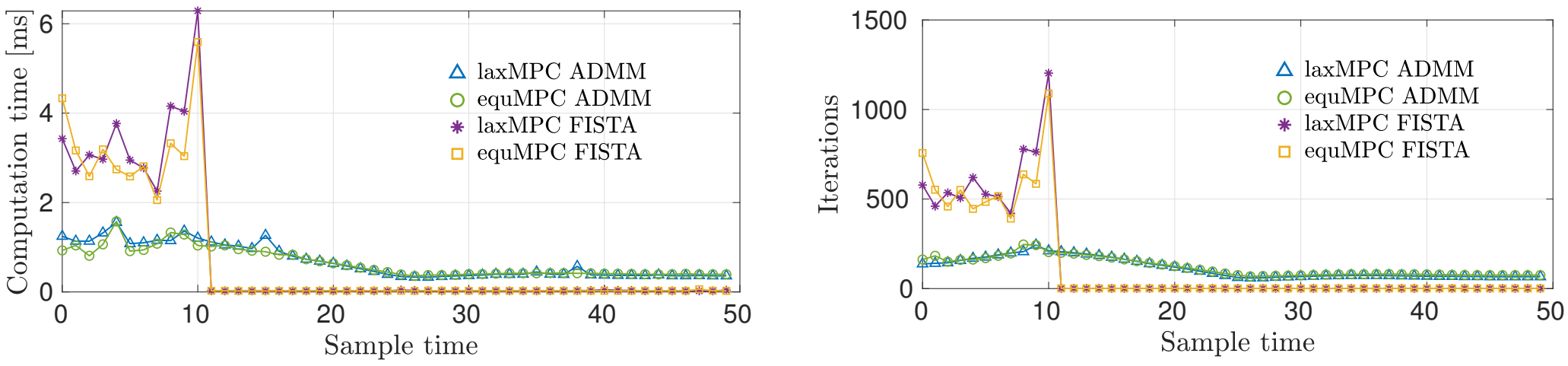

(a) Computation times. (b) Number of iterations.

Figure 5.8: Closed-loop simulation of the ball and plate system: Comparison between the ADMM- and FISTA-based solvers for the standard MPC formulations.

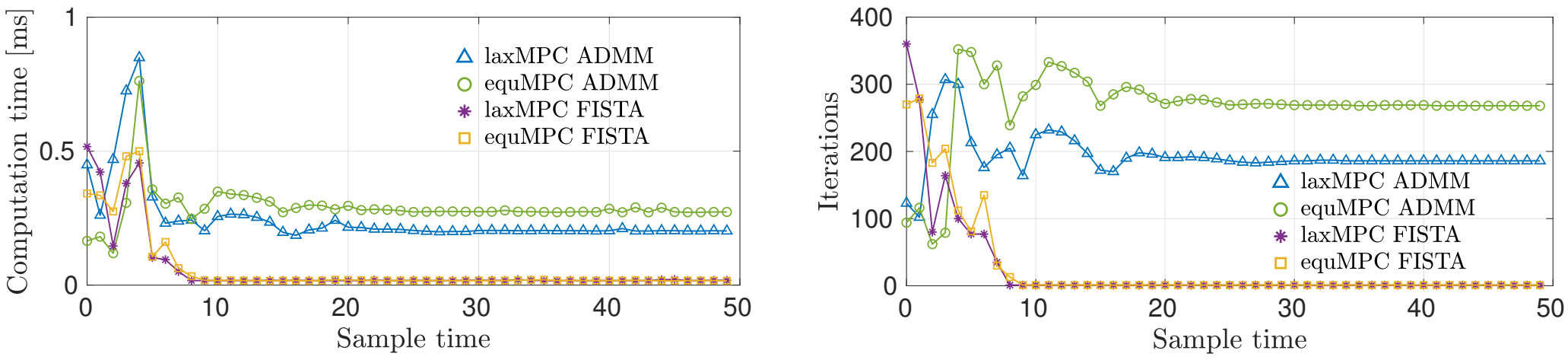

(a) Computation times. (b) Number of iterations.

Figure 5.9: Closed-loop simulation of the oscillating masses system: Comparison between the ADMM- and FISTA-based solvers for the standard MPC formulations.

Figure 5.7b the number of iterations of each solver. The state and input trajectories are very similar to the ones shown in Figure 5.4. Table 5.10 shows an analysis of the number of iterations and computation times of each solver.

Figure 5.8 shows the comparison between the solvers applied to the ball and plate system described in Section 5.7.2. Figure 5.8a shows the computation times and Figure 5.8b the number of iterations of each solver. The state and input trajectories are very similar to the ones shown in Figure 5.5. Table 5.11 shows an analysis of the number of iterations and computation times of each solver.

Figure 5.9 shows the comparison between the solvers applied to the oscillating masses system described in Section 5.7.3. Figure 5.9a shows the computation times and Figure 5.9b the number of iterations of each solver. The state and input trajectories are very similar to the ones shown in Figure 5.6. Table 5.12 shows an analysis of the number of iterations and computation times of each solver.

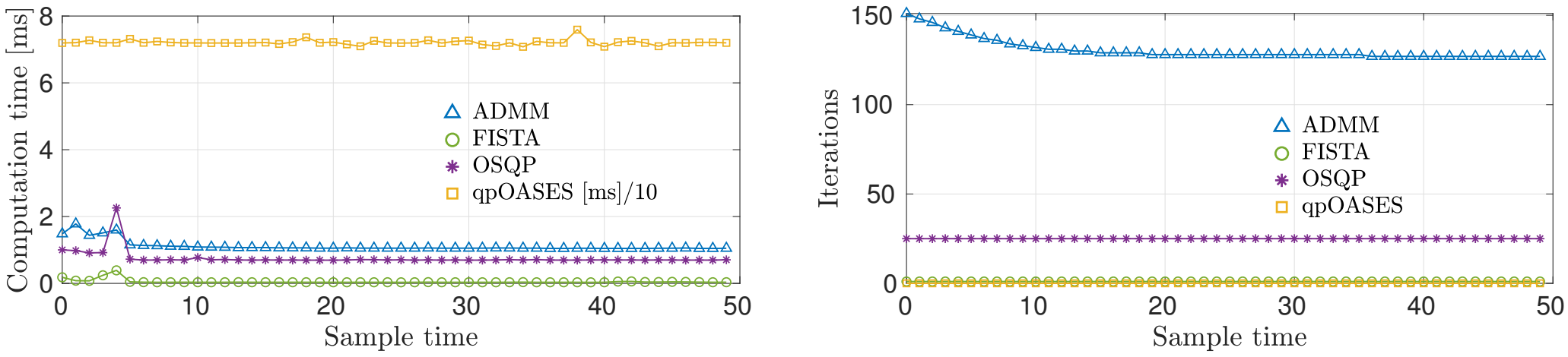

(a) Computation times.

(b) Number of iterations.

Figure 5.10: Closed-loop simulation of the chemical plant: Comparison between solvers applied to the standard MPC formulation with terminal equality constraint. The computation times of qpOASES are shown divided by 10 to be able to appreciate the computation times of the other solvers.

The results indicate that the FISTA-based solver (Algorithm 16) has a good performance. However, as shown in Figure 5.8 and Table 5.11, the number of iterations can increase significantly if there are active constraints in the optimal solution. The comparison between the number of iterations of the FISTA-based and ADMM-based solvers is somewhat unfair due to the difference between their exit conditions. For instance, in Figure 5.9, the suboptimal solutions obtained with both solvers were at a very similar distance to the optimal solution of the problems at each sample time. In Figure 5.8, however, even though the exit tolerances of both algorithms where set to $10^{-4}$, the suboptimal solutions obtained with the FISTA-based solver where up to several orders of magnitude closer to the optimal solutions at each sample time than the ones obtained with the ADMM-based solvers. This may explain the significant increase in the number of iterations of the FISTA-based solvers shown in Figure 5.8 when compared to the mild increase shown in Figure 5.9, even though in both cases there where active constraints in the optimal solutions during the first few sample times.

### 5.8.3 Standard MPC subject to terminal equality constraint

This section compares the ADMM and FISTA-bases solvers for standard MPC formulation (5.16) against other QP solvers from the literature. In particular, we compare Algorithms 16 and 17 with OSQP and qpOASES to solve the QP problem derived from the standard MPC formulation with terminal equality constraint (5.16). The MPC ingredients, solver parameters and references are the same as the ones used in Section 5.8.2.

Figure 5.10 shows the comparison between the solvers applied to the chemical plant described in Section 5.7.1. Figure 5.10a shows the computation times and Figure 5.10b the number of iterations of each solver. The state and input trajectories are very similar to the ones shown in Figure 5.4. Table 5.13 shows an analysis of the number of iterations and computation times of each solver.

| | Solver | ADMM | FISTA | OSQP | qpOASES |
|---|---|---|---|---|---|
| Iterations | Average | 130.64 | 1 | 25 | 0 |
| | Median | 128 | 1 | 25 | 0 |
| | Maximum | 151 | 1 | 25 | 0 |
| | Minimum | 127 | 1 | 25 | 0 |
| Comp. time | Average | 1.116 | 0.04 | 0.764 | 72.09 |
| | Median | 1.067 | 0.033 | 0.704 | 71.99 |
| | Maximum | 1.841 | 0.398 | 2.081 | 75.93 |
| | Minimum | 1.054 | 0.031 | 0.693 | 70.79 |

Table 5.13: Comparison between solvers applied to standard MPC with terminal equality constraint for the chemical plant: number of iterations and computation times [ms].

The active-set method of the qpOASES solver exits after 0 iterations at every sample time due to the optimal solution of the QP problems never having active constraints. Even so, the initial computations still result in a computation time that is much larger than the ones obtained with the solvers based on first order methods. The OSQP solver always converges in 25 iterations (although if smaller exit tolerances are used the number of iterations grows significantly). Even so, its computation times are similar to the ones obtained by the ADMM solver. The FISTA solver outperforms all others in this case because, once again, it converges in a single iteration if there are no active constraints in the optimal solution.

### 5.8.4 Standard MPC without terminal constraint

This section compares the ADMM and FISTA-bases solvers for standard MPC formulation (5.17) against other QP solvers from the literature. In particular, we compare Algorithms 16 and 17 with OSQP and qpOASES to solve the QP problem derived from the standard MPC formulation without terminal constraint (5.17). The MPC ingredients, solver parameters and references are the same as the ones used in Section 5.8.2.

Figure 5.11 shows the comparison between the solvers applied to the ball and plate system described in Section 5.7.2. Figure 5.11a shows the computation times and Figure 5.11b the number of iterations of each solver. The state and input trajectories are very similar to the ones shown in Figure 5.5. Table 5.14 shows an analysis of the number of iterations and computation times of each solver.

In this case the qpOASES solver performs a small number of iterations during the first few sample times, since there are active constraints in the optimal solutions of the resulting QP problems. Once again, the iterations drop to 0 as soon as this is no longer the case. The OSQP solver shows a similar trend to FISTA, although with fewer iterations. Nevertheless, the computation times of FISTA are very similar to the ones obtained with OSQP.

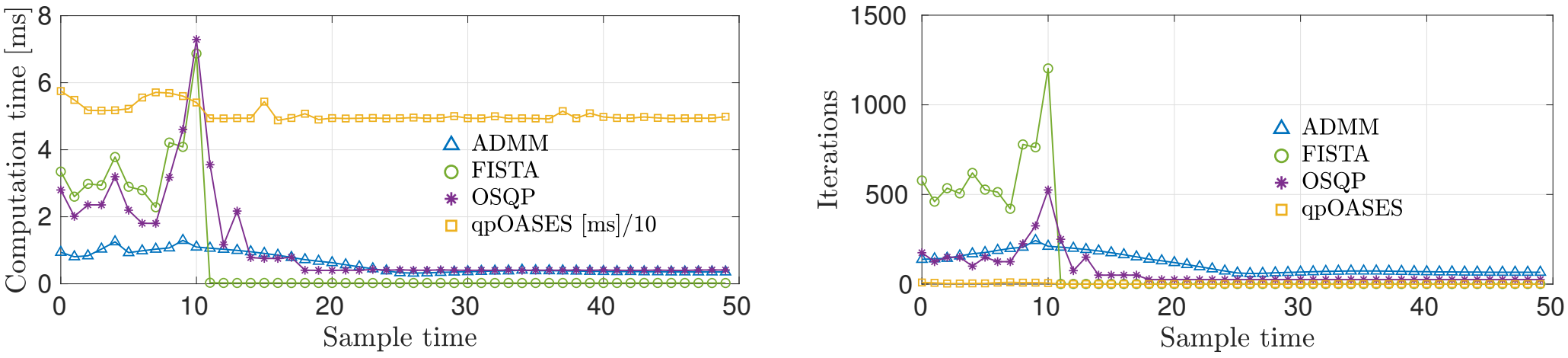

(a) Computation times. (b) Number of iterations.

Figure 5.11: Closed-loop simulation of the ball and plate: Comparison between solvers applied to the standard MPC formulation without terminal constraint. The computation times of qpOASES are shown divided by 10 to be able to appreciate the computation times of the other solvers.

| | Solver | ADMM | FISTA | OSQP | qpOASES |
|---|---|---|---|---|---|
| Iterations | Average | 113.78 | 138.86 | 73 | 1.4 |
| | Median | 73 | 1 | 25 | 0 |
| | Maximum | 243 | 1203 | 525 | 10 |
| | Minimum | 59 | 1 | 25 | 0 |
| Comp. time | Average | 0.618 | 0.795 | 1.131 | 50.75 |
| | Median | 0.399 | 0.024 | 0.41 | 49.45 |
| | Maximum | 1.295 | 6.869 | 7.288 | 57.51 |
| | Minimum | 0.318 | 0.023 | 0.399 | 48.76 |

Table 5.14: Comparison between solvers applied to standard MPC without terminal constraint for the ball and plate: number of iterations and computation times [ms].

### 5.8.5 MPC subject to terminal quadratic constraint

This section compares the ADMM-based solver for the MPC formulation subject to terminal quadratic constraints (5.20) to solving the QP problem that arises from considering the same MPC formulation but substituting the terminal ellipsoidal set with the maximal admissible invariant set of the system, which for controllable linear systems is a polyhedron. In particular, we apply both approaches to the oscillating masses system described in Section 5.7.3.

We compare the solver against two alternatives: OSQP [35] and FalcOpt [80, 135]. OSQP will be implemented using the maximal admissible invariant set of the system, which is a polyhedron, in place of the ellipsoidal terminal set, thus resulting in a QP problem. The comparison with this solver will allow us to evaluate the computational advantages that can be obtained by using the ellipsoidal set instead of a polyhedral one. FalcOpt is a solver for MPC subject to terminal ellipsoidal constraint that generates the code of the solver following a similar philosophy to SPCIES. Therefore, it serves as a comparison of our proposed solver with other state-of-the-art solvers in the literature. However, FalcOpt considers the case of non-linear MPC, does not consider state constraints and the matrix $P$ of its terminal ellipsoidal constraint must be equal to the terminal cost function $T$. That is, its terminal ellipsoidal set is given by $\mathcal{E}(T, x_r, r)$.

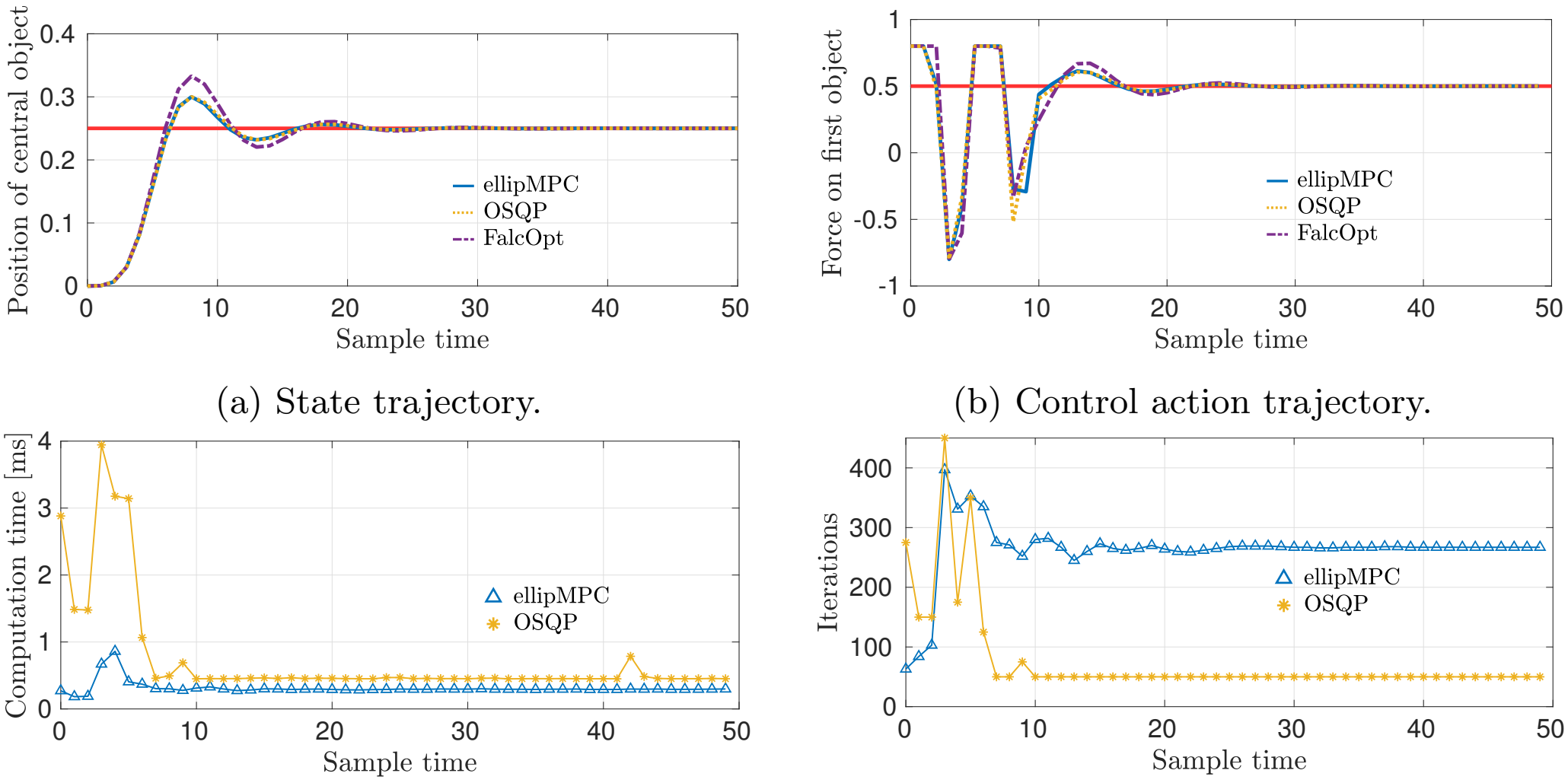

(a) State trajectory.

(b) Control action trajectory.

(c) Computation times.

(d) Number of iterations.

Figure 5.12: Closed-loop simulation of the oscillating masses: Comparison between Algorithm 18 for MPC with terminal quadratic constraints and OSQP applied to MPC using the (polyhedral) maximal admissible invariant set.

The prediction horizon $N$ and the cost function matrices $Q$ and $R$ are shown in Table 5.4. The ingredients $P$, $c$ and $r$ defining the terminal constraint are computed using the procedure described in Section 5.5.2, where the LMI optimization problem is constructed using YALMIP [103] and solved using the SDPT3 solver [145] taking $c = x_r$ and $r = 1$. A feasible solution of problem (5.41) is found for $\lambda = 0.95$. The cost function matrix $T$ is taken as the solution of the Lyapunov function (5.50) for the matrix $K$ obtained from the procedure used to compute $P$. We take the penalty parameter of the ADMM algorithm as $\rho = 15$.

We compute the maximal admissible invariant set of the system using the MPT3 toolbox [148], resulting in a polyhedron $\{ x \in \mathbb{R}^n \,:\, A_t x \leq b_t \}$ described by a matrix $A_t$ with 274 rows, resulting in as many constraints in the QP problem. We solve this QP problem using the OSQP solver.

Vector $r$ of the terminal ellipsoidal set of the FalcOpt solver is obtained using a similar procedure to the one described in Section 5.5.2, but forcing $P = T$ and taking $r$ as a decision variable. We minimize $-r$ to obtain the largest admissible invariant ellipsoid $\mathcal{E}(T, x_r, r)$ of the system, obtaining $r = 3.6552$. We use the FalcOpt solver from commit `5ac104c` of its GitHub repository [80] with an exit tolerance of $10^{-3}$ and its "gradients" option set to "casadi", which uses CasADi [79] (version `v3.4.5`) to compute the gradients.

Figure 5.12 shows the comparison between the two solvers applied to the oscillating masses system described in Section 5.7.3. Figure 5.12c shows the computa-tion times and Figure 5.12d the number of iterations of each solver. Figure 5.12a shows the trajectory of the state $p_2$ and Figure 5.12b the trajectory of the control input $\mathcal{F}_l$. Table 5.15 shows an analysis of the number of iterations and computation times of Algorithm 18 and OSQP. The results using FalcOpt are omitted due to them being several orders of magnitude bigger, either because we were not able to fine-tune the solver correctly and/or because its focus on non-linear MPC makes it less efficient at solving the linear MPC problem.

| | Solver | Algorithm 18 | OSQP |
|---|---|---|---|
| Iterations | Average | 262.52 | 77 |
| | Median | 267 | 50 |
| | Maximum | 397 | 450 |
| | Minimum | 63 | 50 |
| Comp. time | Average | 0.312 | 0.747 |
| | Median | 0.295 | 0.455 |
| | Maximum | 0.861 | 3.944 |
| | Minimum | 0.183 | 0.449 |

Table 5.15: Comparison between solvers applied to MPC with terminal quadratic constraint for the oscillating masses: number of iterations and computation times [ms].

The oscillating masses system is described by 6 states and 2 inputs. Even so, its maximal admissible invariant set results in the inclusion of 274 constraints to the QP problem solved by OSQP. This is a reasonable amount of constraints, but the number can easily become unmanageable, even for average-sized systems, resulting in an inapplicable optimization problem due to the large computation times. In such a case, the computational results shown in Figure 5.12 and Table 5.15 indicate that our proposed solver could result in an applicable controller due to its smaller number of constraints. Furthermore, the computation of an ellipsoidal admissible invariant set is computationally attainable even for large-sized systems (for instance, using the approach described in Section 5.5.2), whereas the computation of a polyhedral admissible invariant set can become very challenging.

The use of the maximal admissible invariant set has benefits in terms of the domain of attraction of the MPC controller [126, §4.2], which would be (possibly) reduced with the use of the quadratic constraint. However, the domain of attraction can be enlarged by increasing the prediction horizon, which, due to the sparse nature of the proposed algorithm, may be an acceptable approach.

### 5.8.6 MPC for tracking

This section compares the proposed EADMM-based solver for the MPCT formulation (5.42) (Algorithm 19) against other alternatives. In particular, we compare it against three other alternatives applied to the ball and plate system described in Section 5.7.2.

The first two alternatives are the OSQP and qpOASES solvers applied to the QP problem that results from (5.42). The reason for using the EADMM algorithm, instead of ADMM, is to be able to recover the inner matrix structures that enable the use of Algorithm 11. However, the QP problem derived from (5.42) is still a sparse QP problem, and thus may be also efficiently solved using sparse QP solvers (such as OSQP) or active-set solvers (such as qpOASES), among other alternatives.

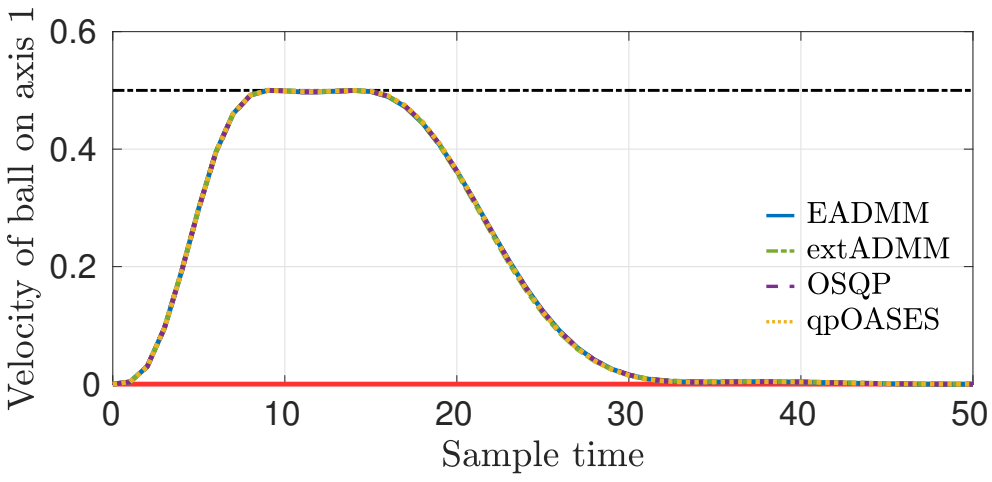

(a) State trajectory.

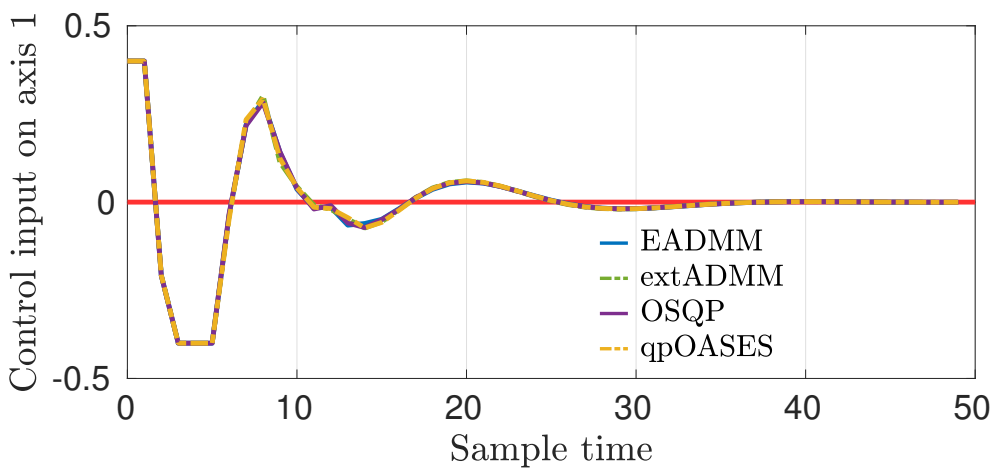

(b) Control action trajectory.

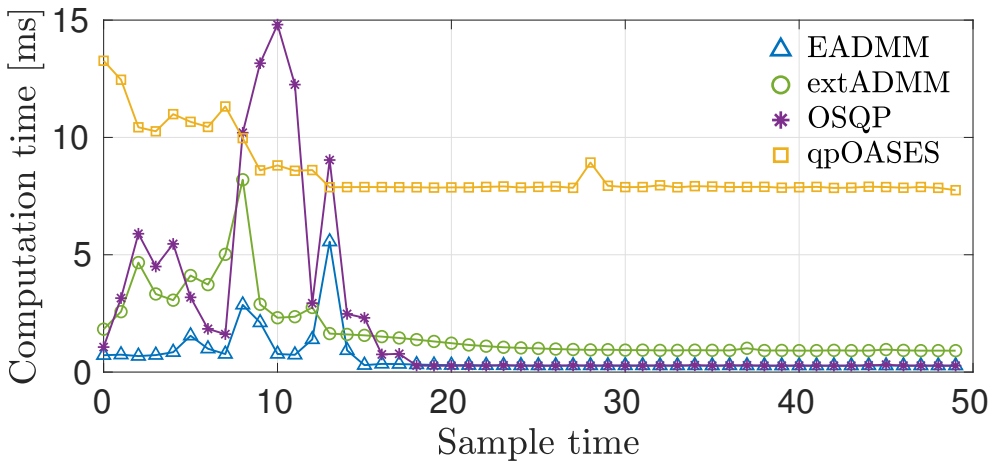

(c) Computation times.

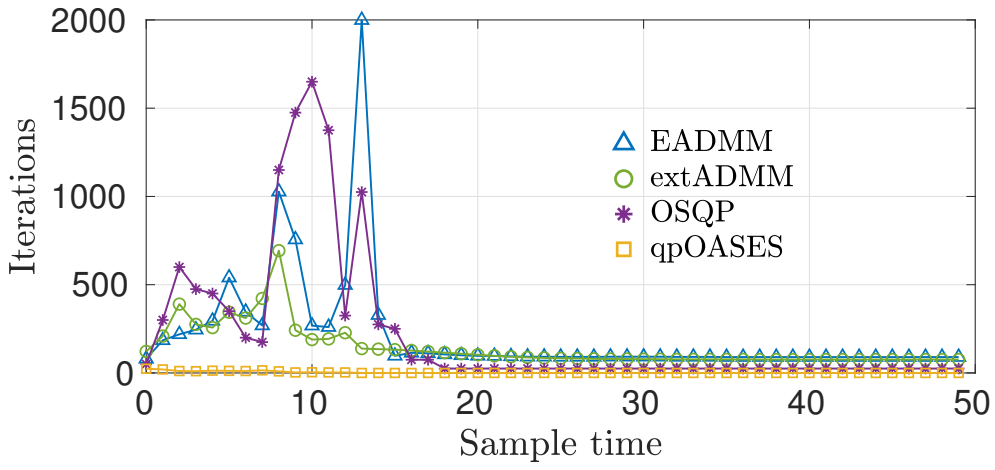

(d) Number of iterations.

Figure 5.13: Closed-loop simulation of the ball and plate: Comparison between Algorithm 19 and other alternatives. The black dashed/dotted line in figure (a) represents the upper bound on the velocity of the ball.

Motivated by the fact that the QP problem derived from (5.42) is sparse, we also compare the proposed Algorithm 19 against a third alternative. Note that the Hessian of (5.42) is not block diagonal due to the inclusion of the artificial reference $(x_s, u_s)$. However, we can recover a block diagonal Hessian by extending the state space model. That is, by defining $\hat{x}_i = (x_i, x_s)$ and $\hat{u}_i = (u_i, u_s)$, $i \in \mathbb{Z}_0^{N-1}$, problem (5.42) can be rewritten so that the resulting QP has a block diagonal Hessian at the expense of a few drawbacks: *(i)* the dimension of the matrices defining the system dynamics are increased in size, which increases the memory requirements; *(ii)* this significantly increases the number of decision variables, which will affect the computational complexity per iteration of the algorithm; *(iii)* additional equality constraints need to be added to ensure that $x_s$ and $u_s$ are the same in all the prediction steps; and *(iv)* the structure of (5.1), while still banded, no longer satisfies Assumtion 5.1.

To solve this problem we use a sparse ADMM algorithm similar to the one presented in Algorithm 17. However, in this case, we cannot employ Algorithm 11 to solve the system of equations. In fact, we find that we obtain better results in this case if the sparse matrices are stored using the *compressed sparse column/row* formats, and the system of equations solved at each iteration of the ADMM algorithm is solved using a sparse LDL representation as done in the QDLDL [35]. This ADMM-based solver for MPCT has also been included in the SPCIES toolbox [4], although we will not explain it in detail since it is out of the scope of this dissertation. In the following, we will refer to this algorithm as *extADMM* to denote that it makes use of the above extended state space representation.

| | Solver | EADMM | extADMM | OSQP | qpOASES |
|---|---|---|---|---|---|
| Iterations | Average | 211.54 | 142.98 | 221.5 | 2.76 |
| | Median | 92 | 84 | 25 | 0 |
| | Maximum | 1999 | 693 | 1650 | 24 |
| | Minimum | 80 | 76 | 25 | 0 |
| Comp. time | Average | 0.627 | 1.726 | 2.085 | 8.541 |
| | Median | 0.281 | 1.022 | 0.278 | 7.891 |
| | Maximum | 5.557 | 8.196 | 14.808 | 13.267 |
| | Minimum | 0.273 | 0.912 | 0.268 | 7.751 |

Table 5.16: Comparison between solvers applied to the MPCT formulation to control the ball and plate: number of iterations and computation times [ms].

The ingredients $Q$ and $R$ are the ones show in Table 5.4. However, we take $N = 15$ in this test. Matrices $T$ and $S$ are selected as:

$$T = \text{diag}(600, 50, 50, 50, 600, 50, 50, 50), \quad S = \text{diag}(0.3, 0.3).$$

We take $\rho = 2000$ as the penalty parameter for the constraints listed in Remark 5.35 and $\rho$=10 for the rest. The penalty parameter of the *extADMM* solver is taken as $\rho = 25$. The exit tolerances of EADMM, *extADMM* and OSQP solvers are all taken as $10^{-4}$. We use the default exit tolerance of the qpOASES solver.

Figure 5.13 shows the comparison between the solvers applied to the ball and plate system described in Section 5.7.2. Figure 5.13a shows the trajectory of the state $\dot{p}_1$ and Figure 5.13b the trajectory of the control input $\ddot{\theta}_1$. Figure 5.13c shows the computation times and Figure 5.13d the number of iterations of each solver. Table 5.16 shows an analysis of the number of iterations and computation times of each solver.

In this test, the number of iterations of the solvers based on FOMs were significantly larger that in the tests shown in the previous sections. This serves to highlight the advantage that can be obtained with the proposed solvers against other sparse QP solvers from the literature thanks to the low computational complexity that we obtain by exploiting the specific structure of the MPC formulations. The *extADMM* solver converges, in general, in fewer iterations than the EADMM solver. However, its computation times are larger due to its higher iteration complexity, mostly due to the large matrices as a result of the extension of the state space. Once again, the computation times of the qpOASES solver are larger than the ones obtained with the sparse FOM-based solvers.

### 5.8.7 Restart methods applied to the FISTA-based solvers

This section shows the application of the restart schemes presented in Part I of this dissertation to FISTA for solving the standard MPC formulations. That is, we show the results of using Algorithms 7, 8 and 10 in Algorithm 16 applied to the MPC formulations (5.17) (labeled *laxMPC*) and (5.16) (labeled *equMPC*). We refer the reader to Sections 3.2.1 to 3.2.3 for the details of these restart schemes. For Algorithm 10 we employ MFISTA (Algorithm 5) instead of FISTA, since it simplifies the implementation of the algorithm, as discussed in Remark 3.28.

Additionally, we also show the results of applying the *objective function value scheme* [45], whose restart condition $E_f$ is given by (3.6); the *gradient alignment scheme* [45], whose restart condition $E_g$ is given by (3.7); and the optimal restart scheme from [46, §5.2.2] using the restart condition $E_f^*$ given by (3.5). We refer the reader to Section 3.1 for the details of these restart schemes.

We apply the above restart schemes to the ball and plate system described in Section 5.7.2 and to the oscillating masses system described in Section 5.7.3, performing the same closed-loop simulations and using the same MPC ingredients that were used in Section 5.8.2. However, we are only interested in the sample times in which there are active constraints in the optimal solutions, since, we recall, Algorithm 16 converges in 1 iteration otherwise. Therefore, we only show the results of the first 11 sample times of the ball and plate system (c.f., Figure 5.8) and first 8 or 9 sample times of the oscillating masses system. For clarity, we include the first sample time without active constraints in the optimal solution. Obviously, there is no point in showing the results for the chemical plant system described in Section 5.7.1, since the non-restarted FISTA algorithm converged in 1 iteration in all the sample times of the closed loop test shown in Section 5.8.2.

We use $\|\mathcal{G}(y_{k-1})\|_{W^{-1}} \leq 10^{-4}$ as the exit tolerance of the FISTA and MFISTA algorithms, instead of taking the exit condition of Algorithm 16. We stop the restart schemes as soon as an iterate $y_{k-1}$ satisfying this condition is found.

Figure 5.14 shows the number of iterations of the restart schemes during the first sample times of the ball and plate system closed-loop test. Figure 5.14a shows the results for the *laxMPC* formulation (5.17) and Figure 5.14b for the *equMPC* formulation (5.16). Table 5.17 show an analysis of the number of iterations of each restart scheme during the sample times shown in Figure 5.14 (excluding the last ones, in which all schemes converged in 1 iteration). Figures 5.15 and 5.16 show the evolution of the iterates of FISTA for each restart scheme during the first sample time of the closed-loop simulations depicted in Figures 5.14a and 5.14b, respectively.

Figure 5.17 shows the number of iterations of the restart schemes during the first sample times of the oscillating masses system closed-loop test. Figure 5.17a shows the results for the *laxMPC* formulation (5.17) and Figure 5.17b for the *equMPC* formulation (5.16). Table 5.18 show an analysis of the number of iterations of each restart scheme during the sample times shown in Figure 5.17 (excluding the last ones, in which all schemes converged in 1 iteration). Fig-

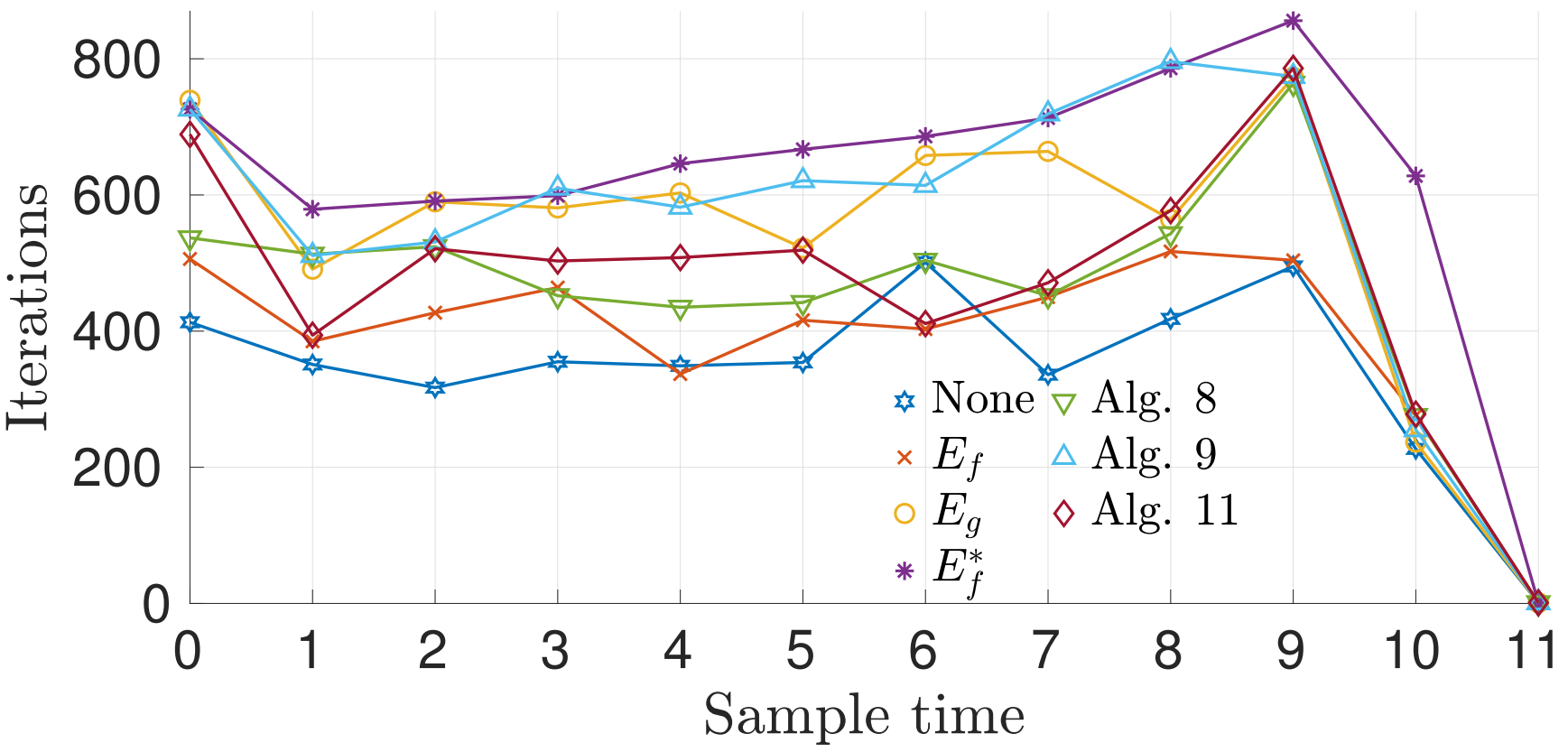

(a) MPC without terminal constraint.

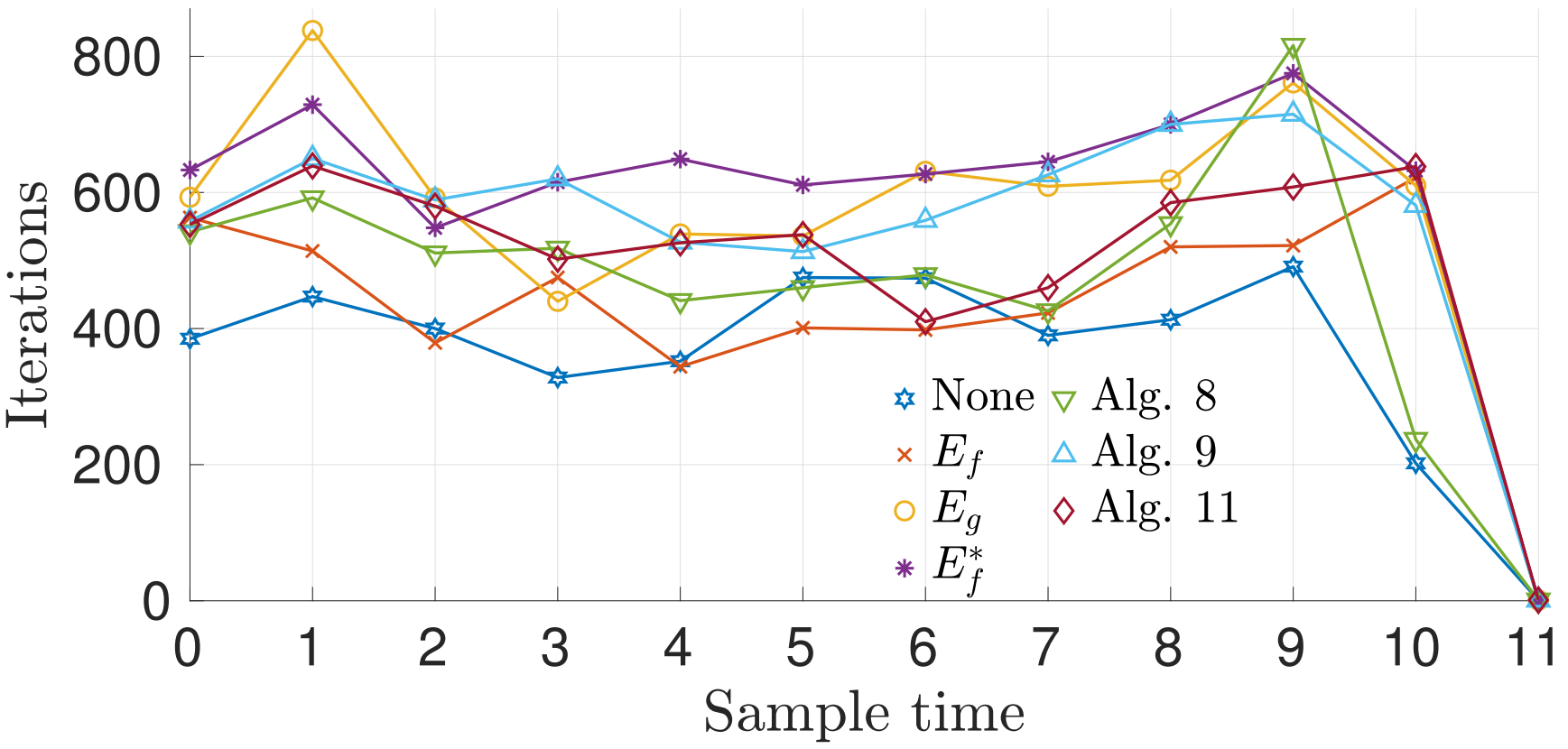

(b) MPC with terminal equality constraint.

Figure 5.14: Number of iterations of FISTA applied to the standard MPC formulations using different restart schemes. Closed loop results on the ball and plate.

| | Restart scheme | None | Alg. 7 | Alg. 8 | Alg. 10 | $E_f$ | $E_g$ | $E_f^*$ |
|---|---|---|---|---|---|---|---|---|
| laxMPC | Avg. Iter. | 374.2 | 494.7 | 612.6 | 514.3 | 425.9 | 583.8 | 679.7 |
| | Med. Iter. | 354 | 504 | 614 | 508 | 427 | 590 | 667 |
| | Max. Iter. | 501 | 764 | 796 | 786 | 517 | 773 | 856 |
| | Min. Iter. | 227 | 276 | 255 | 278 | 276 | 237 | 579 |
| equMPC | Restart scheme | None | Alg. 7 | Alg. 8 | Alg. 10 | $E_f$ | $E_g$ | $E_f^*$ |
| | Avg. Iter. | 396.1 | 506.9 | 603.5 | 549 | 469.2 | 615.2 | 651.3 |
| | Med. Iter. | 400 | 511 | 589 | 553 | 475 | 609 | 633 |
| | Max. Iter. | 491 | 816 | 715 | 639 | 622 | 838 | 775 |
| | Min. Iter. | 202 | 237 | 513 | 410 | 344 | 440 | 548 |

Table 5.17: Analysis of the number of iterations of FISTA with different restart schemes during the sample times shown in Figure 5.14. Application to the ball and plate.

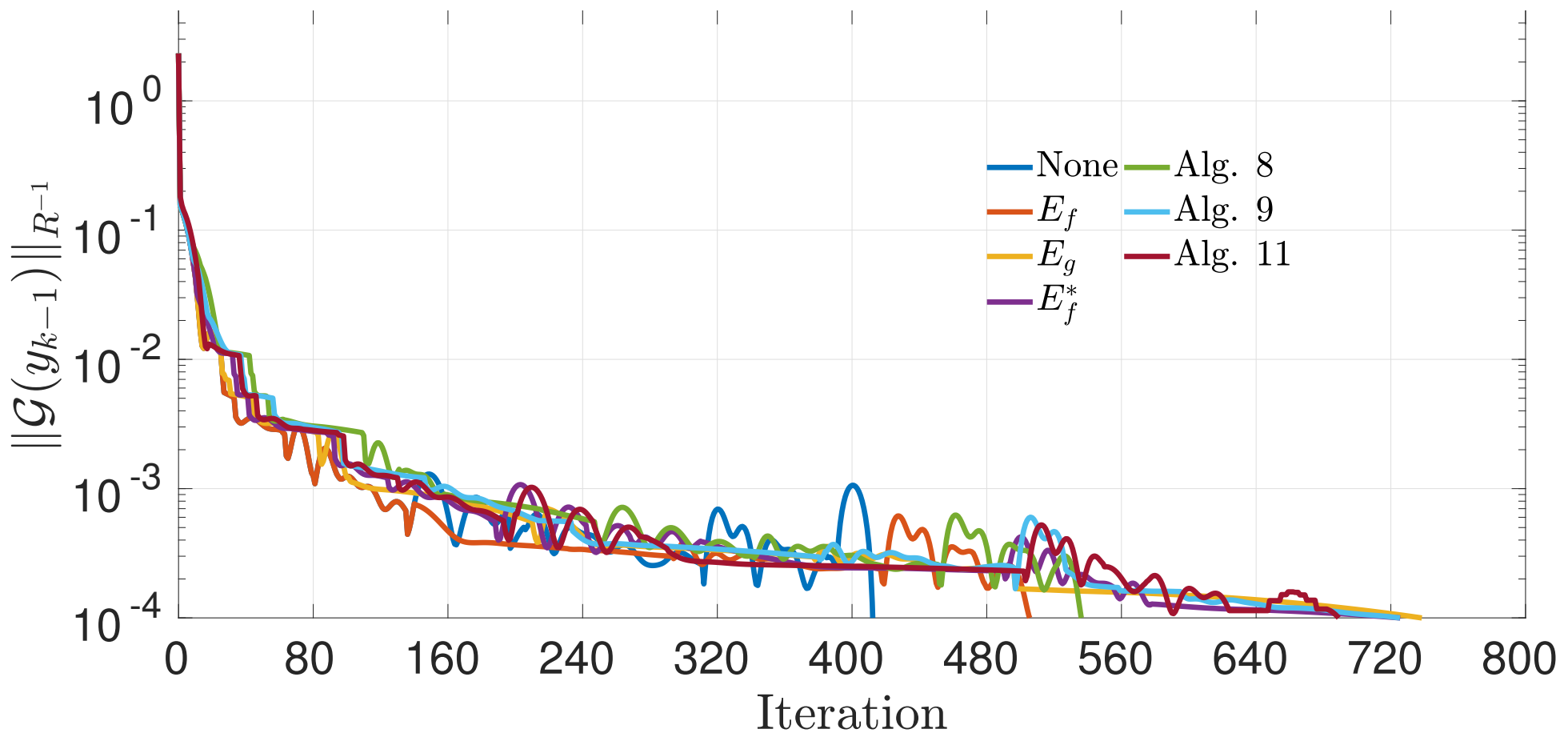

(a) Dual norm of the composite gradient mapping.

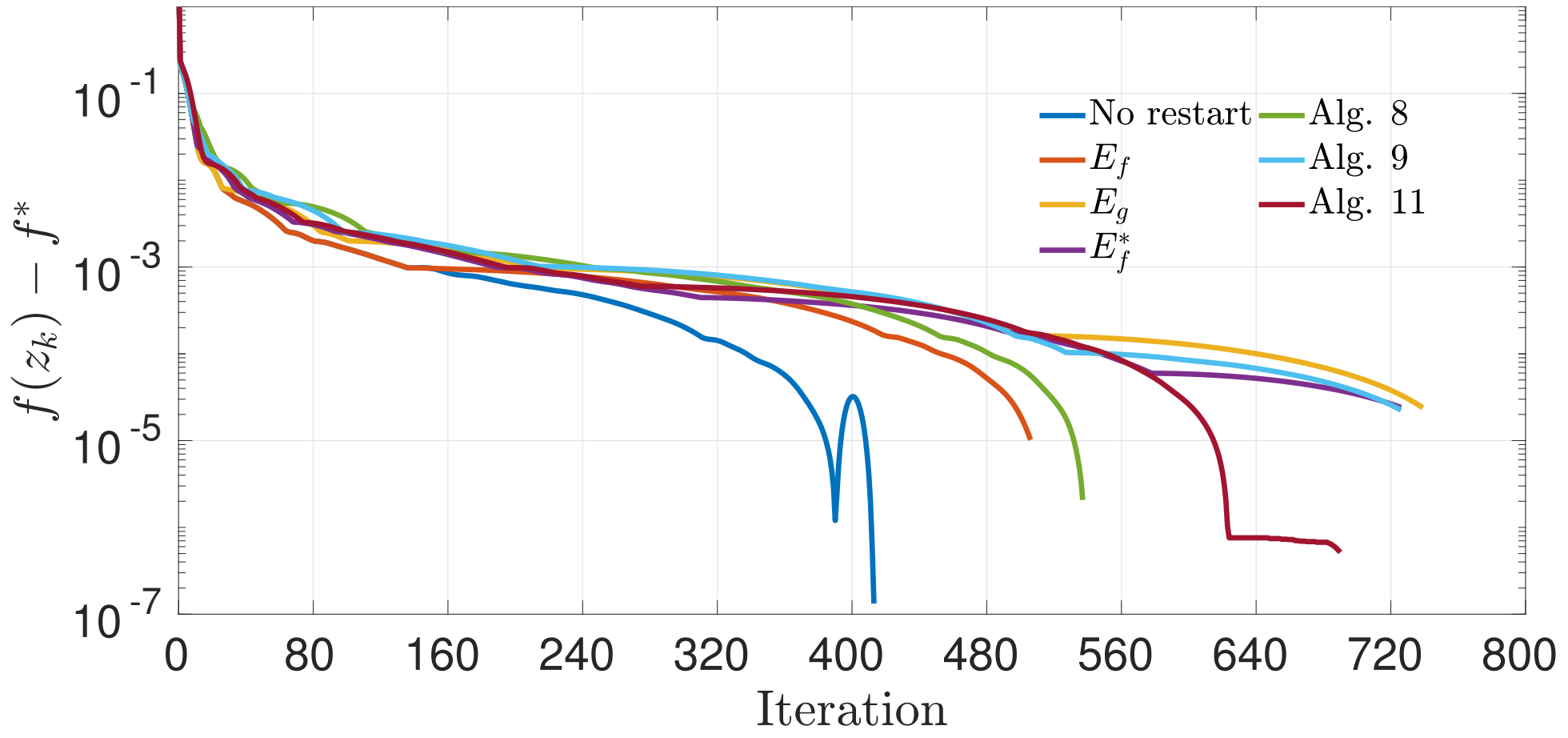

(b) Distance to the optimum in terms of the objective function value.

Figure 5.15: Evolution of the iterates of FISTA for different restart schemes during the first iteration of Figure 5.14a, i.e., applied to *laxMPC* for the ball and plate.

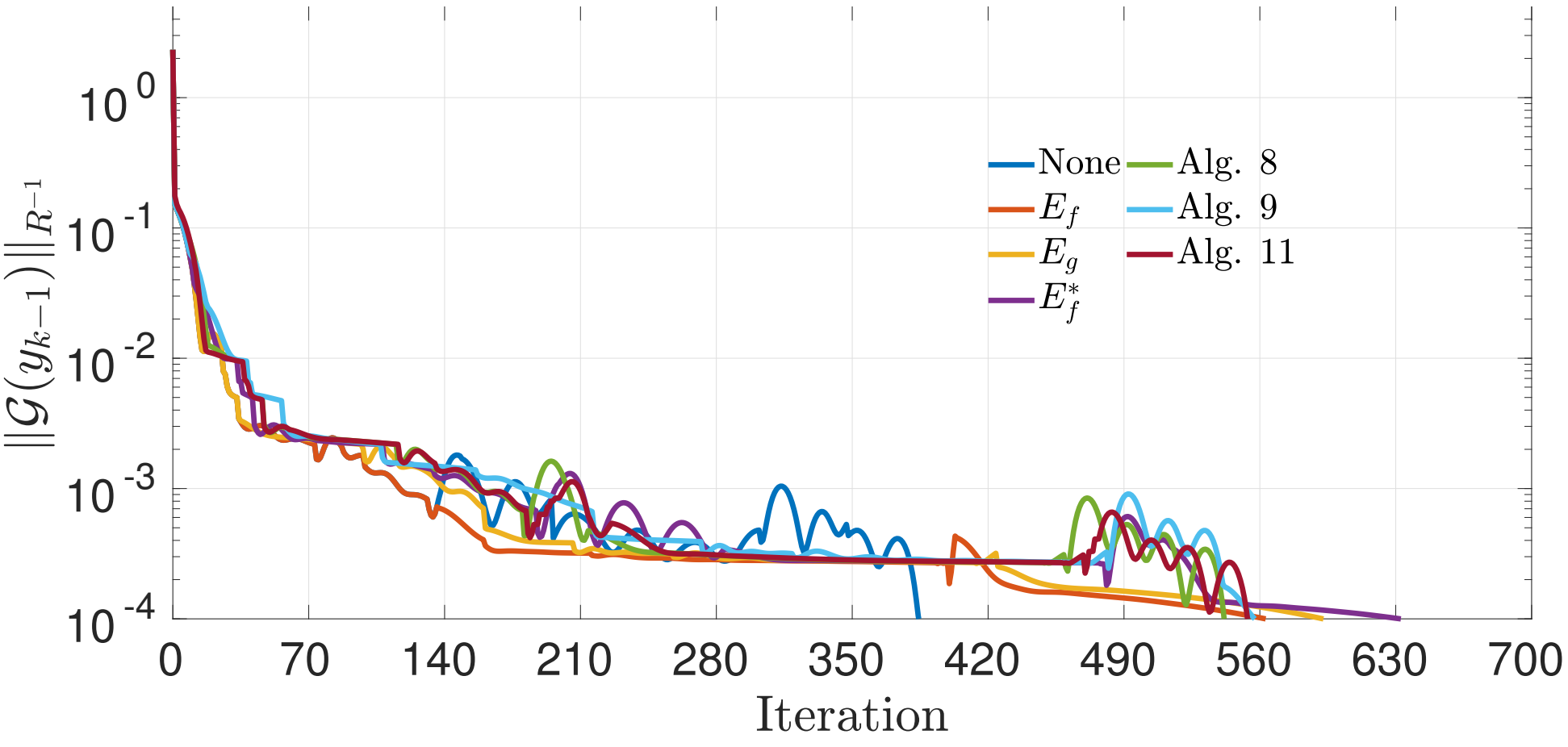

(a) Dual norm of the composite gradient mapping.

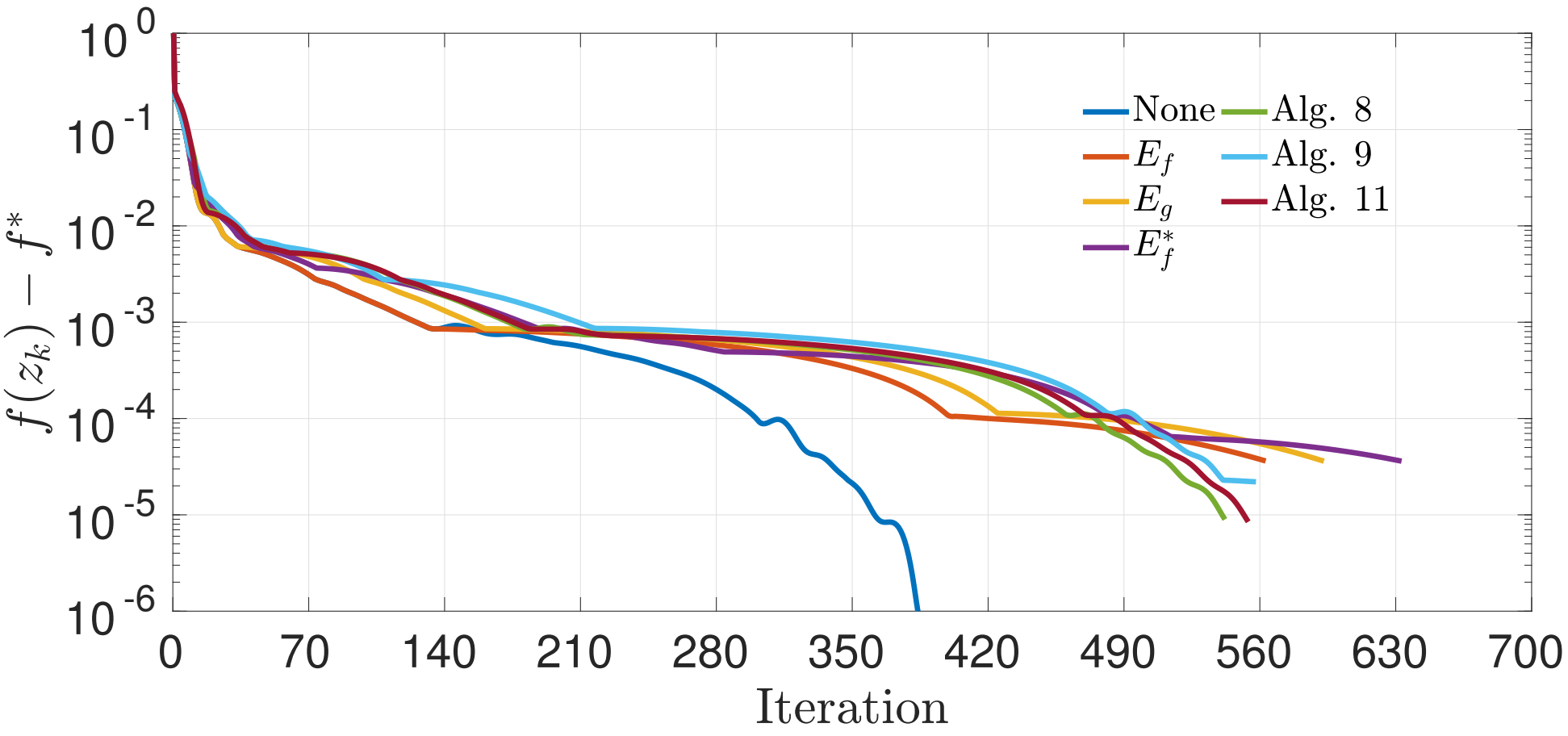

(b) Distance to the optimum in terms of the objective function value.

Figure 5.16: Evolution of the iterates of FISTA for different restart schemes during the first iteration of Figure 5.14b, i.e., applied to *equMPC* for the ball and plate.

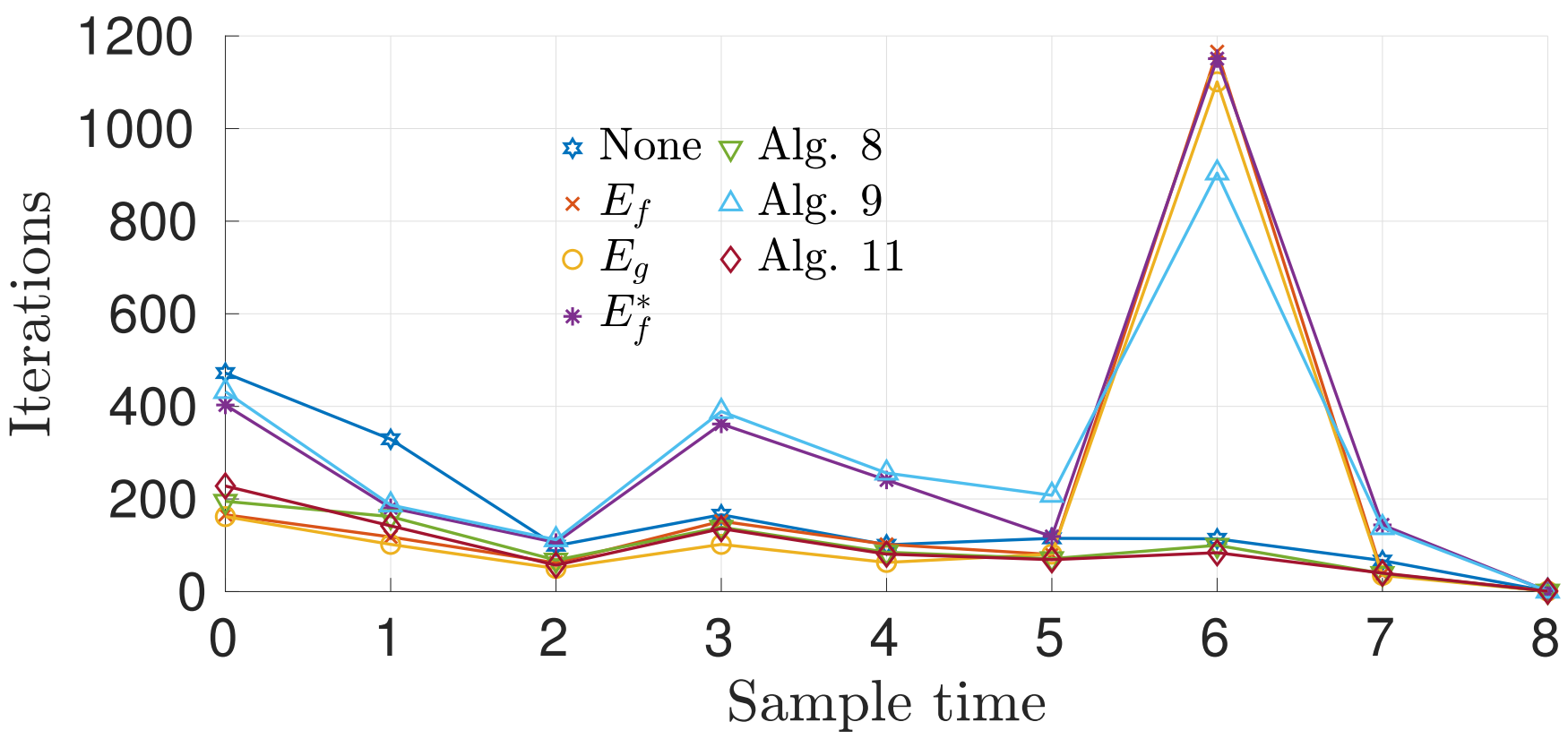

(a) MPC without terminal constraint.

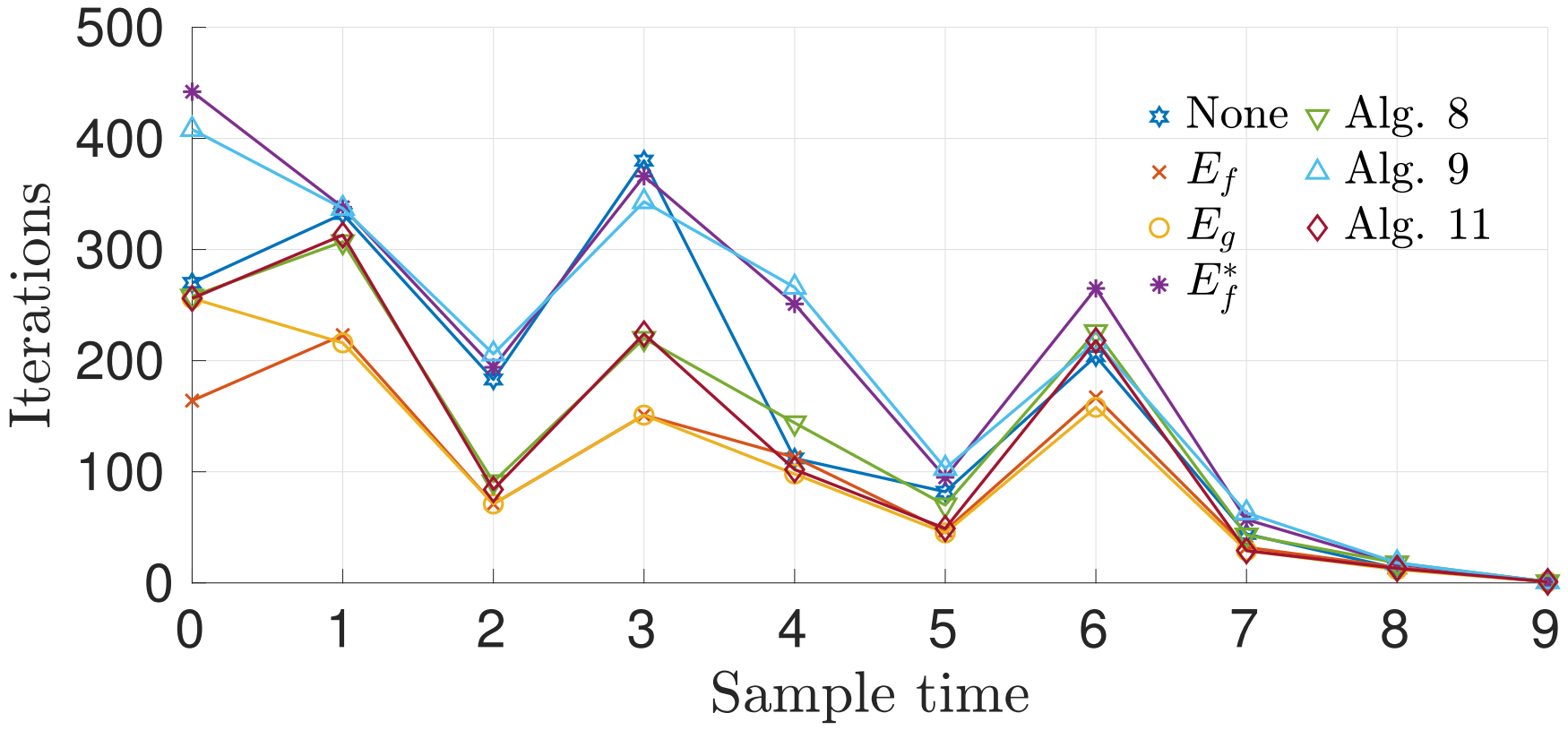

(b) MPC with terminal equality constraint.

Figure 5.17: Number of iterations of FISTA applied to the standard MPC formulations using different restart schemes. Closed loop results on the oscillating masses.

| | Restart scheme | None | Alg. 7 | Alg. 8 | Alg. 10 | $E_f$ | $E_g$ | $E_f^*$ |
|---|---|---|---|---|---|---|---|---|
| laxMPC | Avg. Iter. | 182.9 | 107.4 | 328.3 | 104.6 | 235.5 | 212 | 338.4 |
| | Med. Iter. | 114.5 | 92.5 | 232 | 82.5 | 110 | 91.5 | 211 |
| | Max. Iter. | 472 | 195 | 904 | 228 | 1166 | 1101 | 1151 |
| | Min. Iter. | 67 | 39 | 112 | 40 | 38 | 35 | 106 |
| | Restart scheme | None | Alg. 7 | Alg. 8 | Alg. 10 | $E_f$ | $E_g$ | $E_f^*$ |
| equMPC | Avg. Iter. | 396.1 | 506.9 | 603.5 | 549 | 469.2 | 615.2 | 651.3 |
| | Med. Iter. | 400 | 511 | 589 | 553 | 475 | 609 | 633 |
| | Max. Iter. | 491 | 816 | 715 | 639 | 622 | 838 | 775 |
| | Min. Iter. | 202 | 237 | 513 | 410 | 344 | 440 | 548 |

Table 5.18: Analysis of the number of iterations of FISTA with different restart schemes during the sample times shown in Figure 5.17. Application to the oscillating masses.

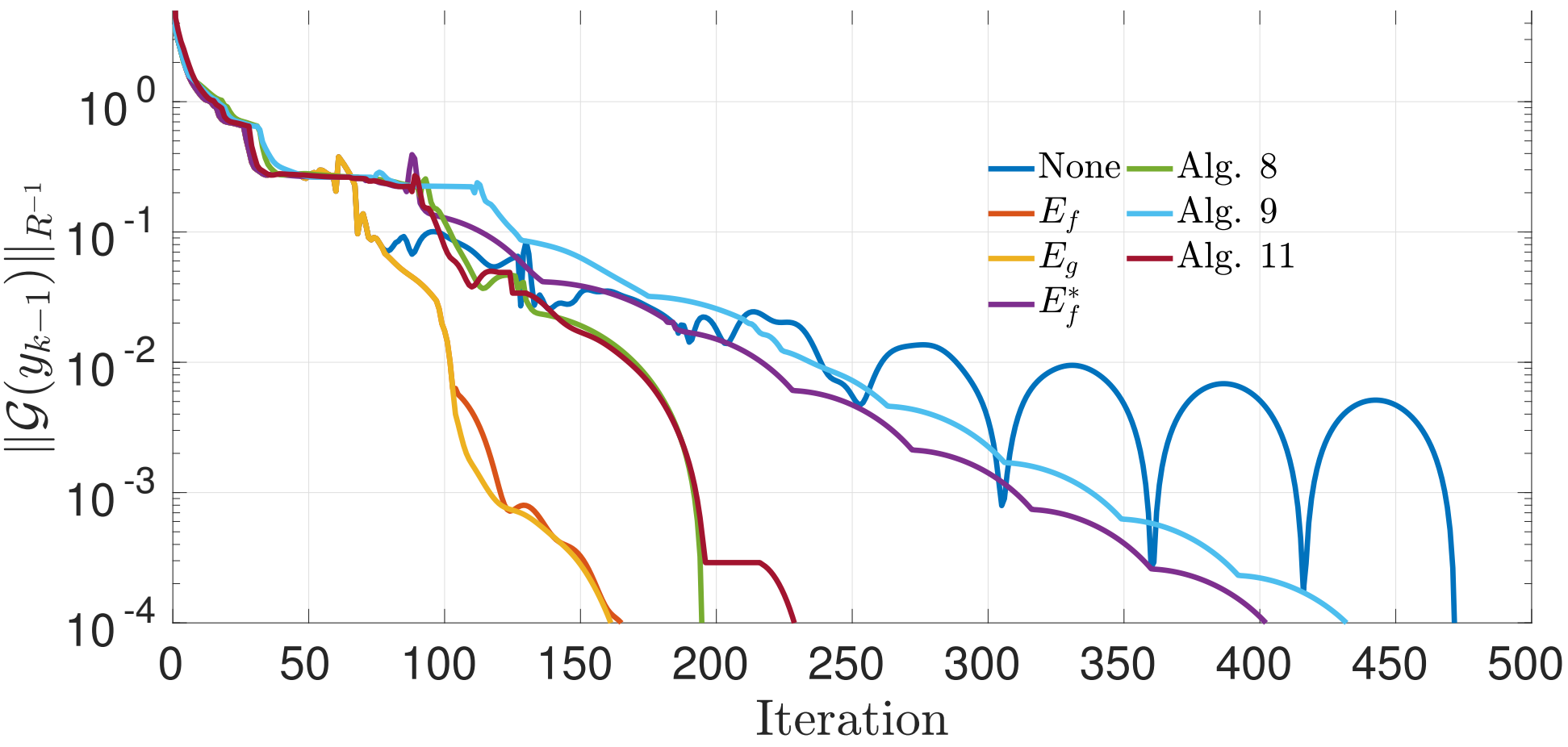

(a) Dual norm of the composite gradient mapping.

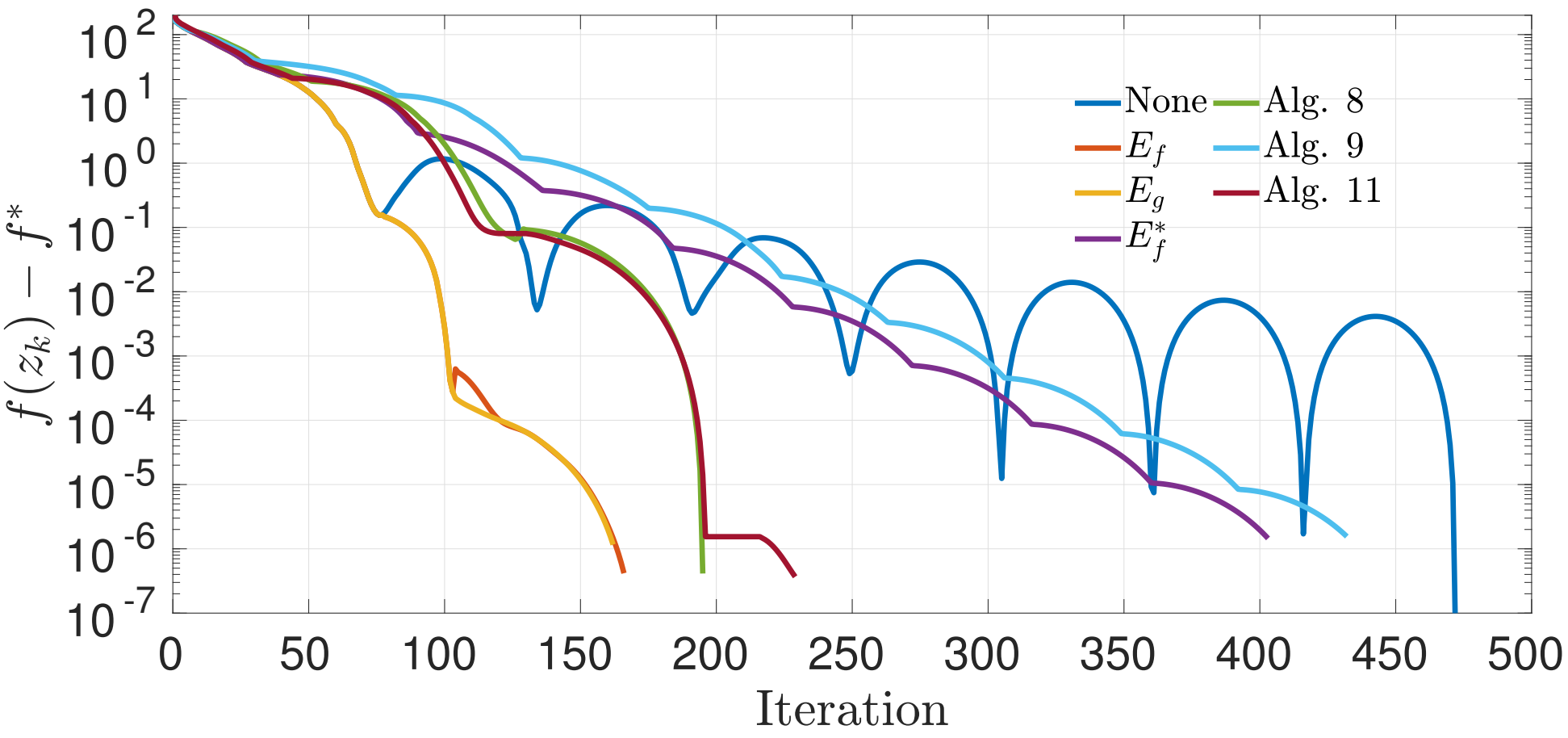

(b) Distance to the optimum in terms of the objective function value.

Figure 5.18: Evolution of the iterates of FISTA for different restart schemes during the first iteration of Figure 5.17a, i.e., applied to *laxMPC* for the oscillating masses.

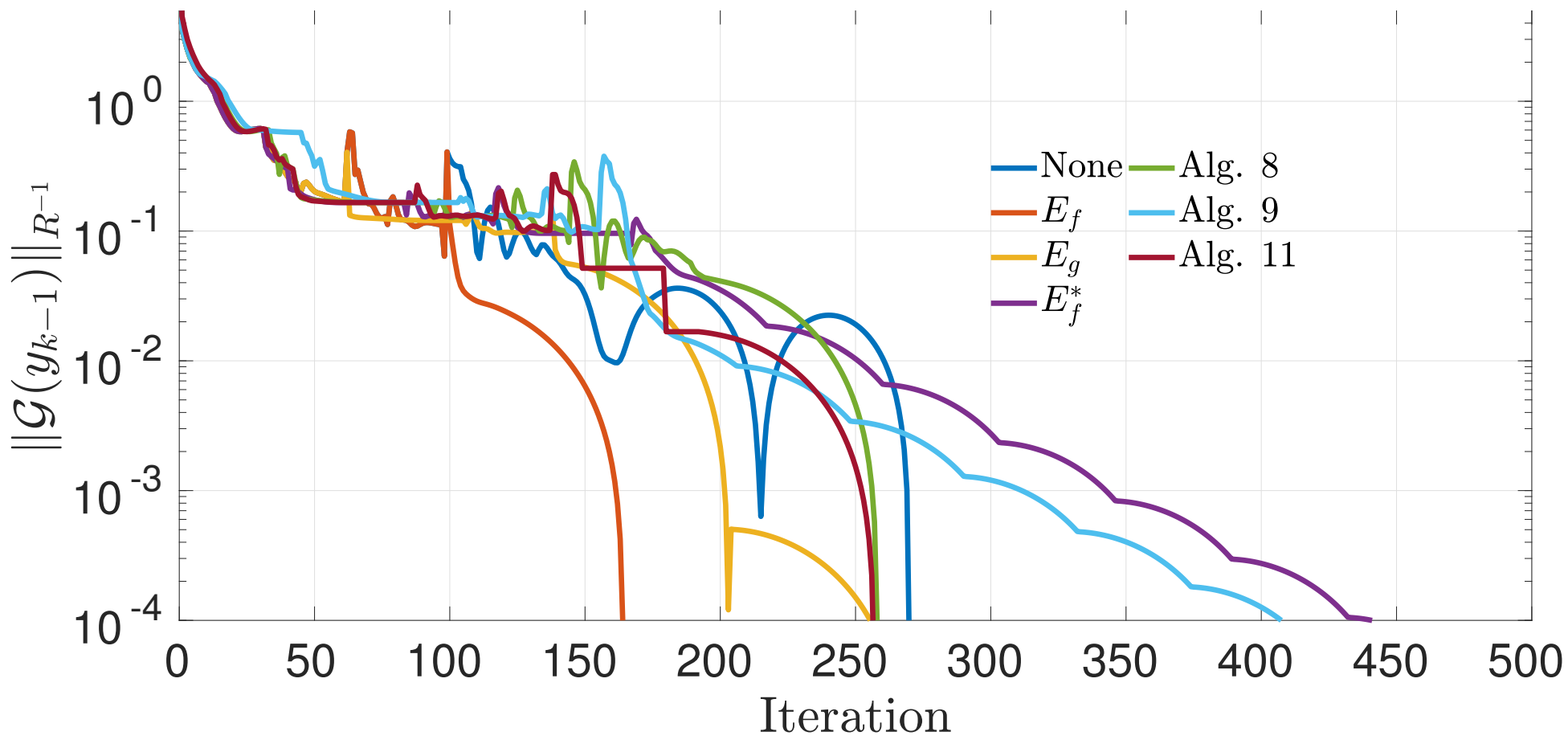

(a) Dual norm of the composite gradient mapping.

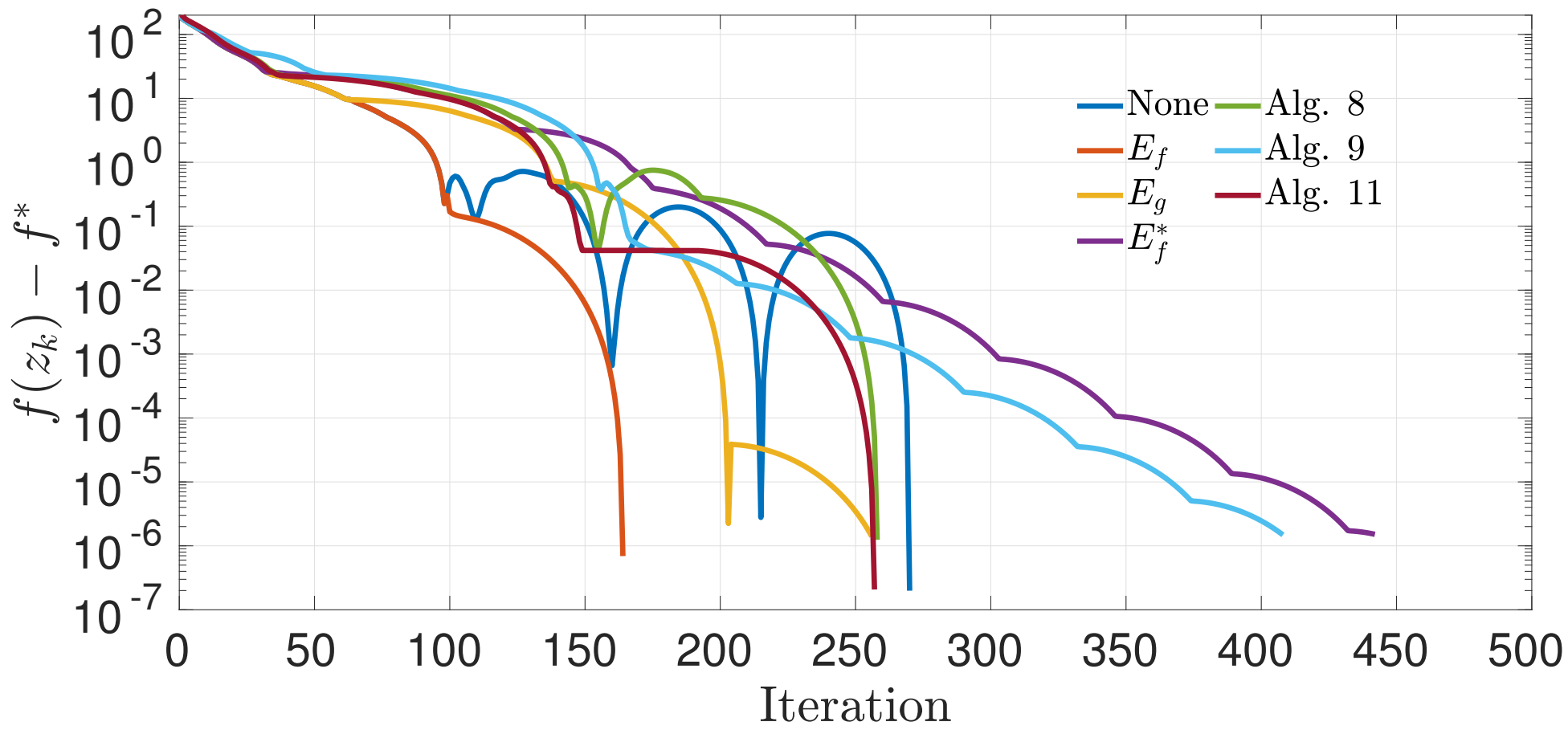

(b) Distance to the optimum in terms of the objective function value.

Figure 5.19: Evolution of the iterates of FISTA for different restart schemes during the first iteration of Figure 5.17b, i.e., applied to *equMPC* for the oscillating masses.

ures 5.18 and 5.19 show the evolution of the iterates of FISTA for each restart scheme during the first sample time of the closed-loop simulations depicted in Figures 5.17a and 5.17b, respectively.

Curiously, the results indicate that the use of restart schemes in the application of FISTA to MPC optimization problems is not always beneficial. In the case of the oscillating masses system, the results show that some of the restart schemes perform better that the non-restarted version, but others provide little to no benefit, with a few moments of significantly worse performance (see Figure 5.17b). The results shown in Table 5.18 show that Algorithms 7 and 10 performed better that the other alternatives, although the extra computational burden required to evaluate their restart conditions may make them less efficient in practice than the non-restarted variant. In the case of the ball and plate system, however, the non-restarted version of FISTA performed better that any of the other alternatives.

## 5.9 Conclusions and future lines of work

The numerical results shown in Section 5.8 indicate that the proposed solvers are suitable for their implementation in embedded systems. Their computation times are in the range of milliseconds in all the case studies shown in this dissertation and in the results shown in [12], which were obtained using a Raspberry Pi model 4B. The computation times obtained in the PLC are larger. However, these results were obtained with a less refined version of the solvers and using a very low resource PLC device.

The result show that the proposed solvers, which are applied to relevant MPC formulations, perform well when compared to other state-of-the-art solvers.

The application of the restart schemes to solving the MPC optimization problems have, for the most part, not provided benefits in terms of the number of iterations of the algorithm. Additional research is required to determine why this is the case. One possibility is that the matrix $W$ that we use leaves little room for improving the non-restarted variant. A future line of work is to study the application of restart schemes to MPC in further detail.

The solvers we present are for the most part only applicable in an academic environment. Their application to a real system would require the inclusion of additional ingredients such as a state estimator or a steady state target optimizer. In [8] we included some of these ingredients in the PLC, achieving an offset-free control of the non-linear system. Additionally, the solvers themselves can be improved further by studying aspects such as: numerical conditioning, the inclusion of soft constraints, certification of the number of iterations, etc. These aspects are all possible future lines of work.

The solvers we present here have been included in the SPCIES toolbox [4]. It takes the model of the system, the parameters of the desired MPC formulation and the options of the solver, and automatically generates the library-free code of the solver for the target embedded system. Currently, the toolbox only considers the generation of plain C-code and MEX files for Matlab. A future line of work is to maintain the development of this toolbox, including additional programming languages (Python, Julia, etc.) and embedded systems (PLCs, FPGAs, etc.), and adding the above mentioned improvements.

# Chapter 6

# Harmonic based model predictive control for tracking

This chapter considers system (4.1) subject to the coupled input-state constraints (4.4).

As discussed in Section 4.2, MPC formulations often make use of a terminal admissible invariant set $\mathcal{X}_t \subseteq \mathbb{R}^n$ (4.7e) to guarantee the asymptotic stability of the closed-loop system to the desired reference $(x_r, u_r) \in \mathbb{R}^n \times \mathbb{R}^m$. However, the use of this terminal set leads to two downsides when the reference to be tracked can change online. The first issue is that the terminal set depends on the value of the reference. If there are a known-before-hand, finite number of references, then a terminal set can be computed offline for each one of them and used online when necessary. Otherwise, the terminal set must be recomputed each time the reference changes, which is typically very computationally demanding. The second issue is that recursive feasibility of the MPC controller can be lost in the event of a reference change, i.e., there may not be a feasible solution of the MPC optimization problem for the current state if the reference has changed since the previous sample time. This second issue is related to the domain of attraction of the MPC controller (Definition 4.6), since feasibility is lost when the current state is no longer within the domain of attraction of the MPC controller due to the reference change. The terminal constraint is the main contributor to this issue when the prediction horizon is not large enough. To see this, note that the predicted state must be able to reach the terminal set within the prediction horizon window, and that systems are typically subject to input constraints.

These issues are particularly relevant when dealing with the online implementation of MPC in embedded systems, since their limited computational and memory resources make them unsuitable for large prediction horizons and the on-line computation of admissible invariant sets.

One possible solution to avoid having to recompute an admissible positive invariant set is to use a terminal equality constraint, i.e., to employ an MPC formulation (5.16). However, as discussed in Section 5.4, this formulation severely suffers from small domains of attraction if the prediction horizon is not large enough.

Another option is to employ the MPC for tracking formulation (5.42) discussed in Section 5.6. This formulation also uses a terminal equality constraint (5.42i), thus not requiring the computation of a new terminal set each time the reference changes. Moreover, its domain of attraction is significantly larger when compared to other standard MPC formulations, especially for small prediction horizons, although it can be further enlarged if soft constraints are considered [149]. Finally, its recursive feasibility and asymptotic stability are guaranteed even in the event of a reference change. However, as we will show further ahead, the closed-loop performance of this formulation may suffer in certain systems if the prediction horizon is too small.

This section presents a novel MPC formulation, originally presented in [13] and [14], which we call *harmonic MPC* and label by HMPC. The intent behind its development is to obtain an MPC formulation that is suitable for its implementation in embedded systems whilst solving the above issues. In particular, we are interested in an MPC formulation that *(i)* does not require the computation of a terminal set each time the reference changes, *(ii)* provides recursive feasibility and asymptotic stability even in the event of a reference change, *(iii)* has a large domain of attraction and good closed-loop performance even when using small prediction horizons, and *(iv)* that the previous three points do not come at the expense of an overly complicated optimization problem, i.e., one that would be more complex to solve online than to simply increase the prediction horizon of the other MPC formulations discussed above.

In short, the objective is to develop an MPC formulation with similar properties to the MPCT formulation (5.42) but with an increased domain of attraction and performance for small prediction horizons. In fact, the HMPC formulation can be viewed as an extension of the MPCT formulation (5.42).

## 6.1 The harmonic MPC formulation

The idea behind the HMPC formulation is to substitute the artificial reference $(x_s, u_s)$ of the MPCT formulation (5.42) by the *artificial harmonic reference* sequences $\{x_{hj}\}$, $\{u_{hj}\}$, $j \in \mathbb{Z}$, given by

$$x_{hj} = x_e + x_s \sin(w(j-N)) + x_c \cos(w(j-N)), \tag{6.1a}$$

$$u_{hj} = u_e + u_s \sin(w(j-N)) + u_c \cos(w(j-N)), \tag{6.1b}$$

where $N \in \mathbb{Z}_{>0}$ is the prediction horizon of the MPC controller and $w \in \mathbb{R}_{\geq 0}$ is the *base frequency*. The harmonic sequences $\{x_{hj}\}$ and $\{u_{hj}\}$ are parameterized by decision variables $x_e \in \mathbb{R}^n$, $x_s \in \mathbb{R}^n$, $x_c \in \mathbb{R}^n$, $u_e \in \mathbb{R}^m$, $u_s \in \mathbb{R}^m$, and $u_c \in \mathbb{R}^m$.

To simplify the text, we use the following notation,

$$\mathbf{x}_H \doteq (x_e, x_s, x_c) \in \mathbb{R}^n \times \mathbb{R}^n \times \mathbb{R}^n,$$
$$\mathbf{u}_H \doteq (u_e, u_s, u_c) \in \mathbb{R}^m \times \mathbb{R}^m \times \mathbb{R}^m,$$

$$y_e \doteq Ex_e + Fu_e, \quad y_s \doteq Ex_s + Fu_s, \quad y_c \doteq Ex_c + Fu_c, \tag{6.2}$$

where we recall that the matrices $E \in \mathbb{R}^{p\times n}$ and $F \in \mathbb{R}^{p\times m}$ define, along with $\underline{y}$ and $\overline{y}$, the constraint set $\mathcal{Y}$ (4.4).

For a given prediction horizon $N$ and base frequency $w$, the HMPC control law for a given state $x(t)$ and reference $(x_r, u_r)$ is derived from the following *second order cone programming* problem labeled by $\mathbb{H}(x; x_r, u_r)$,

$$\mathbb{H}(x(t); x_r, u_r) \doteq \min_{\mathbf{x},\mathbf{u},\mathbf{x}_H,\mathbf{u}_H} J_h(\mathbf{x}, \mathbf{u}, \mathbf{x}_H, \mathbf{u}_H; x_r, u_r) \tag{6.3a}$$
$$s.t.\ x_{j+1} = Ax_j + Bu_j,\ j \in \mathbb{Z}_0^{N-1} \tag{6.3b}$$
$$\underline{y} \leq Ex_j + Fu_j \leq \overline{y},\ j \in \mathbb{Z}_0^{N-1} \tag{6.3c}$$
$$x_0 = x(t) \tag{6.3d}$$
$$x_N = x_e + x_c \tag{6.3e}$$
$$x_e = Ax_e + Bu_e \tag{6.3f}$$
$$x_s \cos(w) - x_c \sin(w) = Ax_s + Bu_s \tag{6.3g}$$
$$x_s \sin(w) + x_c \cos(w) = Ax_c + Bu_c \tag{6.3h}$$
$$\sqrt{y_{s(i)}^2 + y_{c(i)}^2} \leq y_{e(i)} - \underline{y}_{(i)} - \varepsilon_{y(i)},\ i \in \mathbb{Z}_1^p \tag{6.3i}$$
$$\sqrt{y_{s(i)}^2 + y_{c(i)}^2} \leq \overline{y}_{(i)} - \varepsilon_{y(i)} - y_{e(i)},\ i \in \mathbb{Z}_1^p, \tag{6.3j}$$

where $\epsilon_y \in \mathbb{R}^p$ is a vector with arbitrarily small positive components, and the cost function

$$J_h(\mathbf{x}, \mathbf{u}, \mathbf{x}_H, \mathbf{u}_H; x_r, u_r) = \ell_h(\mathbf{x}, \mathbf{u}, \mathbf{x}_H, \mathbf{u}_H) + V_h(\mathbf{x}_H, \mathbf{u}_H; x_r, u_r)$$

is composed of two terms: the summation of stage costs

$$\ell_h(\cdot) = \sum_{j=0}^{N-1} \|x_j - x_{hj}\|_Q^2 + \|u_j - u_{hj}\|_R^2, \tag{6.4}$$

and the offset cost function

$$V_h(\cdot) = \|x_e - x_r\|_{T_e}^2 + \|u_e - u_r\|_{S_e}^2 + \|x_s\|_{T_h}^2 + \|x_c\|_{T_h}^2 + \|u_s\|_{S_h}^2 + \|u_c\|_{S_h}^2, \tag{6.5}$$

where we consider the following assumption.

**Assumption 6.1.** Let Assumtion 4.2 hold and assume that:

*(i)* $\underline{y} + \varepsilon_y < \overline{y} - \varepsilon_y$.

*(ii)* $Q \in \mathbb{S}_{++}^n$, $R \in \mathbb{S}_{++}^m$, $T \in \mathbb{D}_{++}^n$ and $S \in \mathbb{D}_{++}^m$.

We denote the optimal value of optimization problem (6.3) for a state $x(t)$ and a given reference $(x_r, u_r)$ by $\mathbb{H}^*(x(t); x_r, u_r){=}J_h(\mathbf{x}^*, \mathbf{u}^*, \mathbf{x}_H^*, \mathbf{u}_H^*; x_r, u_r)$, where $\mathbf{x}^*$, $\mathbf{u}^*$, $\mathbf{x}_H^*$, $\mathbf{u}_H^*$ are the arguments that minimize (6.3). Furthermore, for every $j \in \mathbb{Z}$, we denote by

$$x_{hj}^* = x_e^* + x_s^* \sin(w(j-N)) + x_c^* \cos(w(j-N)), \tag{6.6a}$$

$$u_{hj}^* = u_e^* + u_s^* \sin(w(j-N)) + u_c^* \cos(w(j-N)), \tag{6.6b}$$

the harmonic signals parameterized by $(\mathbf{x}_H^*, \mathbf{u}_H^*)$. At each discrete time instant $t$, the HMPC control law is given by $u(t) = u_0^*$.

The use of the harmonic artificial reference is heavily influenced by the extensions of the MPCT formulation to the problem of tracking periodic references [150, 151, 152]. However, in this case, even though the artificial harmonic reference (6.1) is periodic, the reference $(x_r, u_r)$ to be tracked is a (piecewise) constant set-point.

Periodic MPC for tracking formulations are used to track a generic periodic signal with period $T \in \mathbb{Z}_{>0}$, i.e., a reference $x_r(t)$, $t \in \mathbb{Z}$, that satisfies $x_r(t) = x_r(t+T)$ $\forall t$, by employing artificial periodic reference signals $\{x_{sj}\}$, $\{u_{sj}\}$, $j \in \mathbb{Z}$, with period $T$. The decision variables $x_{s0}$ to $x_{sT}$ and $u_{s0}$ to $u_{sT}$ must satisfy the system dynamics (equality constraints) and constraints (inequality constraints). Therefore, the amount of decision variables and constraints grows with $T$.

Let us now focus on the problem of tracking a set-point $(x_r, u_r)$, and note that this reference can be viewed as a (constant) periodic reference signal whose period is any arbitrary $T \in \mathbb{Z}_{>0}$. Therefore, a periodic MPC for tracking formulation can be used, where the selection of the period $T$ now becomes a tuning parameter of the controller. It is rather intuitive that the selection of $T$ will have an effect on both the performance and domain of attraction of the controller. Its effect on the performance of the controller is not clear, but it seems reasonable to assume that larger periods $T$ would result in a larger domain of attraction, since we are allowing for more degrees of freedom. However, as discussed above, this would come at the expense of an increase in the number of decision variables and constraints.

The idea behind the HMPC formulation is to use an artificial periodic reference signal whose period does not affect the complexity of the optimization problem, and which can therefore be selected to improve the properties of the controller. In particular, the period of the artificial harmonic reference (6.1) is determined by the design parameter $w$. However, the constraints required to impose the system dynamics and constraints do not depend on the value of $w$. Indeed, constraints (6.3f)-(6.3h) impose the system dynamics (4.1) and constraints (6.3i)-(6.3j) impose the constraints (4.4), as we prove further ahead (see Corollaries 6.5 and 6.8).

Constraint (6.3e) imposes that the predicted state $x_N$ reaches the harmonic artificial reference at the end of the prediction horizon, since $x_{hN} = x_e + x_c$. Then, noting that the artificial harmonic reference satisfies the system dynamics and constraints, it acts as an admissible invariant set of (4.1) under (4.4). That it, the system would be able to remain in a stable admissible trajectory for all $j > N$ by applying the admissible control actions (6.6b).

Problem (6.3) is a second order cone programming problem due to the inclusion of the constraints (6.3i)-(6.3j), instead of the QP problem typically derived from most linear MPC formulations. However, this class of convex optimization problem is a well studied problem in the literature for which several efficient solvers are available, such as COSMO [153] or ECOS [154].

**Remark 6.2.** *Note that the constraints* (6.3b)-(6.3j) *do not depend on the reference. Therefore, the feasibility region (Definition 4.5) of the HMPC controller is independent of the reference. As such, feasibility is never lost in the event of reference changes.*

The following proposition relates constraints (6.3e)-(6.3h) to the satisfaction of the system dynamics. Its proof makes use of the following lemma.

**Lemma 6.3.** Let the elements $v_\ell \in \mathbb{R}^{n_v}$ of a sequence $\{v_\ell\}$ be given by

$$v_\ell = v_e + v_s \sin(w\ell) + v_c \cos(w\ell),\ \forall \ell \in \mathbb{Z},$$

where $w \in \mathbb{R}$ and $v_e \in \mathbb{R}^{n_v}$, $v_s \in \mathbb{R}^{n_v}$, and $v_c \in \mathbb{R}^{n_v}$. Then,

*(i)* $v_{\ell+1} = v_e + v_s^+ \sin(w\ell) + v_c^+ \cos(w\ell),\ \forall \ell \in \mathbb{Z}$,

*(ii)* $\left(v^+_{s(i)}\right)^2 + \left(v^+_{c(i)}\right)^2 = v^2_{s(i)} + v^2_{c(i)},\ i \in \mathbb{Z}_1^{n_v}$,

where $v_s^+ = v_s \cos(w) - v_c \sin(w)$ and $v_c^+ = v_s \sin(w) + v_c \cos(w)$.

***Proof:*** The proof relies on the following well-known trigonometric identities

$$\begin{aligned}\sin(\alpha+\beta) &= \sin(\alpha)\cos(\beta) + \cos(\alpha)\sin(\beta)\\ \cos(\alpha+\beta) &= \cos(\alpha)\cos(\beta) - \sin(\alpha)\sin(\beta).\end{aligned}$$

From these expressions we obtain

$$\begin{aligned}\sin(w(\ell+1)) &= \sin(w)\cos(w\ell) + \cos(w)\sin(w\ell)\\ \cos(w(\ell+1)) &= \cos(w)\cos(w\ell) - \sin(w)\sin(w\ell).\end{aligned}$$

Therefore,

$$\begin{aligned}v_{\ell+1} &= v_e + v_s \sin(w(\ell+1)) + v_c \cos(w(\ell+1))\\ &= v_e + v_s\left[\sin(w)\cos(w\ell) + \cos(w)\sin(w\ell)\right] + v_c\left[\cos(w)\cos(w\ell) - \sin(w)\sin(w\ell)\right]\\ &= v_e + \left[v_s\cos(w) - v_c\sin(w)\right]\sin(w\ell) + \left[v_s\sin(w) + v_c\cos(w)\right]\cos(w\ell)\\ &= v_e + v_s^+\sin(w\ell) + v_c^+\cos(w\ell),\end{aligned}$$

which proves claim *(i)*. Denote now

$$\mathbf{H}_w \doteq \begin{bmatrix} \cos(w) & -\sin(w)\\ \sin(w) & \cos(w)\end{bmatrix}.$$

With this notation,

$$\left(v^+_{s(i)}, v^+_{c(i)}\right) = \mathbf{H}_w \left(v_{s(i)}, v_{c(i)}\right), \; i \in \mathbb{Z}_1^{n_v}.$$

From the identity $\sin^2(w) + \cos^2(w) = 1$ we obtain $\mathbf{H}_w^\top \mathbf{H}_w = \mathbf{I}_2$. We are now in a position to prove claim *(ii)*:

$$\begin{aligned}\left(v^+_{s(i)}\right)^2 + \left(v^+_{c(i)}\right)^2 = \left\|\left(v^+_{s(i)}, v^+_{c(i)}\right)\right\|_2^2 &= \left(v_{s(i)}, v_{c(i)}\right)^\top \mathbf{H}_w^\top \mathbf{H}_w \left(v_{s(i)}, v_{c(i)}\right) \\ &= \left(v_{s(i)}, v_{c(i)}\right)^\top \left(v_{s(i)}, v_{c(i)}\right) = v^2_{s(i)} + v^2_{c(i)}.\end{aligned}$$

■

**Proposition 6.4.** Given the system $x_{j+1} = Ax_j + Bu_j$, suppose that

$$\begin{aligned}&u_{N+\ell} = u_e + u_s \sin(w\ell) + u_c \cos(w\ell), \; \forall \ell \in \mathbb{Z}_{>0}, \\ &x_N = x_e + x_c, \\ &x_e = Ax_e + Bu_e, \\ &x_s \cos(w) - x_c \sin(w) = Ax_s + Bu_s, \\ &x_s \sin(w) + x_c \cos(w) = Ax_c + Bu_c.\end{aligned}$$

Then

$$x_{N+\ell} = x_e + x_s \sin(w\ell) + x_c \cos(w\ell), \; \forall \ell \in \mathbb{Z}_{>0}.$$

***Proof:*** Since $x_N = x_e + x_c$, the claim is trivially satisfied for $\ell = 0$. Suppose now that the claim is satisfied for $\ell \geq 0$, we show that it is also satisfied for $\ell + 1$:

$$\begin{aligned}x_{N+\ell+1} &= Ax_{N+\ell} + Bu_{N+\ell} \\ &= A\left[x_e + x_s \sin(w\ell) + x_c \cos(w\ell)\right] + B\left[u_e + u_s \sin(w\ell) + u_c \cos(w\ell)\right] \\ &= Ax_e + Bu_e + (Ax_s + Bu_s)\sin(w\ell) + (Ax_c + Bu_c)\cos(w\ell) \\ &= x_e + \left[x_s \cos(w) - x_c \sin(w)\right]\sin(w\ell) + \left[x_s \sin(w) + x_c \cos(w)\right]\cos(w\ell) \\ &\overset{(*)}{=} x_e + x_s \sin(w(\ell+1)) + x_c \cos(w(\ell+1)),\end{aligned}$$

where $(*)$ is due to Lemma 6.3. ■

From Proposition 6.4 we derive the following corollary.

**Corollary 6.5.** The artificial harmonic reference $\{x_{hj}\}$ and $\{u_{hj}\}$ (6.1) obtained from any feasible solution of (6.3) satisfies the system dynamics (4.1). That is,

$$x_{hj+1} = Ax_{hj} + Bu_{hj}, \; \forall j \in \mathbb{Z}.$$

The following proposition relates the constraints (6.3i)-(6.3j) to the satisfaction of the system constraints. Its proof makes use of the following lemma.

**Lemma 6.6.** Let the elements $v_\ell \in \mathbb{R}^{n_v}$ of a sequence $\{v_\ell\}$ be given by

$$v_\ell = v_e + v_s \sin(w\ell) + v_c \cos(w\ell), \ \forall \ell \in \mathbb{Z},$$

where $w \in \mathbb{R}$ and $v_e \in \mathbb{R}^{n_v}$, $v_s \in \mathbb{R}^{n_v}$, and $v_c \in \mathbb{R}^{n_v}$. Then, for every $\ell \in \mathbb{Z}$ and $i \in \mathbb{Z}_1^{n_v}$, we have,

$$v_{\ell(i)} \leq v_{e(i)} + \sqrt{v_{s(i)}^2 + v_{c(i)}^2}, \tag{6.7a}$$

$$v_{\ell(i)} \geq v_{e(i)} - \sqrt{v_{s(i)}^2 + v_{c(i)}^2}. \tag{6.7b}$$

***Proof:*** We prove inequality (6.7a). The proof for (6.7b) is similar. For every $\ell \in \mathbb{Z}$ and $i \in \mathbb{Z}_1^{n_v}$, we have that

$$\begin{aligned} v_{\ell(i)} &= v_{e(i)} + v_{s(i)} \sin(w\ell) + v_{c(i)} \cos(w\ell) = v_{e(i)} + \left(v_{s(i)}, v_{c(i)}\right) (\sin(w\ell), \cos(w\ell))^\top \\ &\overset{(*)}{\leq} v_{e(i)} + \left\| \left(v_{s(i)}, v_{c(i)}\right) \right\|_2 \ \|(\sin(w\ell), \cos(w\ell))\|_2 = v_{e(i)} + \sqrt{v_{s(i)}^2 + v_{c(i)}^2}, \end{aligned}$$

where $(*)$ is due to the Cauchy-Schwarz inequality. ■

**Proposition 6.7.** Let Assumtion 6.1 hold. Consider the sequences $\{x_{hj}\}$ and $\{u_{hj}\}$ given by (6.1) and assume that

$$\sqrt{y_{s(i)}^2 + y_{c(i)}^2} \leq y_{e(i)} - \underline{y}_{(i)} - \varepsilon_{y(i)}, \ i \in \mathbb{Z}_1^p, \tag{6.8a}$$

$$\sqrt{y_{s(i)}^2 + y_{c(i)}^2} \leq \overline{y}_{(i)} - \varepsilon_{y(i)} - y_{e(i)}, \ i \in \mathbb{Z}_1^p, \tag{6.8b}$$

where $y_e \in \mathbb{R}^p$, $y_s \in \mathbb{R}^p$ and $y_c \in \mathbb{R}^p$ are given by (6.2). Then,

$$\underline{y} + \varepsilon_y \leq E x_{hj} + F u_{hj} \leq \overline{y} - \varepsilon_y, \ \forall j \in \mathbb{Z}.$$

***Proof:*** Let us define the sequence $\{y_j\}$, $j \in \mathbb{Z}$, where $y_j \in \mathbb{R}^p$ is given by

$$\begin{aligned} y_j &\doteq E x_{hj} + F u_{hj} \\ &= E x_e + F u_e + (E x_s + F u_s) \sin(w(j-N)) + (E x_c + F u_c) \cos(w(j-N)) \\ &\overset{(6.2)}{=} y_e + y_s \sin(w(j-N)) + y_c \cos(w(j-N)). \end{aligned}$$

From Lemma 6.6, we have that, $\forall j \in \mathbb{Z}$,

$$\begin{aligned} y_{j(i)} &\leq y_{e(i)} + \sqrt{y_{s(i)}^2 + y_{c(i)}^2}, \ \forall i \in \\ y_{j(i)} &\geq y_{e(i)} - \sqrt{y_{s(i)}^2 + y_{c(i)}^2}, \ \forall i \in \mathbb{Z}_1^p, \end{aligned}$$

which, taking into account (6.8a) and (6.8b), leads to

$$\begin{aligned} y_{j(i)} &\leq y_{e(i)} + \overline{y}_{(i)} - \varepsilon_{y(i)} - y_{e(i)} = \overline{y}_{(i)} - \varepsilon_{y(i)}, \\ y_{j(i)} &\geq y_{e(i)} + \underline{y}_{(i)} + \varepsilon_{y(i)} - y_{e(i)} = \underline{y}_{(i)} + \varepsilon_{y(i)}, \ \forall j \in \mathbb{Z}, \ \forall i \in \mathbb{Z}_1^p, \end{aligned}$$

which along with Assumtion 6.1.*(i)* lets us conclude that

$$\underline{y} \leq \underline{y} + \varepsilon_y \leq y_j \leq \overline{y} - \varepsilon_y \leq \overline{y}, \; \forall j \in \mathbb{Z}.$$

■

From Proposition 6.7 we derive the following corollary, which states the satisfaction of the system constraints by the artificial harmonic reference (6.1) obtained from any feasible solution of (6.3).

**Corollary 6.8.** The artificial harmonic reference $\{x_{hj}\}$ and $\{u_{hj}\}$ (6.1) obtained from any feasible solution of (6.3) satisfies the system constraints (4.4). That is,

$$\underline{y} \leq E x_{hj} + F u_{hj} \leq \overline{y}, \; \forall j \in \mathbb{Z}.$$

## 6.2 Recursive feasibility of the HMPC formulation

The following theorem states the recursive feasibility of the HMPC formulation (6.3). That is, suppose that a state $x \in \mathbb{R}^n$ belongs to the feasibility region (Definition 4.5) of the HMPC controller. Then, for any feasible solution $\mathbf{x}$, $\mathbf{u}$, $\mathbf{x}_H$ and $\mathbf{u}_H$ of $\mathbb{H}(x(t); x_r, u_r)$, we have that the successor state $Ax(t) + Bu_0$ also belongs to the feasibility region of the HMPC formulation. The proof of the theorem follows the standard approach for proving the recursive feasibility of MPC formulations: it shows that a feasible solution of $\mathbb{H}(Ax(t) + Bu_0; x_r, u_r)$ can always be constructed from a feasible solution of $\mathbb{H}(x(t); x_r, u_r)$. Note that the theorem states that recursive feasibility is maintained even if the reference is changed, as we highlighted in Remark 6.2.

**Theorem 6.9** (Recursive feasibility of the HMPC formulation)**.** Suppose that $x(t) \in \mathbb{R}^n$ belongs to the feasibility region of the HMPC formulation (6.3). Let $\mathbf{x} = \{x_0, x_1, \dots, x_{N-1}\}$, $\mathbf{u} = \{u_0, u_1, \dots, u_{N-1}\}$, $x_e$, $x_s$, $x_c$, $u_e$, $u_s$, $u_c$ be a feasible solution of $\mathbb{H}(x(t); x_r, u_r)$ for a given $(x_r, u_r) \in \mathbb{R}^n \times \mathbb{R}^m$. Then, the successor state $Ax(t) + Bu_0$ belongs to the feasibility region of the HMPC formulation for any $(\hat{x}_r, \hat{u}_r) \in \mathbb{R}^n \times \mathbb{R}^m$.

***Proof:*** In the following, we will show that $\mathbf{x}^+ = \{x_0^+, x_1^+, \dots, x_{N-1}^+\}$, $\mathbf{u}^+ = \{u_0^+, u_1^+, \dots, u_{N-1}^+\}$, $x_e^+$, $x_s^+$, $x_c^+$, $u_e^+$, $u_s^+$, $u_c^+$ given by

$$u_j^+ \doteq u_{j+1}, \; j \in \mathbb{Z}_0^{N-2}, \tag{6.9a}$$

$$u_{N-1}^+ \doteq u_e + u_c, \tag{6.9b}$$

$$x_0^+ \doteq Ax + Bu_0, \tag{6.9c}$$

$$x_{j+1} \doteq Ax_j^+ + Bu_j^+, \; j \in \mathbb{Z}_0^{N-1}, \tag{6.9d}$$

$$u_e^+ \doteq u_e, \tag{6.9e}$$

$$u_s^+ \doteq u_s \cos(w) - u_c \sin(w), \tag{6.9f}$$

$$u_c^+ \doteq u_s \sin(w) + u_c \cos(w), \tag{6.9g}$$
$$x_e^+ \doteq Ax_e + Bu_e, \tag{6.9h}$$
$$x_s^+ \doteq Ax_s + Bu_s, \tag{6.9i}$$
$$x_c^+ \doteq Ax_c + Bu_c, \tag{6.9j}$$

is a feasible solution of $\mathbb{H}(Ax(t) + Bu_0; \hat{x}_r, \hat{u}_r)$ by showing that they satisfy the constraints (6.3b)-(6.3j). That is, we prove in what follows that

$$x_{j+1}^+ = Ax_j^+ + Bu_j^+,\ j \in \mathbb{Z}_0^{N-1}, \tag{6.10a}$$
$$\underline{y} \leq Ex_j^+ + Fu_j^+ \leq \overline{y},\ j \in \mathbb{Z}_0^{N-1}, \tag{6.10b}$$
$$x_0^+ = Ax + Bu_0, \tag{6.10c}$$
$$x_N^+ = x_e^+ + x_c^+, \tag{6.10d}$$
$$x_e^+ = Ax_e^+ + Bu_e^+, \tag{6.10e}$$
$$x_s^+ \cos(w) - x_c^+ \sin(w) = Ax_s^+ + Bu_s^+, \tag{6.10f}$$
$$x_s^+ \sin(w) + x_c^+ \cos(w) = Ax_c^+ + Bu_c^+, \tag{6.10g}$$
$$\sqrt{(y_{s(i)}^+)^2 + (y_{c(i)}^+)^2} \leq y_{e(i)}^+ - \underline{y}_{(i)} - \varepsilon_{y(i)},\ i \in \mathbb{Z}_1^p, \tag{6.10h}$$
$$\sqrt{(y_{s(i)}^+)^2 + (y_{c(i)}^+)^2} \leq \overline{y}_{(i)} - \varepsilon_{y(i)} - y_{e(i)}^+,\ i \in \mathbb{Z}_1^p, \tag{6.10i}$$

where $y_e^+$, $y_s^+$ and $y_c^+$ are given by

$$y_e^+ \doteq Ex_e^+ + Fu_e^+, \quad y_s^+ \doteq Ex_s^+ + Fu_s^+, \quad y_c^+ \doteq Ex_c^+ + Fu_c^+.$$

Equalities (6.10a) and (6.10c) are trivially satisfied by construction, as evident from (6.9c) and (6.9d). Since $x_0^+ = Ax + Bu_0 = x_1$, and $u_j^+ = u_{j1}$, $j \in \mathbb{Z}_0^{N-2}$ (6.9a), we have

$$(x_j^+, u_j^+) = (x_{j+1}, u_{j+1}),\ j \in \mathbb{Z}_0^{N-2}. \tag{6.11}$$

Therefore, from (6.3c) we obtain

$$\underline{y} \leq Ex_j^+ + Fu_j^+ \leq \overline{y},\ j \in \mathbb{Z}_0^{N-2}. \tag{6.12}$$

We now compute the value of $x_{N-1}^+$:

$$x_{N-1}^+ = Ax_{N-2}^+ + Bu_{N-2}^+ = Ax_{N-1} + u_{N-1} = x_N = x_e + x_c = x_{hN}. \tag{6.13}$$

Additionally, $u_{N-1}^+ = u_e + u_c = u_{hN}$. That is, $(x_{N-1}^+, u_{N-1}^+)$ belongs to the artificial harmonic reference obtained from a feasible solution of (6.3). Therefore, from Corollary 6.8, we have that

$$\underline{y} \leq Ex_{N-1}^+ + Fu_{N-1}^+ \leq \overline{y},$$

which along with (6.12) proves (6.10b). The value of $x_N^+$ can be computed from $x_{N-1}^+ \overset{(6.13)}{=} x_N$ and $u_{N-1}^+ \overset{(6.9b)}{=} u_e + u_c$ as follows:

$$\begin{aligned} x_N^+ &= Ax_{N-1}^+ + Bu_{N-1}^+ = Ax_N + B(u_e + u_c) \\ &= A(x_e + x_c) + B(u_e + u_c) \overset{(6.9h),(6.9j)}{=} x_e^+ + x_c^+, \end{aligned}$$

which proves (6.10d). From $x_e^+ \overset{(6.9h)}{=} Ax_e + Bu_e \overset{(6.3f)}{=} x_e$ and $u_e^+ \overset{(6.9e)}{=} u_e$ we have

$$x_e^+ = Ax_e + Bu_e = Ax_e^+ + Bu_e^+,$$

which proves (6.10e). We now prove (6.10f):

$$\begin{aligned} Ax_s^+ + Bu_s^+ &= A(Ax_s + Bu_s) + Bu_s^+ \\ &= A(x_s \cos(w) - x_c \sin(w) + B(u_s \cos(w) - u_c \sin(w)) \\ &= (Ax_s + Bu_s)\cos(w) - (Ax_c + Bu_c)\sin(w) \\ &= x_s^+ \cos(w) - x_c^+ \sin(w). \end{aligned}$$

We prove (6.10g) in a similar way:

$$\begin{aligned} Ax_c^+ + Bu_c^+ &= A(Ax_c + Bu_c) + Bu_c^+ \\ &= A(x_s \sin(w) + x_c \cos(w)) + B(u_s \sin(w) + u_c \cos(w)) \\ &= (Ax_s + Bu_s)\sin(w) + (Ax_c + Bu_c)\cos(w) \\ &= x_s^+ \sin(w) + x_c^+ \cos(w). \end{aligned}$$

Next, we express $y_e^+$, $y_s^+$ and $y_c^+$ in terms of $y_e$, $y_s$ and $y_c$:

$$\begin{aligned} y_e^+ &= Ex_e^+ + Fu_e^+ = Ex_e + Fu_e = y_e, \\ y_s^+ &= Ex_s^+ + Fu_s^+ = E(Ax_s + Bu_s) + Fu_s^+ \\ &= E(x_s \cos(w) - x_c \sin(w)) + F(u_s \cos(w) - u_c \sin(w)) \\ &= (Ex_s + Fu_s)\cos(w) - (Ex_c + Fu_c)\sin(w) \\ &= y_s \cos(w) - y_c \sin(w), \\ y_c^+ &= Ex_c^+ + Fu_c^+ = E(Ax_c + Bu_c) + Fu_c^+ \\ &= E(x_s \sin(w) + x_c \cos(w)) + F(u_s \sin(w) + u_c \cos(w)) \\ &= (Ex_s + Fu_s)\sin(w) + (Ex_c + Fu_c)\cos(w) \\ &= y_s \sin(w) + y_c \cos(w), \end{aligned}$$

which, in view of Lemma 6.3.*(ii)*, leads to

$$\sqrt{\left(y_{s(i)}^+\right)^2 + \left(y_{c(i)}^+\right)^2} = \sqrt{y_{s(i)}^2 + y_{c(i)}^2},\ i \in \mathbb{Z}_1^p,$$

from where we conclude that (6.10h)-(6.10i) are directly inferred from (6.3i)-(6.3j). The fact that the reference $(\hat{x}_r, \hat{u}_r)$ does not have to be the same as $(x_r, u_r)$ is inferred from it not affecting the constraints (6.3b)-(6.3j). ■

## 6.3 Asymptotic stability of the HMPC formulation

In this section we prove the asymptotic stability of the HMPC formulation to the *optimal artificial harmonic reference*, which we now define.

**Definition 6.10** (Optimal artificial harmonic reference)**.** Given a reference $(x_r, u_r) \in \mathbb{R}^n \times \mathbb{R}^m$, we define the *optimal artificial harmonic reference* of the HMPC formulation (6.3) as the harmonic sequences $\{x^\circ_{hj}\}$, $\{u^\circ_{hj}\}$, $j \in \mathbb{Z}$, parameterized by the unique solution $(\mathbf{x}^\circ_H, \mathbf{u}^\circ_H)$ of the strongly convex optimization problem

$$
\begin{aligned}
(\mathbf{x}^\circ_H, \mathbf{u}^\circ_H) = \arg \min_{\mathbf{x}_H, \mathbf{u}_H} \; & V_h(\mathbf{x}_H, \mathbf{u}_H; x_r, u_r) \\
s.t. \; & (6.3f)\text{-}(6.3j).
\end{aligned}
\tag{6.14}
$$

Additionally, we denote by $V^\circ_h(x_r, u_r) \doteq V_h(\mathbf{x}^\circ_H, \mathbf{u}^\circ_H; x_r, u_r)$ the optimal value of problem (6.14).

The following proposition characterizes the optimal artificial harmonic reference by stating that it is an admissible steady state, i.e., that $x^\circ_{hj} = x^\circ_e$ and $u^\circ_{hj} = u^\circ_e$ for all $j \in \mathbb{Z}$.

**Proposition 6.11** (Characterization of the optimal artificial harmonic reference)**.** Consider optimization problem (6.14). Then, for any $(x_r, u_r) \in \mathbb{R}^n \times \mathbb{R}^m$, its optimal solution is the steady state $(x^\circ_e, u^\circ_e)$ of (4.1) satisfying

$$
\underline{y} + \varepsilon_y \leq E x^\circ_e + F u^\circ_e \leq \overline{y} - \varepsilon_y \tag{6.15}
$$

that minimizes the cost $\|x^\circ_e - x_r\|^2_{T_e} + \|u^\circ_e - u_r\|^2_{S_e}$. That is, $\mathbf{x}^\circ_H = (x^\circ_e, \mathbf{0}_n, \mathbf{0}_n)$ and $\mathbf{u}^\circ_H = (u^\circ_e, \mathbf{0}_m, \mathbf{0}_m)$.

***Proof:*** We prove the lemma by contradiction. Assume that $\hat{\mathbf{x}}^\circ_H = (x^\circ_e, x^\circ_s, x^\circ_c)$, $\hat{\mathbf{u}}^\circ_H = (u^\circ_e, u^\circ_s, u^\circ_c)$, is the optimal solution of (6.14) and that at least some (if not all) of $x^\circ_s$, $x^\circ_c$, $u^\circ_s$, $u^\circ_c \neq \mathbf{0}$. First, we show that $\mathbf{x}^\circ_H = (x^\circ_e, \mathbf{0}_n, \mathbf{0}_n)$, $\mathbf{u}^\circ_H = (u^\circ_e, \mathbf{0}_m, \mathbf{0}_m)$ satisfy (6.3f)-(6.3j). Constraints (6.3g) and (6.3h) are trivially satisfied and (6.3f) is satisfied since $(\hat{\mathbf{x}}^\circ_H, \hat{\mathbf{u}}^\circ_H)$ is assumed to be the solution of (6.14). Moreover, since

$$
0 \leq \sqrt{(y^\circ_{s(i)})^2 + (y^\circ_{c(i)})^2}, \; \forall i \in \mathbb{Z}_1^p,
$$

we have that (6.3i) and (6.3j) are also satisfied for $(\mathbf{x}^\circ_H, \mathbf{u}^\circ_H)$. Finally, it is clear from the initial assumption and (6.5) that

$$
V_h(\mathbf{x}^\circ_H, \mathbf{u}^\circ_H; x_r, u_r) < V_h(\hat{\mathbf{x}}^\circ_H, \hat{\mathbf{u}}^\circ_H; x_r, u_r),
$$

contradicting the optimality of $(\hat{\mathbf{x}}^\circ_H, \hat{\mathbf{u}}^\circ_H)$. Finally (6.15) follows from the satisfaction of (6.3i)-(6.3j) along with Proposition 6.7. ■

To prove the asymptotic stability of the HMPC formulation we make use of the following well known asymptotic stability theorem [64, Appendix B.3].

**Theorem 6.12** (Lyapunov asymptotic stability)**.** Consider an autonomous discrete time system $z(t+1) = f(z(t))$, $t \in \mathbb{Z}$, with states $z(t) \in \mathbb{R}^n$ and where the function $f : \mathbb{R}^n \to \mathbb{R}^n$ is continuous and satisfies $f(\mathbf{0}_n) = \mathbf{0}_n$. Let $\mathcal{Z}$ be an invariant set of the system and $\Omega \subseteq \mathcal{Z}$ be a compact set, both including the origin as an interior point. If there exists a continuous function $W : \mathbb{R}^n \to \mathbb{R}_{\geq 0}$ and suitable functions $\alpha_1 \in \mathcal{K}_\infty$ and $\alpha_2 \in \mathcal{K}_\infty$ such that,

*(i)* $W(z(t)) \geq \alpha_1(\|z(t)\|_2)$, $\forall z(t) \in \mathcal{Z}$,

*(ii)* $W(z(t)) \leq \alpha_2(\|z(t)\|_2)$, $\forall z(t) \in \Omega$,

*(iii)* $W(z(t+1)) < W(z(t)), \forall z(t) \in \mathcal{Z} \setminus \{\mathbf{0}_n\}$,
and $W(z(t+1)) = W(z(t))$ if $z(t) = \mathbf{0}_n$,

then $W$ is a Lyapunov function for $z(t+1) = f(z(t))$ in $\mathcal{Z}$ and the origin is asymptotically stable for all initial states in $\mathcal{Z}$.

The following theorem states the asymptotic stability of the HMPC formulation to the optimal artificial harmonic reference $(x_e^\circ, u_e^\circ)$. Its proof is inspired in the proof of the asymptotic stability of the MPCT formulation [155, Theorem 1], and makes use of the following lemma, whose proof is inspired by [155, Lemma 1]. However, in contrast to the proof presented in [155], we directly derive a Lyapunov function that satisfies the conditions of Theorem 6.12.

**Definition 6.13.** Consider a controllable system (4.1). Its *controllability index* is the smallest integer $j > 0$ for which $\text{rank}\big([B, AB, A^2B, \dots, A^{j-1}B]\big) = n$.

**Lemma 6.14.** Consider a system (4.1) subject to (4.4) controlled with the HMPC formulation (6.3). Let Assumtion 6.1 hold and assume that $N$ is greater or equal to the controllability index of (4.1). Let $(x_r, u_r) \in \mathbb{R}^n \times \mathbb{R}^m$ be a given reference and $x \in \mathbb{R}^n$ be a state belonging to the feasibility region of the $\mathbb{H}(x; x_r, u_r)$. Then, $x = x_{h0}^*$ if and only if $x = x_e^\circ$, where $x_{h0}^*$ is given by (6.6a) and $x_e^\circ$ by Proposition 6.11.

***Proof:*** Let $V_h^* \doteq V_h(\mathbf{x}_H^*, \mathbf{u}_H^*; x_r, u_r)$ and $V_h^\circ \doteq V_h(\mathbf{x}_H^\circ, \mathbf{u}_H^\circ; x_r, u_r)$. First, we prove the implication $x = x_{h0}^* \implies x = x_e^\circ$. Assume that $x = x_{h0}^*$. Then, we have that $\mathbb{H}^*(x; x_r, u_r) = V_h^*$, i.e., the optimal solution of $\mathbb{H}(x; x_r, u_r)$ is given by

$$x_j^* = x_{hj}^*, \quad u_j^* = u_{hj}^*, \quad \forall j \in \mathbb{Z}_0^{N-1}, \tag{6.16}$$

where $x_{hj}^*$ and $u_{hj}^*$ are given by (6.6). Indeed, the summation of stage costs (6.4) satisfies $\ell_h(\mathbf{x}^*, \mathbf{u}^*, \mathbf{x}_H^*, \mathbf{u}_H^*) = 0$, which, under Assumtion 6.1.*(ii)* is its smallest possible value. Additionally, from Corollaries 6.5 and 6.8 it is clear that (6.16) is a feasible solution of (6.3b)-(6.3j).

Next, we prove that $V_h^* = V_h^\circ$ by contradiction. Assume that $V_h^* > V_h^\circ$. Under Assumtion 6.1, problem (6.14) is strongly convex. Therefore, $(\mathbf{x}_H^\circ, \mathbf{u}_H^\circ)$ is the unique minimizer of $V_h(\mathbf{x}_H, \mathbf{u}_H; x_r, u_r)$ for all $(\mathbf{x}_H, \mathbf{u}_H)$ that satisfy (6.3f)-(6.3j), which implies that $(\mathbf{x}_H^*, \mathbf{u}_H^*) \neq (\mathbf{x}_H^\circ, \mathbf{u}_H^\circ)$.

Let $\hat{\mathbf{x}}_H$ and $\hat{\mathbf{u}}_H$ be defined as

$$\hat{\mathbf{x}}_H = (\hat{x}_e, \hat{x}_s, \hat{x}_c) = \lambda \mathbf{x}_H^* + (1-\lambda)\mathbf{x}_H^\circ = \lambda(x_e^*, x_s^*, x_c^*) + (1-\lambda)(x_e^\circ, \mathbf{0}_n, \mathbf{0}_n), \ \lambda \in [0,1],$$
$$\hat{\mathbf{u}}_H = (\hat{u}_e, \hat{u}_s, \hat{u}_c) = \lambda \mathbf{u}_H^* + (1-\lambda)\mathbf{u}_H^\circ = \lambda(u_e^*, u_s^*, u_c^*) + (1-\lambda)(u_e^\circ, \mathbf{0}_m, \mathbf{0}_m), \ \lambda \in [0,1].$$

Then, since $N$ is assumed to be greater or equal than to the controllability index of the system, $\mathcal{Y}$ is convex, and $(x_{hj}^*, u_{hj}^*) \in \text{ri}(\mathcal{Y})$ for all $j \in \mathbb{Z}$, there exists a $\hat{\lambda} \in [0,1)$ such that for any $\lambda \in [\hat{\lambda}, 1]$ there is a dead-beat control law $\mathbf{u}^{\text{db}}$ for which the predicted trajectory $\mathbf{x}^{\text{db}}$ satisfying $x_0^{\text{db}} = x_{h0}^*$ and $x_N^{\text{db}} = \hat{x}_{h0}$ is a feasible solution $(\mathbf{x}^{\text{db}}, \mathbf{u}^{\text{db}}, \hat{\mathbf{x}}_H, \hat{\mathbf{u}}_H)$ of problem $\mathbb{H}(x_{h0}^*; x_r, u_r)$.

Taking into account the optimality of (6.16), and noting that there exists a matrix $P \in \mathbb{S}_{++}^n$ such that

$$\sum_{j=0}^{N-1} \|x_j^{\text{db}} - \hat{x}_{hj}\|_Q^2 + \|u_j^{\text{db}} - \hat{u}_{hj}\|_R^2 \leq \|x_0^{\text{db}} - \hat{x}_{h0}\|_P^2,$$

we have that

$$\begin{aligned} V_h^* &= J_h(\mathbf{x}^*, \mathbf{u}^*, \mathbf{x}_H^*, \mathbf{u}_H^*; x_r, u_r) \leq J_h(\mathbf{x}^{\text{db}}, \mathbf{u}^{\text{db}}, \hat{\mathbf{x}}_H, \hat{\mathbf{u}}_H; x_r, u_r) \\ &= \ell_h(\mathbf{x}^{\text{db}}, \mathbf{u}^{\text{db}}, \hat{\mathbf{x}}_H, \hat{\mathbf{u}}_H) + V_h(\hat{\mathbf{x}}_H, \hat{\mathbf{u}}_H; x_r, u_r) \\ &\leq \|x_{h0}^* - \hat{x}_{h0}\|_P^2 + V_h(\hat{\mathbf{x}}_H, \hat{\mathbf{u}}_H; x_r, u_r) \\ &\overset{(*)}{=} (1-\lambda)^2 \|x_{h0}^* - x_e^\circ\|_P^2 + V_h(\hat{\mathbf{x}}_H, \hat{\mathbf{u}}_H; x_r, u_r), \end{aligned} \tag{6.17}$$

where step $(*)$ is using

$$\begin{aligned} x_{h0}^* - \hat{x}_{h0} &= x_{h0}^* - [\lambda x_{h0}^* + (1-\lambda) x_{h0}^\circ] \\ &= (1-\lambda)(x_{h0}^* - x_{h0}^\circ) = (1-\lambda)(x_{h0}^* - x_e^\circ). \end{aligned}$$

From the convexity of $V_h$ we have that

$$V_h(\hat{\mathbf{x}}_H, \hat{\mathbf{u}}_H; x_r, u_r) \leq \lambda V_h^* + (1-\lambda) V_h^\circ, \ \lambda \in [0,1],$$

which combined with (6.17) leads to,

$$V_h^* \leq \theta(\lambda), \ \lambda \in [\hat{\lambda}, 1], \tag{6.18}$$

where

$$\theta(\lambda) \doteq (1-\lambda)^2 \|x_{h0}^* - x_e^\circ\|_P^2 + \lambda(V_h^* - V_h^\circ) + V_h^\circ.$$

The derivative of $\theta(\lambda)$ (w.r.t. $\lambda$) is

$$\theta'(\lambda) = -2(1-\lambda)\|x_{h0}^* - x_e^\circ\|_P^2 + (V_h^* - V_h^\circ).$$

Taking into account the initial assumption $V_h^* - V_h^\circ > 0$, we have that $\theta'(1) > 0$. Therefore, there exists a $\lambda \in [\hat{\lambda}, 1)$ such that $\theta(\lambda) < \theta(1) = V_h^*$, which together with (6.18) leads to the contradiction $V_h^* < V_h^*$. Therefore, we have that $V_h^* \leq V_h^\circ$. Moreover, since $(\mathbf{x}_H^\circ, \mathbf{u}_H^\circ)$ is the unique minimizer of $V_h(\mathbf{x}_H, \mathbf{u}_H; x_r, u_r)$ for all $(\mathbf{x}_H, \mathbf{u}_H)$ that satisfy (6.3f)-(6.3j), we conclude that $x_{h0}^* = x_e^\circ$.

The reverse implication is straightforward. Assume now that $x = x_e^\circ$. Then,

$$\mathbf{x}_H^* = \mathbf{x}_H^\circ,\ \mathbf{u}_H^* = \mathbf{u}_H^\circ,\ x_j^* = x_e^\circ,\ u_j^* = u_e^\circ,\ \forall j \in \mathbb{Z}_0^{N-1}, \tag{6.19}$$

is a feasible solution of $\mathbb{H}(x; x_r, u_r)$, since $(\mathbf{x}_H^\circ, \mathbf{u}_H^\circ)$ satisfies (6.3f)-(6.3j) and $(x_e^\circ, u_e^\circ)$ is a steady state of the system (4.1). Moreover, (6.19) is the optimal solution of $\mathbb{H}(x; x_r, u_r)$. Indeed, note that $\ell_h(\mathbf{x}^*, \mathbf{u}^*, \mathbf{x}_H^*, \mathbf{u}_H^*) = 0$ and that $V_h^* = V_h^\circ$, which, once again, is its minimum value for all $(\mathbf{x}_H, \mathbf{u}_H)$ satisfying (6.3f)-(6.3j). Therefore, due to the strict convexity of $V_h$, we conclude that $\mathbf{x}_H^* = \mathbf{x}_H^\circ$, implying $x_{h0}^* = x_e^\circ$. ■

**Theorem 6.15.** Consider a system (4.1) subject to (4.4) controlled with the HMPC formulation (6.3). Let Assumtion 6.1 hold and assume that $N$ is greater or equal to the controllability index of (4.1). Then, for any reference $(x_r, u_r) \in \mathbb{R}^n \times \mathbb{R}^m$ and any initial state $x(t) \in \mathbb{R}^n$ belonging to the feasibility region of $\mathbb{H}(x(t); x_r, u_r)$, the system controlled by the HMPC formulation's control law is stable, fulfills the system constraints at all future time instants, and asymptotically converges to the optimal artificial harmonic reference $(x_e^\circ, u_e^\circ) \in \mathbb{R}^n \times \mathbb{R}^m$ given by Proposition 6.11.

***Proof:*** The proof is based on finding a function that satisfies the Lyapunov conditions for asymptotic stability given in Theorem 6.12. Let $\mathbf{x}^*$, $\mathbf{u}^*$, $\mathbf{x}_H^*$, $\mathbf{u}_H^*$ be the optimal solution of $\mathbb{H}(x(t); x_r, u_r)$, $\mathbb{H}^*(x(t); x_r, u_r) \doteq J_h(\mathbf{x}^*, \mathbf{u}^*, \mathbf{x}_H^*, \mathbf{u}_H^*)$ be its optimal value, $V_h^* \doteq V_h(\mathbf{x}_H^*, \mathbf{u}_H^*; x_r, u_r)$ and $V_h^\circ \doteq V_h(\mathbf{x}_H^\circ, \mathbf{u}_H^\circ; x_r, u_r)$.

We will now show that the function

$$W(x(t); x_r, u_r) = \mathbb{H}^*(x(t); x_r, u_r) - V_h^\circ$$

is a Lyapunov function for $x(t) - x_e^\circ$ by finding suitable $\alpha_1(\|x(t) - x_e^\circ\|_2) \in \mathcal{K}_\infty$ and $\alpha_2(\|x(t) - x_e^\circ\|_2) \in \mathcal{K}_\infty$ functions satisfying the conditions of Theorem 6.12. For convenience, we will drop the $(x_r, u_r)$ from the notation.

Let $x^+ \doteq Ax(t) + Bu_0^*$ be the successor state and consider the shifted sequence $\mathbf{x}^+$, $\mathbf{u}^+$, $\mathbf{x}_H^+$, $\mathbf{u}_H^+$ defined as in (6.9) but taking $\mathbf{x}^*$, $\mathbf{u}^*$, $\mathbf{x}_H^*$, $\mathbf{u}_H^*$ in the right-hand-side of the equations. It is clear from the proof of Theorem 6.9 that this shifted sequence is a feasible solution of $\mathbb{H}(x^+)$.

The satisfaction of Theorem 6.12.*(i)* for any $x(t)$ belonging to the domain of attraction of the HMPC formulation follows from

$$\begin{aligned} W(x(t)) &= \sum_{j=0}^{N-1} \|x_j^* - x_{hj}^*\|_Q^2 + \|u_j^* - u_{hj}^*\|_R^2 + V_h^* - V_h^\circ \\ &\overset{(*)}{\geq} \|x_0^* - x_{h0}^*\|_Q^2 + \frac{\hat{\sigma}}{2}\|x_{h0}^* - x_e^\circ\|_2^2 \\ &\geq \min\{\lambda_{\min}(Q), \frac{\hat{\sigma}}{2}\} \left(\|x(t) - x_{h0}^*\|_2^2 + \|x_{h0}^* - x_e^\circ\|_2^2\right) \\ &\overset{(**)}{\geq} \frac{1}{2}\min\{\lambda_{\min}(Q), \frac{\hat{\sigma}}{2}\}\|x(t) - x_e^\circ\|_2^2, \end{aligned}$$

where $(**)$ is due to the parallelogram law and $(*)$ follows from the fact that

$$V_h^* - V_h^\circ \geq \frac{\hat{\sigma}}{2}\|x_{h0}^* - x_e^\circ\|_2^2 \tag{6.20}$$

for some $\hat{\sigma} > 0$. To show this, note that, under Assumtion 6.1, $V_h$ is a strongly convex function. Therefore, it satisfies for some $\sigma > 0$ [18, Theorem 5.24], [23, §9.1.2],

$$V_h(z) - V_h(y) \geq \langle \nabla V_h(y), z - y\rangle + \frac{\sigma}{2}\|z - y\|_2^2,$$

for all $z, y \in \mathbb{R}^n \times \mathbb{R}^n \times \mathbb{R}^n \times \mathbb{R}^m \times \mathbb{R}^m \times \mathbb{R}^m$. Particularizing for $z = (\mathbf{x}_H^*, \mathbf{u}_H^*)$ and $y = (\mathbf{x}_H^\circ, \mathbf{u}_H^\circ)$ we have,

$$V_h^* - V_h^\circ \geq \langle \nabla V_h^\circ, (\mathbf{x}_H^*, \mathbf{u}_H^*) - (\mathbf{x}_H^\circ, \mathbf{u}_H^\circ)\rangle + \frac{\sigma}{2}\|(\mathbf{x}_H^*, \mathbf{u}_H^*) - (\mathbf{x}_H^\circ, \mathbf{u}_H^\circ)\|_2^2.$$

From the optimality of $(\mathbf{x}_H^\circ, \mathbf{u}_H^\circ)$ we have that [1, Proposition 5.4.7], [23, §4.2.3],

$$\langle \nabla V_h^\circ, (\mathbf{x}_H, \mathbf{u}_H) - (\mathbf{x}_H^\circ, \mathbf{u}_H^\circ)\rangle \geq 0$$

for all $(\mathbf{x}_H, \mathbf{u}_H)$ satisfying (6.3f)-(6.3j). Since $(\mathbf{x}_H^*, \mathbf{u}_H^*)$ satisfies (6.3f)-(6.3j), this leads to

$$\begin{aligned} V_h^* - V_h^\circ &\geq \frac{\sigma}{2}\|(\mathbf{x}_H^*, \mathbf{u}_H^*) - (\mathbf{x}_H^\circ, \mathbf{u}_H^\circ)\|_2^2 \\ &= \frac{\sigma}{2}\left(\|\mathbf{x}_H^* - \mathbf{x}_H^\circ\|_2^2 + \|\mathbf{u}_H^* - \mathbf{u}_H^\circ\|_2^2\right) \\ &\geq \frac{\sigma}{2}(\|x_e^* - x_e^\circ\|_2^2 + \|x_s^*\|^2 + \|x_c^*\|_2^2) \\ &\geq \frac{\sigma}{2}(\|x_e^* - x_e^\circ\|_2^2 + \|x_s^* \sin(-wN)\|_2^2 + \|x_c^* \cos(-wN)\|_2^2). \end{aligned}$$

Finally, making use of the parallelogram law, inequality (6.20) follows from the fact that there exists a scalar $\hat{\sigma} > 0$ such that

$$\begin{aligned} &\frac{\sigma}{2}(\|x_e^* - x_e^\circ\|_2^2 + \|x_s^* \sin(-wN)\|_2^2 + \|x_c^* \cos(-wN)\|_2^2) \\ &\geq \frac{\hat{\sigma}}{2}\|x_e^* - x_e^\circ + x_s^* \sin(-wN) + x_c^* \cos(-wN)\|_2^2 \\ &= \frac{\hat{\sigma}}{2}\|x_{h0}^* - x_e^\circ\|_2^2. \end{aligned}$$

Since $(x_e^\circ, u_e^\circ) \in \mathrm{ri}(\mathcal{Y})$ (see Proposition 6.11), the system is controllable and $N$ is greater than its controllability index, there exists a sufficiently small compact set $\Omega$ containing the origin in its interior such that, for all states $x(t)$ that satisfy $x(t) - x_e^\circ \in \Omega$, the dead-beat control law

$$u_j^{\mathrm{db}} = K_{\mathrm{db}}(x_j^{\mathrm{db}} - x_e^\circ) + u_e^\circ$$

provides an admissible predicted trajectory $\mathbf{x}^{\mathrm{db}}$ of system (4.1) subject to (4.4), where $x_{j+1}^{\mathrm{db}} = Ax_j^{\mathrm{db}} + Bu_j^{\mathrm{db}}$, $j \in \mathbb{Z}_0^{N-1}$, $x_0^{\mathrm{db}} = x(t)$ and $x_N^{\mathrm{db}} = x_e^\circ$. Then, taking into account the optimality of $\mathbf{x}^*$, $\mathbf{u}^*$, $\mathbf{x}_H^*$, $\mathbf{u}_H^*$, we have that,

$$\begin{aligned} W(x) &= \ell_h(\mathbf{x}^*, \mathbf{u}^*, \mathbf{x}_H^*, \mathbf{u}_H^*) + V_h^* - V_h^\circ \\ &\leq \ell_h(\mathbf{x}^{\mathrm{db}}, \mathbf{u}^{\mathrm{db}}, \mathbf{x}_H^\circ, \mathbf{u}_H^\circ) + V_h^\circ - V_h^\circ \\ &= \sum_{j=0}^{N-1} \|x_j^{\mathrm{db}} - x_e^\circ\|_Q^2 + \|u_j^{\mathrm{db}} - u_e^\circ\|_R^2. \end{aligned}$$

Therefore, there exists a matrix $P \in \mathbb{S}_{++}^n$ such that

$$W(x(t)) \leq \lambda_{\max}(P)\|x(t) - x_e^\circ\|_2^2$$

for any $x - x_e^\circ \in \Omega$, which proves the satisfaction of Theorem 6.12.*(ii)*.

Next, let $\Delta W(x(t)) \doteq W(x^+) - W(x(t))$ and note that, as shown by (6.9b), (6.9h)-(6.9j), (6.11), (6.13) and Lemma 6.3.*(i)*, we have that $x_j^+ = x_{j+1}^*$ for $j \in \mathbb{Z}_0^{N-1}$, $u_j^+ = u_{j+1}^*$ for $j \in \mathbb{Z}_0^{N-1}$, and that $x_{hj}^+ = x_{hj+1}^*$ and $u_{hj}^+ = u_{hj+1}^*$ for $j \in \mathbb{Z}$. Then,

$$\begin{aligned} \Delta W(x(t)) &= \mathbb{H}^*(x^+) - V_h^\circ - \mathbb{H}^*(x(t)) + V_h^\circ \leq \mathbb{H}(x^+) - \mathbb{H}^*(x(t)) \\ &= \sum_{j=0}^{N-1} \Big( \|x_j^+ - x_{hj}^+\|_Q^2 + \|u_j^+ - u_{hj}^+\|_R^2 - \|x_j^* - x_{hj}^*\|_Q^2 - \|u_j^* - u_{hj}^*\|_R^2 \Big) \\ &\quad + V_h(\mathbf{x}_H^+, \mathbf{u}_H^+; x_r, u_r) - V_h^* \\ &\overset{(*)}{=} \sum_{j=1}^{N-1} \big( \|x_j^* - x_{hj}^*\|_Q^2 + \|u_j^* - u_{hj}^*\|_R^2 - \|x_j^* - x_{hj}^*\|_Q^2 - \|u_j^* - u_{hj}^*\|_R^2 \big) \\ &\quad + \|x_N^* - x_{hN}^*\|_Q^2 + \|u_N^* - u_{hN}^*\|_R^2 \\ &= -\|x_0^* - x_{h0}^*\|_Q^2 - \|u_0^* - u_{h0}^*\|_R^2 \\ &\leq -\lambda_{\min}(Q)\|x(t) - x_{h0}^*\|_2^2, \end{aligned}$$

where in step $(*)$ we are making use of the fact that

$$V_h(\mathbf{x}_H^+, \mathbf{u}_H^+; x_r, u_r) = V_h^*.$$

Indeed, note that $x_e^+ = x_e^*$ and $u_e^+ = u_e^*$. Therefore, the first two terms of $V_h(\mathbf{x}_H^+, \mathbf{u}_H^+; x_r, u_r)$ (6.5) are the same as those of $V_h^*$. We now show that, since $T_h$ and $S_h$ are diagonal matrices, the terms $\|x_s\|^2_{T_h} + \|x_c\|^2_{T_h}$ are also the same (terms $\|u_s\|^2_{S_h} + \|u_c\|^2_{S_h}$ follow similarly):

$$\begin{aligned}\|x_s^+\|^2_{T_h} + \|x_c^+\|^2_{T_h} &= \|x_s^* \cos(w) - x_c^* \sin(w)\|^2_{T_h} + \|x_s^* \sin(w) + x_c^* \cos(w)\|^2_{T_h} \\ &= (\sin(w)^2 + \cos(w)^2)\|x_s^*\|^2_{T_h} + (\sin(w)^2 + \cos(w)^2)\|x_c^*\|^2_{T_h} \\ &\quad + 2\cos(w)\sin(w)\langle x_s^*, T_h x_c^*\rangle - 2\cos(w)\sin(w)\langle x_s^*, T_h x_c^*\rangle \\ &= \|x_s^*\|^2_{T_h} + \|x_c^*\|^2_{T_h}.\end{aligned}$$

The satisfaction of Theorem 6.12.*(iii)* now follows from noting that inequality

$$W(x^+) - W(x(t)) \leq -\lambda_{\min}(Q)\|x(t) - x^*_{h0}\|^2_2,$$

along with Lemma 6.14 leads to

$$\begin{aligned}&W(x^+) < W(x(t)), \ \forall x(t) \neq x^\circ_e, \\ &W(x^+) = W(x(t)), \ \text{if}\, x(t) = x^\circ_e.\end{aligned}$$

■

## 6.4 Selection of the base frequency

The ingredients of the HMPC formulation (6.3) are similar to the ones of the MPCT formulation (5.42): the cost function matrices $Q$ and $R$ play the same role, and matrices $T_h$ and $S_h$ have a similar effect that the matrices $T$ and $S$ of MPCT. However, the base frequency $w$ of the artificial harmonic reference (6.1) plays a key role on the performance of the controller. This section discusses and provides an intuitive approach to the selection of this parameter.

It is important to note that the stability and recursive feasibility properties (Theorems 6.9 and 6.15) of the HMPC formulation are satisfied for any value of $w$. However, its proper tuning is still an important aspect of the performance and domain of attraction of the controller. For instance, note that for $w = 2\pi$, the HMPC formulation is equivalent to the MPCT formulation.

There are two main considerations to be made. The first is related to the phenomenon of aliasing and to the selection of the sampling time of continuous-time systems, which will provide an upper bound to $w$. The second relates to the frequency response of linear systems, which will provide some insight into the selection of an initial, and well suited, value of $w$. Subsequent fine tuning may provide better results, but this initial value of $w$ should work well in practice and provide a good starting point.

### Upper bound of the base frequency

The artificial harmonic reference parametrized by any feasible solution $(\mathbf{x}_H, \mathbf{u}_H)$ of (6.3) satisfies the discrete-time system dynamics, as stated in Corollary 6.5. However, $w$ should be selected so that it also correctly describes the underlying continuous-time model. That is, it should be selected small enough to prevent the aliasing phenomenon.

To do so, $w$ must be chosen below the Nyquist frequency for anti-aliasing, i.e., $w < \pi$ [156]. However, since the inputs are applied using a zero-order holder, $w$ should satisfy

$$w \leq \frac{\pi}{2}.$$

### Selection of a suitable base frequency

There are three additional considerations to be made for selecting an adequate $w$: *(i)* high frequencies equate fast system responses, *(ii)* high frequencies tend to have small input-to-state gains, and *(iii)* the presence of state constraints.

At first glance, it would seem that selecting a high value of $w$ would lead to fast system responses. However, this need not be the case, since the gain of the system tends to diminish as the frequency of the input increases, i.e., if $w$ is selected in the *high frequency* band of the system. If the gain is low, then $\{x_{hj}\}$ is very similar to a constant signal of value $x_e$, which results in HMPC controller behaving very similarly to the MPCT controller. Therefore, $w$ should be selected taking into account the gain of the system for that frequency.

A tentative lower bound for $w$ is then the highest frequency of the *low frequency* band of the system. However, a final consideration can be made with regard to the system constraints as follows: the presence of constraints can override the desire for frequencies with large system gains. For instance, take as an example a system with a static gain of 4 with an input $u$ subject to $|u| \leq 1$ and a state $x$ subject to $|x| \leq 2$. Then, selecting a $w$ whose Bode gain is close to the static gain of the system is not desirable because the amplitude of $\{u_{hj}\}$ will be limited by the constraints on $\{x_{hj}\}$. Therefore, we can select a higher frequency. In this case, a proper selection might be to chose $w$ as the frequency whose Bode gain is 2.

**Remark 6.16.** *If the system has multiple states/inputs, then the above considerations should be made extrapolating the idea to the frequency response of MIMO systems. One approach in this case is to focus on the slow dynamics (states) of the system, which are the most restrictive, in that they may require higher prediction horizons in order to be able to reach steady states. Additionally, it is also useful to identify if the system has any integrator states and to take into account their constraints as described in the above discussion.*

## 6.5 Advantages of the HMPC formulation

One of the issues of standard MPC formulations that use a terminal admissible invariant set $\mathcal{X}_t$ (4.7e) is that the system must be able to reach it within the prediction horizon $N$. Moreover, this set must contain $x_r$, or else the closed-loop system cannot converge to the reference. Therefore, the distance to the current reference is one of the mayor determining factors on the domain of attraction. This has an important effect on both the domain of attraction and the performance of the controller that is more pronounced for small prediction horizons and/or if a terminal equality constraint is used.

A similar phenomenon happens in the case of the MPCT formulation. In this case, the effect on the domain of attraction is not as pronounced, since the predicted state need only be able to reach any admissible steady state of the system, which is a less strict imposition than it having to reach the terminal set $\mathcal{X}_t$. However, this may still affect the performance of the controller, particularly so in the presence of integrator states and/or slew-rate input constraints due to them affecting the "speed" at which the system can move.

The HMPC formulation can provide significant advantages in both these aspects. First, the harmonic artificial reference is not necessarily an admissible steady state of the system. Instead, it defines an admissible periodic trajectory of the system that is, effectively, a decision variable of the formulation. Second, the reference does not need to belong to the artificial harmonic reference. Thus, the distance to the reference should have little impact on the domain of attraction of the controller.

These two advantages provide the formulation with performance and domain of attraction advantages when compared to the aforementioned alternatives, as we illustrate in the following with two case studies.

### 6.5.1 Performance advantages of the HMPC formulation

This section presents the potential performance advantages that can be obtained using the HMPC formulation (6.3) by comparing it to the one obtained by the MPCT formulation (5.42) to control the ball and plate system described in Section 5.7.2.

We initialize the system at the origin (i.e., with the ball positioned in the central point of the plate) and set the reference to

$$x_r = (1.8, 0, 0, 0, 1.4, 0, 0, 0), \quad u_r = (0, 0).$$

The HMPC controller is solved using version `v0.7.5` of the COSMO solver [153], while the MPCT controller is solved using version `0.6.0` of the OSQP solver [35]. These two solvers employ the same operator splitting approach, based on the ADMM algorithm [21]. In fact, their algorithms are very similar, with OSQP being particularized to QP problems. The options of both solvers are set to their default values, with the exception of the tolerances `eps_abs`, `eps_rel`, `eps_prim_inf` and `eps_dual_inf`, which are set to $10^{-4}$.

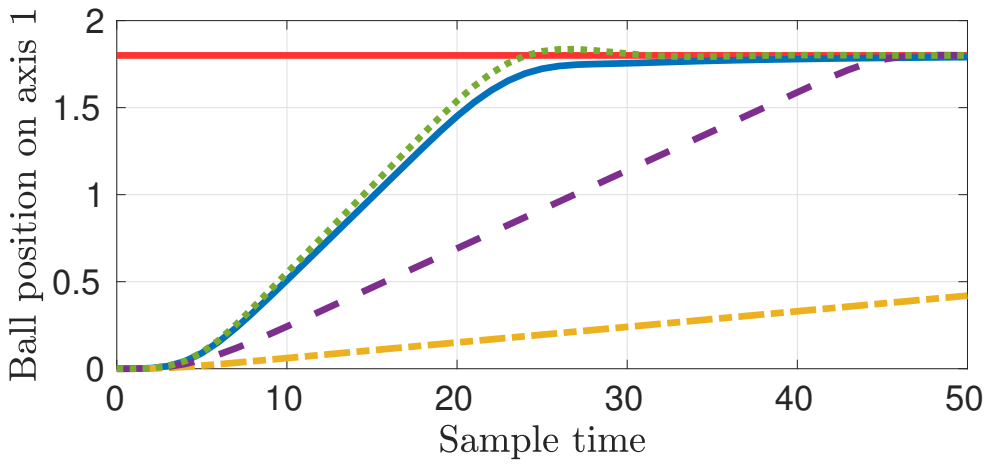

(a) Position of ball on axis 1.

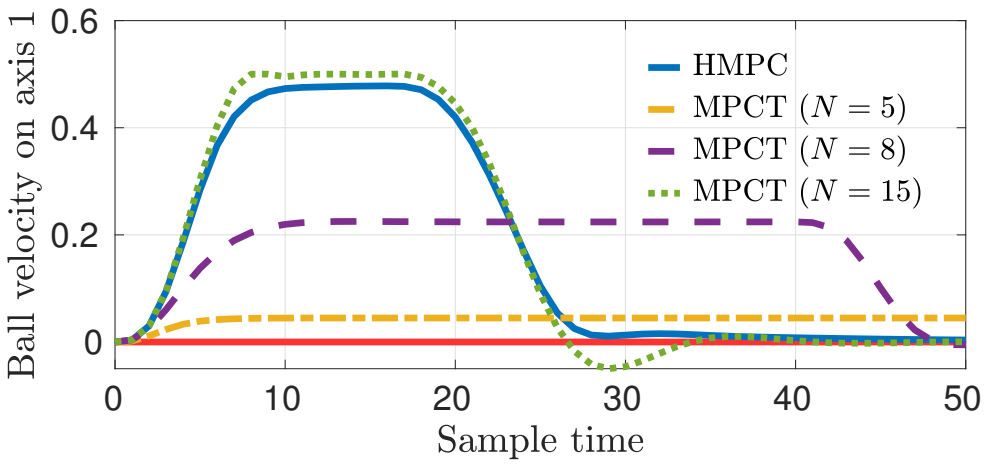

(b) Velocity of ball on axis 1.

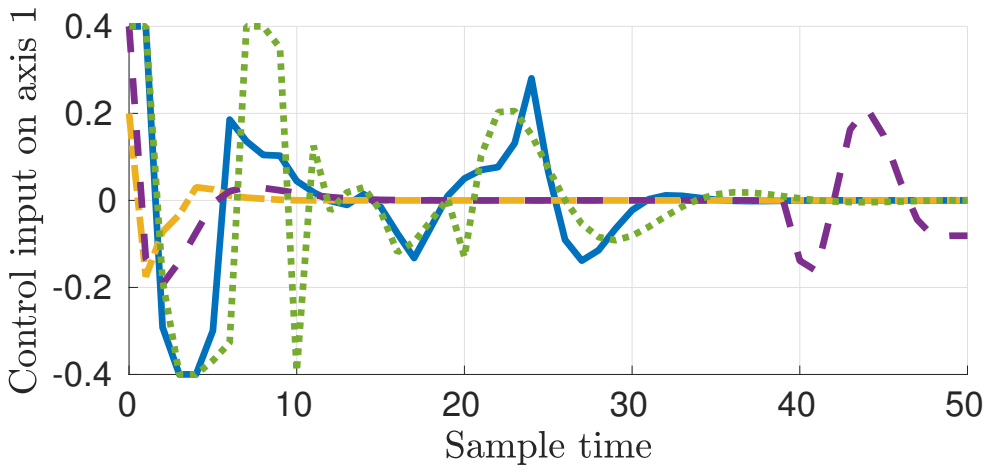

(c) Control input on axis 1.

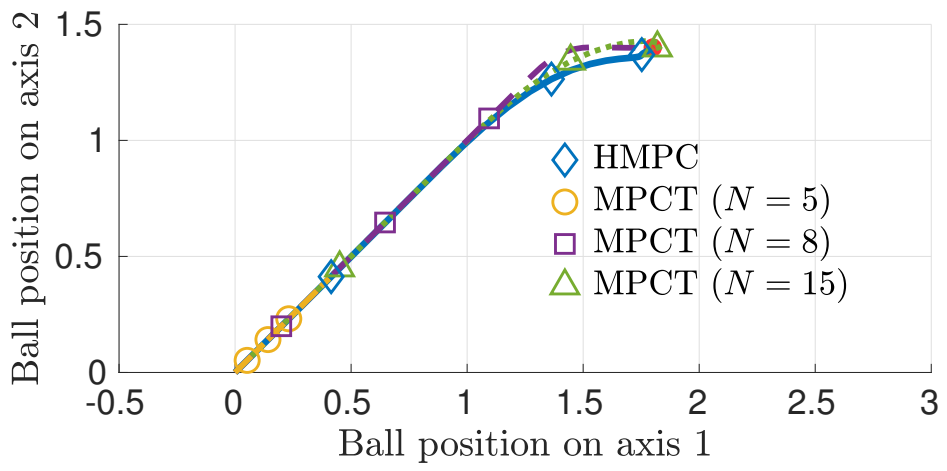

(d) Position of ball on the plate.

Figure 6.1: Closed-loop comparison between HMPC and MPCT.

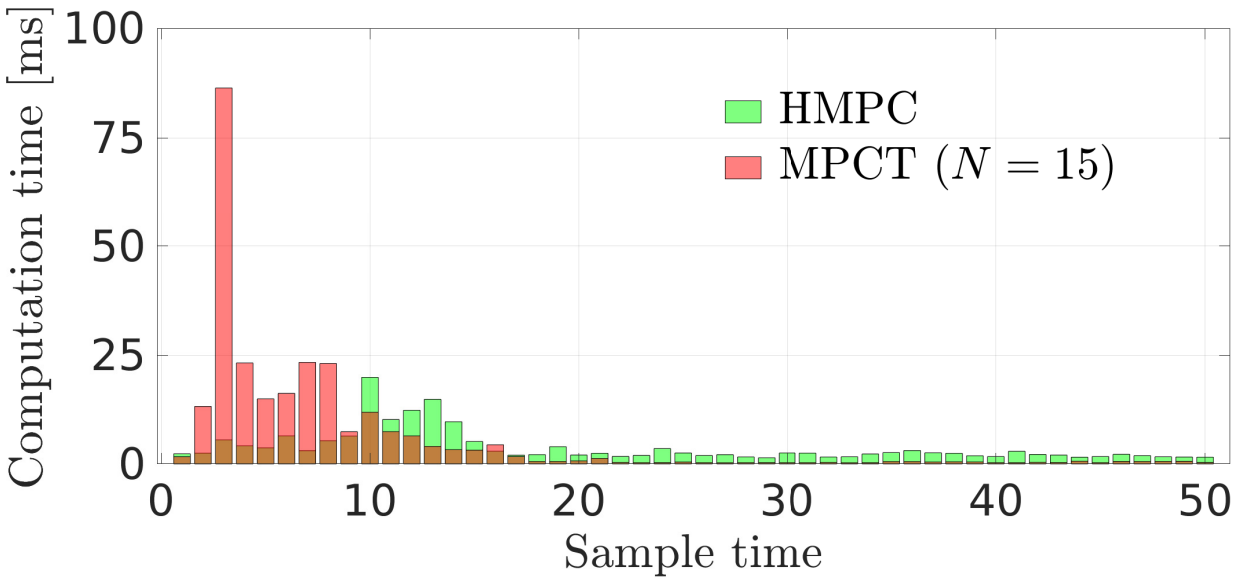

Figure 6.2: Computation times of the COSMO and OSQP solvers.

The parameters of the controllers, which where manually tuned to provide an adequate closed-loop performance, are described in Table 6.1[1]. We compare the HMPC controller with $N = 5$ to three MPCT controllers with prediction horizons $N = 5, 8, 15$. The prediction horizon $N = 15$ was chosen by finding the lowest value for which the MPCT performed well. We measure performance as

$$\Phi \doteq \sum_{k=1}^{N_{\text{iter}}} \|x_k - x_r\|_Q^2 + \|u_k - u_r\|_R^2,$$

where $x_k$, $u_k$ are the states and control actions at each sample time $k$ throughout the simulation and $N_{\text{iter}} = 50$ is the number of sample times. Table 6.2 shows the performance index for each one of the controllers.

[1]Notice that they are the same as the ones used in Section 5.8.6 for the MPCT case study.

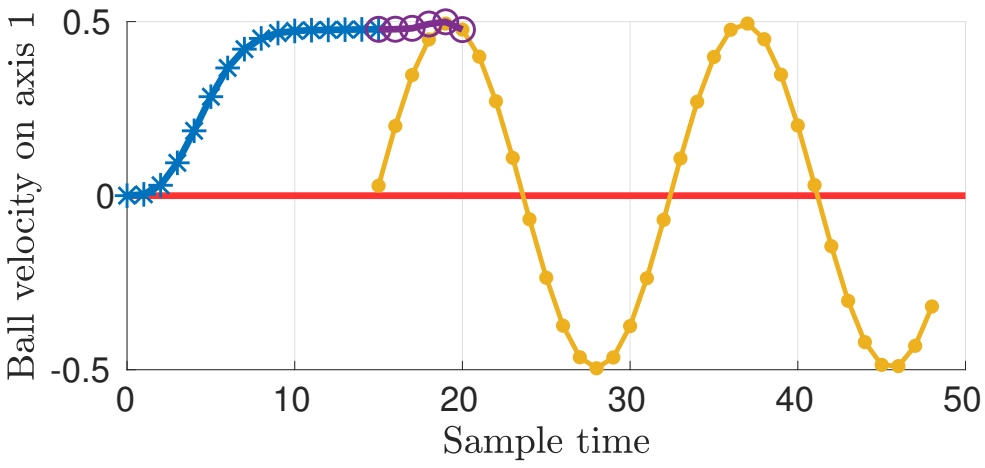

(a) HMPC: Velocity of ball on axis 1.

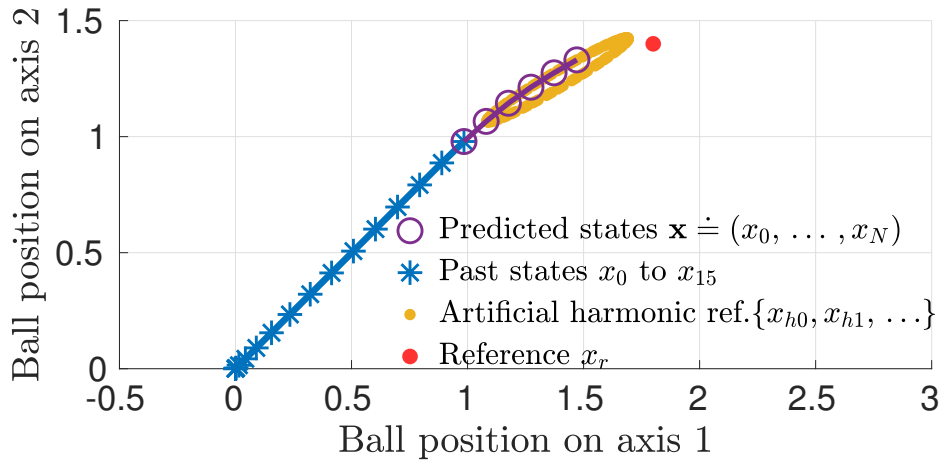

(b) HMPC: Position of ball on plate.

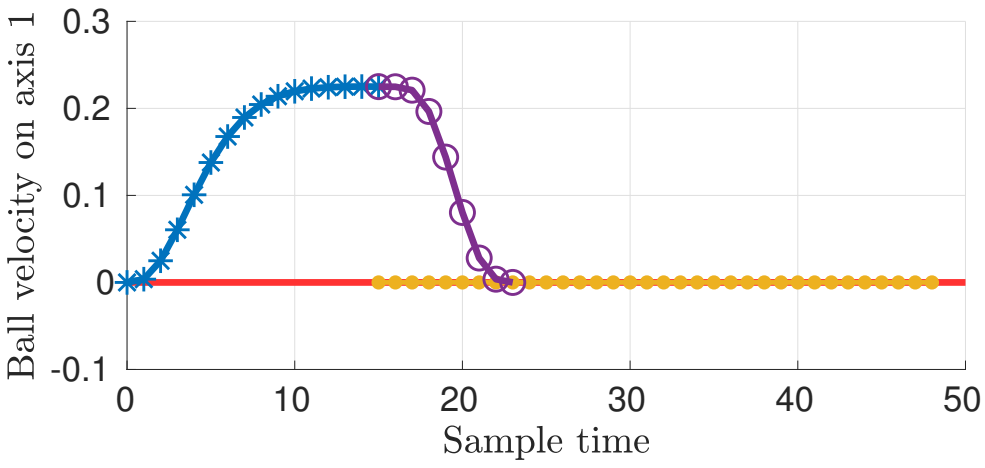

(c) MPCT: Velocity of ball on axis 1.

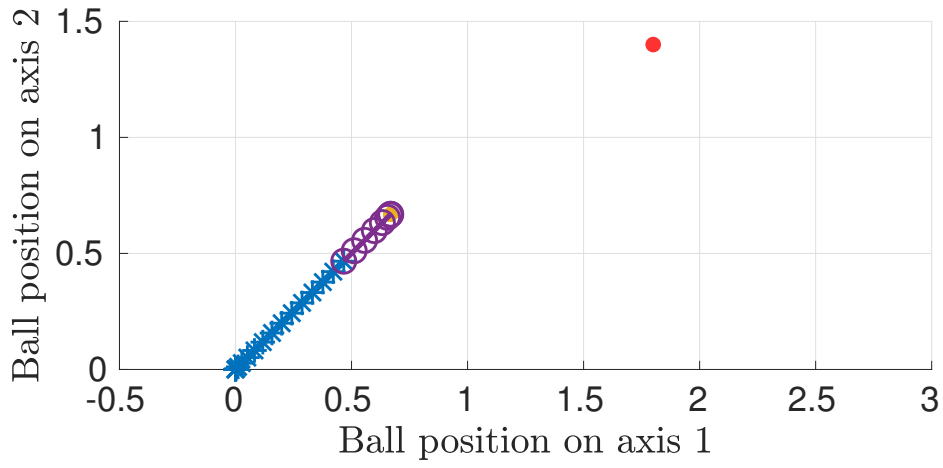

(d) MPCT: Position of ball on plate.

Figure 6.3: Snapshot of HMPC and MPCT at iteration 15.

Figure 6.1 shows the closed-loop simulation results. Figures 6.1a and 6.1b show the position and velocity of the ball on axis 1, i.e., $p_1$ and $\dot{p}_1$, respectively. Figure 6.1c shows the control input on axis 1, i.e., $\ddot{\theta}_1$. Finally, Figure 6.1d shows the trajectory of the ball on the plate. The markers indicate the position of the ball at sample times 10, 20 and 30 for each one of the controllers. The computation times of the HMPC controller and the MPCT controller with the prediction horizon $N = 15$ are shown in Figure 6.2. As can be seen, the COSMO solver applied to the HMPC problem provides computation times that are reasonable when compared to the results of the OSQP solver applied to the MPCT problem, in spite of the fact that the OSQP solver is particularized to QP problems, whereas the COSMO solver is not particularized to the second order cone programming problem (6.3). Our intent with this figure is to show that, even though (6.3) is a more complex problem than (5.42) due to the inclusion of second order cone constraints, it can still be solved in reasonable times using state of the art solvers. We note that COSMO runs on the *Julia* programming language, while OSQP is programmed in *C* and executed using its Matlab interface.

Notice that the velocities obtained with the MPCT controllers with small prediction horizons are far away from the upper bound of 0.5. The HMPC controller, on the other hand, reaches much higher velocities even though its prediction horizon is also small. This results in a much faster convergence of the HMPC controller, as can be seen in Figure 6.1. If the prediction horizon of the MPCT controller is sufficiently large (e.g. $N = 15$), then this issue no longer persists.

To understand why this happens, let us compare the solution of the HMPC controller with the MPCT controller with $N = 8$. Figure 6.3 shows a snapshot of sample time $k = 15$ of the same simulation shown in Figure 6.1. Lines marked with an asterisk are the past states from sample time $k = 0$ to the *current* state at sample time $k = 15$, those marked with circumferences are the predicted states $\{x_j\}$ for $j \in \mathbb{Z}_0^N$, and those marked with dots are the artificial reference. The position of the markers line up with the value of the signals at each sample time, e.g., each asterisk marks the value of the state at each sample time $k \in \mathbb{Z}_0^{15}$. Figures 6.3a and 6.3c show the velocity $\dot{p}_1$ of the ball on axis 1 for the HMPC and MPCT controllers, respectively. Figures 6.3b and 6.3d show the position of the ball on the plate.

| Parameter | Value | | |
|---|---|---|---|
| $Q$ | diag$(10, 0.05, 0.05, 0.05, 10, 0.05, 0.05, 0.05)$ | | |
| $T_e$ | diag$(600, 50, 50, 50, 600, 50, 50, 50)$ | | |
| **Parameter** | **Value** | **Parameter** | **Value** |
| $R$ | diag$(0.5, 0.5)$ | $S_e$ | diag$(0.3, 0.3)$ |
| $T_h$ | $T_e$ | $S_h$ | $0.5S_e$ |
| $T$ | $T_e$ | $S$ | $S_e$ |
| $N$ | 5 (HMPC), 8 and 15 | $\epsilon$ | $(10^{-4}, 10^{-4}, 10^{-4})$ |
| $w$ | 0.3254 | | |

Table 6.1: Parameters of the controllers

| Controller | MPCT | | | HMPC |
|---|---|---|---|---|
| Prediction horizon ($N$) | 5 | 8 | 15 | 5 |
| Performance ($\Phi$) | 2014.03 | 844.16 | 488.88 | 511.09 |

Table 6.2: Performance comparison between controllers

The reason why the velocity does not exceed $\approx$0.2 with the MPCT controller can be seen in Figure 6.3c. The predicted states of the MPCT controller must reach a steady state at $j = N$ due to constraint (5.42i). In our example this translates into the velocity having to be able to reach 0 within a prediction window of length $N = 8$. This is the reason that is limiting the velocity of the ball. A velocity of 0.5 is not attainable with an MPCT controller with a prediction horizon of $N = 8$ because there are no admissible control input sequences $\mathbf{u}$ capable of steering the velocity from 0.5 to 0 in 8 sample times. This issue does not occur with the HMPC controller because it does not have to reach a steady state at the end of the prediction horizon, as can be seen in Figure 6.3a. Instead, it must reach an admissible harmonic reference, which can have a non-zero velocity.

**Remark 6.17.** *It is clear from this discussion, and the results of the MPCT controller with $N = 15$, that this issue will become less pronounced as the prediction horizon is increased. However, for low values of the prediction horizon, the HMPC formulation can provide a significantly better performance than MPCT.*

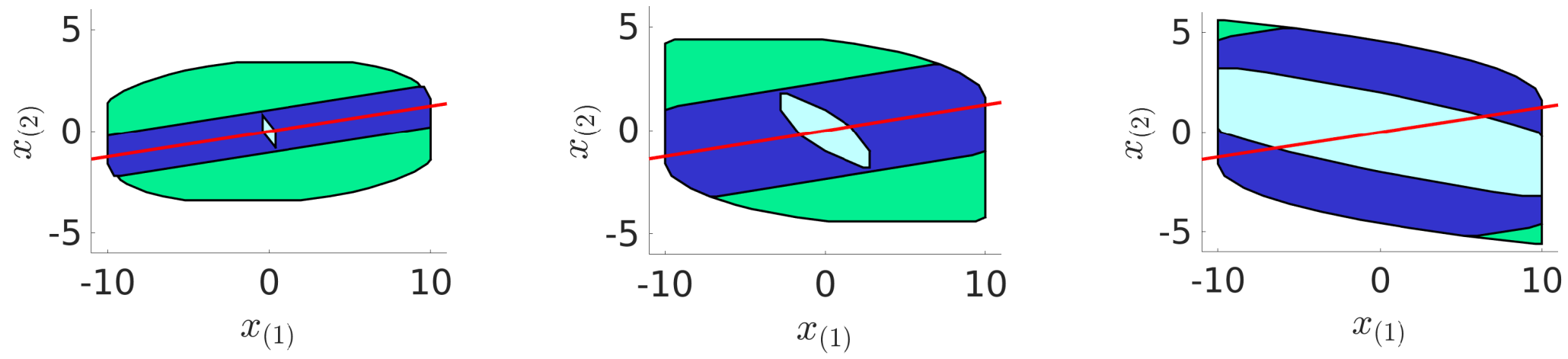

(a) $N = 2$. (b) $N = 4$. (c) $N = 8$.

Figure 6.4: Effect of $N$ on the domain of attraction for $w = \frac{4\pi}{64}$. $\mathcal{V}_{\text{HMPC}}$ in green, $\mathcal{V}_{\text{MPCT}}$ in blue and $\mathcal{V}_{\text{equMPC}}$ in cyan. The red line is the set of equilibrium points of the system.

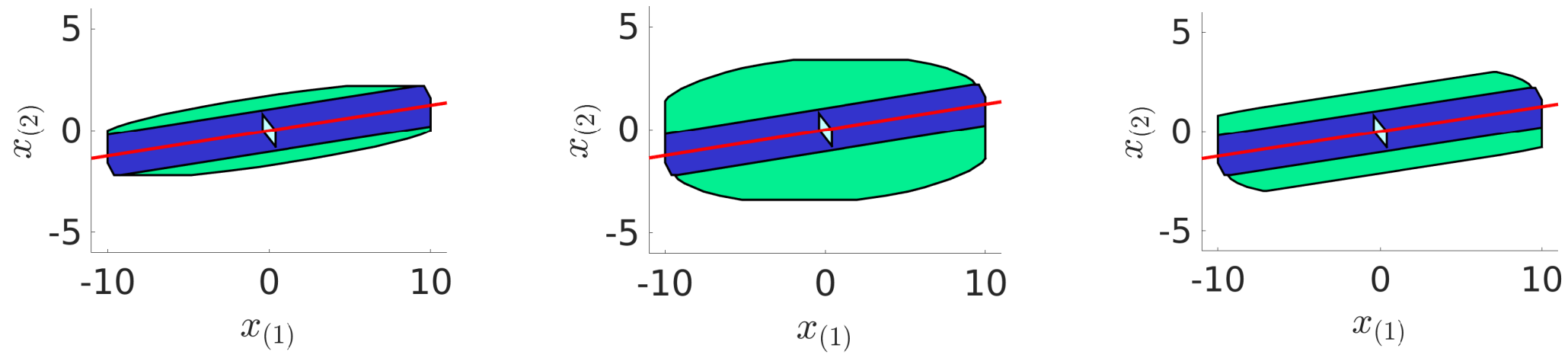

(a) $w = \pi/64$. (b) $w = 4\pi/64$. (c) $w = 11\pi/64$.

Figure 6.5: Effect of $w$ on the domain of attraction for $N = 2$. $\mathcal{V}_{\text{HMPC}}$ in green, $\mathcal{V}_{\text{MPCT}}$ in blue and $\mathcal{V}_{\text{equMPC}}$ in cyan. The red line is the set of equilibrium points of the system.

### 6.5.2 Domain of attraction

This section presents the potential advantages that can be obtained using the HMPC formulation (6.3) in terms of domain of attraction in comparison with the MPCT formulation (5.42) and the standard MPC formulation with terminal constraints (5.16), which we label by *equMPC* in the following. We denote by $\mathcal{V}_{\text{HMPC}}$, $\mathcal{V}_{\text{MPCT}}$ and $\mathcal{V}_{\text{equMPC}}$ the volumes of the domain of attraction of the HMPC, MPCT and equMPC formulations, respectively.

In order to be able to represent the domains of attraction, we consider an academic example where the model (4.1) and constraints (4.4) are given by:

$$A = \begin{bmatrix} 0.9 & 0.8 \\ 0 & 1 \end{bmatrix}, \; B = \begin{bmatrix} 0.8 \\ 1 \end{bmatrix}, \; E = \begin{bmatrix} 1 & 0 \\ 0 & 0 \end{bmatrix}, \; F = \begin{bmatrix} 0 \\ 1 \end{bmatrix},$$

$$\underline{y} = (-10, -0.5), \quad \overline{y} = (10, 0.5).$$

We note that the above constraints can also be expressed as box constraints (4.3).

The optimization problems of the three MPC formulations are programmed using the YALMIP Matlab toolbox [103] and are solved using the CPLEX solver. The values of $\mathcal{V}_{\text{HMPC}}$, $\mathcal{V}_{\text{MPCT}}$ and $\mathcal{V}_{\text{equMPC}}$ (which in this example are areas due to $x \in \mathbb{R}^2$) are computed using the MPT3 toolbox for Matlab [148]. The value of the cost function matrices are $Q = 5\mathbf{I}_2$, $R = 5$, $T = T_e = T_h = 10\mathbf{I}_2$ and $S = S_e = S_h = 0.1$.

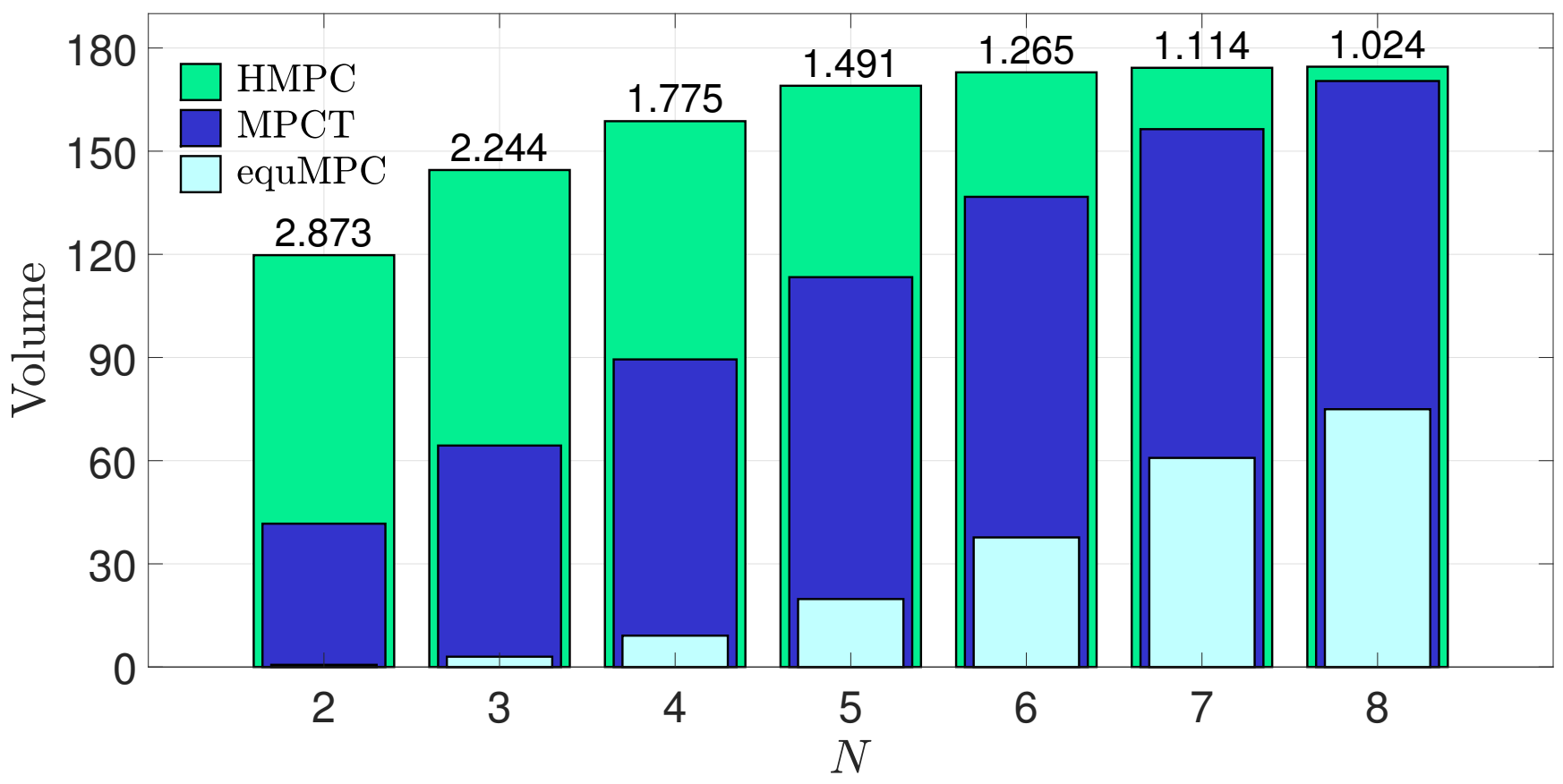

(a) Volumes of the domains of attraction for increasing $N$. $w = 4\pi/64$.

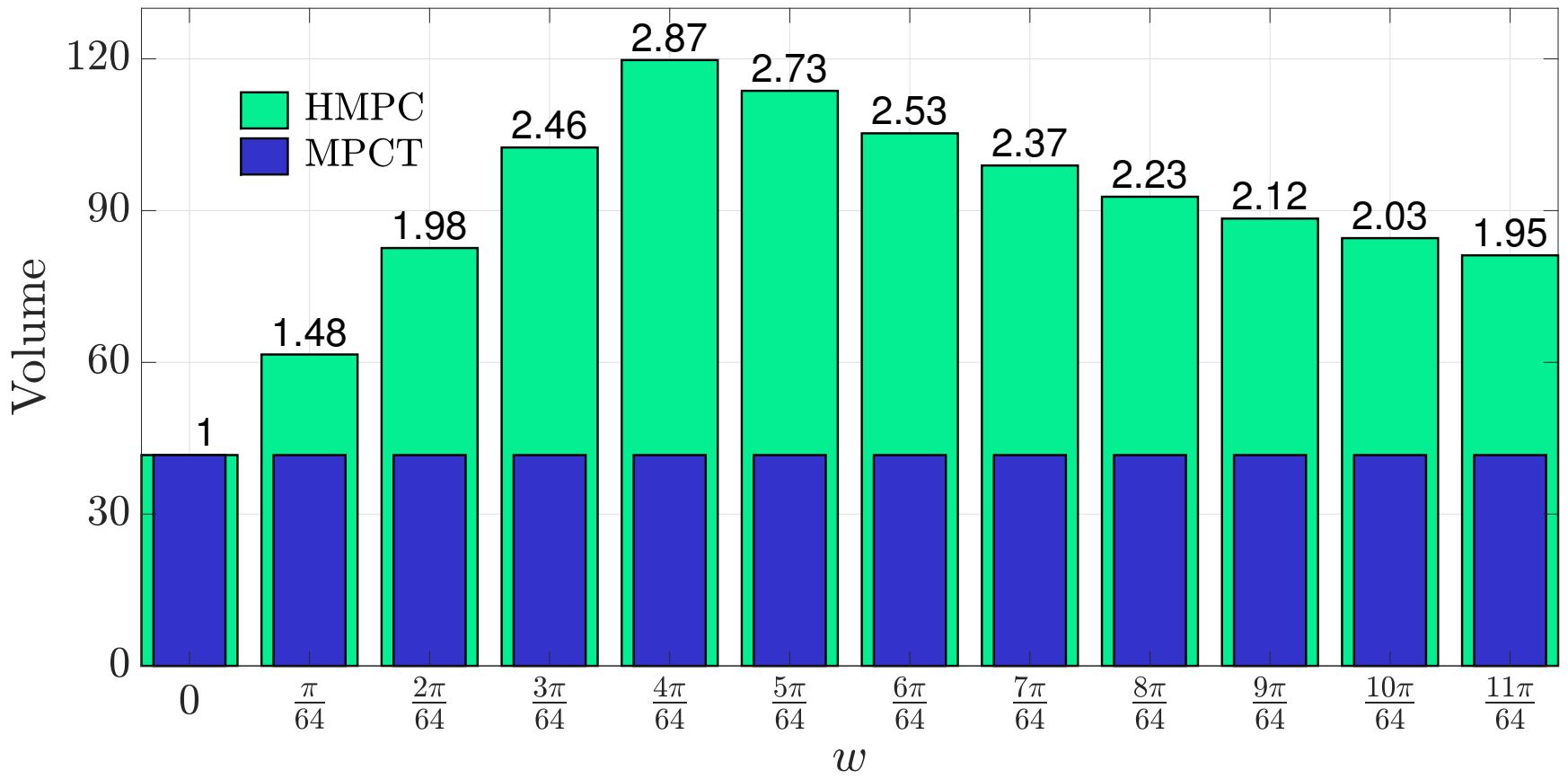

(b) Volumes of the domains of attraction for increasing $w$. $N = 2$.

Figure 6.6: Volumes of the domains of attraction for different values of $N$ and $w$. The numbers above the bars show the value of $\mathcal{V}_{\text{HMPC}}/\mathcal{V}_{\text{MPCT}}$.

Figure 6.4 shows the comparison between the domains of attraction of the MPC formulations for three values of $N$ with a fixed value of $w = 4\pi/64$. Note how the domains of attraction approach a maximum size as $N$ increases. Figure 6.6a shows the volume of the domains of attraction for increasing values of $N$. The numbers above the bars show the value of $\mathcal{V}_{\text{HMPC}}/\mathcal{V}_{\text{MPCT}}$.

Since the number of decision variables of the HMPC controller with prediction horizon $N$ is equal to the number of decision variables of the MPCT controller with prediction horizon $N+2$, a reasonable comparison (in terms of the expected computational complexity per iteration of the solvers) is to compare if $\mathcal{V}_{\text{HMPC}}$ for $N$ is larger than $\mathcal{V}_{\text{MPCT}}$ for $N+2$. Figure 6.6a shows this to be the case for the values of $N \in \mathbb{Z}_2^6$ in this particular example. A similar argument can be made for a comparison between $\mathcal{V}_{\text{HMPC}}$ with $\mathcal{V}_{\text{equMPC}}$.

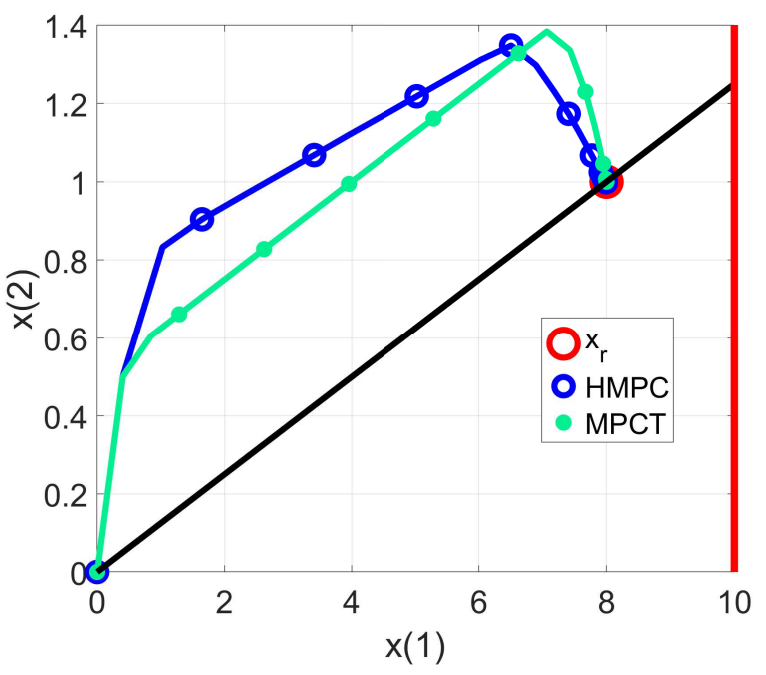

(a) State trajectory for $x_{r1}$

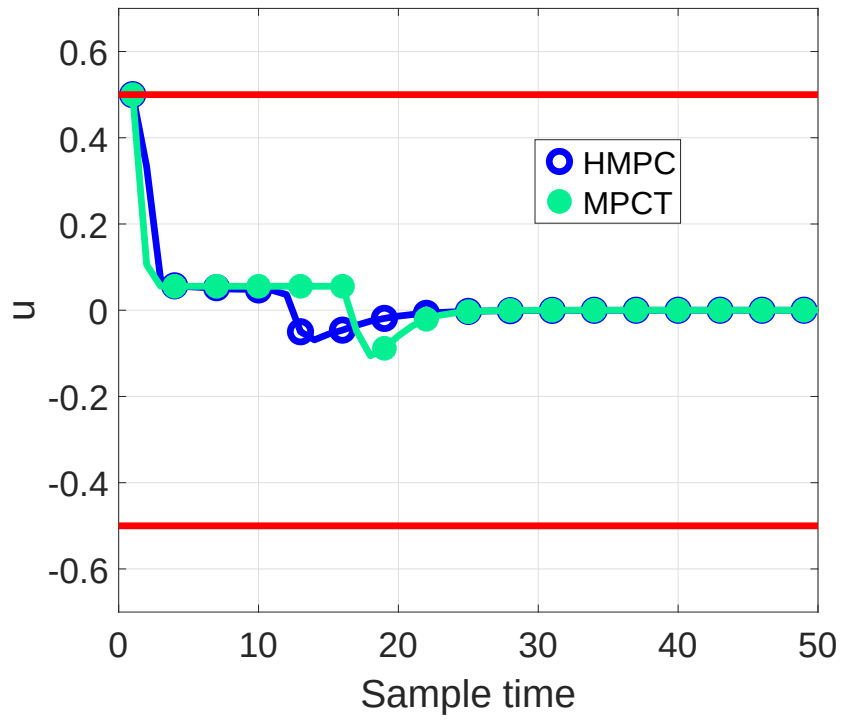

(b) Control input for $x_{r1}$

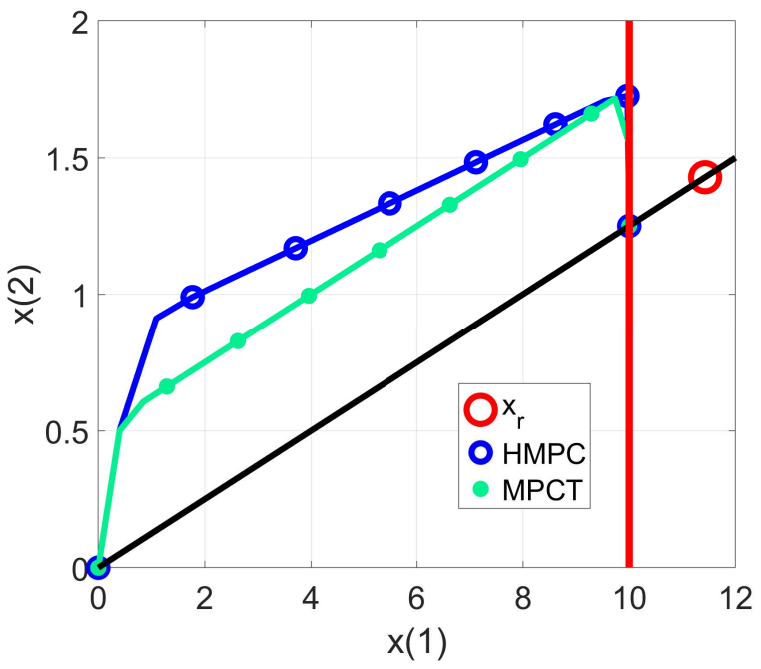

(c) State trajectory for $x_{r2}$

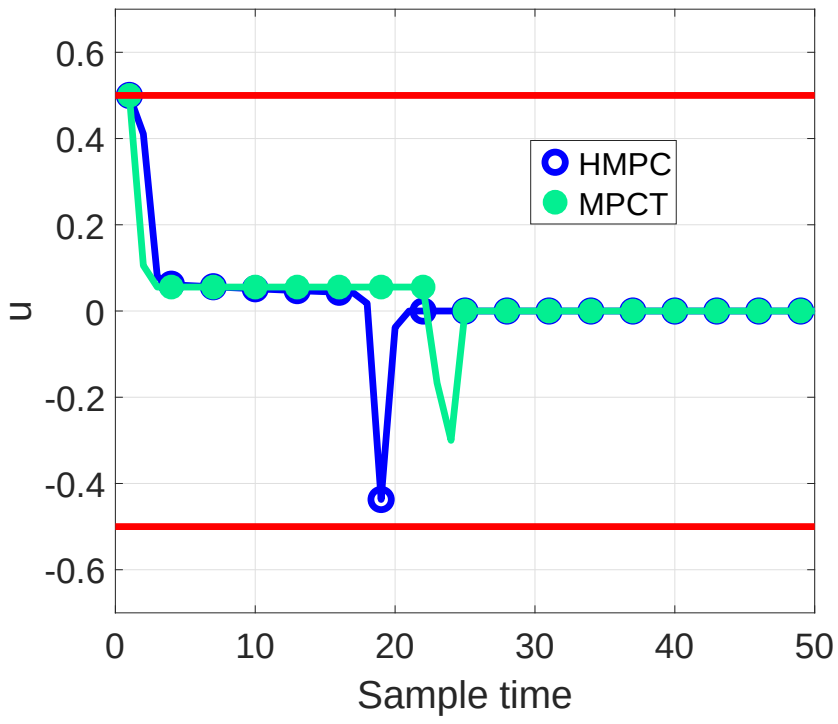

(d) Control input for $x_{r2}$

Figure 6.7: Closed-loop state and control input trajectories. The black line is the set of equilibrium points of the system. The red lines are the state and input bounds.

Figure 6.5 shows a comparison between the domains of attraction of the three controllers for different values of $w$ with a fixed value of $N = 2$. Figure 6.6b shows the value of $\mathcal{V}_{\text{HMPC}}$ for increasing values of $w$. The numbers above the bars show the value of $\mathcal{V}_{\text{HMPC}}/\mathcal{V}_{\text{MPCT}}$. Notice that the HMPC formulation with $w = 0$ has the same domain of attraction as the MPCT formulation, which is not surprising, since in that case the HMPC and MPCT formulations are equivalent.

Finally, we show a comparison between the closed-loop results using the HMPC and MPCT formulations.

We select $w = 4\pi/64$ and $N = 2$. Two tests are shown, both of which are initialized at the origin. The reference of each test is given by,

$$x_{r1} = (8, 1), \quad x_{r2} = (11.4286, 1.4286), \quad u_r = 0.$$

Reference $(x_{r1}, u_r)$ is admissible, whereas $(x_{r2}, u_r)$ is not.

Figures 6.7a and 6.7b show the state and control input trajectories for the reference $(x_{r1}, u_r)$, whereas Figures 6.7c and 6.7d shows them for the reference $(x_{r2}, u_r)$. The red circles indicate the reference; the vertical red line in Figures 6.7a and 6.7c represents the bounds on $x$; and the horizontal red lines in Figures 6.7b and 6.7d represent the bounds on $u$.

As can be seen in Figure 6.7c, if the reference is non-admissible the MPCT and HMPC controllers steer the system to the closest admissible steady state to the reference for the criterion of closeness determined by their terminal costs. In this case, this state happens to be the closest steady state to the reference in the standard Euclidean sense.

## 6.6 Conclusions and future lines of work

The results shown in Section 6.5 indicate that the HMPC formulation can outperform the MPCT and standard MPC formulations in terms of performance and domain of attraction, particularly so for small prediction horizons. We find the performance advantage to be particularly noticeable in systems with integrators and/or slew rate constraints. This is because these systems tend to suffer from the problem we highlight in Figure 6.3.

The fact that these advantages are obtained without requiring the computation of an admissible invariant set, along with the proven recursive feasibility and asymptotic stability of the formulation, make it an interesting candidate for future study on its implementation in embedded systems. The preliminary results we show in Figure 6.2 suggest that the development of a sparse solver for the HMPC formulation could make it suitable for its implementation in these devices. A preliminary solver for the HMPC formulation, loosely based on the COSMO solver [153], shows promising results. However, the development of this solver is still in early stages and was thus omitted in this dissertation.

A possible future line of work is to expand the HMPC formulation to artificial harmonic signals with more harmonics, i.e., taking $\{x_{hj}\}$ and $\{u_{hj}\}$, $j \in \mathbb{Z}$, as

$$x_{hj} = x_e + \sum_{i=1}^{p} x_{si} \sin(w_i(j-N)) + x_{ci} \cos(w_i(j-N)),$$
$$u_{hj} = u_e + \sum_{i=1}^{p} u_{si} \sin(w_i(j-N)) + u_{ci} \cos(w_i(j-N)),$$

for some $p \in \mathbb{Z}_{>0}$, where $w_i = iw$ for $i \in \mathbb{Z}_1^p$, and $w \in \mathbb{R}_{>0}$ is the base frequency. The question is whether any benefits are obtained from the addition of more harmonics compared to the single harmonic version we present in this dissertation.

Finally, another line of future work is to study the use of the HMPC formulation, or a variant of it, to tracking periodic references. The idea would be to determine if we can track generic periodic references in such a way in which the number of constraints added due to the artificial reference does not grow with the period of the reference signal.

# List of Algorithms

# Index